\documentclass[fontsize=11pt,twoside=false,draft=false]{scrartcl}
\usepackage{hyperref}
\hypersetup{plainpages=false,
            linkcolor=black,
            citecolor=black,
            urlcolor=black,
            colorlinks=true}
            
\usepackage[utf8]{inputenc}

\usepackage[automark]{scrlayer-scrpage}

\usepackage[english]{babel}

\usepackage[intlimits]{amsmath}
\usepackage{amsfonts}
\usepackage{amstext}
\usepackage{amsopn}
\usepackage{mathabx}
\usepackage{url}

\usepackage{tikz}
\usetikzlibrary{calc}
\usetikzlibrary{patterns}
\tikzset{
  nw/.style={pattern=north west lines, line width=0.1pt, pattern color=#1},
  nw/.default=black
}
\tikzset{
  ne/.style={pattern=north east lines,line width=0.1pt, pattern color=#1},
  ne/.default=black
}

\tikzset{
  every picture/.append style={
    line join=round,
    line cap=round,
  }
}
\usepackage{longtable}
\usepackage{graphicx}
\usepackage[format=plain]{caption}
\usepackage{subfig}
\usepackage{wrapfig}
\usepackage{enumitem}
\usepackage{nicefrac}

\tikzstyle{every picture}+=[font=\footnotesize]

\usepackage{amsmath,amssymb,mathtools,commath}
\usepackage{cases}
\usepackage[amsmath,hyperref,thmmarks]{ntheorem}
\numberwithin{equation}{section}
\numberwithin{figure}{section}
\numberwithin{table}{section}

\makeatletter
\theoremheaderfont{\bfseries\upshape}
\theoremseparator{.}
\newtheorem{theorem}{Theorem}[section]
\newtheorem{lemma}[theorem]{Lemma}
\newtheorem{proposition}[theorem]{Proposition}

\newtheorem{corollary}[theorem]{Corollary}
\theoremsymbol{}
\theorembodyfont{\upshape}

\theoremseparator{}

\theoremstyle{break}

\theoremstyle{plain}
\theoremheaderfont{\itshape}
\newtheorem{remark}[theorem]{Remark}

\theoremstyle{plain}
\theoremseparator{.}
\newtheorem*{proof2}{Proof}
\theoremsymbol{\ensuremath{\Box}}
\newtheorem*{proof}{Proof}
\newtheorem*{outlineproof}{Outline of the proof}
\theoremseparator{}
\newtheorem*{proofof}{Proof of}

\usepackage{algorithm}
\usepackage{algpseudocode}

\algrenewcommand\algorithmicrequire{\textbf{Input:}}
\algrenewcommand\algorithmicensure{\textbf{Output:}}
\floatname{algorithm}{Algorithm}

\usepackage{dirtree}
\usepackage{zref}

\usepackage{booktabs}

\let\oldbullet\bullet
  \newlength{\raisebulletlen}
  \setbox1=\hbox{$\bullet$}\setbox2=\hbox{\tiny$\bullet$}
  \renewcommand\bullet{\raisebox{\raisebulletlen}{\,\tiny$\oldbullet$}\,}
\pgfdeclarepatterninherentlycolored[\hatchColor]{myLines}{\pgfqpoint{-1pt}{-1pt}}{\pgfqpoint{4pt}{4pt}}{\pgfqpoint{3pt}{3pt}}%
{
    \pgfsetfillcolor{\hatchColor!10}
    \pgfpathrectangle{\pgfqpoint{0pt}{0pt}}{\pgfqpoint{3.1pt}{3.1pt}}
    \pgfusepath{fill}
    \pgfsetlinewidth{0.8pt}
    \pgfsetcolor{\hatchColor}
    \pgfpathmoveto{\pgfqpoint{-0.4pt}{-0.4pt}}
    \pgfpathlineto{\pgfqpoint{3.4pt}{3.4pt}}
    \pgfusepath{stroke}
}

\tikzset{
    hatchColor/.store in=\hatchColor, hatchColor=gray
}
\tikzset{slopetriangle/.style={
  bottom color=black!20,
  middle color=black!5,
  top color=white,
  draw=black
}}

\def\Xint#1{\mathchoice
{\XXint\displaystyle\textstyle{#1}}%
{\XXint\textstyle\scriptstyle{#1}}%
{\XXint\scriptstyle\scriptscriptstyle{#1}}%
{\XXint\scriptscriptstyle\scriptscriptstyle{#1}}%
\!\int}
\def\XXint#1#2#3{{\setbox0=\hbox{$#1{#2#3}{\int}$}
\vcenter{\hbox{$#2#3$}}\kern-.5\wd0}}
\newcommand{\intmean}{\Xint-}

\definecolor{dkgreen}{rgb}{0,0.6,0}
\definecolor{gray}{rgb}{0.5,0.5,0.5}
\definecolor{mauve}{rgb}{0.58,0,0.82}
\usepackage{listings}
\lstdefinestyle{fullsrcsmall}{
  basicstyle=\small\ttfamily,
  xleftmargin=1em,
}
\lstdefinestyle{fullsrcfnsize}{
  basicstyle=\footnotesize\ttfamily,
  xleftmargin=2em,
}
\lstdefinestyle{fullsrcscsize}{
  basicstyle=\scriptsize\ttfamily,
  numbers=none,
  xleftmargin=2em,
}
\lstdefinestyle{inline}{
  basicstyle=\ttfamily,
  numbers=none,
  xleftmargin=0em,
}
\lstdefinestyle{inlinesmall}{
  basicstyle=\small\ttfamily,
  numbers=none,
  xleftmargin=0em,
}
\lstdefinestyle{inlinefnsize}{
  basicstyle=\footnotesize\ttfamily,
  numbers=none,
  xleftmargin=1em,
}

\usepackage{a4wide}
\usepackage[capitalise, english]{cleveref}
\crefname{subsection}{Subsection}{Subsections}
\usepackage{pgfplots}
\usetikzlibrary{intersections, pgfplots.fillbetween,angles,quotes}
\usepackage{placeins}
\usepackage{enumitem}
\usepackage{blindtext}
\usepackage{chngcntr}
\usepackage{bm}
\makeatletter
\newcommand{\mylabel}[2]{#2\def\@currentlabel{#2}\label{#1}}

\newlist{mylist}{enumerate*}{1}
\setlist[mylist]{label=(\alph*), ref=\alph*}
\usepackage{xr}
\makeatother
\usepackage{caption}

\allowdisplaybreaks
\RedeclareSectionCommand[
  beforeskip=.25\baselineskip,
  afterskip=-1em]{paragraph}  
\RedeclareSectionCommand[
  beforeskip=.5\baselineskip,
  afterskip=.5\baselineskip]{subsection}  
  \RedeclareSectionCommand[
  beforeskip=.5\baselineskip,
  afterskip=.5\baselineskip]{subsubsection} 
  \RedeclareSectionCommand[
  beforeskip=.5\baselineskip,
  afterskip=.5\baselineskip]{section}        
\title{Direct guaranteed lower eigenvalue bounds with optimal a priori convergence rates for the bi-Laplacian
}
\author{
Carsten Carstensen\footnote{Department of Mathematics, Humboldt-Universit\"at zu Berlin, Unter den Linden 6, 10099 Berlin, Germany. \href{mailto:cc@math.hu-berlin.de}{cc@math.hu-berlin.de} 
and \href{mailto:puttkams@math.hu-berlin.de}{puttkams@math.hu-berlin.de}}
\and
Sophie Puttkammer$^*$
}
\date{ }

\begin{document}

\maketitle
\begin{abstract}\parindent 0em
	An extra-stabilised Morley finite element method (FEM) directly computes 
	guaranteed lower eigenvalue bounds with optimal a priori convergence rates for the bi-Laplace 
	Dirichlet eigenvalues. 
	The smallness assumption $\min\{\lambda_h,\lambda\}h_{\max}^{4}$ $\le 184.9570$ in $2$D 
	(resp. $\le 21.2912$ in $3$D) on the maximal mesh-size $h_{\max}$
	makes the computed $k$-th discrete eigenvalue  $\lambda_h\le \lambda$ a lower  eigenvalue bound for the $k$-th Dirichlet eigenvalue $\lambda$. 
	This holds for multiple and clusters  of eigenvalues and 
	serves for the localisation of the bi-Laplacian Dirichlet eigenvalues in particular for coarse meshes.
	The analysis requires interpolation error estimates for the Morley FEM with explicit constants in any space dimension $n\ge 2$, 
	which are of independent interest. 
	The convergence analysis in $3$D follows the Babu\v{s}ka-Osborn theory and relies on a companion operator for the Morley
	finite element method. This is based on the Worsey-Farin 3D version of the Hsieh-Clough-Tocher macro element with 
	a careful selection of center points in a further decomposition of each tetrahedron into 12 sub-tetrahedra. 
	Numerical experiments in 2D support the optimal convergence rates of the extra-stabilised Morley FEM 
	and suggest an adaptive algorithm with optimal empirical convergence rates. 
\end{abstract}
\paragraph{Keywords.}{\hspace{.5mm}
				biharmonic eigenvalue problem, direct guaranteed lower eigenvalue bounds, 
				Morley finite element, conforming companion, nonconforming interpolation, 
				Hsieh-Clough-Tocher, Worsey-Farin, a priori error estimates, adaptive mesh-refinement
				}

\section{Introduction}\label{sec:introduction}
	The biharmonic eigenvalue problem $\Delta^2 u =\lambda u$ 
	allows upper bounds from the Rayleigh-Ritz (or) min-max principle for \emph{conforming} finite element 
	methods (FEMs) \cite{BO91,Boffi2010}. 
	Guaranteed \emph{lower} eigenvalue bounds (GLB) can be even more relevant in a safety analysis in 
	computational mechanics, for the detection of spectral gaps, or for valid bounds of the 
	Sobolev embedding $H^2_0(\Omega)\hookrightarrow L^2(\Omega)$. There is a rich literature on lower eigenvalue bounds for the Laplace operator, cf., e.g.,  \cite{VejchodskySebestova2014,Vohraliketal2018,HongXieYueZhang2018} and the references therein. Throughout this paper the focus is on the biharmonic eigenvalue problem with former contribution in  \cite{YLBL12,CGal14,HuHuangLin2014,Liu2015,YangLiBi2016,LiaoShuLiu2019}.  
	
	\subsection{Motivation}
\begin{figure}[htb]
				\centering					
					\scalebox{0.8}{\begin{tikzpicture}

 \definecolor{newteal}{RGB}{27,158,119}
\definecolor{neworange}{RGB}{217,95,2}
\definecolor{newpurple}{RGB}{117,112,179}
\definecolor{newpink}{RGB}{231,41,138}
\definecolor{newgreen}{RGB}{102,166,30}
\definecolor{newyellow}{RGB}{230,171,2}
\definecolor{newred}{RGB}{109, 0, 0}
\definecolor{imblue}{rgb}{0,0.2157,0.4235}
\definecolor{newgreen2}{RGB}{150,200,50}

 
 \colorlet{col0}{imblue}
 \colorlet{col1}{newgreen2}
 \colorlet{col2}{newyellow}
 \colorlet{col3}{neworange}
 \colorlet{col4}{newred}
 
  \pgfplotsset{convergenceplot/.style={
    xlabel=\small{$|\mathcal{T}|$},
    ylabel=\small{error },
    legend style={legend pos=outer north east},
    ymajorgrids=true
  }}
  \pgfplotsset{cases/.style={
    nodes near coords,point meta=explicit symbolic
  }}

\begin{loglogaxis}[convergenceplot]{
			
\addlegendimage{col3, solid, very thick}
			\addlegendentry{$\lambda_k-\lambda_h(k)\phantom{00}$}  
\addlegendimage{col0,solid, very thick}
			\addlegendentry{$\lambda_k-\textup{GLB}(k)$}	
\addlegendimage{col1,solid, very thick}			
			\addlegendentry{$\lambda_k-\lambda_M(k)\phantom{0}$}						
						
\addlegendimage{empty legend} \addlegendentry{}
	\addlegendimage{ mark=*,gray, only marks}
					\addlegendentry{$k=1$}		
	\addlegendimage{ mark=square*,gray, only marks}
					\addlegendentry{$k=4$}	
\addlegendimage{empty legend} \addlegendentry{}
\addlegendimage{ solid, gray}
			\addlegendentry{$\theta=1\phantom{.0}$}	
\addlegendimage{ dashed, gray}
			\addlegendentry{$\theta=0.5$}

 \addplot [forget plot, col1, solid, mark=*, every mark/.append style= {solid,scale=0.9,fill=col1!60!white }]
   			table[x=nrElem, y=lambda-lambda_M]{Pictures/data/BiLaplace_newKappa/BiLaplace_Dumbbell_1_theta_1_lambda_1error.dat};		
    \addplot [forget plot, col1, dashed, mark=*, every mark/.append style= {solid,scale=0.9,fill=col1!60!white,  pattern=myLines,hatchColor=col1 }]
   			table[x=nrElem, y=lambda-lambda_M]{Pictures/data/BiLaplace_newKappa/BiLaplace_Dumbbell_1_theta_0.5_lambda_1error.dat};		
   			 			
 \addplot [forget plot, col0, solid, mark=*, every mark/.append style= {solid,scale=0.9,fill=col0!60!white }]
   			table[x=nrElem, y=lambda-GLB_M]{Pictures/data/BiLaplace_newKappa/BiLaplace_Dumbbell_1_theta_1_lambda_1error.dat};		
    \addplot [forget plot, col0, dashed, mark=*, every mark/.append style= {solid,scale=0.9,fill=col0!60!white,  pattern=myLines,hatchColor=col0 }]
   			table[x=nrElem, y=lambda-GLB_M]{Pictures/data/BiLaplace_newKappa/BiLaplace_Dumbbell_1_theta_0.5_lambda_1error.dat};

\addplot [forget plot, col3, solid, mark=*, every mark/.append style= {solid,scale=0.9,fill=col3!60!white }]
   			table[x=nrElem, y=lambda-GLB_Sk]{Pictures/data/BiLaplace_newKappa/BiLaplace_Dumbbell_1_theta_1_lambda_1error.dat};	
    \addplot [forget plot, col3, dashed, mark=*, every mark/.append style= {solid,scale=0.9,fill=col3!60!white,  pattern=myLines,hatchColor=col3 }]
   			table[x=nrElem, y=lambda-GLB_Sk]{Pictures/data/BiLaplace_newKappa/BiLaplace_Dumbbell_1_theta_0.5_lambda_1error.dat};	
 
 \addplot [forget plot, col1, solid, mark=*, every mark/.append style= {solid,scale=0.9,fill=col1!60!white }]
   			table[x=nrElem, y=lambda-lambda_M]{Pictures/data/BiLaplace_newKappa/BiLaplace_Dumbbell_1_theta_1_lambda_4error.dat};		
    \addplot [forget plot, col1, dashed, mark=*, every mark/.append style= {solid,scale=0.9,fill=col1!60!white,  pattern=myLines,hatchColor=col1 }]
   			table[x=nrElem, y=lambda-lambda_M]{Pictures/data/BiLaplace_newKappa/BiLaplace_Dumbbell_1_theta_0.5_lambda_4error.dat};		
   			
  \addplot [forget plot, col0, solid, mark=square*, every mark/.append style= {solid,scale=0.9,fill=col0!60!white }]
   			table[x=nrElem, y=lambda-GLB_M]{Pictures/data/BiLaplace_newKappa/BiLaplace_Dumbbell_1_theta_1_lambda_4error.dat};		
    \addplot [forget plot, col0, dashed, mark=square*, every mark/.append style= {solid,scale=0.9,fill=col0!60!white,  pattern=myLines,hatchColor=col0 }]
   			table[x=nrElem, y=lambda-GLB_M]{Pictures/data/BiLaplace_newKappa/BiLaplace_Dumbbell_1_theta_0.5_lambda_4error.dat};		

\addplot [forget plot, col3, solid, mark=square*, every mark/.append style= {solid,scale=0.9,fill=col3!60!white }]
   			table[x=nrElem, y=lambda-GLB_Sk]{Pictures/data/BiLaplace_newKappa/BiLaplace_Dumbbell_1_theta_1_lambda_4error.dat};	
    \addplot [forget plot, col3, dashed, mark=square*, every mark/.append style= {solid,scale=0.9,fill=col3!60!white,  pattern=myLines,hatchColor=col3 }]
   			table[x=nrElem, y=lambda-GLB_Sk]{Pictures/data/BiLaplace_newKappa/BiLaplace_Dumbbell_1_theta_0.5_lambda_4error.dat};	
   			
   			
    \shade[top color=black!5,bottom color=black!30]
            (axis cs: 2.00e+04,6.32e-02)
         -- (axis cs: 2.00e+04,2.00e-03)
         -- (axis cs: 6.32e+05,2.00e-03)
         -- cycle;
    \draw   (axis cs: 2.00e+04,6.32e-02)
         -- (axis cs: 2.00e+04,2.00e-03) node [midway,left] {\scriptsize \(1\)}
         -- (axis cs: 6.32e+05,2.00e-03) node [midway,below] {\scriptsize \(1\)}
         -- cycle;  			
 
    \shade[top color=black!5,bottom color=black!30]
            (axis cs: 2.00e+05,3.16e+00)
         -- (axis cs: 2.00e+05,1.00e+00)
         -- (axis cs: 2.00e+06,1.00e+00)
         -- cycle;
    \draw   (axis cs: 2.00e+05,3.16e+00)
         -- (axis cs: 2.00e+05,1.00e+00) node [midway,left] {\scriptsize \(0.5\)}
         -- (axis cs: 2.00e+06,1.00e+00) node [midway,below] {\scriptsize \(1\)}
         -- cycle;

 }
\end{loglogaxis}
\end{tikzpicture}}
					\caption{Convergence history plot for the error in the principal and in the fourth Dirichlet eigenvalue of the 							
								bi-Laplacian on uniform ($\theta=1$, solid) 
							 	and adaptive ($\theta=0.5$, dashed) triangulations of the dumbbell-slit domain from \cref{fig:StartTriangulation}.a.}
							\label{fig:adaptive_Dumbbell_intro}
			\end{figure}
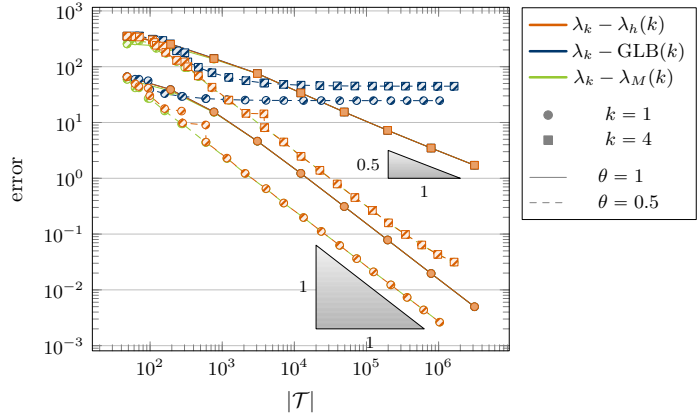

		A post-processing for the nonconforming Morley FEM allows for GLBs for the bi-Laplacian 
		in \cite{CGal14}. 
		The ${k}$-th discrete eigenvalue  $\lambda_{{M}}({k})$ computed from the Morley FEM 
		(displayed in \eqref{eq:disEVP_NC} below) leads to a guaranteed lower bound 
			\begin{align}
				\textup{GLB}({k}):=\frac{\lambda_{{M}}({k})}{1+\lambda_{{M}}({k})\kappa_{2}^2h_{\max}^{4}}
						\le \lambda_k \label{eq:GLB_CR_Morley}
			\end{align}
		for the exact ${k}$-th Dirichlet eigenvalue $\lambda_{k}$ of the  bi-Laplacian. 
		The explicit analytical parameter  $\kappa_2=0.25746$ 
		for ${n}=2$ is known from \cite{CGal14} and 
		$\kappa_2=0.21672$ 
		for ${n}=3$ is provided in \cref{thm:properties_IMorley}.\ref{item:IMorley_kappa} below. 
		The numerical experiments in this paper utilize the improved computational bound 
		$\kappa_2= 0.07353 $ from \cite{LiaoShuLiu2019} for $n=2$. 
		The maximal mesh-size $h_{\max}$ enters as a \emph{global} parameter in \eqref{eq:GLB_CR_Morley} and may cause 
		a significant underestimation for adaptive mesh-refinement, 
		when local mesh-refining leaves a few simplices coarse and $h_{\max}$ large. This leads to a dramatic underestimation in 
		the following motivational example with  
		convergence history plot \cref{fig:adaptive_Dumbbell_intro} 
		and useless post-processed bound \eqref{eq:GLB_CR_Morley}.
		The new method is an extra-stabilised Morley FEM  with an additional piecewise quadratic variable and the fine-tuned 
		parameter $\kappa_2= 0.07353 $. 
		The new method allows  for an optimal empirical convergence rate one  
		(with respect to the number  $|\mathcal{T}|$ of triangles in the triangulation $\mathcal{T}$) with adaptive mesh-refinement. 
		The dumbbell domain with a slit (see \cref{fig:StartTriangulation}.a  below for the initial triangulation) 
		is an extreme example. The adaptive refinement occurs in one of the two cells with minimal coupling, 
		so that $h_{\max}$ is not reduced. 
		The first and fourth Morley eigenvalue $\lambda_M(k)<\lambda_k$ for $k=1,4$ in \cref{fig:adaptive_Dumbbell_intro} are smaller 
		than the approximation  $\lambda_1= 80.93261350$ and $\lambda_4=386.80177939$ of the exact eigenvalues, 
		but this is \emph{not guaranteed} in general, cf. \cite[Sec.~2]{CGal14} for a counter example.
		The eigenvalues of the Morley FEM converge only asymptotically from below \cite{YLBL12,YangLiBi2016} and 
		it remains unclear whether a given triangulation belongs to the asymptotic regime. 
		Since $\lambda_h(k) h_{\max}^{4}\le\kappa_2^{-2}$ holds for all levels in 
		\cref{fig:adaptive_Dumbbell_intro}, \cref{thm:GLB} below implies that the discrete eigenvalue $\lambda_h(k)\le \lambda_k$ for $k=1,4$ is a 
		guaranteed lower eigenvalue bound under the hypothesis of the exact solve of the algebraic eigenvalue problem. (The discussion of interval arithmetic and perturbation analysis for inexact solve in numerical linear algebra is beyond this paper -- the focus here is on the understanding of the discretization error.) 
		The numerical results for $\textup{GLB}\le\lambda_h\le\lambda_M$ are almost indistinguishable for uniform mesh-refinement in 
		\cref{fig:adaptive_Dumbbell_intro} and 
		result in one line with empirical convergence rate $1$ for the principal and $1/2$ for the fourth eigenvalue.  
		
		In short, if the GLB relies exclusively on \eqref{eq:GLB_CR_Morley}, a naive adaptive mesh-refinement appears useless in this 
		example, while the new bound displays optimal empirical convergence rates. 
		\subsection{Eigenvalue problems and main results}
			The continuous eigenvalue problem seeks eigenpairs $(\lambda, u) \in \mathbb{R}^+\times V$ with
			 \begin{align}
				 	a(u,v)= \lambda\, b(u,v)\quad \text{for all } v\in V\quad\text{and}\quad 
				 	\Vert u\Vert_{L^2(\Omega)}=1\label{eq:contEVP}
			 \end{align}
			 in the Hilbert space $V:=H^{2}_0(\Omega)$ with the energy scalar product 
			 $a(\bullet,\bullet):=(D^2\bullet, D^2\bullet)_{L^2(\Omega)}$ for the Hessian $D^2$ and 
			 the $L^2$ scalar product $b(\bullet,\bullet):=(\bullet, \bullet)_{L^2(\Omega)}$;
			 the infinite but countable many eigenvalues $0<\lambda_1\le\lambda_2\le\dots$ with 
			 $\lim_{j\to \infty}\lambda_j=\infty$ in \eqref{eq:contEVP} are enumerated in ascending order counting 
			 multiplicities. 
			 For any shape-regular triangulation $\mathcal{T}$ of $\Omega\subset\mathbb{R}^n$ into simplices, 
			 the piecewise constant mesh-size function $h_{\mathcal{T}}\in P_0(\mathcal{T})$ is defined by 
			$h_{\mathcal{T}}|_T=h_T:=\textup{diam}(T)$ in each simplex $T\in\mathcal{T}$ and 
			 $h_{\max}:=\max_{T\in\mathcal{T}}h_T$. 
			 The discrete space 
			 $\boldsymbol{V_h}:= P_{2}(\mathcal{T})\times M({\mathcal{T}})\subset  P_{2}(\mathcal{T})\times 
			 P_{2}(\mathcal{T})$ 
			 consists of piecewise quadratic polynomials. The Morley space $M(\mathcal{T})$ 
			 is well established for two-dimensional plate problems \cite{Mor68} and generalized in \cite{MX06} 
			 for any space dimension (cf. \cref{sec:Morley} below for details).   
			 The algebraic eigenvalue problem of the extra-stabilised method seeks discrete eigenpairs 
			 $(\lambda_h, \boldsymbol{u_h}) \in \mathbb{R}^+\times (\boldsymbol{V_h}\setminus\{0\})$ with 
			\begin{align}
					\boldsymbol{a_h}(\boldsymbol{u_h},\boldsymbol{v_h})
					=\lambda_h \boldsymbol{b_h}(\boldsymbol{u_h},\boldsymbol{v_h})
					\quad\text{ for all }\boldsymbol{v_h}\in\boldsymbol{V_h}.  \label{eq:dis_EVP_alt}
			\end{align}
			The discrete scalar product $\boldsymbol{a_h}$ contains the  scalar product 
			$a_{\mathrm{pw}}(\bullet, \bullet):=(D^2_\mathrm{pw}\bullet,D^2_\mathrm{pw}\bullet)_{L^2(\Omega)}$ for the piecewise Hessian of
			the Morley functions in $M(\mathcal{T})$ and some stabilisation; the bilinear form 
			$\boldsymbol{b_h}$ is the $L^2$ scalar product of the piecewise quadratic components in 
			$P_{2}(\mathcal{T})$, 
			\begin{align*}
				\boldsymbol{a_h}(\boldsymbol{v_h},\boldsymbol{w_h})
						&:=a_{\mathrm{pw}}(v_{M}, w_{M})+\kappa_{2}^{-2}(h_{\mathcal{T}}^{-4}
						(v_{\mathrm{pw}}-v_{M}), w_{\mathrm{pw}}-w_{M})_{L^2(\Omega)},\\
				\boldsymbol{b_h}(\boldsymbol{v_h},\boldsymbol{w_h})&:=(v_{\mathrm{pw}},w_{\mathrm{pw}})_{L^2(\Omega)}
						\quad\text{for all }
						\boldsymbol{v_h}=(v_{\mathrm{pw}},v_{M}),\, 
						\boldsymbol{w_h}=(w_{\mathrm{pw}},w_{M})\in\boldsymbol{V_h}.
			\end{align*} 
			Since $(\boldsymbol{V_h},\boldsymbol{a_h})$ is a Hilbert space  and  $\boldsymbol{b_h}$ is a 
			semi-scalar product with kernel $\{0\}\times M(\mathcal{T})\subset \boldsymbol{V_h}$, 
			the algebraic eigenvalue problem \eqref{eq:dis_EVP_alt} has 
			$M:=\textup{dim}(P_{2}(\mathcal{T}))=\binom{{2}+n}{n}|\mathcal{T}|$ 
			finite and positive algebraic eigenvalues $0< \lambda_h(1)\le\dots\le\lambda_h(M)<\infty$ enumerated in ascending order counting 
			multiplicities. 		
			The new method \eqref{eq:dis_EVP_alt} directly computes  guaranteed lower eigenvalue bounds for any space 
			dimension 
			$n\ge2$. The a priori and a posteriori smallness assumption is explicit in terms of the maximal mesh-size 
			$h_{\max}$, but 
			surprisingly robust with respect to the shapes of the simplices in the triangulation $\mathcal{T}$. 
			The interpolation estimates of \cref{thm:properties_IMorley} below define the global parameter 
				\begin{align}					 		
					 	\kappa_1:=\sqrt{\frac{1}{\pi^2}+\frac{1}{2n(n+1)(n+2)}}
					 \quad \text{and}\quad
						\kappa_2:= \frac{\kappa_{1}}{\pi}	+\sqrt{\frac{n\kappa_{1}^2+2\kappa_{1} }{2(n-1)(n+1)(n+2)}}.	\label{eq:def_kappa}					
				\end{align}			
			\begin{theorem}[GLB]\label{thm:GLB}
			For any  $k=1,\dots, M$, the $k$-th eigenvalue $\lambda_k$ from \eqref{eq:contEVP} 
			and the $k$-th eigenvalue $\lambda_h(k)$ from \eqref{eq:dis_EVP_alt} satisfy that 
			$\min\{\lambda_h({k}),\lambda_k\} \kappa^2_2 h_{\max}^{4}\le 1 $ implies $\lambda_h({k})\le \lambda_{k}$.
			\end{theorem}
			Notice that $\min\{\lambda_h({k}),\lambda_k\} \kappa^2_2 h_{\max}^{4}\le 1 $ in \cref{thm:GLB} means that each of the conditions 
			\textup{(i)} $\lambda_k\kappa^2_2 h_{\max}^{4}\le 1$ (a priori) or 
			\textup{(ii)} $\lambda_h({k})\kappa^2_2 h_{\max}^{4}\le 1$ (a posteriori) implies the GLB property 
			$\lambda_h({k})\le \lambda_{k}$. Remarks~\ref{rem:kappa2prime_1}--\ref{rem:kappa2prime_2}  below explain that 
			the choice $\kappa_2= 0.07353 $ in $2$D (resp. $\kappa_2=0.21672$ in $3$D) is possible in the discrete system \eqref{eq:dis_EVP_alt} and in \cref{thm:GLB} and leads to the condition 
			$\min\{\lambda_h,\lambda\}h_{\max}^{4}\le 184.9570$ in $2$D 
			(resp. $\le 21.2912$ in $3$D) sufficient for $\lambda\le\lambda_h$.{ \phantom{xx} }
			
			\cref{thm:GLB} leads under some condition at least a posteriori to GLB and hence the next question is the quality of those. To describe optimal a priori convergence rates, let  $\mathbb{T}$ denote the set of uniformly shape-regular triangulation of a fixed bounded polyhedral 
			Lipschitz domain 
			$\Omega\subset\mathbb{R}^3$ into tetrahedra with respect to a global shape-regularity constant $C_{\mathrm{sr}}>0$: 
			Any tetrahedron $T\in\mathcal{T}\in\mathbb{T}$ with diameter $h_T$ and volume $|T|$ satisfies
			$|T|^{1/3}\le h_T\le  C_{\mathrm{sr}}|T|^{1/3}$.
			The subset $\mathbb{T}(\delta)\subset \mathbb{T}$ denotes the triangulations with maximal mesh-size 
			$h_{\max}\le \delta$. Let $\sigma:=\min\{1,\sigma_{\mathrm{reg}}\}$ denote the minimum of one 
			and the index of elliptic regularity $\sigma_{\mathrm{reg}}>0$ from \eqref{eq:def_sigma} below. 
			\begin{theorem}[a priori convergence]\label{thm:BabuskaOsborn} 
				Suppose $\lambda$ is an eigenvalue of \eqref{eq:contEVP}  of multiplicity $\mu$ with eigenspace 
				$E(\lambda)\subset H^{2+t}(\Omega)\cap V$  
				for some $t$ with 
				$\sigma \le t\le 1$. Then 
				there exist $\delta,C>0$ such that any triangulation $\mathcal{T}\in\mathbb{T}(\delta)$ and the discrete space 
				$\boldsymbol{V_h}:=P_2(\mathcal{T})\times M(\mathcal{T})$ lead in \eqref{eq:dis_EVP_alt} to 
				exactly $\mu$ algebraic eigenvalues  
				$\lambda_{h,1},\dots, \lambda_{h,{\mu}}$ of \eqref{eq:dis_EVP_alt} (counting multiplicities), that converge to 
				$\lambda$ as $h_{\max}\to 0$. 
				Let $E_h:=\textup{span}\{\boldsymbol{u_h}\in E_h(\lambda_{h,k}): k=1,\dots,\mu\}$ 
				abbreviate the span of the discrete eigenspaces 
				$E_h(\lambda_{h,k})\subset \boldsymbol{V_h}$			
				of $\lambda_{h,k}$ for $k=1,\dots, \mu$. 
				Then 
				\begin{align*}
						h_{\max}^{-t}\max_{k=1,\dots,\mu}|\lambda-\lambda_{h,k}|
						+h_{\max}^{-\sigma}\max _{\substack{u\in E(\lambda)\\\Vert u\Vert_{L^2(\Omega)}=1}}
						\min_{\substack{\boldsymbol{u_h}=(u_{\mathrm{pw}},u_{M})\in E_h\\
						\Vert u_{\mathrm{pw}}\Vert_{L^2(\Omega)}=1}}\Vert u-u_{\mathrm{pw}}\Vert_{L^2(\Omega)}&\\
						+h_{\max}^{-\sigma}
						\max_{\substack{\boldsymbol{u_h}=(u_{\mathrm{pw}},u_{M})\in E_h\\
						\Vert u_{\mathrm{pw}}\Vert_{L^2(\Omega)}=1}}
						\min _{\substack{u\in E(\lambda)\\\Vert u\Vert_{L^2(\Omega)}=1}}						
						\Vert u-u_{\mathrm{pw}}\Vert_{L^2(\Omega)}
						&\le C h_{\max}^{t}. 
				\end{align*}
			\end{theorem}
			The results of \cref{thm:GLB} and \ref{thm:BabuskaOsborn} assume exact solve of the algebraic eigenvalue problem \eqref{eq:dis_EVP_alt}, 
			but standard perturbation results in numerical linear algebra \cite{Par98} can be added to obtain rigorous bounds in practical applications. 		
							 		
		\subsection{Outline}
			\cref{sec:eigensolver} analyses the discrete eigenvalue problem \eqref{eq:dis_EVP_alt} and proves \cref{thm:GLB} 
			in any space dimension $n\ge 2$. 
			\cref{sec:Morley} recalls the Morley finite element (FE) and 
			presents interpolation error estimates in \cref{thm:properties_IMorley}. 
			The  interpolation constant $\kappa_2$ in 
			\eqref{eq:def_kappa} 
			leads to the guaranteed lower bound property from \cref{thm:GLB} in \cref{sec:GLB} 
			and a generalization of \cite{CGal14} in \cref{thm:GLB_NC}. 
		    \cref{sec:comments} introduces a reduced formulation for the new method, remarks on the relation to the standard 
			Morley eigenvalue problem, and introduces a related extra-stabilised Crouzeix-Raviart method. 
			The a priori convergence analysis in  $3$D of 
			\cref{sec:convergence} is based on a conforming companion operator, i.e., a right-inverse of the interpolation operator in $M(\mathcal{T})$ 
			with the extra properties in \cref{thm:J_Morley}. 
			This operator relies on the conforming Hsieh-Clough-Tocher finite element in $3$D suggested by Worsey-Farin (WF) 
			in \cite{WF87}  
			and allows for $L^2$ error estimates of separate interest. 
			The analysis of the conforming companion contains some technical details like the correct scaling of the 
			$\textit{WF}$ 
			basis functions, 
			which extends \cite[\S~6.1, p.340ff]{Ciarlet78} 
			to $n=3$ and is explained in the self-contained supplement to this paper.  
			The preparations in Subsections~\ref{sec:companion}--\ref{sec:convergence_source} 
			allow the proof of  \cref{thm:BabuskaOsborn} in \cref{sec:convergence_ev}.
			Since the method is new in any space dimension, the $2$D numerical experiments in \cref{sec:experiments} 			
			confirm the theoretical results, present details on \cref{fig:adaptive_Dumbbell_intro}, and provide  
			striking numerical evidence for the superiority of adaptive mesh-refinement for the bi-Laplace Dirichlet eigenvalue problem.
			
		\subsection{Notation}
			Standard notation on Lebesgue and Sobolev spaces applies throughout this paper; $(\bullet,\bullet)_{L^2(\Omega)}$ abbreviates
			the $L^2$ scalar product and $H^{2}(T)$ 
			abbreviates $H^{2}(\textup{int}(T ))$ for a compact set $T$ with non-void interior $\textup{int}(T )$.
			The vector space $H^{2}(\mathcal{T}):=\{v\in L^2(\Omega):\, v|_T\in H^{2}(T)\}$ consists of  
			piecewise $H^{2}$ functions and is equipped with the 
			semi-norm 
			$\vvvert\bullet\vvvert_{\mathrm{pw}}^2
			:=(D^2_{\mathrm{pw}}\bullet,D^2_{\mathrm{pw}}\bullet)_{L^2(\Omega)}$. The piecewise Hessian $D^2_{\mathrm{pw}}$ 
			is understood with respect to 
			the non-displayed regular triangulation $\mathcal{T}$ of $\Omega\subset\mathbb{R}^n$ into simplices. 	
			The context-depending notation $|\bullet|$ denotes the euclidean length of a vector, the cardinality of 
			a finite set, as well as the non-trivial $n$-,$(n-1)$-, or $(n-2)$- dimensional Lebesgue measure of a subset of  $\mathbb{R}^n$. 
			Let $P_2(T)$ denote the space of quadratic functions on $T\in\mathcal{T}$ and 
			$P_2(\mathcal{T}):=\{v\in L^2(\Omega):\, v|_T\in P_2(T)\text{ for all }T\in\mathcal{T}\}$ the space of piecewise quadratic functions. 
			Given a function $v\in L^2(\omega)$, define the integral mean 
			$\intmean_\omega v\,\textup{d}x:= 1/|\omega|\,\int_\omega v\,\textup{d}x$. 
			The $L^2$ projection $\Pi_0$ onto the piecewise constant functions $ P_0(\mathcal{T})$ reads 
			$\Pi_0(f)|_T:=\intmean_T f\,\textup{d}x$ for all 
			$f\in L^2(\Omega)$ and  $T\in\mathcal{T}$.  
			For any $A\in P_0(\mathcal{T};\mathbb{R}^{\ell\times \ell})$ SPD,  
			$(\bullet, \bullet)_A:=(A\bullet,\bullet)_{L^2(\Omega)}$  abbreviates  the weighted $L^2$ 
			scalar product with induced $A$-weighted $L^2$ norm 
			$\Vert \bullet \Vert_A:=\Vert A^{1/2} \bullet\Vert_{L^2(\Omega)}$. 
			Let $\sigma:=\min\{1,\sigma_{\mathrm{reg}}\}$ denote the minimum of one for the approximation property and 
			the positive index of elliptic regularity $\sigma_{\mathrm{reg}}>0$ for the source problem of 
			the bi-Laplacian $\Delta^2$ in $H^2_0(\Omega)$ on  the bounded polyhedral Lipschitz domain 
			$\Omega\subset \mathbb{R}^n$: Given any right-hand side $f\in L^2(\Omega)$, the weak
			solution $u\in V$ to $\Delta^2 u=f $ satisfies 
			\begin{align}
					u\in H^{{2}+\sigma}(\Omega)
					\text{ and } \Vert u\Vert_{H^{2+\sigma}(\Omega)}\le C(\sigma)\Vert f\Vert_{L^2(\Omega)}. 
					\label{eq:def_sigma}
			\end{align}
			The Sobolev space $H^{2+s}(\Omega)$ is defined for $0<s<1$ by complex interpolation of $H^2(\Omega)$ and $H^3(\Omega)$. 
			Notice $E(\lambda)\subset H^{2+\sigma}(\Omega)$ in \cref{thm:BabuskaOsborn} follows from \eqref{eq:def_sigma} but 
			$E(\lambda)\subset H^{2+t}(\Omega)$ is possible for $t\ge \sigma$ for some eigenvalues $\lambda$. 
			Throughout this paper, $a \lesssim b$ abbreviates $a\le Cb$ with a generic constant $C$ only 
			dependent on $\sigma$ in \eqref{eq:def_sigma} and the shape-regularity constant $C_{\mathrm{sr}}$ 
			of $\mathcal{T}\in\mathbb{T}$; 
			$a \approx b$ stands for $a\lesssim  b\lesssim   a$.  
					
\section{Eigensolver for guaranteed lower bounds in any dimension}\label{sec:eigensolver}
\subsection{The Morley finite element}\label{sec:Morley}
			Given a shape-regular triangulation $\mathcal{T}$ of a bounded polyhedral Lipschitz domain 
			$\Omega\subset\mathbb{R}^n$ into $n$-simplices (tetrahedra in $3$D) in the sense of Ciarlet \cite{BS08,Bra13,BBF13}, let
			$\mathcal{V}$ (resp. $\mathcal{V}(\Omega)$ or $\mathcal{V}(\partial\Omega)$) denote the set 
			of all (resp. interior or boundary)  vertices, let $\mathcal{F}$ (resp. $\mathcal{F}(\Omega)$ or $\mathcal{F}(\partial\Omega)$) denote the set 
			of all (resp. interior or boundary) $(n-1)$-subsimplices (faces in $3$D), and 
			let $\mathcal{E}$ 
			(resp. $\mathcal{E}(\Omega)$ or $\mathcal{E}(\partial\Omega)$) denote the set of all (resp. interior or boundary) 
			$(n-2)$-subsimplices (edges in $3$D and vertices in $2$D) in $\mathcal{T}$.
		The degrees of freedom for the Morley element 
		\cite[Def.~1]{MX06} on an  $n$-simplex $T\in\mathcal{T}$ are 
		the integral means of the function $f$ along any $(n-2)$-subsimplex $E\in\mathcal{E}(T)$ of $T$ and
		of the normal derivative ${\partial f}/{\partial\nu}$ for each $(n-1)$-subsimplex $F\in\mathcal{F}(T)$ of $T$. 
		Let the integral mean over a node $z\in\mathcal{V}$ be the point evaluation, to see that  
		this reduces to the classical definition \cite{Mor68} for $n=2$. 
		The $m:=|\mathcal{F}|+|\mathcal{E}|$ global degrees of freedom are labelled, for any $f\in H^2(\Omega)$, by 
		\begin{align*}
			L_E(f):= \intmean_{E} f\,\textup{d}s \text{ for  any }E\in \mathcal{E} 
			\quad\text{and}\quad  L_F(f):=\intmean_{F} \nabla f\cdot\nu_F\,\textup{d}\sigma \text{ for  any } F\in \mathcal{F},
		\end{align*}
		where $\nu_F$  denotes the unit normal for any side $F\in\mathcal{F}$ with a fixed orientation. 				
		Section~6 in \cite{CGH14} introduces the dual basis for  the Morley finite element in $2$D and specifies 
		an implementation in $30$ lines of matlab. 
		The nodal basis in \cite[Thm.~1]{MX06} reads for $n\ge 3$ as follows.

		Given any side $F\in\mathcal{F}$ with unit normal $\nu_F$ of fixed orientation, 
		the support $\textup{supp}(\phi_F)=\overline{\omega(F)}$ 
		of the basis function $\Phi_F\in P_2(\mathcal{T})$ dual to $L_F$ consists of all					
		adjacent $n$-simplices. For any interior side 
		$F=\partial T_+\cap \partial T_-\in\mathcal{F}(\Omega)$ the side-patch $\overline{\omega(F)}:=T_+\cup T_-$ consists of 
		the neighbouring simplices ${T}_{\pm}\in\mathcal{T}$ 
		with $\nu_F=\nu_{T_+}|_F=-\nu_{T_-}|_F$.  
		For any boundary side $F=\partial T_+\cap \partial \Omega\in \mathcal{F}(\partial\Omega)$ set $\overline{\omega(F)}:=T_+$ with  $\nu_F=\nu_{T_+}|_F$.
		Suppose that $F=F_j$ is the side opposite the vertex $P_j$ with barycentric coordinate $\lambda_j$ in 
		$T_\pm\subset \overline{\omega(F)}$,  then 
			\begin{align}
				\phi_{F}|_{T_\pm}:=
						\pm \big(\lambda_{j}	
						(n\lambda_{j}-2)\big)/\big(2\vert\nabla\lambda_{j}\vert\big).
						\label{eq:Morley_basis_F}
			\end{align}
		Let $E_{jk}:=\textup{conv}\{P_1,\dots,P_{j-1},P_{j+1},\dots,P_{k-1}, P_{k+1},\dots, P_{n+1}\}\in\mathcal{E}(T)$ denote the 			
		$(n-2)$-subsimplex of $T:=\textup{conv}\{P_1,\dots,P_{n+1}\}\in\mathcal{T}$ 
		in the intersection $E\in \partial F_j\cap\partial F_k$ of the sides $F_j,F_k\in\mathcal{F}(T)$.
		Given a $(n-2)$-subsimplex $E\in\mathcal{E}$, the support  
		$\textup{supp}(\phi_E)=\bigcup \mathcal{T}(E)$ of the basis function $\Phi_E\in P_2(\mathcal{T})$ dual to $L_E$ 
		consists of all adjacent $n$-simplices $\mathcal{T}(E):=\{T\in\mathcal{T}:\, E\in\mathcal{E}(T)\}$.  Suppose  that 
		$E=E_{jk}$ is the $(n-2)$-subsimplex in the intersection $E\in \partial F_j\cap\partial F_k$ of the sides 
		$F_j,F_k\in\mathcal{F}(T)$ in $T\in\mathcal{T}(E)$ with barycentric coordinates $\lambda_j$ and $\lambda_k$ in $T$ 
		(associated with the opposite vertices $P_j$ and $P_k$). Then
			\begin{align}
				\phi_E|_T= 1-(n-1)(\lambda_j+\lambda_k)+n(n-1)\lambda_j\lambda_k-(n-1)\nabla \lambda_j
							\cdot\nabla\lambda_k
							\sum_{\ell\in \{j,k\}}\frac{\lambda_\ell(n\lambda_\ell-2)}{2\vert\nabla\lambda_\ell\vert^2}.
							\label{eq:Morley_basis_E}
			\end{align}
		This defines the nodal basis functions: 
		$L_G(\phi_F|_T)=\delta_{FG}$, $L_D(\phi_E|_T)=\delta_{DE}$,
		and $L_F(\phi_E|_T)=0=L_E(\phi_F|_T)$ follows for any $n$-simplex $T\in\mathcal{T}$, 
		any $(n-1)$-simplices $F,G\in\mathcal{F}(T)$, and any $(n-2)$-simplices $D,E\in \mathcal{E}(T)$.  
		The Morley finite element space with  homogeneous boundary conditions reads 
			\begin{align*}
				M(\mathcal{T})&:=\textup{span}\{\phi_F :\,F\in\mathcal{F}(\Omega)\}\oplus 
														\textup{span}\{\phi_E:\,E\in\mathcal{E}(\Omega)\}\subset P_2(\mathcal{T}).
			\end{align*}
				Given the dual basis for the Morley FEM in 
				\eqref{eq:Morley_basis_F}--\eqref{eq:Morley_basis_E}, define the 
				interpolation operator
				$I_{M}:\, V\to{M}(\mathcal{T})$ for any $v\in V:= H^2_0(\Omega)$ by 
					\begin{align}
						I_M(v):=\sum_{F\in\mathcal{F}(\Omega)}\intmean_F\nabla v\cdot\nu_F\,\textup{d}\sigma\, \phi_F
								+\sum_{E\in\mathcal{E}(\Omega)}\intmean_E v\,\textup{d}s\,\phi_E. \label{eq:def_IM}
					\end{align}
				This interpolation operator has the following important properties with the explicit constants $\kappa_1$ and $\kappa_2$ 
				from \eqref{eq:def_kappa}, which are not quantified in \cite{MX06}.
				\begin{theorem}[properties of $I_M$]\label{thm:properties_IMorley}
					\begin{mylist}
						\item\label{item:Pi0Morley} 
							Any $v\in V$ satisfies $\Pi_0 D^2v=D_{\mathrm{pw}}^2 I_Mv$, in particular 
							$a_{\mathrm{pw}}(v-I_M v,w_{M})=0${ for all } $w_{M}\in M(\mathcal{T})${ and }$v\in V$; 
							$I_M$ is the $a_{\mathrm{pw}}$-orthogonal projection onto  $M(\mathcal{T})$ with 
				 			$\vvvert v-I_Mv\vvvert_{\mathrm{pw}}
				 			=\min_{v_{M}\in M(\mathcal{T})}\vvvert v-v_{M}\vvvert_{\mathrm{pw}}$
				 			for any $v\in V$.					
						\item\label{item:IMorley_kappa}
						 Any $v\in H^2(T)$ in $T\in\mathcal{T}$ satisfies 
								$
									\vert v-I_Mv\vert_{H^{2-\ell}(T)}\le \kappa_\ell h_T^\ell \vert v-I_M v\vert_{H^2(T)}
								$ for $\ell=1,2$. 
					\end{mylist}
				\end{theorem}
				\begin{proof2}
				This is known for $n=2$ from \cite{CGal14}, so let $n\ge 3$ in the sequel. 
				\begin{proofof}\textit{(\ref{item:Pi0Morley}).}
				The definition of  $I_M$ implies $\intmean_F {\nabla I_Mv\cdot\nu_F}\,\textup{d}\sigma
					= \intmean_F{\nabla v\cdot \nu_F}\,\textup{d}\sigma$
					for any $F\in\mathcal{F}$. 
					This and an integration by parts 
					prove $\Pi_0 D^2v=D_{\mathrm{pw}}^2 I_Mv$.
				  Since $w_M\in M(\mathcal{T})\subset P_2(\mathcal{T})$, 
				 this	concludes the proof of  
							$a_{\mathrm{pw}}( v-I_Mv, w_M)=((1-\Pi_0)D^2v, D^2_{\mathrm{pw}}w_M)_{L^2(\Omega)}=0$.
				\end{proofof} 
				\begin{proofof}\textit{(\ref{item:IMorley_kappa}) for $\ell=1$.} 
					A key observation is that the piecewise gradient $\nabla_{\mathrm{pw}} v_M\in \textit{CR}^1_0(\mathcal{T})^n$ 
					of any $v_M\in M(\mathcal{T})$ is a 
					Crouzeix-Raviart function in $n$ components.
					(This follows from $\intmean_F\nabla I_Mv\,\textup{d}\sigma=\intmean_F\nabla v\,\textup{d}\sigma$, 
					since for any $v\in C^1(T)$ in $T\in\mathcal{T}$ the Morley degrees of freedom uniquely determine  
					$\intmean_F\nabla v\,\textup{d}\sigma$ for any $F\in\mathcal{F}(T)$ \cite[Lem.~1]{MX06}.)
					Moreover, the Crouzeix-Raviart interpolation operator $I_{CR}$ \cite{CR73} 
					(applied component-wise) satisfies 
					$\nabla_{\mathrm{pw}} I_M v=I_{CR}\nabla v$ for any $v\in H^2(T)$ in $T\in\mathcal{T}$. 
					Lemma A.1 in \cite{CZZ18} shows 
					$\Vert f-I_{CR}f\Vert_{L^2(T)}\le \kappa_1 h_T \vert f-I_{CR} f\vert_{H^1(T)}$ for $f\in H^1(T)$ and any $n\ge2$. 
					The choice $f=\partial v/\partial x_j$ for $j=1,\dots,n$ concludes the proof of 
					$\vert v-I_Mv\vert_{H^1(T)}\le \kappa_1 h_T \vert v-I_M v\vert_{H^2(T)}$ for $v\in H^2(T)$.
				\end{proofof}
				\begin{proofof}\textit{(\ref{item:IMorley_kappa}) for $\ell=2$.}
			Let $g:= v-I_Mv \in H^2(T)$ and set 
				$
					I_{CR}(g):=\sum_{F\in\mathcal{F}}\big(\intmean_F g\,\textup{d}\sigma\big) \psi_F 
				$			 
			with the side-oriented Crouzeix-Raviart basis function $\psi_F\in \textit{CR}^1(\mathcal{T})$ with 
			$\psi_F(\textup{mid}(G))=\delta_{FG}$ for all $F,G\in\mathcal{F}$. 
			The local mass matrix $M(T)\in \mathbb{R}^{(n+1)\times (n+1)}$ for $\mathcal{F}(T):=\{F_1,\dots,F_{n+1}\}$ 
			reads
				\begin{align*}
					M(T):=\Big(\int_T\psi_{F_j}\psi_{F_k}\,\textup{d}\sigma\Big)_{j,k=1,\dots,n+1}
					=\Bigg( \frac{|T|(2-n+n^2\delta_{jk})}{(n+1)(n+2)}\Bigg)_{j,k=1,\dots,n+1}.
				\end{align*}
		 	The eigenvalue $|T|/(n+1)$ of $M(T)$ has the eigenvector $(1,\dots,1)\in \mathbb{R}^{n+1}$. 
		 	The eigenvalue $|T|n^2/((n+1)(n+2))$ has the
	    		$n$-dimensional eigenspace of vectors in $\mathbb{R}^{n+1}$ perpendicular to $(1,\dots,1)$.
			Hence the coefficient vector 
			$x:=(\intmean_F g\,\textup{d}\sigma:\, F \in \mathcal{F}(T))\in\mathbb{R}^{n+1}$ of $I_{CR}g$ satisfies 
				\begin{align}
					\Vert I_{CR}g\Vert_{L^2(T)}^2 
												&=\int_T \Big(\sum_{F\in\mathcal{F}(T)}\intmean_F g\,\textup{d}\sigma\,\psi_F\Big)^2
														\,\textup{d}x
												  =x\cdot M(T)x\notag
												\\&   \le \frac{|T|n^2}{(n+1)(n+2)} |x|^2
												  = \frac{|T|n^2}{(n+1)(n+2)} \sum_{F\in\mathcal{F}(T)}
												  					\Big(\intmean_F g\,\textup{d}\sigma\Big)^2.\label{eq:proof_kappa2_ICR}
				\end{align}
			If the $m$-simplex $F$ has the $(m-1)$-subsimplex $E$ opposite to the vertex $P$ in $F$ for $m\ge 2$, then any $v\in H^1(F)$ satisfies 
			the trace identity  
			\begin{align}
				\intmean_{E} v \,\textup{d}s=\intmean_{F} v\,\textup{d}\sigma +{m^{-1}}\intmean_{F} (x-P)\cdot\nabla v\,\textup{d}\sigma. 
				\label{eq:traceIdentity} 
			\end{align}		
			(This follows from an integration by parts \cite{CGR12,CH17}.)  
			In each $(n-1)$-subsimplex $F\in\mathcal{F}(T)$ 
			with midpoint $\textup{mid}(F)$, the set of $(n-2)$-subsimplices (edges in $3$D) 
			$\mathcal{E}(F):=\{E\in\mathcal{E}(T):E\subset \partial F\}=\{E_1,\dots,E_n\}$
			defines the sub-triangulation of $F$ into $F_j:=\textup{conv}(\textup{mid}(F),E_j)\subset F$ for $j=1,\dots,n$.  
			Since the function $g|_{F}:=(v-I_Mv)|_{F}\in H^1(F)$ satisfies 
			$\intmean_{E_j} g\,\textup{d}s=0$ for any  $j=1,\dots,n$,		
			the trace identity \eqref{eq:traceIdentity} 
			for each $(n-1)$-simplex $F_j\subset F$ proves
				\begin{align*}
					\intmean_{F}g\,\textup{d}\sigma
							&=-\frac{1}{|F|(n-1)}\int_{F}(x-\textup{mid}(F))\cdot\nabla g\,\textup{d}\sigma
							\le \frac{1}{|F|(n-1)} \Vert \bullet -\textup{mid}(F)\Vert_{L^2(F)}\Vert \nabla g\Vert_{L^2(F)}
				\end{align*} with the Cauchy-Schwarz inequality in the last step. 
			For the $(n-1)$-dimensional simplex $F=\textup{conv}\{P_1,\dots, P_n\} $ assume without loss of generality 
			that  
			$\textup{mid}(F)=\frac{1}{n}\sum_{j=1}^nP_j=0$. 
			An estimation of the  mass-matrix  for the Courant basis 
			functions 
			associated 	with the vertices $P_1,\dots, P_n$ of $F$  leads to 
				\begin{align*}
					\Vert \bullet -\textup{mid}(F)\Vert_{L^2(F)}^2
					=\int_F|x|^2\,\textup{d}\sigma=\frac{|F|}{n(n+1)}\sum_{\ell=1}^n|P_\ell|^2.
				\end{align*}
			 Similar to \cite[Lem.~A.1]{CZZ18}, elementary algebra 
			 with $\frac{1}{n}\sum_{j=1}^nP_j=0$ and $|P_j-P_k|\le h_F$ for all $j,k=1,\dots, n$ with 
			$j\not =k$ lead to 
			\begin{align*}
			\sum_{\ell=1}^n|P_\ell|^2
				={1}/{(2n)}\sum_{j,k=1}^n|P_j-P_k|^2\le {h_F^2(n-1)}/{2}.
			\end{align*}
			The combination of the last three displayed estimates reads
			\begin{align*}
				\Bigg(\intmean_{F}g\,\textup{d}\sigma\Bigg)^2
						&\le \frac{1}{(n-1)^2|F|^2} \frac{h_F^2{|F|(n-1)}}{{2n(n+1)}}\Vert \nabla g\Vert_{L^2(F)}^2
						\le\frac{h_T^2}{|F|} \big({2n(n-1)(n+1)}\big)^{-1}\Vert \nabla g\Vert_{L^2(F)}^2.
			\end{align*}
			The trace inequality 
			$\Vert v\Vert_{L^2(F)}^2\le ({|F|}/{|T|})\Vert v\Vert_{L^2(T)}(\Vert v\Vert_{L^2(T)} 
					+{2h_T}/{n}\Vert \nabla v\Vert_{L^2(T)})$ for $v\in H^1(T)$ and $F\in \mathcal{F}(T)$
			 is a direct consequence of the trace identity \eqref{eq:traceIdentity} for the $n$-simplex $T\in\mathcal{T}$ with $(n-1)$-subsimplex 
			 $F\in\mathcal{F}(T)$.
			This and Young's inequality show
			\begin{align*}
				\frac{|T|}{|F|}\Vert \nabla g\Vert_{L^2(F)}^2
						\le 
							(1+(\kappa_1n)^{-1})\Vert \nabla g \Vert_{L^2(T)}^2 
											+h_T^2 \kappa_1 n^{-1}\Vert D^2 g\Vert_{L^2(T)}^2.
			\end{align*}
			The proven \cref{thm:properties_IMorley}.\ref{item:IMorley_kappa} for $\ell=1$ shows 
			$
				\Vert \nabla g \Vert_{L^2(T)}\le \kappa_{1}h_T \Vert D^2 g \Vert_{L^2(T)}.
			$
			In combination with the last two displayed inequalities and \eqref{eq:proof_kappa2_ICR}, this reads 
			\begin{align}
				\Vert I_{CR}g\Vert_{L^2(T)}^2  &
				\le  \frac{n\kappa_{1}^2+2\kappa_1}{2(n-1)(n+1)(n+2)}h_T^4\Vert D^2 g\Vert_{L^2(T)}^2.  \label{eq:INCg_norm}
			\end{align}
			On the other hand, the Crouzeix-Raviart interpolation operator $I_{CR}$ satisfies
			$\Vert (1-I_{CR})g\Vert_{L^2(T)}$
			$\le \kappa_{1}h_T\Vert \nabla (1-I_{CR})g\Vert _{L^2(T)}$ and 
			$\nabla I_{CR}g=\Pi_0 \nabla g$ 
			as in \cite{CGal14,CGed14} for $n=2$ and in \cite{CZZ18} for $n\ge 3$.
			Hence, the Poincar\'e inequality with Payne-Weinberger constant \cite{PW60,Bebendorf2003} shows 
			\begin{align*}
				\Vert (1-I_{CR})g\Vert_{L^2(T)}
				&\le \kappa_{1}h_T\Vert \nabla (1-I_{CR})g\Vert _{L^2(T)}
				= \kappa_{1}h_T\Vert (1-\Pi_0)\nabla g\Vert _{L^2(T)}
				\le\frac{\kappa_{1}}{\pi}h_T^2\Vert D^2 g\Vert _{L^2(T)}.
			\end{align*}
			The combination of this with \eqref{eq:INCg_norm} and a triangle inequality concludes the proof of 
			$\Vert v-I_Mv\Vert_{L^2(T)}\le \kappa_2 h_T^2 \vert v-I_M v\vert_{H^2(T)}$ for any $v\in H^2(T)$. 
		\end{proofof}
			The constant $\kappa_2=0.25746$ 
			for $n=2$ from \cite[Thm.~3]{CGal14} is recovered, 
			if the Poincar\'e constant $1/\pi$ is replaced by the optimal $1/j_{1,1}$ in $2$D \cite{LS10}. 
			The computational bounds  $\kappa_1\le 0.1893$ from \cite{Liu2015} and $\kappa_2\le 0.07353 $ \cite{LiaoShuLiu2019} improve the analytical bounds from \cite{CGal14} in $2$D.    
			\end{proof2}
			\begin{corollary}[further properties]\label{lem:InterpolationOperator} 
					\begin{enumerate}[wide, label=(\alph*), ref=\alph* ]	
						\item\label{item:cor_I_alpha} 
							Any $v\in H^{2+s}(\Omega)$ with $0\le s\le 1$ satisfies 
							\begin{align*}
								\vvvert (1-I_M)v\vvvert_{\mathrm{pw}}
								\le (h_{\max}/\pi)^{s} \Vert v\Vert_{H^{2+s}(\Omega)}. 
							\end{align*}
						\item\label{item:cor_I_a} Any $v,\, w\in V$ and $v_{M}\in M(\mathcal{T})$ satisfy
									$
							 			a_{\mathrm{pw}}(v, v_{M})=a_{\mathrm{pw}}(I_Mv,v_{M})		
									$
								and
								\begin{align*}
							 		a_{\mathrm{pw}}(v, (1-I_M)w)&=a_{\mathrm{pw}}((1-I_M)v, (1-I_M)w)
							 			\\&\le \min_{v_{M}\in M(\mathcal{T})} \vvvert v-v_{M}\vvvert_{\mathrm{pw}}
							 			\min_{w_{M}\in M(\mathcal{T})} \vvvert w-w_{M}\vvvert_{\mathrm{pw}}	. 		 
								\end{align*}
					\end{enumerate}
			\end{corollary}
			\begin{proof}
					\cref{thm:properties_IMorley}.\ref{item:Pi0Morley}  and a piecewise Poincar\'e inequality (as above from \cite{PW60,Bebendorf2003}) 
					 show for $s=0$ and $s=1$ that 
					\begin{align}
						\vvvert (1-I_M)v\vvvert_{\mathrm{pw}}=\Vert D^{2} v-\Pi_0D^{2}v\Vert_{L^2(\Omega)}
									\le (h_{\max}/\pi)^{s}\Vert v\Vert_{H^{{2}+s}(\Omega)}. \label{eq:proof_Ialpha}
					\end{align}
					Since $H^{2+s}(\Omega)$ is defined by complex interpolation of $H^2(\Omega)$ and $H^3(\Omega)$, 
					each component of the Hessian $D^2v$ belongs to $H^s(\Omega)\in [L^2(\Omega),H^1(\Omega)]_s$ 
					in the complex interpolation space  between  $L^2(\Omega)$ and $H^1(\Omega)$ \cite{Tatar2007}. The interpolation of 
					\eqref{eq:proof_Ialpha} concludes the proof of (\ref{item:cor_I_alpha}). 
					\cref{thm:properties_IMorley}.\ref{item:Pi0Morley} implies the first claim in (\ref{item:cor_I_a}). 
					The combination with the Cauchy-Schwarz inequality implies the second. 
			\end{proof}	
			
\subsection{Guaranteed lower bounds}\label{sec:GLB}
		This section proves  that the discrete method \eqref{eq:dis_EVP_alt} indeed provides GLBs for the continuous 
		eigenvalues in \emph{any} space dimension $n\ge 2$. 
		\begin{proofof}\textit{\cref{thm:GLB}.\,}
			Abbreviate $\lambda=\lambda_{k}$ from \eqref{eq:contEVP} and $\lambda_h=\lambda_h({k})$ from 
			\eqref{eq:dis_EVP_alt}. 
			Let $\phi_1,\dots,\phi_{k}$ denote the first ${k}$ $b$-orthonormal eigenfunctions of \eqref{eq:contEVP}; 
			the min-max principle \cite{StrangFix2008, Boffi2010} guarantees 
			$\vvvert \phi\vvvert^2\le \lambda$ for any $\phi\in \textup{span}\{\phi_1,\dots,\phi_{k}\}$ with $\Vert\phi\Vert_{L^2(\Omega)}=1$. 
			Let $\Pi_2$ denote $L^2$ projection onto $P_2(\mathcal{T})$. 
			\begin{enumerate}[wide, label=\textit{Case \arabic*.}]
			\item
				Assume the $L^2$-projections $\Pi_2\phi_1,\dots, \Pi_2\phi_{k}$ are linear dependent. 
				Then there exists some 			
				$\phi \in \textup{span}\{\phi_1,\dots,\phi_{k}\}$ with $\Vert \phi\Vert_{L^2(\Omega)}=1$ and $\Pi_2\phi=0$. 
				Let $\kappa_2^\prime$ denote the best possible  constant in 
				\begin{align}
					\Vert (1-\Pi_2)\psi\Vert_{L^2(\Omega)}
						\le \kappa_2^\prime {h_{\max}^2} \vvvert (1-I_{M})\psi\vvvert_{\mathrm{pw}}\quad\text{for all }\psi\in H^2_0(\Omega).
						\label{eq:kappa2prime}
				\end{align}
				The approximation property of $\Pi_2$ and $M(\mathcal{T})\subset P_2(\mathcal{T})$ 
				imply \eqref{eq:kappa2prime} with $\kappa_2^\prime\le \kappa_2$. 
				The above $\phi$ therefore satisfies 
				\begin{align*}
					1&=\Vert \phi\Vert_{L^2(\Omega)}
						=\Vert (1-\Pi_2)\phi\Vert_{L^2(\Omega)}
						\le \kappa_2^\prime {h_{\max}^2} \vvvert (1-I_{M})\phi\vvvert_{\mathrm{pw}}.
				\end{align*}
				The Pythagoras theorem from \cref{thm:properties_IMorley}.\ref{item:Pi0Morley} and 
				$\vvvert \phi\vvvert^2\le \lambda$ from the min-max principle \cite{StrangFix2008, Boffi2010} for \eqref{eq:contEVP} show
				\begin{align*}
					\vvvert (1-I_{M})\phi\vvvert_{\mathrm{pw}}^2+\vvvert I_{M}\phi\vvvert_{\mathrm{pw}}^2
						= \vvvert \phi\vvvert^2
						\le \lambda.
				\end{align*}
				The combination of the last two displayed inequalities reads $1\le \lambda (\kappa_2^\prime)^2 h_{\max}^4$. 
				Throughout this  paper, the values used for $\kappa_2$ satisfy $\kappa_2^\prime <\kappa_2$ 
				(see Remarks~\ref{rem:kappa2prime_1}--\ref{rem:kappa2prime_2} below), whence 
				$\lambda \kappa_2^2 h_{\max}^4>1$. In other words, the a priori condition  $\lambda \kappa_2^2 h_{\max}^4\le 1$ fails 
		 		and (the remaining hypothesis of \cref{thm:GLB}) $\lambda_h \kappa_2^2h_{\max}^4\le 1$ holds.
				This proves $\lambda_h \kappa_2^2h_{\max}^4\le 1<\lambda \kappa_2^2 h_{\max}^4$ and concludes the proof of  
		 		$\lambda_h\le \lambda$ in the first case. 
			\item 
				Assume that the projections $\Pi_2\phi_1,\dots, \Pi_2\phi_{k}$ are linear independent. 
				Set $\boldsymbol{S}_k:=\textup{span}\big\{\big(\Pi_2\phi_1,I_{M}\phi_1\big),\dots,\allowbreak
				\big(\Pi_2\phi_{k},I_{M}\phi_{k}\big)\}\allowbreak
				\subset\boldsymbol{V_h}$ with $\textup{dim}(\boldsymbol{S}_k)=k$. 
				Since $\boldsymbol{b_h}$ is positive definite on ${\boldsymbol{S}_k\times \boldsymbol{S}_k}$,  
				the min-max principle \cite{StrangFix2008, Boffi2010} 
				for \eqref{eq:dis_EVP_alt} shows
				\begin{align}
		 			\lambda_h\le \max_{\boldsymbol{v_h}\in \boldsymbol{S}_{k}\setminus\{0\}}
		 			\frac{\boldsymbol{a_h}(\boldsymbol{v_h},\boldsymbol{v_h})}
		 			{\boldsymbol{b_h}(\boldsymbol{v_h},\boldsymbol{v_h})}.\label{eq:min_max_dis}
		 		\end{align}
				Let $\boldsymbol{v_h}=(\Pi_2\phi,I_{M}\phi) \in \boldsymbol{S}_{k}\setminus\{0\}$  be some maximizer in the upper bound
				of \eqref{eq:min_max_dis} with 
				$\Vert\phi\Vert_{L^2(\Omega)}=1$ and deduce $\vvvert \phi\vvvert^2\le \lambda$ 
				from the min-max principle \cite{StrangFix2008, Boffi2010} for \eqref{eq:contEVP}. 
				The inequality \eqref{eq:min_max_dis} ensures for $\boldsymbol{v_h}=(\Pi_2\phi,I_{M}\phi) \in \boldsymbol{S}_{k}\setminus\{0\}$ that 
				\begin{align*}
					\lambda_h\Vert \Pi_2\phi\Vert_{L^2(\Omega)}^2
					=\lambda_h \boldsymbol{b_h}(\boldsymbol{v_h},\boldsymbol{v_h})
					\le \boldsymbol{a_h}(\boldsymbol{v_h},\boldsymbol{v_h})
					=\vvvert I_{M} \phi\vvvert_{\mathrm{pw}}^2	 			
		 			+\kappa_{2}^{-2}\Vert h_{\mathcal{T}}^{-{2}}(\Pi_2- I_{M})\phi\Vert_{L^2(\Omega)}^2.
				\end{align*}
				Since the  piecewise constant mesh-size  $h_{\mathcal{T}}$ does not interact with the piecewise $L^2$ 
				projections, the Pythagoras theorem  shows
				\begin{align*}
				 	\Vert h_{\mathcal{T}}^{-{2}}(\Pi_2- I_{M})\phi\Vert_{L^2(\Omega)}^2
				 	&=\Vert h_{\mathcal{T}}^{-{2}}(1- I_{M})\phi\Vert_{L^2(\Omega)}^2
				 			-\Vert h_{\mathcal{T}}^{-{2}}(1- \Pi_2)\phi\Vert_{L^2(\Omega)}^2
				 	\\& \le \kappa_{2}^2\vvvert (1-I_{M})\phi\vvvert_{\mathrm{pw}}^2 
				 				-h_{\max}^{-{4}}+h_{\max}^{-{4}}\Vert \Pi_2\phi\Vert_{L^2(\Omega)}^2
				\end{align*}
				with \cref{thm:properties_IMorley}.\ref{item:IMorley_kappa} for $\ell=2$ and 
				$1-\Vert \Pi_2\phi\Vert_{L^2(\Omega)}^2\le h_{\max}^{4}\Vert h_{\mathcal{T}}^{-{2}}(1-\Pi_2)\phi\Vert_{L^2(\Omega)}^2$  in the last step. 
				The combination of the previous two displayed estimates leads to 
				\begin{align}
					(\lambda_h - \kappa_2^{-2}h_{\max}^{-{4}})\Vert \Pi_2\phi\Vert_{L^2(\Omega)}^2
					+ \kappa_2^{-2}h_{\max}^{-{4}}
					\le 
					\vvvert I_{M} \phi\vvvert_{\mathrm{pw}}^2+\vvvert (1-I_{M})\phi\vvvert_{\mathrm{pw}}^2
					=\vvvert\phi\vvvert^2\le \lambda	\label{eq:proof_GLB}
				\end{align}						
				with the Pythagoras theorem from \cref{thm:properties_IMorley}.\ref{item:Pi0Morley} in the equality and 
				$\vvvert \phi\vvvert^2\le \lambda$ from above in the last step. 
				The inequality \eqref{eq:proof_GLB} implies 
				\begin{align}
					1-\lambda\kappa_2^{2}h_{\max}^{{4}}	\le (1-\lambda_h\kappa_2^{2}h_{\max}^{{4}})\Vert \Pi_2\phi\Vert_{L^2(\Omega)}^2.		
					\label{eq:proof_GLB2}
				\end{align}
				Without loss of generality assume $\lambda\kappa_2^{2}h_{\max}^{{4}}\le 1$, 
				otherwise the remaining hypothesis $\lambda_h\kappa_2^{2}h_{\max}^{{4}}\allowbreak\le 1<\lambda\kappa_2^{2}h_{\max}^{{4}}$ 
				proves the claim. 
				Since the linear independence of $\Pi_2\phi_1,\dots, \Pi_2\phi_{k}$ shows  $\Vert\Pi_2\phi\Vert_{L^2(\Omega)}>0$, 
				\eqref{eq:proof_GLB2} implies $\lambda_h\kappa_2^{2}h_{\max}^{{4}}\le 1$.  
				Hence $\Vert\Pi_2\phi\Vert_{L^2(\Omega)}\le\Vert\phi\Vert_{L^2(\Omega)}=1$ and \eqref{eq:proof_GLB2} imply 
				$1-\lambda\kappa_2^{2}h_{\max}^{{4}}	\le 1-\lambda_h\kappa_2^{2}h_{\max}^{{4}}$. This proves $\lambda_h\le \lambda$ 
				in the second case.
		 	\end{enumerate}	
		\end{proofof}
		
			\begin{remark}[choice of $\kappa_2=0.07353$ in $2$D]\label{rem:kappa2prime_1}
		 	The computational bound  $\kappa_2= 0.07353>\kappa_2^\star$  in $2$D 	from \cite{LiaoShuLiu2019}
			is a guaranteed upper bound for some 
			optimal $\kappa_2^\star$ with  $\Vert v-I_Mv\Vert_{L^2(T)}\le \kappa_2^\star h_T^2 \vert v-I_M v\vert_{H^2(T)}$ 
			for all $v\in H^2(T)$, $T\in\mathcal{T}$.	
			Hence $\kappa_2^\prime:= \kappa_2^\star<\kappa_2$ and \cref{thm:GLB} holds for this improved choice of $\kappa_2$.
			\end{remark}
			\begin{remark}[choice of $\kappa_2$ in  \eqref{eq:def_kappa}]\label{rem:kappa2prime_2}
			Standard arguments including successive (piecewise) Poin-car\'e inequalities \cite{PW60,Bebendorf2003}  eventually 
			imply $\kappa_2^\prime\le 1/\pi^2$ 
			and the analytical bound in \eqref{eq:def_kappa} satisfies $1/\pi^2<\kappa_2$, hence 
			the claim $\kappa_2^\prime<\kappa_2$ follows for $n\ge 3$ as well. 
			\end{remark}
			The standard Morley eigenvalue problem seeks  
		$(\lambda_{M}, \phi_{M})\in \mathbb{R}^+\times \big(M({\mathcal{T}})\setminus\{0\}\big)$ with 
					\begin{align}
							a_{\mathrm{pw}}(\phi_{M},v_{M})
							=\lambda_{M} (\phi_{M},v_{M})_{L^2(\Omega)} 
														\quad\text{for all } v_{M}\in M({\mathcal{T}}).  \label{eq:disEVP_NC}
					\end{align}
				Assume the $N:=\textup{dim}(M(\mathcal{T}))$ algebraic eigenvalues 
				$0<\lambda_{M} (1)\le\lambda_{M}(2)\le \dots\le \lambda_{M} (N)<\infty$ of \eqref{eq:disEVP_NC} 
				are enumerated in ascending order counting 		
				multiplicities. 
		The following \cref{thm:GLB_NC} refines  \cite{CGal14} for $n=2$ by removing an unnecessary separation condition 
		 \cite[Rem.~2.2]{Liu2015} and generalizes it to $n\ge 3$.
		\begin{theorem}[GLB in \eqref{eq:GLB_CR_Morley}]\label{thm:GLB_NC}
			For any $k=1,\dots, N$ the ${k}$-th eigenvalue $\lambda_{M}({k})$ of \eqref{eq:disEVP_NC} 
			leads in \eqref{eq:GLB_CR_Morley} 
			to 	a guaranteed lower bound $\textup{GLB}({k})\le \lambda_k$ for  the ${k}$-th eigenvalue $\lambda_{k}$ of 
			\eqref{eq:contEVP}. 
		\end{theorem}
		\begin{proof}
		Let $\phi_1,\dots,\phi_k$ denote the first $k$ $L^2$-orthonormal eigenfunctions in \eqref{eq:contEVP}.				
		Assume first, that  the interpolations $I_{M}\phi_1,\dots, I_{M}\phi_{k}$ are linear dependent. Let  		 	
				 $\phi \in \textup{span}\{\phi_1,\dots,\phi_{k}\}$ with $\Vert \phi\Vert_{L^2(\Omega)}=1$ and $I_{M}\phi=0$. 
				\cref{thm:properties_IMorley}.\ref{item:IMorley_kappa} and 
				the min-max principle for \eqref{eq:contEVP} show
				\begin{align*}
					1=\Vert \phi\Vert_{L^2(\Omega)}^2
						=\Vert (1-I_{M})\phi\Vert_{L^2(\Omega)}^2
						\le \kappa_{2}^2 h_{\max}^{4}\vvvert (1-I_{M})\phi\vvvert_{\mathrm{pw}}^2
						= \kappa_{2}^2 h_{\max}^{4}\vvvert \phi\vvvert^2\le \lambda_k\kappa_{2}^2 h_{\max}^{4}.
				\end{align*} 
				It follows $\kappa_{2}^{-2}h_{\max}^{-4}\le\lambda_k$ and 
				$\textup{GLB}({k})= {\lambda_{M}(k)}/\big({1+\lambda_{M}(k) 
				\kappa_{2}^2h_{\max}^{4}}\big)\le 
				\kappa_{2}^{-2}h_{\max}^{-4}\le \lambda_k$ concludes the proof. 
				Theorem~2 in \cite{CGal14} states \eqref{eq:GLB_CR_Morley} under a separation condition  
				$h_{\max}^4\kappa_2^2<(\sqrt{1+k^{-1}}-1)/\sqrt{\lambda_k}$ (with adapted notation) 
				in $2$D. 
				If $I_M\phi_1,\dots, I_M\phi_k$ are linearly independent, the proof  in \cite{CGal14} does not need 
				the separation condition to show \eqref{eq:GLB_CR_Morley} and holds with  \cref{thm:properties_IMorley}.\ref{item:IMorley_kappa}  
				 for $n\ge 3$. 
				The remaining details in \cite{CGal14} apply verbatim and are therefore omitted. 
		\end{proof}			
\subsection{Comments}\label{sec:comments} 
 	This section  introduces a reduced problem \eqref{eq:disEVP} as a disturbed nonconforming eigenvalue problem 	
 	and compares the new method with the standard Morley formulation.
			\subsubsection{Reduced eigenvalue problem}
				Under the condition $\lambda_h\kappa_{2}^2h^{{4}}_{\max}< 1$, the algebraic eigenvalue problem 
				\eqref{eq:dis_EVP_alt} is equivalent to a reduced form that seeks 
				$(\lambda_h,u_{{M}})\in \mathbb{R}^+\times \big(M(\mathcal{T})\setminus\{0\}\big)$ with 
				\begin{align}
					a_{\mathrm{pw}}(u_{{M}},v_{{M}})
					=\lambda_h\, b\Big(\frac{u_{{M}}}{1-\lambda_h\kappa_{2}^2 h_{\mathcal{T}}^{{4}}},v_{{M}}\Big)
							\qquad\text{for all }v_{{M}}\in M(\mathcal{T}). \label{eq:disEVP}
				\end{align}
				The formulation \eqref{eq:disEVP}  is reduced in that the additional variables $u_{\mathrm{pw}}$ and $v_{\mathrm{pw}}$ in 
				\eqref{eq:dis_EVP_alt} are condensed out; but 
				\eqref{eq:disEVP} is a (simple) rational eigenvalue problem with the same dimension and sparsity as the Morley eigenvalue problem 
				\eqref{eq:disEVP_NC}.  A solution $(\lambda_h, u_M)$ to \eqref{eq:disEVP} is also called an eigenpair and the (geometric) multiplicity $\ge 1$ is the 
				dimension of the eigenspace of all $u_M\in M(\mathcal{T})\setminus\{0\}$ so that $(\lambda_h,u_M)$ solves \eqref{eq:disEVP}. 
				The numerical treatment of the rational eigenvalue problem \eqref{eq:disEVP} is left as a topic for future research in numerical linear algebra. 
				 
			\begin{proposition}[equivalence]\label{lem:AquivProb}
			\begin{enumerate}[label=(\alph*),wide]
			\item\label{item:AquivProp_1} If the eigenpair $(\lambda_h,\boldsymbol{u_h})$ 
				 	of  \eqref{eq:dis_EVP_alt} satisfies $\lambda_h\hspace{-.5mm}<\hspace{-.5mm}\kappa_{2}^{\hspace{-.5mm}-2}h_{\max}^{-4}$ with 
				 	 $\boldsymbol{u_h}=(u_{\mathrm{pw}},u_{{M}})\in \boldsymbol{V_h}\setminus\{0\}$, 
				 	then   $(\lambda_h, u_{{M}})$ solves \eqref{eq:disEVP} and 
				 	$u_{\mathrm{pw}}=(1-\lambda_h\kappa_{2}^2h_{\mathcal{T}}^{{4}})^{-1}u_{{M}}$.
			\item\label{item:AquivProp_2}  If  
					$(\lambda_h,u_{{M}})$ is a solution to \eqref{eq:disEVP}  with $0<\lambda_h<\kappa_{2}^{-2}h_{\hspace{-.5mm}\max}^{-4}$, 
					then 
					$\lambda_h$, $u_{{M}}$, and 
					$u_{\mathrm{pw}}=(1-\lambda_h\kappa_{2}^2h_{\mathcal{T}}^{{4}})^{-1}u_{{M}}$ in  
					$ \boldsymbol{u_h}=(u_{\mathrm{pw}},u_{{M}})$ 
					form an eigenpair $(\lambda_h,\boldsymbol{u_h})$ of \eqref{eq:dis_EVP_alt}. 
			\item\label{item:AquivProp_3}  
					The set of eigenvalues of \eqref{eq:dis_EVP_alt} in the open interval $(0, \kappa_{2}^{-2}h_{\max}^{-4})$ is equal to the set of 
					solutions $\lambda_h$ to \eqref{eq:disEVP} in $(0, \kappa_{2}^{-2}h_{\max}^{-4})$ counting (geometric) multiplicities. 
			\end{enumerate}
			\end{proposition}
			\begin{proof}
				Throughout this proof with a specific and fixed $\lambda_h$ with $0<\lambda_h<\kappa_{2}^{-2}h_{\max}^{-4}$, abbreviate 
					\begin{align}
						\delta:=\frac{1}{1-\lambda_h\kappa_{2}^2h_{\mathcal{T}}^{{4}}}-1
							=\frac{\lambda_h\kappa_{2}^2h_{\mathcal{T}}^{{4}}}{1-\lambda_h\kappa_{2}^2h_{\mathcal{T}}^{{4}}}
							=h_{\mathcal{T}}^{{4}}\kappa_{2}^2\lambda_h(1+\delta)\in P_0(\mathcal{T}). \label{eq:def_delta}
					\end{align}
				\begin{proofof}\textit{\ref{item:AquivProp_1}.} 	
					Suppose that  $(\lambda_h, \boldsymbol{u_h}) \in \mathbb{R}^+\times \boldsymbol{V_h}$ is an eigenpair 
					of \eqref{eq:dis_EVP_alt}. 
					For $v_\mathrm{pw}\in P_2(\mathcal{T})$, the test function $(v_{\mathrm{pw}},0)\in \boldsymbol{V_h}$ in 
					\eqref{eq:dis_EVP_alt} shows
					$
						\kappa_{2}^{-2}h_{\mathcal{T}}^{-{4}}(u_{\mathrm{pw}}-u_{{M}})=\lambda_h u_{\mathrm{pw}}. 
					$
					This is equivalent to $u_{\mathrm{pw}}=(1+\delta)u_{{M}}$. The test function 
					$\boldsymbol{v_h}=(v_{{M}},v_{{M}})\in M({\mathcal{T}})\times M({\mathcal{T}})\subset\boldsymbol{V_h}$ 
					in \eqref{eq:dis_EVP_alt} leads to 
					\begin{align*}
						a_{\mathrm{pw}}(u_{{M}},v_{{M}})=\lambda_h (u_{\mathrm{pw}}, v_{{M}})_{L^2(\Omega)}
							= \lambda_h ((1+\delta)u_{{M}}, v_{{M}})_{L^2(\Omega)}.
					\end{align*}
					Since this holds for all $v_{{M}}\in M({\mathcal{T}})$, $(\lambda_h, u_{{M}})$ is a solution to
					\eqref{eq:disEVP}. 
				\end{proofof}
				\begin{proofof}\textit{\ref{item:AquivProp_2}.} 
					Suppose $(\lambda_h, u_{{M}})\in \mathbb{R}^+\times M({\mathcal{T}})$ is a solution to 
					\eqref{eq:disEVP} with $0<\lambda_h<\kappa_2^{-2}h_{\max}^{-4}$ and $\delta$ 
					in \eqref{eq:def_delta}. 
					Then $u_{\mathrm{pw}}:=(1+\delta)u_{{M}}$ and $\boldsymbol{u_h}:=(u_{\mathrm{pw}},u_{{M}})$ in
					\eqref{eq:disEVP} lead to
					\begin{align*}
						\boldsymbol{a_h}(\boldsymbol{u_h},\boldsymbol{v_h})&
							= \lambda_h  ((1+\delta)u_{{M}},v_{{M}})_{L^2(\Omega)}
								+ (\kappa_{2}^{-2}h_{\mathcal{T}}^{-{4}}\delta u_{{M}}, v_{\mathrm{pw}}-v_{{M}})_{L^2(\Omega)}
					\end{align*}
					for all $\boldsymbol{v_h}=(v_{\mathrm{pw}},v_{{M}})\in \boldsymbol{V_h}$.
					Recall $\kappa_{2}^{-2}h_{\mathcal{T}}^{-{4}}\delta=\lambda_h(1+\delta)$ from 
					\eqref{eq:def_delta} to verify 
					\begin{align*}
						\boldsymbol{a_h}(\boldsymbol{u_h},\boldsymbol{v_h})&
							=\lambda_h ((1+\delta)u_{{M}},v_{\mathrm{pw}})_{L^2(\Omega)}
							=\lambda_h \boldsymbol{b_h}(\boldsymbol{u_h},\boldsymbol{v_h}). 
					\end{align*}
					Hence $(\lambda_h,\boldsymbol{u_h})$ is an eigenpair of \eqref{eq:dis_EVP_alt}.  
					\end{proofof}
					\begin{proofof}\textit{\ref{item:AquivProp_3}.} 
					The combination of \ref{item:AquivProp_1}--\ref{item:AquivProp_2} proves the equality of the solution sets for $\lambda_h$ 
					in \eqref{eq:dis_EVP_alt} resp. \eqref{eq:disEVP} in the open interval $(0,\kappa_2^{-2}h_{\max}^{-4})$. 
					It remains to prove  equality of the multiplicities. Let $p$ resp. $q$ denote the multiplicities 
					of the eigenvalue $\lambda_h$ in \eqref{eq:dis_EVP_alt} resp. \eqref{eq:disEVP}. 
					The point is that $u_{M,1},\dots, u_{M,\gamma}\in M(\mathcal{T})$ are linear independent if and only if 
					$\big((1+\delta)u_{M,1},u_{M,1})\big),\dots, \big((1+\delta)u_{M,\gamma},u_{M,\gamma})\big)$ are linear independent in $\boldsymbol{V_h}$ 
					for the fixed $\delta$ from \eqref{eq:def_delta} for $\lambda_h$ with $0<\lambda_h<\kappa_2^{-2}h_{\max}^{-4}$. 
					This and \ref{item:AquivProp_1} resp. \ref{item:AquivProp_2} imply $p\le q$ resp. $q\le p$. Consequently $p=q$ 
					concludes the proof. 
					\end{proofof}
			\end{proof}
		\subsubsection{Comparison with Morley eigenvalues}			
			For $k\le N:=\textup{dim}(M(\mathcal{T}))$,
			the following \cref{lem:compare} allows the placement of $\lambda_h(k)$ from \eqref{eq:dis_EVP_alt} in 
			$\textup{GLB}(k)\le \lambda_h(k)\le \lambda_M(k)$ between $\textup{GLB}(k)$ from \eqref{eq:GLB_CR_Morley} and 
		 $\lambda_M(k)$ from \eqref{eq:disEVP_NC}. 
			\begin{lemma}[comparison]\label{lem:compare}
				For any  ${k}=1,\dots, N$,  the ${k}$-th algebraic eigenvalues $\lambda_h({k})$ of \eqref{eq:disEVP} and  
				$\lambda_{M} ({k})$ of \eqref{eq:disEVP_NC} satisfy  $\lambda_h({k})\le \lambda_{M} ({k})$.
				If $\lambda_h({k})$ satisfies $\lambda_h({k})\kappa_{2}^2h^{4}_{\max}< 1$, then
				$\textup{GLB}({k})\le \lambda_h({k})$ holds.  For a 
				uniform triangulation $\mathcal{T}$ with $h_{\max}=h_{\mathcal{T}}$ a.e. in $\Omega$ follows equality 
				$\textup{GLB}({k})= \lambda_h({k})$.
			\end{lemma}
			\begin{proof}
					The first result is a straightforward modification of \cite[Thm.~6.2]{CZZ18}. 
					Since the Morley eigenfunctions are linearly independent, the pairs $(\phi_M(1),\phi_M(1)),\dots,(\phi_M(k),\phi_M(k))$ 
					form a $k$-dimensional subspace of $\boldsymbol{V_h}$. Hence the min-max principle proves the claim. 
					The test functions $(v_{M},v_{M})\in M(\mathcal{T})\times M(\mathcal{T})\subset\boldsymbol{V_h}$ and  
					$(v_{\mathrm{pw}},0)\in \boldsymbol{V_h}$  in \eqref{eq:dis_EVP_alt} show for the first $k$ eigenpairs 
					$(\lambda_h,\boldsymbol{u_h})\in\mathbb{R}_+\times\boldsymbol{V_h}$ with 
					$\boldsymbol{u_h}=(u_{\mathrm{pw}},u_{M})\in P_2(\mathcal{T})\times M(\mathcal{T})$ 
					of \eqref{eq:dis_EVP_alt},
					 that $a_{\mathrm{pw}}(u_{M},v_{M})=\lambda_h b(u_{\mathrm{pw}},v_{M})$ for all $v_M\in M(\mathcal{T})$
				    \text{ and }
					$u_{M}=(1-\lambda_h\kappa_{2}^2h_{\mathcal{T}}^{4})u_{\mathrm{pw}}$. 
					Hence the arguments for the second inequality are analogue to \cite[Thm.~6.4]{CZZ18} with $\varepsilon$ 
					replaced by $\kappa_2^2 h_{\max}^4$ 
					and therefore further details are omitted. 
					On a uniform mesh with $h_{\max}=h_{\mathcal{T}}$ a.e. in $\Omega$ the scaling on the 
					right-hand side of \eqref{eq:disEVP} is constant, thus 
					 \eqref{eq:disEVP} and 
					\eqref{eq:disEVP_NC} are equivalent with 
					$\lambda_h({k})={\lambda_{M} ({k})}/({1+\kappa_{2}^2\lambda_{M} ({k})h_{\max}^{4}})=GLB({k})$.
				\end{proof}
			\begin{remark}[verification of the mesh-size condition]
				\cref{lem:compare} and $\lambda_{M} (k)<\kappa_{2}^{-2}h_{\max}^{-4}$ guarantee that the 
				discrete eigenvalue $\lambda_h(k)$ satisfies $\lambda_h(k) <\kappa_{2}^{-2}h_{\max}^{-4}$. 
				This is an a priori test sufficient for the applicability of \cref{lem:AquivProb}.\hfill $\Box$
			\end{remark}
			\subsubsection{An extra-stabilized Crouzeix-Raviart FEM}
			The arguments of this paper allow for an eigenvalue solver of the Dirichlet eigenvalue of the 
			Laplacian with guaranteed lower eigenvalue bound. 
			For the Laplace eigenvalue problem $-\Delta u=\lambda u$ in $H^1_0(\Omega)$ an extra-stabilised 
			Crouzeix-Raviart FEM comparable to \eqref{eq:dis_EVP_alt} 
			seeks 
			$\big(\lambda_h, (u_{\mathrm{pw}},u_{CR})\big)
					\in \mathbb{R}_+\times (P_1(\mathcal{T})\times \textit{CR}^1_0(\mathcal{T}))\setminus\{0\}$, 
			such that 
			\begin{align}
				(\nabla_{\mathrm{pw}}u_{CR},\nabla_{\mathrm{pw}}v_{CR})_{L^2(\Omega)}+\kappa_{1}^{-2}(h_{\mathcal{T}}^{-2}
						(u_{\mathrm{pw}}-u_{CR}), v_{\mathrm{pw}}-v_{CR})_{L^2(\Omega)}
						=\lambda_h b(u_{\mathrm{pw}},v_{\mathrm{pw}})\label{eq:disEVP_CR}						
			\end{align}
			for any $(v_{\mathrm{pw}},v_{CR})\in P_1(\mathcal{T})\times \textit{CR}^1_0(\mathcal{T})$. 
			The eigenvalue problem \eqref{eq:disEVP_CR} is for $n=2$ the lowest-order skeleton method in \cite{CZZ18}; 
			for $n\ge 3$ it is a completely different method.
			The standard interpolation operator (see e.g. \cite{CGed14,CP18}) 
			for the Crouzeix-Raviart FE $\textit{CR}^1_0(\mathcal{T})\subset P_1(\mathcal{T})$ \cite{CR73} 
			satisfies the conditions in \cref{thm:properties_IMorley}.\ref{item:Pi0Morley}--\ref{item:IMorley_kappa}
			(with the Hessian $D^2$ replaced by the gradient $\nabla$, $V=H^2$ replaced by $H^1$, $h_T^2$ by $h_T$, 
			and $P_2(\mathcal{T})$ by $P_1(\mathcal{T})$). 
			Hence the results analogue to those of  \cref{sec:GLB}--\ref{sec:comments} hold 
			for the Dirichlet eigenvalues of the Laplacian and the discrete eigenpairs of \eqref{eq:disEVP_CR}.
			(A conforming companion with the properties in \cref{thm:J_Morley} (again $D^2$ 
			replaced by $\nabla$, $V=H^2$ by $H^1$, $h_{\mathcal{T}}^2$ by $h_{\mathcal{T}}$, 
			and $P_2(\mathcal{T})$ by $P_1(\mathcal{T})$) is designed in \cite[Prop.~2.3]{CGS15} for $n=2$;  
			a generalization for $n\ge 3$ is straight-forward.)
			\FloatBarrier	
	
\section{Convergence rates in $3$D}\label{sec:convergence}	
This section presents a conforming companion in $3$D to apply  the Babu\v{s}ka-Osborn convergence analysis  \cite{BO91} for the discrete eigenvalue problem \eqref{eq:dis_EVP_alt} and the standard Morley eigenvalue problem \eqref{eq:disEVP_NC}. 
For the latter the paper \cite{YangLiBi2016} for $n\ge 2$ follows \cite{Rannacher1979} for $n=2$ and utilizes the trace inequality for second order derivatives $\partial^\alpha u/\partial x_\alpha$ for $|\alpha|=2$ under the regularity assumption $u\in W^{3,p}(\Omega)$ for $4/3<p\le 2$. Those terms arise in an integration by parts in the classical a priori error analysis of the Morley FEM. The present paper circumvents this by using the companion operator $J_M$ following  \cite{CGS13,Gal15,CCNN21}.  
This allows results for 
a general $u\in H^{2+\sigma}(\Omega)$ even for small $\sigma$ with $0<\sigma\le 1$.
\subsection{Conforming companion}\label{sec:companion}
		The  conforming companion operator $J_M$ in this paper is seen as a right-inverse of the Morley interpolation operator
		 $I_M: V\to  M(\mathcal{T})$ from \eqref{eq:def_IM} 
		 with an additional $L^2$ orthogonality.
		
			\begin{theorem}[properties of $J_M$]\label{thm:J_Morley}
			 	There exists a constant $M_2\approx 1$ (that exclusively depends on $\mathbb{T}$) 
			 	and a conforming companion $J_Mv_M\in V:=H^2_0(\Omega)$ for any $v_M\in M(\mathcal{T})$ with 
				 \begin{enumerate}[label=(\alph*), ref=\alph* ,wide]
			 		\item\label{item:lem_J_rightInverse} 
			 				$J_M$ is a right inverse to the interpolation $I_M$ in that
			 					$
			 						I_M\circ J_M =\textup{id} \text{ in  } {M(\mathcal{T})}, 
			 					$
					\item\label{item:lem_J_approx} 
						$
							\Vert h_{\mathcal{T}}^{-{2}} (1-J_M)v_{M}\Vert_{L^2(\Omega)}
							+\vvvert (1-J_M)v_{M}\vvvert_{\mathrm{pw}}
							\le {M_2} \min_{v\in V}\vvvert v_{M}-v\vvvert_{\mathrm{pw}}
						$,			
			 		\item\label{item:lem_J_orthogonality}
						the orthogonality	$(1-J_M)(M(\mathcal{T}))\perp P_{2}(\mathcal{T})$  
						holds in $L^2(\Omega)$.		
			 	\end{enumerate}
			\end{theorem}
			\begin{outlineproof}	
				The design can follow the $2$D discussions in \cite{Gal15, VZ19} in the spirit of \cite{CGS15}: 
				one subtle issue is the scaling of the nodal basis functions for the $\textit{WF}$ FEM 
				as a generalization of 				
				\cite[Thm.~6.1.3]{Ciarlet78} to $3$D. While the technical details of the proof are provided in the supplement, an outline of 
				the design will follow here. 
				\\[1ex]
				\textit{$\textit{WF}$ 
				partition.}
				 Unlike the $HCT$ partition of each triangle in $3$ subtriangles, 
				the subdivision in the $3$D  $\textit{WF}$ 
				finite element scheme  \cite{WF87,Sor09} 
				depends on the triangulation $\mathcal{T}\in\mathbb{T}$. Each tetrahedron $T\in\mathcal{T}$ is 
				divided into $12$ sub-tetrahedra with respect to a careful selection of center points $c_F$
				on each facet $F\in\mathcal{F}(T)$ of $T$ and $c_T$ inside the tetrahedron $T\in\mathcal{T}$:   
			    $c_T$ is the midpoint of the incircle of $T$ and 
				$c_F$ is the intersection of $F\in\mathcal{F}(\Omega)$ with the straight line through 
				$c_{T_+}$ and $c_{T_-}$ for $T_{\pm}\in\mathcal{T}$ aligned to $F=\partial T_+\cap\partial T_-$, while 
				  $c_F:=\textup{mid}(F)$ is simply the center of gravity for a triangle $F\in\mathcal{F}(\partial\Omega)$ on the boundary.
				Theorem~A.3 
				of the supplement guarantees the (uniform) shape-regularity of the resulting 
				subtriangulation $\widehat{\mathcal{T}}$ 
				for $\mathcal{T}\in\mathbb{T}$ and that 
				the distance of each center point $c_T$ (resp. $c_F$) to the boundary $\partial T$ (resp. the relative boundary $\partial F$) 
				is bounded from below
				by some global constant times $h_T$ (resp. $h_F$). 
				\\[1ex]
				\textit{$\textit{WF}$ 
				finite element. }
				The $28$ local degrees of freedom for any $K=\textup{conv}\{Q_1,$ $Q_2,$ $Q_3,$ $Q_4\}$ $\in\mathcal{T}$  are
				the evaluation of the function $f\in H^2(K)$ and its gradient $\nabla f$ at the vertices 
				$Q_1,\dots, Q_4$ of $K$ and 
				the evaluation of the gradient $\tau_{E}\times\nabla f(\textup{mid}(E))$ 
				at each edge midpoint $\textup{mid}(E)$ for $E\in\mathcal{E}(T)$ with unit tangent vector $\tau_E$.
				This determines $\nabla f(\textup{mid}(E))$ in the non-tangential directions $\nu_{E,1},\nu_{E,2}$ 		
				with $\textup{span}\{\tau_{E},\nu_{E,1},\nu_{E,2}\}=\mathbb{R}^3$. 
				Those $28$ degrees of freedom define a finite element $\big(K, C^1(K)\cap P_3(\widehat{\mathcal{T}}(K)),\{L_1,\dots,L_{28}\}\big)$ 
				in the sense of Ciarlet. 
				Since the explicit proof of this is not included in \cite{WF87},    
				Theorem~A.1 
				provides it in the supplement.
				The facet center point 
				$c_F$ of an interior  facet $F=\partial T_+\cap \partial T_-\in\mathcal{F}(\Omega)$  shared by the two neighbouring tetrahedra 
				$T_\pm\in\mathcal{T}$ belongs to the same straight line as their center points $c_{T_+}$ and $c_{T_-}$ and then \cite{WF87} 
				implies $C^1$~conformity of  $\textit{WF}(\mathcal{T})
				:=P_3(\widehat{\mathcal{T}})\cap V$.  
				Theorem~A.2 
				in the  supplement provides a comprehensive proof, that is supposed to be readable without 
				profound a priori knowledge of Bernstein polynomials \cite{Boor1987} in multivariate $C^1$ splines. 
				\\[1ex]
				\textit{Scaling of the $\textit{WF}$ 
				basis functions. }
				Let $\varphi_{z,1}, \dots \varphi_{z,4}$ and $\varphi_{E,1},\varphi_{E,2}$ for  $z\in\mathcal{V}(\Omega)$ and $E\in\mathcal{E}(\Omega)$ 
				denote the nodal $\textit{WF}$ 
				basis functions dual to the global degrees of freedom for  $z\in\mathcal{V}$, $j=1,2,3$, $E\in\mathcal{E}$, and $\mu=1,2$, 
				\begin{align*}
						L_{z,1}f:=f(z), \quad L_{z,j+1}:=\frac{\partial f}{\partial x_j}(z),\quad \text{and}\quad
						L_{E,\mu}f:=\frac{\partial f}{\partial \nu_{E,\mu}}(\text{mid}(E))
				\end{align*}
				(such that $L_{z,j}(\varphi_{a,k})=\delta_{za}\delta_{jk}$,  
				$L_{E,\mu}(\varphi_{F,\kappa})=\delta_{EF}\delta_{\mu\kappa}$, 
				and $L_{z,j}(\varphi_{E,\mu})=0=L_{E,\mu}(\varphi_{z,j})$ 
				for any $a,z\in\mathcal{V}$, $E,F\in\mathcal{E}$, $j,k=1,\dots,4$, and $\mu,\kappa=1,2$). 
				Theorem~A.4 
				in the supplement  generalizes a conclusion of \cite[Thm.~6.1.3]{Ciarlet78} to $\textit{WF}$ 
				and asserts 
				the expected scaling of the nodal basis functions. 
				For $s=0,1,2$ 
				\begin{align}
					h_\ell \Vert \varphi_{z,1}\Vert_{H^s(\Omega)}+
					\Vert \varphi_{z,j+1}\Vert_{H^s(\Omega)}+
					\Vert \varphi_{E,\mu}\Vert_{H^s(\Omega)}\lesssim h_\ell^{5/2-s} \label{eq:Scaling}
				\end{align}							 
				 holds with the volume $h_\ell^{3}:=|\textup{supp}(\varphi_\ell)|$ of the nodal patch 
				$\overline{\omega(z)}:=\bigcup\mathcal{T}(z)$, $\mathcal{T}(z):=\{T\in\mathcal{T}:\, z\in\mathcal{V}(T)\}$ 
				for $\varphi_\ell= \varphi_{z,1}$ or $\varphi_\ell= \varphi_{z,j+1}$,  and 
				of the edge patch $\overline{\omega(E)}:=\bigcup\mathcal{T}(E)$, $\mathcal{T}(E):=\{T\in\mathcal{T}:\, E\in\mathcal{E}(T)\}$ 
				for $\varphi_\ell=\varphi_{E,\mu}$. The point is that the constants in \eqref{eq:Scaling} are uniformly bounded in terms of the 
				uniform shape-regularity of $\mathbb{T}$. 
				\\[1ex]
				The $\textit{WF}$ 
				allows the four-step design of $J_M\equiv J_4$ below. Details of the proofs are provided in Supplement~B. 
				\\[1ex]
				\textit{Definition of $J_1$. }
				The enrichment operator $J_1:M(\mathcal{T})\to \textit{WF}(\mathcal{T})$ 
				with homogeneous boundary 
				conditions is defined by averaging of degrees of freedom of $\textit{WF}(\mathcal{T})$ 
				The scaling of the $\textit{WF}$ 
				basis function \eqref{eq:Scaling} is a key argument in the proof of the local approximation property 
				\begin{align}
					h_T^{-4}\Vert v_M-J_1v_M\Vert_{L^2(T)}^2	
		 				\lesssim \sum_{z\in\mathcal{V}(T)}\sum_{F\in\mathcal{F}(z)}h_F\Vert [D^2 v_M]_F\times 
		 				\nu_F\Vert_{L^2(F)}^2\label{eq:proof_J1_approx}
				\end{align}
		 		for any $T\in\mathcal{T}$;  $\mathcal{F}(z):=\{F\in\mathcal{F}:z\in\partial F\}$ denotes the set of faces 
		 		with vertex $z\in\mathcal{V}$ in \eqref{eq:proof_J1_approx} and $[D^2v_M]_F\times \nu_F$ denotes the tangential components  of the jump $[D^2v_M]_F$ 
				across a side $F\in\mathcal{F}(z)$ with the row-wise cross product $[D^2v_M]_F\times\nu_F$ with the unit normal $\nu_F\in\mathbb{R}^3$. 
				 The estimate \eqref{eq:proof_J1_approx} and the $2$D arguments from \cite[Prop.~2.3]{Gal15} modified with the $\textup{Curl}$ operator 
				in $3$D as in \cite{CBJ01} lead to 
				$$\Vert h_{\mathcal{T}}^{-{2}} (1-J_1)v_{M}\Vert_{L^2(\Omega)}\lesssim \min_{v\in V}\vvvert v_{M}-v\vvvert_{\mathrm{pw}}.$$ 
				\\[1ex]
				\textit{Definition of $J_2$. }
				For each edge $E\in\mathcal{E}(\Omega)$ define  below a function 
				$\xi_E\in H^2_0(\widehat{\omega}(E))\subset H^2_0(\Omega)$ with $\intmean_G\xi_E\textup{d}s=\delta_{GE}$ 
				for all edges $G\in \mathcal{E}$, such that the  support 
				$\textup{supp}(\xi_E)\subset \widehat{\omega}(E)$ 
				is contained in the edge patch ${\widehat{\omega}}(E):=\textup{int}\big(\bigcup \widehat{\mathcal{T}}(E)\big)$ of $E$ 
				in the $\textit{WF}$ 
				partition $\widehat{\mathcal{T}}$. Then 
				$$  {J}_2(v_M):= J_1 v_M+
												\sum_{E\in\mathcal{E}(\Omega)}\Big(\intmean_E (v_M-J_1v_M)\,\textup{d}s\Big)\xi_E\in V.$$
				The shape-regularity of $\widehat{\mathcal{T}}$ (from Theorem~A.3) 
				allows the choice of a ball 
			 	$B:=B(\textup{mid}(E),R_E)$ $\subset\widehat{\omega}(E)$ with midpoint $\textup{mid}(E)$ and radius 
				$R_E$ such that  $h_T\approx R_E\approx h_E$ in the definition of $\xi_E\in C^1(\mathbb{R}^3)\cap H^2_0(B)$ by 
				 	\begin{align*}
				 		\xi_E(x):=
				 				\frac{|E|}{R_E}\Big(1-3\frac{|y|^2}{R_E^2}+2\frac{|y|^3}{R_E^3}\Big) \quad &\text{ for }  x\in \overline{B} 
				 				\text{ and }y:=x-\textup{mid}(E).
				 	\end{align*}
				\textit{Definition of $J_3$. }	
				For each side ${F}\in{\mathcal{F}}(\Omega)$ define below a function $\zeta_{{F}}\in H^2_0(\omega(F))\subset H^2_0(\Omega)$ with  
			 	$\intmean_{{G}} \nabla \zeta_{{F}}\cdot\nu_{G}\,\textup{d}s =\delta_{{G}{F}}$ for all sides 
			 	${G}\in{\mathcal{F}}$ and support 
			 	$\textup{supp}(\zeta_{{F}})\subset \overline{\omega(F)}$ in the face 
			 	patch ${\omega(F)}:=\textup{int}(T_+\cup T_-)$ of the neighbouring tetrahedra $T_\pm\in \mathcal{T}(F)$ 
				with $F=\partial T_+\cap \partial T_-$ and with unit normal vectors of a fixed orientation $\nu_F=\nu_{T_+}|_F=-\nu_{T_-}|_F$
			 	of $F$ in $\mathcal{T}$. Then  
			 		$$ J_3(v_M):= {J}_2 		
							 				v_M+\sum_{F\in\mathcal{F}(\Omega)}\Big(\intmean_{{F}}\nabla(v_M-{J}_2v_M)
							 					\cdot\nu_{{F}}\,\textup{d}s\Big)\zeta_{{F}}
							 	\in V.
					$$
			 	Suppose ${F}=\textup{conv}\{z_1,z_2, z_3\}\in{\mathcal{F}}(T_\pm)$ is  the common face  of 
				$T_\pm=\textup{conv}\{z_1,\dots,z_{4}^{\pm}\}\in\mathcal{T}$ opposite the vertex $z_4^{\pm}\in\mathcal{V}(T_\pm)$. 
				Let ${\lambda}_k^{\pm}$ denote the barycentric coordinate in $T_\pm$ associated with the vertex 
				$z_k^{\pm}\in\mathcal{V}(T_\pm)$ for $k=1,\dots, 4$. Then $\zeta_F\in P_{7}(\mathcal{T})\cap C^1(\Omega)\cap H^2_0(\omega(F))$ 
				reads 
					\begin{align*}
						\zeta_{F}|_{T_\pm}
							:=\pm\frac{7!}{2} \textup{dist}(z_{4}^{\pm},{F}) ({\lambda}_1\lambda_2{\lambda}_3)^2\, {\lambda}_{4}^{\pm} 
							\in P_{7}(T_\pm) 
					\end{align*}
				in $T_\pm\in\mathcal{T}(F)$ (and vanishes outside the face patch ${\omega(F)}$). 
				The integral mean corrections guarantee that $J_3v_M$ satisfies  (\ref{item:lem_J_rightInverse}).
				\\[1ex]
				\textit{Definition of $J_4$. }	
				The correction $J_4v_M\in V$ is designed such  that its $L^2$ projection $\Pi_2(J_4v_M)$ onto $P_2(\mathcal{T})$ 
				coincides with $v_M\in M(\mathcal{T})\subset P_2(\mathcal{T})$, i.e., 	
				$J_4\equiv J_M$ satisfies (\ref{item:lem_J_orthogonality}). 
				For any $T\in\mathcal{T}$, recall the barycentric coordinate $\lambda_z$ associated with the vertex $z\in \mathcal{V}(T)$ in $T$, and 	
				define the scaled squared volume-bubble function 
				$b_T:=4^8\prod_{z\in\mathcal{V}(T)}\lambda_z^2\in P_{8}(T)\cap H^2_0(T)\subset H^2_0(\Omega)$ with 
				$\Vert b_T\Vert_{L^\infty(T)}=1$.  Let $v_T\in P_2(T)$ denote the Riesz representation of the linear functional 
				$w_T\mapsto\int_T (v_M-J_3v_M) w_T\,\textup{d}x$ in the Hilbert space $P_2(T)$ endowed with the weighted $L^2$ scalar product 
				$(b_T\bullet,\bullet)_{L^2(\Omega)}$, such that 
				$
					(v_M-J_3v_M,w_T)_{L^2(T)}=(b_Tv_T,w_T)_{L^2(T)}\quad\text{for all }w_T\in P_2(T).  
				$
				 	Set 
					$$
						J_4 v_M:=J_3v_{M}+\sum_{T\in\mathcal{T}}v_Tb_T.
					$$
				\\[1ex]
				\textit{Outline of the proof of (\ref{item:lem_J_rightInverse})--(\ref{item:lem_J_orthogonality}). } 
				Since  $(v_Tb_T)|_{\partial T}\equiv 0 \equiv (\nabla(v_Tb_T))|_{\partial T}$  vanishes along the boundary 
				$\partial T$ of $T\in\mathcal{T}$ and $\zeta_F|_{\partial F}=0$ vanishes along the boundary $\partial F$ of 
				$F\in\mathcal{F}(\Omega)$, 
				$J_M\equiv J_4$ satisfies (\ref{item:lem_J_rightInverse}) and (\ref{item:lem_J_orthogonality}). 					
				The  above correction functions satisfy
				$\Vert\xi_E\Vert_{L^2(T)}\approx h_T^{3/2}$, $\Vert \zeta_F\Vert_{L^2(T)}\approx h_T^{5/2}$, and 
				$\Vert v_Tb_T\Vert_{L^2(T)} \lesssim \Vert v_M-J_3v_M\Vert_{L^2(T)}$ for 
				any tetrahedron $T\in\mathcal{T}$ with edge $E\in\mathcal{E}(T)$ and face $F\in\mathcal{F}(T)$.  
				This and 
				a combination of  inverse estimates \cite{BS08}, 
				Cauchy-Schwarz, and discrete trace inequalities ensure 
					$
						\Vert (1-J_M)v_{M}\Vert_{L^2(T)}\lesssim \Vert (1-J_1)v_{M}\Vert_{L^2(T)}.
					$ 
				Hence (\ref{item:lem_J_approx}) follows from the  local analysis of the averaging operator $J_1$  above.
				The details on the universal constant $M_2\approx 1$ are provided in Supplement~B.
			\end{outlineproof}	

			\begin{corollary}[further properties]\label{lem:ConformingCompanion} 
				Any $w\in V$ and  any $v_{M}\in M(\mathcal{T})$, $\mathcal{T}\in\mathbb{T}$, satisfy
		 		\begin{enumerate}[label=(\alph*), ref=\alph* ]
		 			\item\label{item:cor_J_b_vh}
		 			$\displaystyle b(v_{M}-J_Mv_{M},w)=b(v_{M}-J_Mv_{M},w-I_Mw)
		 			\le \Vert v_{M}-J_Mv_{M}\Vert_{L^2(\Omega)}\Vert w-I_Mw\Vert_{L^2(\Omega)}$\\
		 			$\displaystyle\phantom{xxxxxxxxxxxxxi}\le h_{\max}^{4} \kappa_{2}^2{M_2} 
		 			\min_{v\in V}\vvvert v_{M}-v\vvvert_{\mathrm{pw}}\min_{w_{M}\in M(\mathcal{T})}\vvvert w-w_{M}\vvvert_{\mathrm{pw}};$
		 			\item\label{item:cor_J_a_vh}
		 			$\displaystyle a_{\mathrm{pw}}( v_{M}-J_Mv_{M},w)
		 			=a_{\mathrm{pw}}(v_{M}-J_Mv_{M},w-I_Mw)\le
		 			\vvvert v_{M}-J_Mv_{M}\vvvert_{\mathrm{pw}}\vvvert w-I_Mw\vvvert_{\mathrm{pw}}$\\
		 			$\displaystyle\phantom{xxxxxxxxxxxxxxxi}	\le {M_2}\min_{v\in V}\vvvert v-v_{M}\vvvert_{\mathrm{pw}} \min_{w_{M}\in M(\mathcal{T})}
		 			\vvvert w-w_{M}\vvvert_{\mathrm{pw}}.$
		 		\end{enumerate}
		 	\end{corollary}
		 	\begin{proofof}{\textit{(\ref{item:cor_J_b_vh})}.}
		 		The identity follows from the orthogonality in \cref{thm:J_Morley}.\ref{item:lem_J_orthogonality} and 
		 		 the Cauchy-Schwarz inequality implies the first inequality.  
		 		The first term is controlled by \cref{thm:J_Morley}.\ref{item:lem_J_approx}
		 		and the second term by 
		 		\cref{thm:properties_IMorley}  for the Morley interpolation ((\ref{item:IMorley_kappa}) for $\ell=2$).
		 	\end{proofof}
		 	\begin{proofof}{\textit{(\ref{item:cor_J_a_vh})}.}
		 		The identity follows from the orthogonality in \cref{thm:properties_IMorley}.\ref{item:Pi0Morley}, 
		 		which also allows to bound the second term resulting from the Cauchy-Schwarz inequality. This and 
		 		\cref{thm:J_Morley}.\ref{item:lem_J_approx} conclude the proof of (\ref{item:cor_J_a_vh}). 
		 	\end{proofof}
		 	\begin{remark}[Guaranteed upper eigenvalue bounds]
		 		The companion operator $J_M$ can be employed in a postprocessing for guaranteed upper eigenvalue bounds as follows.
		 		Given $m\in\mathbb{N}$, let $(\lambda_{h}(j),\boldsymbol{u_{h,j}})$ with 
		 		$\boldsymbol{u_{h,j}}=(u_{\mathrm{pw},j},u_{M,j})\in V_h\setminus\{0\}$ denote the $j$-th eigenpair of 
		 		\eqref{eq:dis_EVP_alt} (or alternatively $(\lambda_{M}(j),u_{M,j})$ the $j$-th eigenpair of \eqref{eq:disEVP_NC}).  
		 		If $u_{M,1},\dots, u_{M,m}$ are linearly independent,  
		 		then $J_Mu_{M,1},\dots, J_Mu_{M,m}$ are linear independent vectors in $V$ as well, because $I_MJ_Mu_{M,j}=u_{M,j}$ from 
		 		\cref{thm:J_Morley}.\ref{item:lem_J_rightInverse}. 
		 		For the linear independence $u_{M,1},\dots, u_{M,m}$ in \eqref{eq:disEVP} the mesh-size condition $\lambda_h(m)\kappa_{2}^2h^{{4}}_{\max}< 1$ 
		 		is sufficient according to \cref{lem:AquivProb} (in \eqref{eq:disEVP_NC} the condition $m\le\textup{dim}(M(\mathcal{T}))$ is sufficient). 
		 		Then an $m\times m$ generalized algebraic eigenvalue problem with $A:=(a(J_Mu_{M,j},J_Mu_{M,k}): j,k=1,\dots,m)$ and 
		 		$B:=(b(J_Mu_{M,j},J_Mu_{M,k}): j,k=1,\dots,m)$ leads to algebraic eigenvalues $\mu_1\le\mu_2\le\dots\le\mu_m$. 
		 		The exact eigenvalue $\lambda_j\le\mu_j$ of \eqref{eq:contEVP} has the guaranteed upper bound $\mu_j$ by the min-max principle 
		 		\cite{StrangFix2008,Boffi2010}. The same strategy applies to the $\textit{CR}$-eigenvalue problem as well \cite{CGed14}; cf.  [LLX12] for an alternative post-processing.
			\end{remark}
		\subsection{Convergence analysis for the source problem }\label{sec:convergence_source}
			Recall $a(\bullet,\bullet):=(D^2\bullet,D^2\bullet)_{L^2(\Omega)}$ and its piecewise version
			$a_{\mathrm{pw}}(\bullet,\bullet):=(D^2_{\mathrm{pw}}\bullet,D^2_{\mathrm{pw}}\bullet)_{L^2(\Omega)}$. 
			Given $f\in L^2(\Omega)$, let $u\in V\equiv H^{2}_0(\Omega)$ solve  
			\begin{align}
				a(u, v)= (f,v)_{L^2(\Omega)}\quad\text{for all }v\in V. \label{eq:def_T}
			\end{align}
			Let $u_{M}\in M(\mathcal{T})$ denote the discrete solution to the Morley source problem 
			\begin{align}
				a_{\mathrm{pw}}(u_{M},v_{M})=(f, v_{M})_{L^2(\Omega)}
																			\quad\text{ for all }v_{M}\in M(\mathcal{T}).				
					\label{eq:nonconf_sourceproblem}
			\end{align}			
			Given the $L^2$ projection $\Pi_2$ onto $P_2(\mathcal{T})$, the second order  oscillation 
			of   $f\in L^2(\Omega)$ reads 
			$\textup{osc}_2(f,\mathcal{T}):=\Vert h_{\mathcal{T}}^2(1-\Pi_2)f\Vert_{L^2(\Omega)}$. Recall $0<\sigma\le1$ from \eqref{eq:def_sigma}.
		\begin{lemma}[discrete error estimate in $M(\mathcal{T})$]\label{cor:convergence_u_nc}
				There exist constants $C_{1},C_2\approx 1$ (that exclusively depend on $\mathbb{T}$) 
				such that given any $f\in L^2(\Omega)$, the exact solution $u\in H^{{2}+\sigma}(\Omega)\cap V$ to \eqref{eq:def_T} 
				 and the discrete solution $u_{M}\in M(\mathcal{T})$ to \eqref{eq:nonconf_sourceproblem} for $\mathcal{T}\in\mathbb{T}$ with maximal mesh-size 
				 $h_{\max}$ satisfy (\ref{item:norm_estimate_source_nc})--(\ref{item:L2_skp_evp_nc}). 
				\begin{enumerate}[label=(\alph*), ref =\alph*, wide]	
				\item \label{item:norm_estimate_source_nc} 
				 If  $u\in H^{{2}+s}(\Omega)$ for some $s$ with $\sigma\le s\le 1$, then 
								$$\displaystyle \vvvert u-u_{M}\vvvert_{\mathrm{pw}}+h_{\max}^{-\sigma}\Vert u-u_{M}\Vert_{L^2(\Omega)}
									\le C_{1} \big(h_{\max}^{s} \Vert u\Vert_{H^{2+s}(\Omega)}
									+\textup{osc}_{2}(f,\mathcal{T})\big).
								$$ 
				\item\label{item:L2_skp_evp_nc}
				Given any eigenvalue $\lambda$ of \eqref{eq:contEVP} with eigenspace 
				$E(\lambda)\subset H^{2+t}(\Omega)\cap V$ 
				for some $t$ with 
				$\sigma\le t\le {1}$, suppose $f, g\in E(\lambda)$. 
				Then 
							$$\displaystyle |(u-u_{M},g)_{L^2(\Omega)}|
									\le C_2(\lambda^{-1}+\kappa_2h_{\max}^4)^2 h_{\max}^{2t}\Vert f\Vert_{H^{2+t}(\Omega)}\Vert g\Vert_{H^{2+t}(\Omega)}.$$ 
				\end{enumerate}		  
			\end{lemma}
			\begin{proofof}{\textit{(\ref{item:norm_estimate_source_nc}).}}
			 \cref{thm:properties_IMorley} 
			shows that the interpolation operator $I_M$ satisfies the utilized properties \cite[Eqn.~(2.3),(2.5)]{Gal15}. 
			The $3$D companion operator $J_M$ in \cref{thm:J_Morley} satisfies in particular 
			the orthogonalities $\Pi_0(v_M-J_Mv_M)=0$ and $\Pi_0(D^2_{\mathrm{pw}}(v_M-Jv_M))=0$ for all $v_M\in M(\mathcal{T})$
			 and the approximation property \cite[Eqn.~(2.7),(2.8)]{Gal15}.
			Hence the arguments in \cite[Prop.~2.9--2.10]{Gal15} apply verbatim to $n=3$ and further details are omitted.
			\end{proofof}
				\begin{proofof}{\textit{(\ref{item:L2_skp_evp_nc}).}}
						Since $u\in V$ solves \eqref{eq:def_T} with the right-hand side 
						$f\in E(\lambda)$,  $u=f/\lambda\in H^{2+t}(\Omega)\cap V$. 
						\cref{rem:kappa2prime_2} and \cref{lem:InterpolationOperator}.\ref{item:cor_I_alpha} 
						prove $\textup{osc}_2(f,\mathcal{T})\le h_{\max}^{4+t}/\pi^{2+t}\,\Vert f\Vert_{H^{2+t}(\Omega)}$. 
						This and (\ref{item:norm_estimate_source_nc}) show 
						\begin{align}
						\vvvert u-u_M\vvvert_{\mathrm{pw}}
									\le C_1h_{\max}^t \big({1}/{\lambda} + {h_{\max}^{4}}/{\pi^{2+t}}\big)\Vert f\Vert_{H^{2+t}(\Omega)}.
									\label{eq:proof_error_EVP}
						\end{align}
						Since $g\in E({\lambda})$ is an eigenvector in \eqref{eq:contEVP}, 
						$\lambda (u-J_{M}u_{M},g)_{L^2(\Omega)}= a(u-J_{M}u_{M},g)$.
						\cref{lem:InterpolationOperator}.\ref{item:cor_I_a} implies 
						$a_{\mathrm{pw}}(u_{M}, g)=a_{\mathrm{pw}}(u_{M}, I_{M}g)$.
						This shows the first equality in 
						\begin{align*}
							\lambda (u-J_{M}u_{M},g)_{L^2(\Omega)}
									&= a(u,g)-a_{\mathrm{pw}}(u_{M}, I_{M}g)
									+a_{\mathrm{pw}}(u_{M}-J_{M}u_{M},g)\notag\\
									& = (f,g-I_{M}g)_{L^2(\Omega)}+a_{\mathrm{pw}}(u_{M}-J_{M}u_{M},g).\label{eq:0}
						\end{align*}
						The second equality follows because $u$ solves \eqref{eq:def_T} and $u_{M}\in M(\mathcal{T})$ solves 
						\eqref{eq:nonconf_sourceproblem} with right-hand side 
						$f$. The term $(f,g-I_{M}g)_{L^2(\Omega)}=(f,J_{M}I_{M}g-I_{M}g)_{L^2(\Omega)}+(f,g-J_{M}I_{M}g)_{L^2(\Omega)}$ 
						is split into two. 
						\cref{lem:ConformingCompanion}.\ref{item:cor_J_b_vh} controls the first contribution 
						\begin{align*}
							(f,J_{M}I_{M}g-I_{M}g)_{L^2(\Omega)}&
							\le {M_{2}}\ h_{\max}^{4}\kappa_{2}^{2}\vvvert f-I_{M}f\vvvert_{\mathrm{pw}}
								\vvvert g-I_{M}g\vvvert_{\mathrm{pw}}.\label{eq:1}
						\end{align*}
						\cref{lem:ConformingCompanion}.\ref{item:cor_J_a_vh} ensures 
						$$ 
							a_{\mathrm{pw}}(u_{M}-J_{M}u_{M},g)
								\le {M_{2}}\vvvert u-u_{M}\vvvert_{\mathrm{pw}}\vvvert g-I_{M}g\vvvert_{\mathrm{pw}}.
						$$ 
						Since $f\in E(\lambda)$ is an eigenvector in \eqref{eq:contEVP}, 
						$\lambda(f,g-J_{M}I_{M}g)_{L^2(\Omega)}=a(f, g-J_{M}I_{M}g)$. 
						The right-inverse property \cref{thm:J_Morley}.\ref{item:lem_J_rightInverse} and 
						the orthogonality in \cref{thm:properties_IMorley}.\ref{item:Pi0Morley} show
							$a(f, g-J_{M}I_{M}g)=a_{\mathrm{pw}}(f-I_{M}f, g-J_{M}I_{M}g)$. 
						A triangle inequality   and \cref{thm:J_Morley}.\ref{item:lem_J_approx} ensure 
						$\vvvert g-J_{M}I_{M}g\vvvert_{\mathrm{pw}}
							\le \vvvert g-I_Mg\vvvert_{\mathrm{pw}}  +\vvvert g-J_{M}I_{M}g\vvvert_{\mathrm{pw}}
							\le (1+{M_{2}})\vvvert g-I_{M}g\vvvert_{\mathrm{pw}}. $
						This and a Cauchy-Schwarz inequality verify 
						\begin{align*}
							\lambda(f,&g-J_{M}I_{M}g)_{L^2(\Omega)}\hspace{-.5mm}\le\hspace{-.5mm} 	\vvvert f-I_{M}f\vvvert_{\mathrm{pw}}
							\vvvert g-J_{M}I_{M}g\vvvert_{\mathrm{pw}}					
							\hspace{-1mm}
							\le \hspace{-1mm}(1+{M_{2}}) \vvvert f-I_{M}f\vvvert_{\mathrm{pw}}\vvvert g-I_{M}g\vvvert_{\mathrm{pw}}. 
						\end{align*}	
						The combination of the last four displayed estimates reads 
						\begin{align*}
							\lambda&(u-J_{M}u_{M},g)_{L^2(\Omega)}
								\hspace{-.5mm}\le \hspace{-.5mm}\vvvert g-I_{M}g\vvvert_{\mathrm{pw}} M_2
									\bigg(\hspace{-1mm}\bigg(\hspace{-.5mm}\frac{1+M_2^{-1}}{\lambda}+h_{\max}^{4}\kappa_{2}^{2}\hspace{-1mm}\bigg) 
									\vvvert f-I_{M}f\vvvert_{\mathrm{pw}}\hspace{-1mm}+\hspace{-1mm}\vvvert u-u_{M}\vvvert_{\mathrm{pw}}\hspace{-1mm}\bigg). 
						\end{align*}
						This,  \cref{lem:InterpolationOperator}.\ref{item:cor_I_alpha}, and 
						\eqref{eq:proof_error_EVP}  verify that 
						\begin{align*}
								\lambda&(u-J_{M}u_{M},g)_{L^2(\Omega)}\hspace{-.5mm}
								\le\hspace{-.5mm} h_{\max}^{2t} \Vert f\Vert_{H^{2+t}(\Omega)}\Vert g\Vert_{H^{2+t}(\Omega)} 
								\frac{M_2}{\pi^{t}} \hspace{-.5mm}
								\bigg(\hspace{-.5mm}\frac{1+M_2^{-1}}{\lambda\pi^t}+\frac{h_{\max}^{4}\kappa_2}{\pi^{t}}
								+C_1\hspace{-.5mm}\bigg(\hspace{-.5mm}\frac{1}{\lambda} + \frac{h_{\max}^{4}}{\pi^{2+t}}\hspace{-.5mm}\bigg)\hspace{-1mm} \bigg).
						\end{align*}
						For the term $(J_{M}u_{M}-u_{M},g)_{L^2(\Omega)}$, \cref{lem:ConformingCompanion}.\ref{item:cor_J_b_vh} 
						followed by  \cref{lem:InterpolationOperator}.\ref{item:cor_I_alpha} and \eqref{eq:proof_error_EVP}  
						show 
						\begin{align*}
							(J_{M}u_{M}-u_{M},g)_{L^2(\Omega)}
							 &\le M_2\kappa_2 h_{\max}^4\vvvert g-I_Mg\vvvert_{\mathrm{pw}}
							 	\vvvert u-u_M\vvvert_{\mathrm{pw}}
							\\&
							\le h_{\max}^{2t} \Vert f\Vert_{H^{2+t}(\Omega)}\Vert g\Vert_{H^{2+t}(\Omega)} 
								\frac{C_1M_2\kappa_2h_{\max}^{4}}{\pi^t} \bigg(\frac{1}{\lambda} + \frac{h_{\max}^{4}}{\pi^{2+t}}\bigg).
						\end{align*}  
						The combination of the last two displayed inequalities proves that 
						$(u-u_{M},g)_{L^2(\Omega)}
								=(J_{M}u_{M}-u_{M},g)_{L^2(\Omega)}+(u-J_{M}u_{M},g)_{L^2(\Omega)}
								\lesssim  h_{\max}^{2t} \Vert f\Vert_{H^{2+t}(\Omega)}\Vert g\Vert_{H^{2+t}(\Omega)}$. 
						The bookkeeping of the multiplicative constants concludes the proof. 								
					\end{proofof}
			Given any right-hand side $f\in L^2(\Omega)$,  let 
			$\boldsymbol{u_h}=(u_{\mathrm{pw}},u_{M})\in\boldsymbol{V_h}$ 
			denote the 
			discrete solution to the extra-stabilised source problem 
			\begin{align}
				\boldsymbol{a_h}(\boldsymbol{u_h},\boldsymbol{v_h})=(f, v_{\mathrm{pw}})_{L^2(\Omega)}
						\quad \text{for all }\boldsymbol{v_h}=(v_{\mathrm{pw}},v_{M})\in\boldsymbol{V_h}.  
						\label{eq:def_Th}
			\end{align}
			The analysis of \eqref{eq:def_Th} reduces to that of \cref{cor:convergence_u_nc} plus perturbation terms. 
		\begin{lemma}[discrete error estimate in $\boldsymbol{V_h}$]\label{lem:ThtoT}
				There exists a constant $C_{\mathrm{pw}}>0$ (that exclusively depends on $\mathbb{T}$),  such that given any 
				 $f\in L^2(\Omega)$, the exact solution $u\in H^{{2}+\sigma}(\Omega)\cap V$ to \eqref{eq:def_T} and 
				the discrete 
				solution $\boldsymbol{u_h}=(u_{\mathrm{pw}},u_{M})\in \boldsymbol{V_h}$ 
				to \eqref{eq:def_Th} for $\mathcal{T}\in\mathbb{T}$ with maximal mesh-size $h_{\max}$ satisfy 
				(\ref{item:norm_estimate_source_h})--(\ref{item:L2_skp_evp_h}). 	 
				\begin{enumerate}[label=(\alph*), ref =\alph*,wide]	
				
				\item\label{item:norm_estimate_source_h} 				 
				If  $u\in H^{{2}+s}(\Omega)$ for some $s$ with  $\sigma\le s\le 1$, then										
								$$\displaystyle \vvvert u-u_{\mathrm{pw}}\vvvert_{\mathrm{pw}}
								+h_{\max}^{-\sigma}\Vert u-u_{\mathrm{pw}}\Vert_{L^2(\Omega)}
									\le C_{\mathrm{pw}} \big(h_{\max}^{s} \Vert u\Vert_{H^{2+s}(\Omega)}
										+\textup{osc}_{2}(f,\mathcal{T})\big).$$
				\item\label{item:L2_skp_evp_h}
				Given any eigenvalue $\lambda$ of \eqref{eq:contEVP} with eigenspace 
				$E(\lambda)\subset H^{2+t}(\Omega)\cap V$ 
				for some $t$ with 
				$\sigma\le t\le {1}$, suppose $f, g\in E(\lambda)$  and $C_2$ from \cref{cor:convergence_u_nc}. Then 
							$$\displaystyle |(u-u_{\mathrm{pw}},g)_{L^2(\Omega)}|
							\le \big(C_2(\lambda^{-1}+\kappa_2h_{\max}^4)^2 +\kappa^2_2h_{\max}^2\big) 
							h_{\max}^{2t}\Vert f\Vert_{H^{2+t}(\Omega)}\Vert g\Vert_{H^{2+t}(\Omega)}.$$
				\end{enumerate}		  
			\end{lemma}
			\begin{proofof}\textit{(\ref{item:norm_estimate_source_h}). }
				For $v_{\mathrm{pw}}\in  P_2(\mathcal{T})$,
				the test-function $(v_{\mathrm{pw}},0)\in \boldsymbol{V_h}$  in \eqref{eq:def_Th} leads to 
				$\kappa_{2}^{-2}h_{\mathcal{T}}^{-{4}}(u_{\mathrm{pw}}-u_{M})=\Pi_{2} f$. 
				Thus $u_{\mathrm{pw}}=u_M+\kappa_{2}^{2}h_{\mathcal{T}}^{{4}}\Pi_{2} f$ and a triangle inequality shows
				\begin{align}
					&\vvvert u-u_{\mathrm{pw}}\vvvert_{\mathrm{pw}}+ h_{\max}^{-\sigma}\Vert u-u_{\mathrm{pw}}\Vert_{L^2(\Omega)}
					\notag\\
						 &\le \big(\vvvert u-u_{M}\vvvert_{\textup{pw}}+h_{\max}^{-\sigma}\Vert u-u_M \Vert_{L^2(\Omega)}\big)
							+ \kappa_{2}^2\big(\vvvert h_{\mathcal{T}}^{4} \Pi_{2} f\vvvert_{\mathrm{pw}}
							+h_{\max}^{-\sigma}\Vert h_{\mathcal{T}}^{4} \Pi_{2} f\Vert_{L^2(\Omega)}\big).\label{eq:norm_estimatesource_1}
				\end{align}
				The test-function $(v_{M},v_{M})\in M(\mathcal{T})\times M(\mathcal{T})\subset\boldsymbol{V_h}$ shows  
				$a_{\mathrm{pw}}(u_{M},v_{M})=(f, v_{M})_{L^2(\Omega)}$. 
				In other words, the solution component $u_{M}\in M({\mathcal{T}})$ solves the Morley source problem 	
				\eqref{eq:nonconf_sourceproblem}. 
				\cref{cor:convergence_u_nc}.\ref{item:norm_estimate_source_nc}  
				controls the first term in \eqref{eq:norm_estimatesource_1}, 
				$$\vvvert u-u_{M}\vvvert_{\textup{pw}}+h_{\max}^{-\sigma}\Vert u-u_M \Vert_{L^2(\Omega)}
					\le C_1  \big(h_{\max}^{s} \Vert u\Vert_{H^{2+s}(\Omega)}+\textup{osc}_{2}(f,\mathcal{T})\big). $$
				An 	inverse estimate for $P_{2}(\mathcal{T})$ with constant $c_{\mathrm{inv}}>0$  and the boundedness of $\Pi_2$    
				show 
				$$\vvvert h_{\mathcal{T}}^{4} \Pi_{2} f\vvvert_{\mathrm{pw}}
					+h_{\max}^{-\sigma}\Vert h_{\mathcal{T}}^{4} \Pi_{2} f\Vert_{L^2(\Omega)}
				\le (c_{\mathrm{inv}}+h_{\max}^{2-\sigma})\Vert h_{\mathcal{T}}^{2} \Pi_{2} f\Vert_{L^2(\Omega)}
				\le (c_{\mathrm{inv}}+h_{\max}^{2-\sigma})\Vert h_{\mathcal{T}}^{2}  f\Vert_{L^2(\Omega)}
				. $$			
				The efficiency  estimate $\Vert h^2_{\mathcal{T}}f\Vert_{L^2(\Omega)}\lesssim \vvvert u-u_{M}\vvvert_{\mathrm{pw}}
				+\textup{osc}_2(f,\mathcal{T})$ follows from the bubble-function methodology due to  \cite{Verf2013}, 
				cf. \cite[Thm.~2]{BdVNS07}. The combination with \cref{cor:convergence_u_nc}.\ref{item:norm_estimate_source_nc} 
				shows 
				$$\Vert h_{\mathcal{T}}^{2}f\Vert_{L^2(\Omega)} 
				\lesssim h_{\max}^{s} \Vert u\Vert_{H^{2+s}(\Omega)}
				+\textup{osc}_{2}(f,\mathcal{T}).$$
				The last two displayed inequalities bound the second term in \eqref{eq:norm_estimatesource_1} 
				and that concludes the proof of (\ref{item:norm_estimate_source_h}).
				\end{proofof}
				\begin{proofof}\textit{(\ref{item:L2_skp_evp_h}).}
				The substitution of $u_{\mathrm{pw}}=u_{M}+\kappa_{2}^2h_{\mathcal{T}}^{4}\Pi_{2} f$  from part (\ref{item:norm_estimate_source_h}) leads to 
						\begin{align*}
								(u-u_{\mathrm{pw}} ,g )_{L^2(\Omega)}
									= (u-u_{M},g)_{L^2(\Omega)}-\kappa_{2}^2(h_{\mathcal{T}}^4\Pi_{2} f,g)_{L^2(\Omega)}.
						\end{align*}
				Since $u_{M}$ solves the Morley source problem 	
				\eqref{eq:nonconf_sourceproblem},  \cref{cor:convergence_u_nc}.\ref{item:L2_skp_evp_nc}  
				controls $(u-u_{M},g)_{L^2(\Omega)}$.
				This and  
						$
							(\Pi_{2} f,g)_{L^2(\Omega)}\le \Vert f\Vert_{L^2(\Omega)}\Vert g\Vert _{L^2(\Omega)}
						$
			  	conclude the proof.
			\end{proofof}
		\subsection{Convergence rates for the eigenvalue problem }\label{sec:convergence_ev}
				The preparations in \cref{sec:companion}--\ref{sec:convergence_source} 
				allow the proof of the optimal a priori convergence rates in \cref{thm:BabuskaOsborn} 
				with fundamental arguments from \cite{BO91}. 

					\begin{proofof}\textit{\cref{thm:BabuskaOsborn}.}  
						Given any right-hand side $f\in L^2(\Omega)$ let $S(f):=u\in V$ denote the continuous 
						solution to \eqref{eq:def_T} 
						and let  $S_{h}(f):=u_{\mathrm{pw}}\in P_{2}(\mathcal{T})$ 
						denote the 
						first component of the solution 
						$\boldsymbol{u_h}=(u_{\mathrm{pw}},u_{M})\in\boldsymbol{V_h}$ to  
						\eqref{eq:nonconf_sourceproblem}.  
						This defines solution operators $S:L^2(\Omega)\to L^2(\Omega)$ and $S_{h}:L^2(\Omega)\to L^2(\Omega)$. 	
						\cref{lem:ThtoT}.\ref{item:norm_estimate_source_h} implies the convergence ${{S}}_h\to {{S}}$ 
						in the operator norm of $\mathcal{L}(L^2(\Omega))$ as 
						$h_{\max}\to 0$. 
						Suppose $(\lambda,\phi)\in\mathbb{R}^+\times V$ denotes an eigenpair of  \eqref{eq:contEVP} and 
						$(\lambda_h,\boldsymbol{u_{h}})\in\mathbb{R}^+\times \boldsymbol{V_h}$ with 
						$\boldsymbol{u_{h}}=(u_{\mathrm{pw}},u_{M})$ 
						denotes an eigenpair of  \eqref{eq:dis_EVP_alt}, then $(1/\lambda,\phi)$ is an eigenpair of 
						$S$ and $(1/\lambda_h, u_{\mathrm{pw}})$ is an eigenpair 
						of $S_{h}$.
						Hence, the Babu\v{s}ka-Osborn theory \cite{BO91}  (see also  \cite[Sec.~9]{Boffi2010} or 
						\cite[Sec.~1.4.2]{SZ17}) 
						implies for any non-zero eigenvalue $1/{\lambda}$  of ${{S}}$ with eigenspace 
						$E({\lambda})=\ker(\lambda^{-1}-{{S}})$ 
						of dimension $\mu=\textup{dim}(E(\lambda))\in \mathbb{N}$, that there exist exactly $\mu$ eigenvalues 
						$1/\lambda_{h,1},\dots, 1/\lambda_{h,\mu}$ of 
						${{S}}_h$,  which converge to $1/\lambda$ as $h_{\max}\to 0$. 
						The error estimates for the selfadjoint operator $S$ in \cite[Rem.~7.5]{BO91} 
						read (with a generic constant which depends on $\lambda$) 
						$
							\Vert u-u_{\mathrm{pw},k}\Vert_{L^2(\Omega)}
							\lesssim\Vert (S-S_{h})|_{E(\lambda)}\Vert_{\mathcal{L}(E(\lambda);L^2(\Omega))}
						$ and
					 	\begin{align*}
							\max_{1\le k\le \mu}\big|{\lambda}-{{\lambda}_{h,k}}\big|&\lesssim 			
							\max_{1\le k\le \mu}\big|\lambda^{-1}-{{\lambda}_{h,k}}^{-1}\big|\\
											&\lesssim \sup_{\phi,\psi\in E({\lambda})\setminus\{0\}}
												\frac{|((S-S_{h})\phi,\psi)_{L^2(\Omega)}|}
																{\Vert \phi\Vert_{L^2(\Omega)}\Vert\psi\Vert_{L^2(\Omega)}}
												+\Vert (S-S_{h})|_{E({{\lambda}})}\Vert_{\mathcal{L}(E({\lambda});{L^2(\Omega)})}^2.
						\end{align*}
						 For the finite-dimensional eigenspace $E(\lambda)\subset H^{2+t}(\Omega)$ there  exists 
						 $C_t:=\sup_{\phi\in E(\lambda)} \frac{\Vert \phi\Vert_{H^{2+t}(\Omega)}}{\Vert \phi\Vert_{L^2(\Omega)}}$ $<\infty$ so that 
						$
							\Vert \phi\Vert_{H^{2+t}(\Omega)}\le C_t\Vert \phi\Vert_{L^2(\Omega)}\text{ for all }\phi\in E(\lambda).		
						$
						The combination of  \cref{lem:ThtoT}.\ref{item:norm_estimate_source_h} with \cref{rem:kappa2prime_2} and 
						\cref{lem:InterpolationOperator}.\ref{item:cor_I_alpha} shows (as an analog to \eqref{eq:proof_error_EVP}) that 
						\begin{align*}
							\Vert (S-S_{h})|_{E(\lambda)}\Vert_{\mathcal{L}(E(\lambda);L^2(\Omega))}
									\le  h_{\max}^{t+\sigma}C_tC_{\mathrm{pw}} \big({1}/{\lambda} + {h_{\max}^{4}}/{\pi^{2+t}}\big).
						\end{align*}
						\cref{lem:ThtoT}.\ref{item:L2_skp_evp_h} bounds the remaining term 
						\begin{align*}
							\sup_{\phi,\psi\in E({\lambda})\setminus\{0\}}
												\frac{|((S-S_{h})\phi,\psi)_{L^2(\Omega)}|}
																{\Vert \phi\Vert_{L^2(\Omega)}\Vert\psi\Vert_{L^2(\Omega)}}
									\le h_{\max}^{2t} C_t^2\big(C_2(\lambda^{-1}+\kappa_2h_{\max}^4)^2 +\kappa^2_2h_{\max}^2\big) .
						\end{align*}
						This concludes the proof.
						\end{proofof}
					Unlike \cite{Rannacher1979,YangLiBi2016}, 
					the following \cref{thm:BabuskaOsborn_nonconforming} specifies the convergence rates for \eqref{eq:disEVP_NC} 
					directly in terms of $\sigma=\min\{1,\sigma_{\mathrm{reg}}\}$ for the index of  elliptic regularity $\sigma_{\mathrm{reg}}$
					 from  \eqref{eq:def_sigma} and the Sobolev regularity $t$ of $E(\lambda)$. 
				\begin{theorem}[a priori convergence for \eqref{eq:disEVP_NC}]\label{thm:BabuskaOsborn_nonconforming}
						Given a non-zero eigenvalue  $\lambda$ of \eqref{eq:contEVP} of multiplicity $\mu$, suppose that    		
						$E(\lambda)\subset H^{2+t}(\Omega)\cap V$ holds for some $t$ with 
						$\sigma \le t\le 1$. Then there exist $\delta,C>0$ such that any $\mathcal{T}\in\mathbb{T}(\delta)$ 
						and the discrete space $M(\mathcal{T})$ 
						lead in \eqref{eq:disEVP_NC} to exactly $\mu$ algebraic eigenvalues  
						$\lambda_{M,1},\dots, \lambda_{M,{\mu}}$ that converge to $\lambda$ as 
						$h_{\max}\to 0$. Let $E_M:=\textup{span}\{u_M\in E_M(\lambda_{M,k}): k=1,\dots,\mu\}$ 
						denote the union of the discrete eigenspaces $E_M(\lambda_{M,k})\subset M(\mathcal{T})$ for 
						$\lambda_{M,1},\dots, \lambda_{M,{\mu}}$. Then the convergence results in \cref{thm:BabuskaOsborn} hold with 
						$\lambda_{h,k}$ replaced by $\lambda_{M,k}$,  and $\boldsymbol{u_h}=(u_{\mathrm{pw}},u_{M})\in E_h$ with 
						$\Vert u_{\mathrm{pw}}\Vert_{L^2(\Omega)}=1$ replaced by $\phi_M\in E_M$ with $\Vert \phi_M\Vert_{L^2(\Omega)}=1$.
				\end{theorem}
					\begin{proof}
						Define the solution operator $S_{M}:L^2(\Omega)\to  L^2(\Omega)$ with $S_{M}(f):=u_{M}\in M(\mathcal{T})$ 
						the Morley finite element solution of \eqref{eq:nonconf_sourceproblem} for right-hand side $f\in L^2(\Omega)$.
						The proof is verbatim to the proof of \cref{thm:BabuskaOsborn} with the results in \cref{cor:convergence_u_nc} 
						instead of \cref{lem:ThtoT}  and so the details are omitted here. 
					\end{proof}

\section{Numerical experiments in $2$D}\label{sec:experiments}
	This section presents numerical evidence for the superiority of the new GLB in \cref{thm:GLB} over 
	  \eqref{eq:GLB_CR_Morley} and the asymptotic convergence rates from \cref{thm:BabuskaOsborn} in $2$D.
		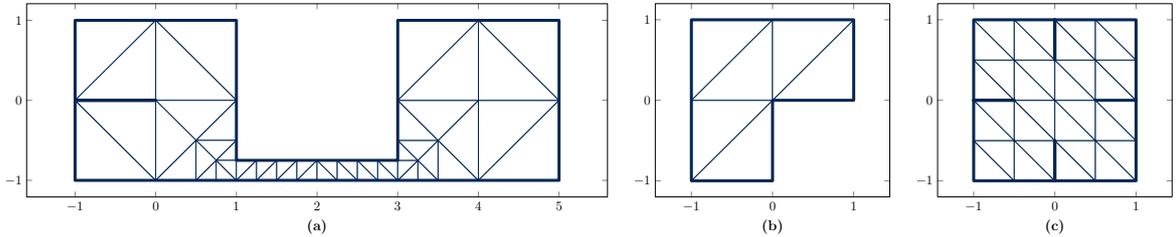
\begin{figure}[htb]
		\begin{center}
			\scalebox{0.45}{\begin{tikzpicture}
\definecolor{newteal}{RGB}{27,158,119}
\definecolor{neworange}{RGB}{217,95,2}
\definecolor{newpurple}{RGB}{117,112,179}
\definecolor{newpink}{RGB}{231,41,138}
\definecolor{newgreen}{RGB}{102,166,30}
\definecolor{newyellow}{RGB}{230,171,2}
\definecolor{newred}{RGB}{163, 0, 0}
\definecolor{imblue}{rgb}{0,0.2157,0.4235}
 \definecolor{newblue}{rgb}{0,0.15,0.33}
 
\colorlet{col0}{newblue}
 \colorlet{col1}{newgreen}
 \colorlet{col2}{newyellow}
 \colorlet{col3}{neworange}
 \colorlet{col4}{newred}
		  \definecolor{imblue}{RGB}{49,104,142}
		  \pgfplotsset{triangulationplot/.style={
		    xlabel=\large \textbf{(a)},
		    axis equal,
		    tick label style={font=\normalsize}
		  }}
		  \pgfplotsset{patchstyle/.style={
		    patch,
		    white,
		    faceted color=newblue,
		    line width=.4pt
		  }}
		  \begin{axis}[triangulationplot,xtick={-1,0,...,5},ytick={-1,0,1}, height=\axisdefaultheight , width=2.2*\axisdefaultwidth]
		    \addplot[patchstyle,
		                patch table={Pictures/data/Dumbbell_1/n4e_0.dat}]
		    				table{Pictures/data/Dumbbell_1/c4n_0.dat};
		   \draw (axis cs: 0,0) node[shape=coordinate] (P1) {P1};
		   	 \draw (axis cs: -1,0) node[shape=coordinate] (P2) {P2};
		   	 \draw (axis cs: -1,-1) node[shape=coordinate] (P3) {P3};
			  \draw (axis cs: 5,-1) node[shape=coordinate] (P5) {P5};
			   \draw (axis cs: 5,1) node[shape=coordinate] (P6) {P6};
			    \draw (axis cs: 3,1) node[shape=coordinate] (P7) {P7};
			    \draw (axis cs: 3,-0.75) node[shape=coordinate] (P8) {P8};
			    \draw (axis cs: 1,-0.75) node[shape=coordinate] (P9) {P9};
			    \draw (axis cs: 1,1) node[shape=coordinate] (P10) {P10};
			    \draw (axis cs: -1,1) node[shape=coordinate] (P11) {P11};
			   \draw[color=newblue, line width=2.5pt]  (P1)--(P2)--(P3)--(P5)--(P6)--(P7)--(P8)--(P9)--(P10)--(P11)--(P2)--(P1)--cycle; 
		 \end{axis}
	
\end{tikzpicture}}\hspace{1mm}
			\scalebox{0.45}{\begin{tikzpicture}

 \definecolor{newblue}{rgb}{0,0.15,0.33}


		  \pgfplotsset{triangulationplot/.style={
		     xlabel=\large \textbf{(b)},
		    axis equal,
		    tick label style={font=\normalsize}
		  }}
		  \pgfplotsset{patchstyle/.style={
		    patch,
		    white,
		    faceted color=newblue,
		    line width=.4pt
		  }}
		  \begin{axis}[triangulationplot,xtick={-3,-2,...,3},ytick={-3,-2,...,3}]
		    \addplot[patchstyle,
		             patch table={Pictures/data/Lshape/n4e_0.dat}]
		    				table{Pictures/data/Lshape/c4n_0.dat};
		   	 \draw (axis cs: 0,0) node[shape=coordinate] (P1) {P1};
		   	 \draw (axis cs: 1,0) node[shape=coordinate] (P2) {P2};
		   	 \draw (axis cs: 1,1) node[shape=coordinate] (P3) {P3};
			  \draw (axis cs: -1,1) node[shape=coordinate] (P4) {P4};
			  \draw (axis cs: -1,-1) node[shape=coordinate] (P5) {P5};
			   \draw (axis cs: 0,-1) node[shape=coordinate] (P6) {P6};
			   \draw[color=newblue, line width=2.5pt]  (P1)--(P2)--(P3)--(P4)--(P5)--(P6)--(P1)--cycle; 
		  \end{axis}  
\end{tikzpicture}}\hspace{1mm}			
			\scalebox{0.45}{\begin{tikzpicture}
 	\definecolor{newblue}{rgb}{0,0.15,0.33}
 
		  \pgfplotsset{triangulationplot/.style={
		    xlabel=\large \textbf{(c)},
		    axis equal,
		    tick label style={font=\normalsize}
		  }}
		  \pgfplotsset{patchstyle/.style={
		    patch,
		    white,
		    faceted color=newblue,
		    line width=.4pt
		  }}
		  \begin{axis}[triangulationplot,xtick={-3,-2,...,3},ytick={-3,-2,...,3}]
		    \addplot[patchstyle,
		                patch table={Pictures/data/PerturbSlit/n4e_0.dat}]
		    				table{Pictures/data/PerturbSlit/c4n_0.dat};
		    				
		    	 \draw  (axis cs: 0.5,0) node[shape=coordinate] (P1) {P1};			
		   	 \draw (axis cs: 1,0) node[shape=coordinate] (P2) {P2};
		   	 \draw (axis cs: 1,1) node[shape=coordinate] (P3) {P3};
		   	  \draw  (axis cs: 0,1) node[shape=coordinate] (P4) {P4};
		   	  \draw  (axis cs: 0,0.5) node[shape=coordinate] (P5) {P5};
			  \draw (axis cs: -1,1) node[shape=coordinate] (P6) {P6};
			   \draw  (axis cs: -1,0) node[shape=coordinate] (P7) {P7};
			     \draw  (axis cs: -0.5,0) node[shape=coordinate] (P8) {P8};
			  \draw (axis cs: -1,-1) node[shape=coordinate] (P9) {P9};
			   \draw  (axis cs: 0,-1) node[shape=coordinate] (P10) {P10};
			  \draw  (axis cs: 0,-0.5) node[shape=coordinate] (P11) {P11};
			   \draw (axis cs: 1,-1) node[shape=coordinate] (P12) {P12};			   
		   \draw[color=newblue, line width=2.5pt] (P1)--(P2)--(P3)--(P4)--(P5)--(P4)--(P6)--(P7)--(P8)--(P7)--(P9)--(P10)--(P11)--(P10)--(P12)--(P2);
		 \end{axis}
	
\end{tikzpicture}}
			\caption{Initial triangulation $\mathcal{T}_0$ of dumbbell-slit (a), L-shaped (b),  
							and four-slit domain (c).}\label{fig:StartTriangulation}
		\end{center}
		\end{figure}
		\subsection{Implementation}
			The implementation in this paper is realized in MATLAB based on the data structure and assembling from  
			\cite[Sec.~7.8]{CB18}. 
			\cref{fig:StartTriangulation} displays the initial triangulations $\mathcal{T}_0$ for the numerical experiments 
			below. 
			The $k$-th eigenpair $(\lambda_h(k),\boldsymbol{u_h}(k))\in\mathbb{R}^+\times \boldsymbol{V_h}$ of 
			\eqref{eq:dis_EVP_alt} with 
			$\boldsymbol{u_h}(k)=(u_{\mathrm{pw}},u_{M}(k))\in P_2(\mathcal{T})\times M(\mathcal{T})$ and for 
			comparison the post-processed Morley bound $\textup{GLB}(k)$ in \eqref{eq:GLB_CR_Morley} from 
			\cite{CGed14} 
			are computed with the MATLAB routine \texttt{eigs} exactly; the 
			termination and round-off errors are small and neglected for simplicity. 	
			Under the condition $\kappa_2^2h_{\max}^4\lambda_h(k)\le 1$, \cref{thm:GLB} guarantees 
			 $\lambda_h(k)\le\lambda_k$ for the $k$-th eigenvalue $\lambda_k$ of \eqref{eq:contEVP}. Otherwise 
			(if $1< \kappa_2^2h_{\max}^4\lambda_h(k)$ on a coarser mesh) the 
			 value $\lambda_h(k)$ is set zero, but the point is that this never occurs in all the examples displayed in this paper.
			The adaptive algorithm \cite{Doerfler1996,MorinNochettoSiebert2002,CFPP14,CR16} is based on the 
			refinement indicator $\eta(T)$ defined in \eqref{eq:def_eta} below for any triangle 
			$T\in\mathcal{T}$.
			Given the discrete solution $\big(\lambda_h, \boldsymbol{u_h}\big)\in \mathbb{R}^+\times \boldsymbol{V_h}$ 
			of \eqref{eq:dis_EVP_alt} of number ${k}$, $\lambda_h:=\lambda_h(k)$,  
			the local contribution $\eta^2(T)=(\eta(T))^2$ for any $T\in\mathcal{T}$ with area $|T|$ and set of edges
			$\mathcal{F}(T)$ solely depends on the Morley component 
			$u_{M}\in M(\mathcal{T})$ of 
			$\boldsymbol{u_h}=(u_{\mathrm{pw}},u_{M})\in  \boldsymbol{V_h}$ and reads
			\begin{align}
				\eta^2(T) = |T|^{2}\Vert \lambda_h u_{M} \Vert^2_{L^2(T)}
									+|T|^{1/2}\sum_{F\in\mathcal{F}(T)}\Vert [{D}^{2} u_{M}]_F\times \nu_F\Vert^2_{L^2(F)} 
									\label{eq:def_eta}
			\end{align}
			with the tangential components $[D^2v]_F\times \nu_F$ of the jump $[D^2v]_F$ along any edge $F\in\mathcal{F}$ and  
			the (piecewise) Hessian $D^2$.
			The respective convergence history plots in \cref{fig:adaptive_Dumbbell_intro} and \cref{fig:adaptive_Lshape_Slit} 
			display the difference $\lambda_k-\lambda_h(k)$ and $\lambda_k-\textup{GLB}(k)$ of the exact eigenvalue $\lambda_k$
			and guaranteed lower bounds $\lambda_h(k)$ and $\textup{GLB}(k)$ for uniform  red-mesh-refinement 
			$\theta=1$ (solid line and filled markers) and 
			adaptive mesh-refinement  with a bulk parameter $\theta=0.5$  in the D\"orfler marking algorithm and newest vertex bisection (dashed line and striped markers) 
			plotted against the number of triangles $|\mathcal{T}|$. 
			The computational bound  $\kappa_2= 0.07353 $ from \cite{LiaoShuLiu2019} 
			improves the analytical bound from \cite{CGal14} and the effect is investigated in \cref{fig:adaptive_Lshape_Slit}.a  with a comparison 
			between the bounds computed with 
			$\kappa_2=0.07353$ from \cite{LiaoShuLiu2019} (line color orange/blue)  
			and $\kappa_2=0.25746$ from \cite{CGal14} (line color red/green). 
			On uniform meshes $\textup{GLB}(k)$ (line color blue) and $\lambda_h(k)$ (line color orange) coincide by 
			\cref{lem:compare} and are visible in orange only in \cref{fig:adaptive_Dumbbell_intro} in the introduction 
			and in \cref{fig:adaptive_Lshape_Slit} below. 			
				\subsection{Dumbbell-slit domain}\label{sec:Dumbbellslit}
				The principal and fourth eigenvalue  $\lambda_1= 80.93261350$ and $\lambda_4=386.80177939$ 
				on the non-convex dumbbell domain with a slit $\Omega:=(-1,1)\times(-1,5)\setminus ([0,1)\times \{0\}\cup [1,3]\times [-0.75,1])$ 
				of \cref{fig:StartTriangulation}.a  are approximated with Bogner-Fox-Schmidt and Aitken extrapolation as in \cite{CGal14}. 
				\cref{fig:adaptive_Dumbbell_intro} has been discussed in the introduction as an example where uniform mesh-refining leads to a better 
				GLB from \eqref{eq:GLB_CR_Morley} than adaptive refinement. 
				The complex geometry suggests a large computational pre-asymptotic regime, but 
				the new method converges systematically even for course triangulations.  
				The example provides striking numerical evidence for the superiority of the adaptive version of 
				the extra-stabilized Morley eigensolver.
				
				\subsection{L-shaped domain}\label{sec:Lshape}
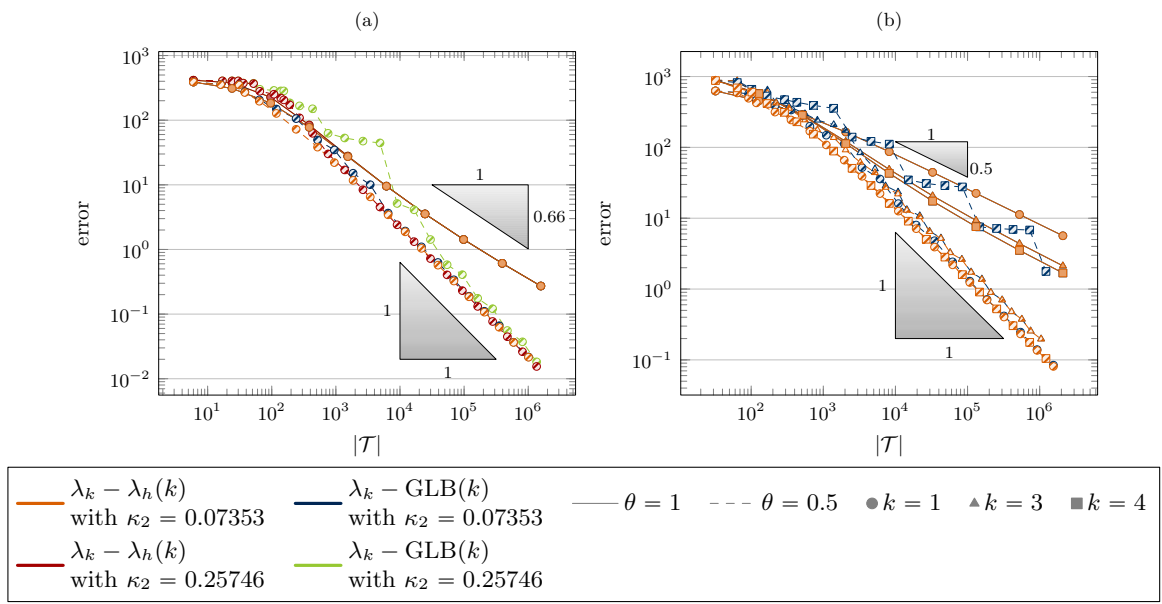
\begin{figure}[htb]
				\centering
					\scalebox{0.8}{\begin{tikzpicture}


 \definecolor{newteal}{RGB}{27,158,119}
\definecolor{neworange}{RGB}{217,95,2}
\definecolor{newpurple}{RGB}{117,112,179}
\definecolor{newpink}{RGB}{231,41,138}
\definecolor{newgreen}{RGB}{102,166,30}
\definecolor{newyellow}{RGB}{230,171,2}
\definecolor{newred}{RGB}{163, 0, 0}
\definecolor{imblue}{rgb}{0,0.2157,0.4235}
 \definecolor{newblue}{rgb}{0,0.15,0.33}
 \definecolor{newgreen2}{RGB}{150,200,50}
 
\colorlet{col0}{newblue}
 \colorlet{col1}{newgreen2}
 \colorlet{col2}{newyellow}
 \colorlet{col3}{neworange}
 \colorlet{col4}{newred} 
 
  \pgfplotsset{convergenceplot/.style={
    title ={(a)},
    xlabel=\small{$|\mathcal{T}|$},
    ylabel=\small{error },
    legend style={legend pos=outer north east},
    ymajorgrids=true
  }}
  \pgfplotsset{cases/.style={
    nodes near coords,point meta=explicit symbolic
  }}

\begin{loglogaxis}[convergenceplot]{
			
%

%
 \addplot [forget plot, col1, solid, mark=*, every mark/.append style= {solid,scale=0.9,fill=col1!60!white }]
   			table[x=nElem, y=lambda-GLB_M]{Pictures/data/BiLaplace_Adap/BiLaplace_Lshape_theta_1_lambda_1_error.dat};		
    \addplot [forget plot, col1, dashed, mark=*, every mark/.append style= {solid,scale=0.9,fill=col1!60!white,  pattern=myLines,hatchColor=col1 }]
   			table[x=nElem, y=lambda-GLB_M]{Pictures/data/BiLaplace_Adap/BiLaplace_Lshape_theta_0.5_lambda_1_error.dat};		

\addplot [forget plot, col4, solid, mark=*, every mark/.append style= {solid,scale=0.9,fill=col4!60!white }]
   			table[x=nElem, y=lambda-GLB_Sk]{Pictures/data/BiLaplace_Adap/BiLaplace_Lshape_theta_1_lambda_1_error.dat};	
    \addplot [forget plot, col4, dashed, mark=*, every mark/.append style= {solid,scale=0.9,fill=col4!60!white,  pattern=myLines,hatchColor=col4 }]
   			table[x=nElem, y=lambda-GLB_Sk]{Pictures/data/BiLaplace_Adap/BiLaplace_Lshape_theta_0.5_lambda_1_error.dat};	
   				
 \addplot [forget plot, col0, solid, mark=*, every mark/.append style= {solid,scale=0.9,fill=col0!60!white }]
   			table[x=nrElem, y=lambda-GLB_M]{Pictures/data/BiLaplace_newKappa/BiLaplace_Lshape_theta_1_lambda_1error.dat};		
    \addplot [forget plot, col0, dashed, mark=*, every mark/.append style= {solid,scale=0.9,fill=col0!60!white,  pattern=myLines,hatchColor=col0 }]
   			table[x=nrElem, y=lambda-GLB_M]{Pictures/data/BiLaplace_newKappa/BiLaplace_Lshape_theta_0.5_lambda_1error.dat};		

\addplot [forget plot, col3, solid, mark=*, every mark/.append style= {solid,scale=0.9,fill=col3!60!white }]
   			table[x=nrElem, y=lambda-GLB_Sk]{Pictures/data/BiLaplace_newKappa/BiLaplace_Lshape_theta_1_lambda_1error.dat};	
    \addplot [forget plot, col3, dashed, mark=*, every mark/.append style= {solid,scale=0.9,fill=col3!60!white,  pattern=myLines,hatchColor=col3 }]
   			table[x=nrElem, y=lambda-GLB_Sk]{Pictures/data/BiLaplace_newKappa/BiLaplace_Lshape_theta_0.5_lambda_1error.dat};

%
    \shade[top color=black!5,bottom color=black!30]
            (axis cs: 1.00e+04,6.32e-01)
         -- (axis cs: 1.00e+04,2.00e-02)
         -- (axis cs: 3.16e+05,2.00e-02)
         -- cycle;
    \draw   (axis cs: 1.00e+04,6.32e-01)
         -- (axis cs: 1.00e+04,2.00e-02) node [midway,left] {\scriptsize \(1\)}
         -- (axis cs: 3.16e+05,2.00e-02) node [midway,below] {\scriptsize \(1\)}
         -- cycle;
   
    \shade[top color=black!5,bottom color=black!30]
            (axis cs: 1.00e+06,1.02e+00)
         -- (axis cs: 1.00e+06,1.00e+01)
         -- (axis cs: 3.16e+04,1.00e+01)
         -- cycle;
    \draw   (axis cs: 1.00e+06,1.02e+00)
         -- (axis cs: 1.00e+06,1.00e+01) node [midway,right,xshift=-.5mm] {\scriptsize \(0.66\)}
         -- (axis cs: 3.16e+04,1.00e+01) node [midway,above,yshift=-.5mm] {\scriptsize \(1\)}
         -- cycle;
 }
\end{loglogaxis}
\end{tikzpicture}}
					\scalebox{0.8}{\begin{tikzpicture}

 \definecolor{newteal}{RGB}{27,158,119}
\definecolor{neworange}{RGB}{217,95,2}
\definecolor{newpurple}{RGB}{117,112,179}
\definecolor{newpink}{RGB}{231,41,138}
\definecolor{newgreen}{RGB}{102,166,30}
\definecolor{newyellow}{RGB}{230,171,2}
\definecolor{newred}{RGB}{109, 0, 0}
\definecolor{imblue}{rgb}{0,0.2157,0.4235}
\definecolor{newgreen2}{RGB}{150,200,50}

 
 \colorlet{col0}{imblue}
 \colorlet{col1}{newgreen2}
 \colorlet{col2}{newyellow}
 \colorlet{col3}{neworange}
 \colorlet{col4}{newred}
 
  \pgfplotsset{convergenceplot/.style={
  	title ={(b)},
    xlabel=\small{$|\mathcal{T}|$},
    ylabel=\small{error },
    legend style={legend pos=outer north east},
    ymajorgrids=true
  }}
  \pgfplotsset{cases/.style={
    nodes near coords,point meta=explicit symbolic
  }}

\begin{loglogaxis}[convergenceplot]{

    \addplot [forget plot, col0, dashed, mark=*, every mark/.append style= {solid,scale=0.9,fill=col0!60!white,  pattern=myLines,hatchColor=col0 }]
   			table[x=nrElem, y=lambda-GLB_M]{Pictures/data/BiLaplace_newKappa/BiLaplace_PerturbSlit2_theta_0.5_lambda_1error.dat};		
	\addplot [forget plot, col0, solid, mark=*, every mark/.append style= {solid,scale=0.9,fill=col4!60!white }]
   			table[x=nrElem, y=lambda-GLB_M]{Pictures/data/BiLaplace_newKappa/BiLaplace_PerturbSlit2_theta_1_lambda_1error.dat};		

\addplot [forget plot, col3, solid, mark=*, every mark/.append style= {solid,scale=0.9,fill=col3!60!white }]
   			table[x=nrElem, y=lambda-GLB_Sk]{Pictures/data/BiLaplace_newKappa/BiLaplace_PerturbSlit2_theta_1_lambda_1error.dat};	
    \addplot [forget plot, col3, dashed, mark=*, every mark/.append style= {solid,scale=0.9,fill=col3!60!white,  pattern=myLines,hatchColor=col3 }]
   			table[x=nrElem, y=lambda-GLB_Sk]{Pictures/data/BiLaplace_newKappa/BiLaplace_PerturbSlit2_theta_0.5_lambda_1error.dat};	
    
  \addplot [forget plot, col0, solid, mark=triangle*, every mark/.append style= {solid,scale=0.9,fill=col0!60!white }]
   			table[x=nrElem, y=lambda-GLB_M]{Pictures/data/BiLaplace_newKappa/BiLaplace_PerturbSlit2_theta_1_lambda_3error.dat};		
    \addplot [forget plot, col0, dashed, mark=triangle*, every mark/.append style= {solid,scale=0.9,fill=col0!60!white,  pattern=myLines,hatchColor=col0 }]
   			table[x=nrElem, y=lambda-GLB_M]{Pictures/data/BiLaplace_newKappa/BiLaplace_PerturbSlit2_theta_0.5_lambda_3error.dat};		

\addplot [forget plot, col3, solid, mark=triangle*, every mark/.append style= {solid,scale=0.9,fill=col3!60!white }]
   			table[x=nrElem, y=lambda-GLB_Sk]{Pictures/data/BiLaplace_newKappa/BiLaplace_PerturbSlit2_theta_1_lambda_3error.dat};	
    \addplot [forget plot, col3, dashed, mark=triangle*, every mark/.append style= {solid,scale=0.9,fill=col3!60!white,  pattern=myLines,hatchColor=col3 }]
   			table[x=nrElem, y=lambda-GLB_Sk]{Pictures/data/BiLaplace_newKappa/BiLaplace_PerturbSlit2_theta_0.5_lambda_3error.dat};	
   			
   \addplot [forget plot, col0, solid, mark=square*, every mark/.append style= {solid,scale=0.9,fill=col0!60!white }]
   			table[x=nrElem, y=lambda-GLB_M]{Pictures/data/BiLaplace_newKappa/BiLaplace_PerturbSlit2_theta_1_lambda_4error.dat};		
    \addplot [forget plot, col0, dashed, mark=square*, every mark/.append style= {solid,scale=0.9,fill=col0!60!white,  pattern=myLines,hatchColor=col0 }]
   			table[x=nrElem, y=lambda-GLB_M]{Pictures/data/BiLaplace_newKappa/BiLaplace_PerturbSlit2_theta_0.5_lambda_4error.dat};		

\addplot [forget plot, col3, solid, mark=square*, every mark/.append style= {solid,scale=0.9,fill=col3!60!white }]
   			table[x=nrElem, y=lambda-GLB_Sk]{Pictures/data/BiLaplace_newKappa/BiLaplace_PerturbSlit2_theta_1_lambda_4error.dat};	
    \addplot [forget plot, col3, dashed, mark=square*, every mark/.append style= {solid,scale=0.9,fill=col3!60!white,  pattern=myLines,hatchColor=col3 }]
   			table[x=nrElem, y=lambda-GLB_Sk]{Pictures/data/BiLaplace_newKappa/BiLaplace_PerturbSlit2_theta_0.5_lambda_4error.dat};	  			

    \shade[top color=black!5,bottom color=black!30]
            (axis cs: 1.00e+04,6.32e+00)
         -- (axis cs: 1.00e+04,2.00e-01)
         -- (axis cs: 3.16e+05,2.00e-01)
         -- cycle;
    \draw   (axis cs: 1.00e+04,6.32e+00)
         -- (axis cs: 1.00e+04,2.00e-01) node [midway,left] {\scriptsize \(1\)}
         -- (axis cs: 3.16e+05,2.00e-01) node [midway,below] {\scriptsize \(1\)}
         -- cycle;

    \shade[top color=black!5,bottom color=black!30]
            (axis cs: 1.00e+05,3.79e+01)
         -- (axis cs: 1.00e+05,1.20e+02)
         -- (axis cs: 1.00e+04,1.20e+02)
         -- cycle;
    \draw   (axis cs: 1.00e+05,3.79e+01)
         -- (axis cs: 1.00e+05,1.20e+02) node [midway,right,yshift=-1.5mm,xshift=-1mm] {\scriptsize \(0.5\)}
         -- (axis cs: 1.00e+04,1.20e+02) node [midway,above,yshift=-1mm] {\scriptsize \(1\)}
         -- cycle;       
 }
\end{loglogaxis}
\end{tikzpicture}}
					\begin{tikzpicture}
\definecolor{newteal}{RGB}{27,158,119}
\definecolor{neworange}{RGB}{217,95,2}
\definecolor{newpurple}{RGB}{117,112,179}
\definecolor{newpink}{RGB}{231,41,138}
\definecolor{newgreen}{RGB}{102,166,30}
\definecolor{newyellow}{RGB}{230,171,2}
\definecolor{newred}{RGB}{163, 0, 0}
\definecolor{imblue}{rgb}{0,0.2157,0.4235}
 \definecolor{newblue}{rgb}{0,0.15,0.33}
 \definecolor{newgreen2}{RGB}{150,200,50}
 
\colorlet{col0}{newblue}
 \colorlet{col1}{newgreen2}
 \colorlet{col2}{newyellow}
 \colorlet{col3}{neworange}
 \colorlet{col4}{newred}
 
\begin{axis}[%
    legend columns=7,
    scale only axis,width=1mm, 
    hide axis,
	xmin=1,
    legend style={legend cell align=left, cells={align=left}},
	legend style={/tikz/every even column/.append style={column sep=0.3cm}}
    ]
\addlegendimage{col3, solid, very thick}
			\addlegendentry{$\lambda_k-\lambda_h(k)$ \\ with $\kappa_2=0.07353$  }  
\addlegendimage{col0,solid, very thick}
			\addlegendentry{$\lambda_k-\textup{GLB}(k)$ \\ with $\kappa_2=0.07353$ }	
\addlegendimage{ solid, gray}
			\addlegendentry{$\theta=1$}	
\addlegendimage{ dashed, gray}
			\addlegendentry{$\theta=0.5$}
	\addlegendimage{ mark=*,gray, only marks}
					\addlegendentry{$k=1$}
\addlegendimage{ mark=triangle*,gray, only marks}
					\addlegendentry{$k=3$}	
\addlegendimage{ mark=square*,gray, only marks}
					\addlegendentry{$k=4$}					
\addlegendimage{col4, solid, very thick}
			\addlegendentry{$\lambda_k-\lambda_h(k)$ \\ with $\kappa_2=0.25746$}  
\addlegendimage{col1,solid, very thick}
			\addlegendentry{$\lambda_k-\textup{GLB}(k)$\\ with $\kappa_2=0.25746$ }	
 \addplot [forget plot,  solid]
   			table[x=a, y=b]{Pictures/data/dummy.dat};				
  		
\end{axis}
\end{tikzpicture}
					\caption{Comparison of the distance between $\lambda_k$ and $\lambda_h(k)$ 
								(resp. $\textup{GLB}(k)$) computed on uniform ($\theta=1$, solid) and adaptive 
								($\theta=0.5$, dashed) meshes of the L-shaped domain for $k=1$ in (a) 
								and the four-slit domain for $k=1,3,4$ in (b).}
								\label{fig:adaptive_Lshape_Slit}
				\end{figure}

				The principal eigenvalue $\lambda_1= 418.97504246688220$ 
				on the  non-convex L-shaped domain  $\Omega:=(-1,1)^2\setminus [0,1)\times (-1,0]$ of \cref{fig:StartTriangulation}.b 
				is approximated in \cite{CGal14}. 
				The associated eigenfunction in $H^2_0(\Omega)\setminus H^3(\Omega)$ results  in the reduced empirical 
				convergence rate $0.66$ for uniform mesh-refinement in \cref{fig:adaptive_Lshape_Slit}.a. 
				The adaptive mesh-refinement  with \eqref{eq:def_eta} 
				allows to recover the optimal convergence rate one (with respect to the number of triangles $|\mathcal{T}|$ in the 	
				triangulation $\mathcal{T}$) for $\lambda_h(1)$.
				The choice of $\kappa_2=0.07353$ (line color blue) instead of $\kappa_2=0.25746$ (line color green) 
				improves the  guaranteed lower bound $\textup{GLB}(k)$ significantly. The bound computed with $\kappa_2=0.25746$   	
				suffers from the involvement of $h_{\max}$ visible in form of steps, 
				while the choice $\kappa_2=0.07353$ leads to a straight line in the convergence history plot. 	
				Undisplayed experiments on graded meshes \cite{CB18} of the L-shaped domain, e.g., with grading parameter  
				$\beta=10/7$, recover the optimal convergence rates  and confirm  
				\cref{lem:compare} as well.		
				\subsection{Four-slit domain}\label{sec:fourslit}
				The principal eigenvalue  $\lambda(1)=830.21478777$ and 
				the double eigenvalue  
				$\lambda(3)=1125.1279=\lambda(4)$ on the four-slit domain 
				$\Omega:=(-1,1)^2\setminus \big([0,0.5)\times \{0\}\cup [0,-0.5)\times \{0\} \cup \{0\}\times [0,0.5)\cup \{0\}
				\times[0,-0.5)\big)$ of \cref{fig:StartTriangulation}.c
				are approximated as in  \cite{CGal14}. 
				The associated eigenfunctions on the non-convex  domain seem to belong to  
				$H^2_0(\Omega)\setminus H^3(\Omega)$ because uniform mesh-refinement leads to the reduced 
				convergence rates $0.5$ for the first and $0.55$ for the third and fourth in \cref{fig:adaptive_Lshape_Slit}.b. 
				The AFEM algorithm with bulk parameter $\theta=0.5$ driven by the estimator \eqref{eq:def_eta} 
				allows to recover the optimal convergence rate one. 
				 The $\textup{GLB}$ in \cref{fig:adaptive_Lshape_Slit}.b  are computed with $\kappa_2=0.07353$. 
				Undisplayed comparison with  $\kappa_2=0.25746$ lead to worse $\textup{GLB}$. 
				A clustering adaptive algorithm as in \cite{Gal15_cluster} was not necessary for the double eigenvalue $\lambda_3=\lambda_4$.
				
			\subsection{Comments and Conclusions}\label{sec:conclusion}
				The empirical observations of the numerical experiments 
				in \cref{sec:Lshape}--\ref{sec:fourslit} show:			
				\begin{enumerate}[wide,label=(\roman*),ref=\roman*,itemsep=0.1em]
				\item  
					All experiments confirm the a priori convergence rates of  \cref{thm:BabuskaOsborn}.
					The empirical convergence rate depends only on the smoothness of the approximated eigenfunction. For instance 
					\cref{fig:adaptive_Dumbbell_intro} displays for uniform refinement the optimal convergence rate one for the principal eigenvalue 
					 despite the reduced empirical convergence rate for the 
					fourth eigenvalue. 		
				\item 
					\cref{thm:BabuskaOsborn} predicts a convergence for  a sufficiently  fine initial mesh.  
					In all examples the convergence rate is visible even for  moderately fine triangulations, so this 
					restriction does not affect the numerical examples much. 
				\item 
					If the condition on the mesh-size is satisfied, the method \eqref{eq:dis_EVP_alt} 
					provides indeed guaranteed lower eigenvalue bounds in all 
					numerical experiments and so confirms \cref{thm:GLB}. 	
				 \item\label{item:constant} 
				The constant $\kappa_2=0.07353$ from \cite{LiaoShuLiu2019} leads to a significant improvement of the 
				 known bound \eqref{eq:GLB_CR_Morley} in examples with adaptive mesh-refinement.
				 \item\label{item:ratio}
				 The (undisplayed) improvement factor 
				 $q:=(\lambda_k-\lambda_h(k))/(\lambda_k-\textup{GLB}(k))$  was computed with 
				 $\kappa_2=0.07353$ on the adaptive triangulations. 
				 For the principal eigenvalue of the L-shaped and four-slit domain the improvement with the new method is marginal and 
				 the ratio $q$ oscillates between $0.6$ and $1$. 
				  In the remaining examples the improvement is more significant. For the fourth eigenvalue of the four-slit domain 
				  the ratio $q$ oscillates between $0.05$ and $0.25$ for triangulations with more than $4500$ triangles. 
				  For the first (resp. fourth) eigenvalue of the dumbbell-slit domain the ratio $q$ decreases from $0.16$ (resp. $0.15$) 
				  to $0.0008$ (resp. $0.005$) for triangulations with more than $600$ (resp. $3800$) triangles. 
				 \item The new method increases the number of degrees of freedom by a factor four in the $2$D numerical benchmarks. The equivalent rational problem \eqref{eq:disEVP} from \cref{lem:AquivProb} could be efficiently addressed by a Newton scheme, so the final comparison is beyond this paper. As it stands, the new method is favourable at least for the examples in 
				  \cref{sec:Dumbbellslit} and the fourth eigenvalue in \cref{sec:fourslit}. 
				 \item 
				 	The adaptive algorithm driven by the estimator \eqref{eq:def_eta} recovers the optimal empirical convergence 
				 	rates in all examples for the extra stabilised method. 
				 	The analysis of optimal convergence rates and further details of the proposed adaptive algorithm shall appear 
				 	 in \cite{CP_Part2}.
				 \item
				 	The overall conclusion of the numerical experiments is that there exist examples, where the post-processed GLB 
				 	\eqref{eq:GLB_CR_Morley} may fail completely for localized triangulations. In contrast, the new scheme is compatible with adaptive 
				 	mesh-refining and leads to GLB that cannot be reached with \eqref{eq:GLB_CR_Morley}.
				\end{enumerate} 
\FloatBarrier
	
\paragraph{Acknowledgements.}
The authors thank the referees, e.g., for the suggestion of the constant $\kappa_2< 0.07353 $ from \cite{LiaoShuLiu2019} in the numerical experiments and the improvement factor $q$ that led to the comparison in \cref{sec:Lshape} and  \cref{sec:conclusion}.\ref{item:constant}--\ref{item:ratio}. This work has been supported by the Deutsche Forschungsgemeinschaft (DFG) in the Priority Program 1748 ‘Reliable simulation techniques in solid mechanics. Development of non-standard discretization methods, mechanical and mathematical analysis’ under the project CA 151/22-2. The second author is supported by the Berlin Mathematical School.
\footnotesize{
\bibliographystyle{alpha}
\bibliography{BibGLB}
}
\end{document}


\maketitle
	\appendix
	
\section{The Worsey-Farin FEM in $3$D}\label{sec:Appendix_CT3D} 
	This section aims at an elementary self-contained introduction to the Worsey-Farin ($\text{WF}$)  FEM and 
	the analysis of the associated nodal basis functions to generalize \cite[Thm.~6.1.3]{Ciarlet78}. 
	The Clough-Tocher finite element has been proposed for any space dimension $n\ge 3$ in \cite{WF87}. 
	More than two decades later, \cite{Sor09} observed  additional constraints on the choice of the subtriangulation 
	of the macro element for $n\ge 4$ necessary for $C^1$ conformity, i.e., \cite{WF87} is wrong for $n\ge4$. 
	It appears still an open problem, whether all conditions for $n=4$ can be satisfied simultaneously \cite[p.42, l.40ff]{Sor09}.
	This underlines that the details are technical and explains the restriction to $n=3$ in this paper.	
	This supplement summarizes the necessities on Bernstein polynomials  \cite{Boor1987} 
	and the 
	\cref{alg:CT3D} from \cite{WF87} for $n=3$ for a general audience.

	\subsection{Bernstein polynomials in a simplex}
			Given $m\in\{1,2,3,4\}$ points $P_1,P_2,\dots, P_m\in\mathbb{R}^3$ with linearly independent differences $P_2-P_1,\dots,P_m-P_1$, 
			their convex hull 
			$T=\textup{conv}\{P_1,\dots,P_m\}$ is an $(m-1)$-simplex (of 
			positive $(m-1)$-dimensional Hausdorff measure). 
			The set of vertices $\mathcal{V}(T)=\{ P_1,\dots, P_m\}$ of $T$ is the set of extremal 
			points and so 
			uniquely defined. We fix an ordering and 
			so identify the $(m-1)$-simplex $T$ with an ordered list $(P_1,\dots,P_m)$ of vertices. 
			We will encounter a finite union of 
			simplices and 
			decouplings of $(m-1)$-simplices and need to 
			keep track of the vertex set. Below $(P_1,\dots,P_m)$ will be replaced by $(P_{\sigma(1)},\dots, P_{\sigma(m)})$ for global 
			numbers $\sigma(1),\dots,\sigma(m)$ in a finite list of all 
			vertices $\mathcal{V}$ in a regular triangulation $\mathcal{T}$ of $\Omega\subset\mathbb{R}^3$ 
			into tetrahedra in the sense of Ciarlet. 
			Let $ \langle T\rangle:=P_1+\textup{span}\{P_2-P_1,\dots,P_m-P_1\}$ denote the unique 
			$(m-1)$-dimensional affine subspace $\langle T\rangle$ of  $\mathbb{R}^3$ that contains $T\equiv(P_1,\dots,P_m)$.  
			The barycentric coordinates $(\lambda_1,\dots,\lambda_m)$ of a point 
			$x\in \langle T\rangle$ solve the $4\times m$ linear system of equations
				\begin{align}
					\begin{pmatrix}
					1&1&\dots &1 \\P_1&P_2&\dots&P_m
					\end{pmatrix}
					\begin{pmatrix}
					\lambda_1\\\vdots\\\lambda_m
					\end{pmatrix}
					=
					\begin{pmatrix}
					1\\x
					\end{pmatrix}.\label{eq:BarycentricCoordinates4x}
				\end{align}
			The dependence on $x$ in $\lambda_j(x)=\lambda_j$ from \eqref{eq:BarycentricCoordinates4x} is 
			typically not written out explicitly for brevity. 
			It holds $x\in T$ if and only if $x\in\mathbb{R}^3$ satisfies \eqref{eq:BarycentricCoordinates4x} with 
			$\lambda_1,\dots,\lambda_m\ge 0$. 
			The relative interior of $T$ reads
			 $$\mathrm{relint}(T):=\Bigg\{\sum_{\mu=1}^m\lambda_\mu P_\mu\,\Big|\, \sum_{\mu=1}^m\lambda_\mu=1\text{ and }
			 				\lambda_1,\dots,\lambda_m> 0\Bigg\}.$$ 
			 Hence $T=\overline{\mathrm{relint}(T)}$ and the relative boundary 
			 $\mathrm{relbdy}(T)=\partial T:=T\setminus \mathrm{relint}(T)$  of 
			 $T$ reads 
			$$\partial T=\Bigg\{\sum_{\mu=1}^m\lambda_\mu P_\mu\,\Big|\, \sum_{\mu=1}^m\lambda_\mu=1\text{, }
							\lambda_1,\dots,\lambda_m\ge 0	
								\text{, and }\lambda_j=0 \text{ for at least one }1\le j\le m\Bigg\}.$$

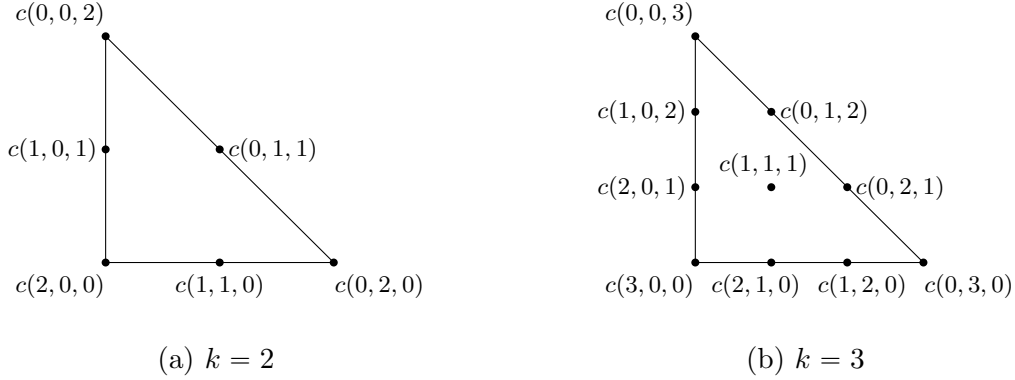
\begin{figure}
						\begin{center}
						\begin{minipage}{0.49\textwidth}
							\begin{center}
								\begin{tikzpicture}[x=30mm,y=30mm]
									\phantom{\draw[->,line width=0.7pt] (1,0)--(1,-0.3);}
						  					\draw (0,0)--(1,0)--(0,1)--cycle;
						  					\fill (0,0) circle (1.5pt);
						  					\draw(0,0)node[shape=coordinate,label={[xshift=-6mm, yshift=-6mm] $c(2,0,0)$}] (P1) {P1};
						  					\draw(0.5,0)node[shape=coordinate,label={[xshift=0mm, yshift=-6mm] $c(1,1,0)$}] (P1P2) {P1P2};
						  					\fill (1,0) circle (1.5pt);
						  					\draw(1,0)node[shape=coordinate,label={[xshift=6mm, yshift=-6mm] $c(0,2,0)$}] (P2) {P2};
						  					\fill (0.5,0.5) circle (1.5pt);
						  					\draw(0.5,0.5)node[shape=coordinate,label={[xshift=7mm, yshift=-3mm] $c(0,1,1)$}] (P2P3) {P2P3};
						  					\fill (0,1) circle (1.5pt);
						  					\draw(0,1)node[shape=coordinate,label={[xshift=-6mm, yshift=0mm] $c(0,0,2)$}] (P3) {P3};
						  					\fill (0.5,0) circle (1.5pt);
						  					\fill (0,0.5) circle (1.5pt);
						  					\draw(0,0.5)node[shape=coordinate,label={[xshift=-7mm, yshift=-3mm] $c(1,0,1)$}] (P1P3) {P1P3};
									\end{tikzpicture} \\ {(a)} $k=2$
								\end{center}
						\end{minipage}
							\begin{minipage}{0.49\textwidth}
								\begin{center}
									\begin{tikzpicture}[x=30mm,y=30mm]
										\phantom{\draw[->,line width=0.7pt] (1,0)--(1,-0.3);}
						  					\draw (0,0)--(1,0)--(0,1)--cycle;
						  					\fill (0,0) circle (1.5pt);
						  					\draw(0,0)node[shape=coordinate,label={[xshift=-6mm, yshift=-6mm] $c(3,0,0)$}] (P1) {P1};
						  					\fill (0.333,0) circle (1.5pt);
						  					\draw(0.333,0)node[shape=coordinate,label={[xshift=-2mm, yshift=-6mm] $c(2,1,0)$}] (P1P2_1) {P1P21};
						  					\fill (0.666,0) circle (1.5pt);
						  					\draw(0.666,0)node[shape=coordinate,label={[xshift=+2mm, yshift=-6mm] $c(1,2,0)$}] (P1P2) {P1P2};
						  					\fill (1,0) circle (1.5pt);
						  					\draw(1,0)node[shape=coordinate,label={[xshift=6mm, yshift=-6mm] $c(0,3,0)$}] (P2) {P2};
						  					\fill (0.666,0.333) circle (1.5pt);
						  					\draw(0.666,0.333)node[shape=coordinate,label={[xshift=7mm, yshift=-3mm] $c(0,2,1)$}] (P2P3_1) {P2P31};
						  					\fill (0.333,0.666) circle (1.5pt);
						  					\draw(0.333,0.666)node[shape=coordinate,label={[xshift=7mm, yshift=-3mm] $c(0,1,2)$}] (P2P3) {P2P3};
						  					\fill (0,1) circle (1.5pt);
						  					\draw(0,1)node[shape=coordinate,label={[xshift=-6mm, yshift=0mm] $c(0,0,3)$}] (P3) {P3};
						  					\fill (0,0.333) circle (1.5pt);
						  					\draw(0,0.333)node[shape=coordinate,label={[xshift=-7mm, yshift=-3mm] $c(2,0,1)$}] (P1P3_1) {P1P31};
						  					\fill (0,0.666) circle (1.5pt);
						  					\draw(0,0.666)node[shape=coordinate,label={[xshift=-7mm, yshift=-3mm] $c(1,0,2)$}] (P1P3) {P1P3};
						  					\draw(0.333,0.333)node[shape=coordinate,label={[xshift=-1mm, yshift=0mm] $c(1,1,1)$}] (MT) {MT};
						  					\fill (MT) circle (1.5pt);
									\end{tikzpicture}\\ {(b)} $k=3$
								\end{center}	
						\end{minipage}	
						\end{center}
							\caption{Illustration of the position associated with the ordinates $c(\alpha)$ for $\alpha\in A_k$ with $m=3$ and $k=2$ 
											(a) and $k=3$ (b).}
							\label{fig:illustration_alpha_position}	
				\end{figure}

			The space $P_k(T)$ of algebraic polynomials of total degree at most $k\in\mathbb{N}_0$ is seen as a subspace of 
			$C^\infty(T)$ and allows many representations. 
			The Bernstein polynomials form a natural basis of $P_k(T)$ with standard multi-index notation  for 
			$\alpha=(\alpha_1,\dots,\alpha_m)\in\mathbb{N}_0^m$ 
			with length $|\alpha|=\alpha_1+\dots+\alpha_m$ 
			and
				\begin{align}
					\frac{\lambda^\alpha}{\alpha!}:=\frac{\lambda_1^{\alpha_1}\lambda_2^{\alpha_2}\cdots \lambda_m^{\alpha_m}}
														{\alpha_1!\alpha_2!\cdots \alpha_m!} 
					\quad\text{at}\quad x=\sum_{\mu=1}^m\lambda_\mu P_\mu\in T.
					\label{eq:MultiIndex}
				\end{align} 
			(We follow the convention ${\lambda^\alpha}/{\alpha!}\equiv 0$ if at least one of the indices 
			$\alpha_1,\dots, \alpha_m$ of $\alpha$ is negative.)  
			If the length 
			$m\in\{1,2,3,4\}$ is fixed, abbreviate $A_k:=\{\alpha\in\mathbb{N}_0^m|\, |\alpha|=k\}\subset \mathbb{N}_0^m$.
			The Bernstein polynomials of degree $k\in\mathbb{N}_0$ are all $\lambda^\alpha/\alpha!$ from \eqref{eq:MultiIndex} for $\alpha \in A_k$.
			In fact, any (real) polynomial $f\in P_k(T)$ has unique (real) coefficients 
			$(c(\alpha)|\,\alpha\in A_k)$ called \emph{ordinates} such that 
				\begin{align}
					f
					=k!\sum_{\alpha\in A_k}c(\alpha)\frac{\lambda^\alpha}{\alpha!}\quad
					\text{ at }\quad x=\sum_{\mu=1}^m\lambda_\mu P_\mu \in T.\label{eq:BernsteinPolynomial}
				\end{align} 
			The definition \eqref{eq:BernsteinPolynomial} makes sense also in the entire affine subspace 
			$\langle T\rangle:=P_1+\textup{span}\{P_2-P_1,\dots,P_m-P_1\}\ni x$ 
			and immediately extends any $f\in P_k(T)$ to $f\in P_k(\langle T\rangle)$. 
			The ordinates of $f$ form a family  $(c(\alpha)|\,\alpha\in A_k)$ and can be arranged in an $m$-dimensional array. 
			Any ordinate $c(\alpha)$ for $\alpha\in A_k$ is associated with the position $x=\sum_{\mu=1}^m\frac{\alpha_\mu}{k}P_\mu\in T$ 
			as illustrated in \cref{fig:illustration_alpha_position} for $m=3$ and $k=2,3$.  
			If the length of $m$ is fixed (and there is no need to highlight the length $m$ of the multi-indices from the context), then 
			$e_j:=(0,\dots,0,1,0,\dots,0)\in\mathbb{N}_0^m$ denotes the $j$-th 
			canonical unit vector with the coefficients $e_j(\ell)=\delta_{j\ell}$ for $j,\ell\in\{1,\dots,m\}$.
			It cannot be overemphasized that the ordinates $(c(\alpha)|\,\alpha\in A_k)$ of the polynomial $f$ in 
			\eqref{eq:BernsteinPolynomial} are neither coefficients of a monomic tensor basis 
			$x_1^{\alpha_1}x_2^{\alpha_2}\dots x_n^{\alpha_n}$ nor the values of the polynomial in general. However, at the $m$ vertices, 
			the ordinates are the values of the polynomial $f$: 
			For any $j=1,\dots,m$, the vertex $P_j$ is identified with the multi-index $ke_j$ and $c(ke_j)=f(P_j)$ 
			follows from \eqref{eq:BernsteinPolynomial} and the barycentric coordinates $(\lambda_1,\dots,\lambda_m)=e_j$ of $P_j$.  
			
				\begin{remark}[gradient $\nabla f$]\label{rem:derivativeBernsteinBasis}
					The derivative of $f$ from \eqref{eq:BernsteinPolynomial} reads at $x=\sum_{\mu=1}^m\lambda_\mu P_\mu \in T$ 
						\begin{align}
							\nabla f 
							=k!\sum_{\alpha\in A_k}\sum_{\mu=1}^m c_j(\alpha)\frac{\lambda^{(\alpha-e_\mu)}}{(\alpha-e_\mu)!}\nabla \lambda_\mu				
							=k!\sum_{\beta\in A_{k-1}}\frac{\lambda^\beta}{\beta!}\sum_{\mu=1}^mc(\beta+e_{\mu})\nabla\lambda_\mu .
							\label{eq:derivativeBernstein}
						\end{align}
				\end{remark}
			Some examples on the evaluation of $f$ and $\nabla f$ on the boundary $\partial T$ of a triangle ($m=3$) and 
			tetrahedron ($m=4$) conclude this subsection for the sake of an illustration as well as future reference. The reader may skip the proof at 
			first reading,
			but requires details in the analysis of \cref{alg:CT3D} later in \cref{sec:HCT3D=FE}.
			\begin{example}[triangle, $k=2$, $m=3$]\label{exp:m=2=k}
				For $m=3$ consider a triangle $T$ identified with its vertices $(P_1,P_2{,}P_3)$. 
						The first edge $E_1=\textup{conv}\{P_2,P_3\}$  of the triangle $T$ lies opposite to the first 
						vertex $P_1$ and has midpoint $\textup{mid}(E_1):=P_{23}:=(P_2+P_3)/2$. 
						In 
						the two-dimensional affine space
						$\langle T\rangle = P_1+\textup{span}\{P_2-P_1,P_3-P_1\}$, the edge $E_1$ belongs to the zero set 
						$\{\lambda_1=0\}\supset E_1$ of the first barycentric 
						coordinate $\lambda_1$, where $\lambda_2+\lambda_3=1$ 
						for $\lambda_2,\lambda_3\ge 0$.
						The polynomial $\lambda^\alpha$ vanishes at $x=\lambda_2P_2+\lambda_3P_3\in E_1$,  if $\alpha_1\ge 1$ and so 
						\eqref{eq:BernsteinPolynomial} 
						reduces to 
							\begin{align}
								f=k!\sum_{\substack{\alpha\in A_k\\\alpha_1=0}}c(\alpha){\lambda^\alpha}/{\alpha!}\quad
								\text{at}\quad x=\lambda_2P_2+\lambda_3P_3\in E_1.\label{eq:exp_m=2=k}
							\end{align}  
						For $k=2$ and with $c(0,2,0)=f(P_2)$ and $c(0,0,2)=f(P_3)$, the formula \eqref{eq:exp_m=2=k} shows at 
						$ P_{23}=\textup{mid}(E_1)$ 
						that 
							\begin{align*}
								4f(P_{23}) = f(P_2)+2c(0,1,1)+f(P_3).
							\end{align*}
						The derivative $\nabla f$ of $f$ along the edge $E_1\subset \{\lambda_1=0\}$ follows from \eqref{eq:derivativeBernstein}, 
							\begin{align}
								\nabla f=k!\sum_{\substack{\alpha\in A_{k-1}\\\alpha_1=0}}\Bigg(\sum_{\mu=1}^3c(\alpha+e_{\mu})\nabla\lambda_{\mu}\Bigg)
								\frac{\lambda^\alpha}{\alpha!}\quad\text{ at }x=\lambda_2P_2+\lambda_3P_3\in E_1.\label{eq:exp_m=2=k_grad}
							\end{align}
						Since $k=2$, there are solely two summands in the sum of the formula \eqref{eq:exp_m=2=k_grad} 
						for  $(\alpha_2,\alpha_3)=(1,0)$ and $(\alpha_2,\alpha_3)=(0,1)$. Hence  
							\begin{align*}
								\nabla f=&2\big(c(1,1,0)\nabla\lambda_1+c(0,2,0)\nabla\lambda_2+c(0,1,1)\nabla\lambda_3\big)\lambda_2
								\\&+2\big(c(1,0,1)\nabla\lambda_1+c(0,1,1)\nabla\lambda_2+c(0,0,2)\nabla\lambda_3\big)\lambda_3\quad\text{ at }x\in E_1.
							\end{align*}
						This defines an affine vector function along $E_1$. 
						In conclusion, the data $f|_{E_1}$ and $\nabla f|_{E_1}$ determine the ordinates
						 $c(0,0,2),\,c(0,2,0),\,c(1,1,0),\,c(1,0,1),\,c(0,1,1)$ 
						and vice 
						versa.
			 \end{example}
			 \begin{example}[tetrahedron, $k=2$, $m=4$]\label{exp:k=2_m=4} 
				For $m=4$ consider a tetrahedron $T$ identified with its vertices $(P_1,P_2{,}P_3{,}P_4)$ in $\mathbb{R}^3$
				with edge midpoints $P_{k\ell}=(P_k+P_\ell)/2$ for $k,\ell=1,\dots, 4$ and $k\not =\ell$. 
				 The ordinates  $(c(\alpha)|\alpha\in A_2)$ define the quadratic polynomial 
					$f=2\sum_{\alpha\in A_2}c(\alpha)\lambda^\alpha/\alpha!\in P_2(T)$, which assumes,  for all $j,k,\ell=1,\dots,4$ 
					and $k\not =\ell$, the following values.
					\begin{enumerate}[label=(\alph*), ref=\alph*]
						\item\label{item:gPj} $\displaystyle f(P_j)=c(2e_j),$ 
						\item\label{item:gPkell} 
						$\displaystyle 4f(P_{k\ell})=c(2e_k)+2c(e_k+e_\ell)+c(2e_\ell)=f(P_k)+2c(e_k+e_\ell)+f(P_\ell).$
					\end{enumerate}
					\begin{proofof}\textit{(\ref{item:gPj}).}
						The first formula (\ref{item:gPj}) follows with $(\lambda_1,\dots,\lambda_4)=e_j$ at $x=P_j$  and the evaluation of the non-zero 
						$\lambda^\alpha/\alpha!$ for $\alpha=2e_j$.
					\end{proofof}
					\begin{proofof}\textit{(\ref{item:gPkell}).}
						At the edge midpoint $x=P_{k\ell}$ in (\ref{item:gPkell}), we have $(\lambda_1,\dots,\lambda_4)=(e_k+e_\ell)/2$ and 
						$\lambda^\alpha/\alpha!$ is possibly non-zero only for 
						 $\alpha=2e_k$, $e_k+e_\ell$, and $2e_\ell$.
						The polynomial $\lambda^\alpha/\alpha!$  assumes the respective values $\lambda^\alpha/\alpha!=1/8, 1/4,1/8$. 
						This leads to the formulas in (\ref{item:gPkell}). 
					\end{proofof}
			\end{example}
			\begin{example}[tetrahedron, $k=3$, $m=4$]\label{exp:k=3_m=4}
				In the notation of \cref{exp:k=3_m=4}, the ordinates $(c(\alpha)|\alpha\in A_3)$ 
					 define the cubic polynomial  
					$f=6\sum_{\alpha\in A_3}c(\alpha)\lambda^\alpha/\alpha!\in P_3(T)$ and its gradient $\nabla f\in P_2(T;\mathbb{R}^3)$.
					For all $\{j,k,\ell,m\}=\{1,2,3,4\}$ it follows 
					\begin{enumerate}[label=(\alph*), ref=\alph*]
						\item\label{item:fPj} $\displaystyle f(P_j)=c(3e_j),$ 
						\item\label{item:fPkell} $\displaystyle 8 f(P_{k\ell})=f(P_k)+3c(2e_k+e_\ell)+3c(e_k+2e_\ell)+f(P_\ell),$
						\item\label{item:nablafPj_1}$\displaystyle(P_j-P_k)\cdot\nabla f(P_j)=3\big(c(3e_j)-c(2e_j+e_k)\big),$
						\item\label{item:nablafPkell} 
							$\displaystyle \frac{4}{3}(P_j-P_{k\ell})\cdot\nabla f(P_{k\ell})
								=c(2e_\ell+e_j)+2c(e_\ell+e_k+e_j)+c(2e_k+e_j)-4f(P_{k\ell})$,
						\item\label{item:nablafPkell_2}
							$\displaystyle\frac{8}{3}(P_k-P_{k\ell})\cdot\nabla f(P_{k\ell})=f(P_k)+c(2e_k+e_\ell)-c(2e_\ell+e_k)-f(P_\ell)$.
					\end{enumerate}
					\begin{proofof}\textit{(\ref{item:fPj})--(\ref{item:fPkell}).}
						The formulas (\ref{item:fPj})--(\ref{item:fPkell}) follow with $(\lambda_1,\dots,\lambda_4)=e_j$ at $x=P_j$ and 
						the only possible $\alpha=3e_j$ 
						and  $(\lambda_1,\dots,\lambda_4)=(e_k+e_\ell)/2$ 
						at the midpoint $x=P_{k\ell}$ with possible 
						$\alpha=3e_k,\, 2e_k+e_\ell, 2e_\ell+e_k,\, 3e_\ell$
						and respective values $\lambda^\alpha/\alpha!=1/48, 1/16,1/16,1/48$.
						\end{proofof}
						\begin{proofof}\textit{(\ref{item:nablafPj_1}).}
						The polynomial $\lambda^\beta/\beta!$ vanishes at $x=P_j$ for all $\beta\not=2e_j$, while 
						$\lambda^{2e_j}/(2e_j)!=1/2$ at $x=P_j$. Consequently, 
						\begin{align}
							\nabla f(P_j)=3\sum_{\mu=1}^4c(2e_j+e_\mu)\nabla\lambda_\mu. \label{eq:nablaf(Pj)}
						\end{align} 
						The barycentric coordinates satisfy 
						$(P_j-P_k)\cdot\nabla \lambda_\mu=\lambda_\mu(P_j)-\lambda_\mu(P_k)=\delta_{j\mu}-\delta_{k\mu}$  
						and so (\ref{item:nablafPj_1}) follows from \eqref{eq:nablaf(Pj)}. 
						\end{proofof}
						\begin{proofof}\textit{(\ref{item:nablafPkell}).}
						At the edge midpoint $x=P_{k\ell}$, \cref{exp:k=2_m=4} shows $\lambda^\beta/\beta!=1/8,1/4,1/8$ for 
						$\beta=2e_k,e_k+e_\ell,2e_\ell$. Hence 
						\begin{align}
							\nabla f(P_{k\ell})=\frac{3}{4}\sum_{\mu=1}^4\big(c(2e_\ell+e_\mu)+2c(e_\ell+e_k+e_\mu)
						 		+c(2e_k+e_\mu)\big)\nabla\lambda_\mu.\label{eq:nablaf(Pjk)}
						\end{align}
						Since
						$(P_j-P_{k\ell})\cdot \nabla\lambda_\mu=\delta_{j\mu}-\frac{1}{2}(\delta_{k\mu}+\delta_{\ell\mu})$ 
						is one for $\mu=j$, $-1/2$ for $\mu=k,\ell$, and vanishes otherwise, (\ref{item:nablafPkell}) follows from (\ref{item:fPkell}).
						\end{proofof}
						\begin{proofof}\textit{(\ref{item:nablafPkell_2}).}
						Since 
						$(P_k-P_{k\ell})\cdot \nabla\lambda_\mu=\delta_{k\mu}-\frac{1}{2}(\delta_{k\mu}+\delta_{\ell\mu})$ is $1/2$ for $\mu=k$ and 
						$-1/2$ for 
						$\mu=\ell$ and vanishes otherwise, 		
						\eqref{eq:nablaf(Pjk)} leads to (\ref{item:nablafPkell_2}).
					\end{proofof}
			\end{example}

	\subsection{Smoothness across an interface}
			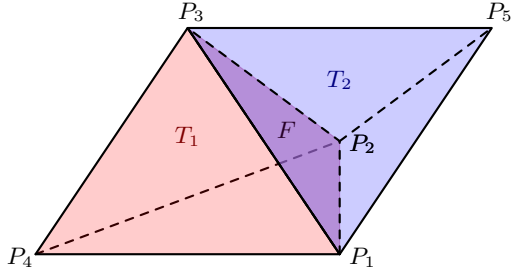
\begin{figure}
					\begin{center}
					\begin{tikzpicture}[x=20mm,y=30mm]
							\draw(0,0)node[shape=coordinate, label={[xshift=-2mm, yshift=-3mm] $P_4$}] (P4) {P4};
							\draw(2,0)node[shape=coordinate,label={[xshift=3mm, yshift=-3mm] $P_1$}] (P1) {P1};
		  					\draw(1,1)node[shape=coordinate,label={[xshift=0.5mm, yshift=-0.5mm] $P_3$}] (P3) {P3};
		  					\draw(3,1)node[shape=coordinate,label={[xshift=1mm, yshift=-0.5mm] $P_5$}] (P5) {P5};
		  					\draw(2,0.5)node[shape=coordinate,label={[xshift=3mm, yshift=-3mm] $P_2$}] (P2) {P2};
		  					
		  					\draw(1.6,0.4)node[shape=coordinate] (Q1Q4) {Q1Q4};			  					
		  					
		  					\draw[dashed, thick,gray] (Q1Q4)--(P2);
		  					\draw[dashed, thick] (P4)--(Q1Q4);  
		  					\filldraw[red, opacity=0.2](P1)--(P4)--(P3)--cycle;	
		  					\filldraw[blue, opacity=0.2](P1)--(P3)--(P5)--cycle;	
		  					\filldraw[violet, opacity=0.3](P1)--(P2)--(P3)--cycle;
		  					\draw [thick] (P1)--(P3)--(P4)--cycle;
		  					\draw [thick] (P1)--(P3)--(P5)--cycle;
		  					\draw [thick,dashed] (P1)--(P2)--(P3)--cycle;
		  					\draw [thick,dashed] (P2)--(P5);	
		  					
		  					\draw(1.7,0.5)node[shape=coordinate,label={[xshift=-1mm, yshift=-1mm] \textcolor{black!70!violet}{$F$}}] (c1) {c1};
		  					\draw(1.0,0.5)node[shape=coordinate,label={[xshift=0mm, yshift=-2mm]\textcolor{black!50!red}{ $T_1$}}] (cK) {cK};
		  					\draw(2.0,0.75)node[shape=coordinate,label={[xshift=0mm, yshift=-2mm] \textcolor{black!50!blue}{$T_2$}}] (cK) {cK};
		  					\draw(2,0.5)node[shape=coordinate,label={[xshift=3mm, yshift=-3mm] $P_2$}] (P2) {P2};
					\end{tikzpicture}\end{center}
		
				\caption{Illustration of two neighbouring tetrahedra $T_1=\textup{conv}\{P_1,P_2,P_3,P_4\}$ 
							in red and $T_2=\textup{conv}\{P_1,P_2,P_3,P_5\}$ in blue with common 
							face $F=\textup{conv}\{P_1,P_2,P_3\}$ in violet. 
							}\label{fig:Appendix_neighbouringTetrahedra}
			\end{figure} 
			This subsection characterizes the function space $C^1(T_1\cup T_2)\cap P_k(\{T_1,T_2\})$ for two tetrahedra $T_1$, $T_2$
			with a common face $F=\partial T_1\cap \partial T_2$. 
			\cref{fig:Appendix_neighbouringTetrahedra} illustrates the following situation.
			Suppose $T_1=\textup{conv}\{P_1,P_2,P_3,P_4\}$ and $T_2=\textup{conv}\{P_1,P_2,P_3,P_5\}$ 
			 share the face 
			$F=\partial T_1\cap \partial T_2=\textup{conv}\{P_1,P_2,P_3\}$ with unit normal $\nu_F$.
			Recall  $A_k:=\{\alpha\in\mathbb{N}_0^4|\, |\alpha|=k\}\subset \mathbb{N}_0^4$ and set 
			$B_k:=\{\alpha\in\mathbb{N}_0^{3}|\, |\alpha|=k\}\subset \mathbb{N}_0^{3}$. 
			In this context the notation $c_j(\beta;0)$ in \eqref{eq:con1_C1T1T2}--\eqref{eq:con2_C1T1T2} below abbreviates 
			 $c_j(\beta;0):=c_j((\beta_1,\beta_2,\beta_3,0))$ for some $\beta\in B_k$ or $\beta\in B_{k-1}$ and $j=1,2$. 
			 Identify $T_1$ with  $(P_1,P_2,P_3,P_4)$ and $T_2$ with $(P_1,P_2,P_3,P_5)$.
			The function $g\in P_k(\{T_1,T_2\})$ is given in $T_1$ resp. $T_2$ as the polynomial \eqref{eq:BernsteinPolynomial} and with ordinates 
			$(c_1(\alpha)|\,\alpha\in A_k)$ resp. $(c_2(\alpha)|\,\alpha\in A_k)$,
			\begin{align}
				g_j=g|_{T_j}=k!\sum_{{\alpha\in A_k}}c_j(\alpha)\frac{\lambda^\alpha}{\alpha!}\quad\text{ at }\quad 
				x=\lambda_1P_1+\lambda_2P_2+\lambda_3P_{3}+\lambda_4P_{\ell} \in T_j \label{eq:gonT1T2}
			\end{align}
			for $j=1,2$ with $\ell=\ell(1)=4$ for $j=1$ and $\ell=\ell(2)=5$ for $j=2$. 	

			\begin{lemma}\label{lem:g_C1_onT1T2}
				The polynomials $g_j\in P_k(T_j)$ with \eqref{eq:gonT1T2} for $j=1,2$ form a function $g\in P_k(\{T_1, T_2\})$ by $g|_{T_j}=g_j$ for 	
				$j=1,2$. 
				Then 
				\begin{enumerate}[label=(\alph*), ref=\alph*, wide]
				  \item\label{item:con1_C1T1T2} $g\in C^0(T_1\cup T_2)$ if and only if 
					\begin{align}
						&c_1(\beta;0)=c_2(\beta;0)\quad\text{for all }\beta \in B_k;\label{eq:con1_C1T1T2}
				    \end{align}
				  \item\label{item:con2_C1T1T2} $g\in C^1(T_1\cup T_2)$ if and only if \eqref{eq:con1_C1T1T2} and 
					\begin{align}
						&\sum_{\mu=1}^4c_1((\beta;0)+e_\mu)(\nu_F\cdot\nabla\lambda_\mu|_{T_1})
						 =\sum_{\mu=1}^4c_2((\beta;0)+e_\mu)(\nu_F\cdot\nabla\lambda_\mu|_{T_2})
						\text{ for all }\beta \in B_{k-1}.\label{eq:con2_C1T1T2}
					\end{align} 
				\end{enumerate}
			\end{lemma}
			\begin{proofof}\textit{(\ref{item:con1_C1T1T2}).}
					Recall that the common face $F=\partial T_1\cap\partial T_2$ satisfies $F\subset \{\lambda_4=0\}$. 
					Hence \eqref{eq:gonT1T2} implies 
						\begin{align*}
							&g_j|_F=k!\sum_{\beta\in B_k}c_j(\beta;0)\frac{\lambda^\beta}{\beta!}\quad
							\text{at }x=\lambda_1P_1+\lambda_2P_2+\lambda_3P_3\in F \text{ and }j=1,2
						\end{align*}
					with the abbreviation $c_j(\beta;0):=c_j((\beta_1,\beta_2,\beta_3,0))$ 
					and ${\lambda^\beta}/{\beta!}=\frac{\lambda_1^{\beta_1}\lambda_2^{\beta_2}\lambda_3^{\beta_3}}{\beta_1!\beta_2!\beta_3!}$.
					Therefore, $g_1|_F=g_2|_F$ for $g\in C^0(T_1\cup T_2)$ is equivalent to 
						\begin{align*}
							0=\sum_{\beta\in B_k}\big(c_1(\beta;0)-c_2(\beta;0)\big){\lambda^\beta}/{\beta!}\quad 
							\text{at }x=\lambda_1P_1+\lambda_2P_2+\lambda_3P_3\in F. 
						\end{align*}
					Since the Bernstein polynomials $\big({\lambda^\beta}/{\beta!}|\, \beta\in B_k\big)$ form a basis of $P_k(F)$ the last statement is 
					equivalent to \eqref{eq:con1_C1T1T2}. This proves (\ref{item:con1_C1T1T2}). 
					\end{proofof}
				\begin{proofof}\textit{(\ref{item:con2_C1T1T2}).}
					Since $g\in C^1(T_1\cup T_2)$ includes $g\in C^0(T_1\cup T_2)$, \eqref{eq:con1_C1T1T2} has to hold. 
					Additionally we need continuity of the gradient at $F$. 					 
					The gradient $\nabla g_j|_F$ from \eqref{eq:derivativeBernstein} reads 
						\begin{align*}
							\nabla g_j|_F &
							=k!\sum_{\beta\in A_{k-1}}\frac{\lambda^\beta}{\beta!}\bigg|_F\sum_{\mu=1}^4c_j(\beta+e_\mu)\nabla\lambda_\mu|_{T_j}
							=k!\sum_{\beta \in B_{k-1}}\frac{\lambda^\beta}{\beta!}\sum_{\mu=1}^4c_j((\beta;0)+e_\mu)\nabla\lambda_\mu|_{T_j}
						\end{align*}
						at $x=\lambda_1P_1+\lambda_2P_2+\lambda_3P_3\in F$. 
					Hence the continuity of the gradients $\nabla g_1|_F=\nabla g_2|_F$ is equivalent to 
						\begin{align*}
						 0=\sum_{\beta\in B_{k-1}}\frac{\lambda^\beta}{\beta!}\sum_{\mu=1}^4\Big( c_1\big((\beta;0)+e_\mu\big)
						 \nabla\lambda_\mu|_{T_1}-c_2\big((\beta;0)+e_\mu\big)\nabla\lambda_\mu|_{T_2}\Big)
						 \quad\text{at }x=\sum_{\mu=1}^3\lambda_\mu P_\mu\in F. 
						\end{align*}
					Since the Bernstein polynomials $\big({\lambda^\beta}/{\beta!}|\, \beta\in B_{k-1}\big)$ form a basis of $P_{k-1}(F)$ 
					the last statement is equivalent to
						\begin{align}
							 0=\sum_{\mu=1}^4\Big( c_1\big((\beta;0)+e_\mu\big)\nabla\lambda_\mu|_{T_1}-c_2\big((\beta;0)+e_\mu\big)
							 \nabla\lambda_\mu|_{T_2}\Big)\in \mathbb{R}^3\quad\text{for all }\beta\in B_{k-1}. \label{eq:proof_cond2_gC1}
						\end{align}
					It remains to show, that \eqref{eq:proof_cond2_gC1} 
					is equivalent to \eqref{eq:con2_C1T1T2}. In other words, to verify that solely the normal derivatives 
					$\nu_F\cdot \nabla\lambda_\mu|_{T_j}$ for $j=1,2$ and $\mu=1,2,3$ are of interest.  
					Fix two tangential directions $\tau_1,\,\tau_2$ of unit length parallel to $\langle F\rangle$, 
					such that $\tau_1,\tau_2, \nu_F$ form an orthonormal basis of $\mathbb{R}^3$. Then 
					$$\nabla \lambda_\mu|_{T_j}
						= (\nabla \lambda_\mu|_{T_j}\cdot \tau_1)\tau_1+(\nabla \lambda_\mu|_{T_j}\cdot \tau_2)\tau_2
							+(\nabla \lambda_\mu|_{T_j}\cdot \nu_F)\nu_F.$$
					The components of the difference $\nabla\lambda_\mu|_{T_1}-\nabla\lambda_\mu|_{T_2}$ for $\mu=1,2,3$ in 
					$\mathbb{R}^3$ in the direction $\tau_1,\tau_2\in\textup{span}\{P_2-P_1,P_3-P_1\}$ parallel to $\langle F\rangle$ 
					vanish owing to the Hadamard  jump condition: 
					  The piecewise gradient $\nabla g$ of $g\in C^0(T_1\cup T_2)\cap P_k(\{T_1,T_2\})$ jumps across the interface 
					 $F$ 
					  and the jump $[\nabla g]_F\parallel \nu_F$ exclusively  points in the normal direction. 
					  For $\mu=1,2,3$ and the tangential direction $\tau_j\in \textup{span}\{P_2-P_1,P_3-P_1\}$ for $j=1,2$, we have 
						   \begin{align*}
							   \tau_j\cdot \nabla\lambda_\mu|_{T_1}
							   =\frac{\partial}{\partial s}\lambda_\mu(x_0+s\tau)\big|_{s=0}=\tau_j\cdot \nabla\lambda_\mu|_{T_2}
							   =: \tau_j\cdot \nabla\lambda_\mu
							   \quad\text{at any }x_0\in \text{relint}(F).
						   \end{align*}
					For $\mu=4$, notice that $P_r-P_1\perp \nu_F\parallel\nabla\lambda_4|_{T_j}$ for $r=2,3$ and $j=1,2$. 
					In fact, suppose $\varrho_\ell>0$  denotes  the height  of $P_\ell$ 
					  over the plane $\langle F\rangle$, then 
					   $\varrho_\ell \nabla\lambda_4|_{T_j}=(-1)^j\nu_F$ for $\ell=4,5$, $j=1,2$  
					   provided the orientation of $\nu_F$ is fixed from $P_4 $ to $P_5$ in that $\nu_F$ is the outer unit normal in $F\subset\partial T_1$.  
					   In particular this means $\tau_j\cdot\nabla\lambda_4|_{T_1}=0=\tau_j\cdot\nabla\lambda_4|_{T_2}$ for $j=1,2$.
					Hence, \eqref{eq:proof_cond2_gC1} is recast for all $\beta\in B_{k-1}$ as the first identity in 
					\begin{align*}
						-\nu_F \sum_{\mu=1}^4&\Big( c_1\big((\beta;0)+e_\mu\big)(\nabla\lambda_\mu|_{T_1}\cdot\nu_F)-c_2\big((\beta;0)+e_\mu\big)
							 (\nabla\lambda_\mu|_{T_2}\cdot\nu_F)\Big)
							\\& =	\sum_{j=1}^2\tau_j
							\sum_{\mu=1}^3\Big( c_1\big((\beta;0)+e_\mu\big)-c_2\big((\beta;0)+e_\mu\big)\Big)
							 \nabla\lambda_\mu\cdot\tau_j\Big)
							 =	 0						 
					\end{align*}					 					
					with the equality of ordinates from \eqref{eq:con1_C1T1T2} in the last step. 
					 The last identity is equivalent to \eqref{eq:con2_C1T1T2}. 
				\end{proofof}
			The coefficient relations \eqref{eq:con1_C1T1T2}--\eqref{eq:con2_C1T1T2} concern the factors one in 
			\eqref{eq:con1_C1T1T2} or 
			geometric factors $\nu_F\cdot\nabla\lambda|_{T_j}$ for $j=1,2$ in \eqref{eq:con2_C1T1T2}. 
			However they do \emph{not} explicitly display $k$  and this gives rise to a factorization (cf. \cite[\S14]{Boor1987} in terms of B-forms). 
			Suppose $(c_j(\alpha)|\alpha\in A_k)$ defines $g_j$ in \eqref{eq:gonT1T2} for 
			$F=\partial T_1\cap \partial T_2=\textup{conv}\{P_1,P_2,P_3\}$ and $j=1,2$. Define for a fixed index 
			$\xi\in\{1,2,3\}$ suppressed in the notation
				\begin{align}
				 	g^\prime|_{T_j}
				 		:=g_j^\prime:=(k-1)!\sum_{\beta \in A_{k-1}}c_j(e_\xi+\beta){\lambda^\beta}/{\beta!}\quad\text{for }j=1,2.\label{eq:gprimeon T1T2} 
				\end{align}
			Then the composition $g^\prime\in P_{k-1}(\{T_1,T_2\})$ from \eqref{eq:gprimeon T1T2} is a piecewise polynomial of degree at most $k-1$.
				\begin{example}[$g^\prime$ for $k=0,\dots,3$]\label{ex:gprime} 
				Abbreviate $\sum_{\mu\le \ell}:=\sum_{\substack{\mu,\ell=1\\\mu\le\ell}}^4$.
				\begin{enumerate}[label=(\alph*)]
					\item If $k=0$ and $g_j=c_j((0,0,0,0))$, then $g^\prime_j=0$ by definition.
					\item If $k=1$ and $g_j=\sum_{\mu=1}^4 c_j(e_\mu){\lambda^{e_\mu}}$, then $g^\prime_j=c_j(e_\xi)$.  
					\item If $k=2$ and $g_j=2\sum_{\mu\le\ell} c_j(e_\mu+e_\ell)\frac{\lambda^{e_\mu+e_\ell}}{(e_\mu+e_\ell)!}$, then 
								$g^\prime_j=\sum_{\mu=1}^4 c_j(e_\xi+e_\mu){\lambda_\mu}$; e.g., 	
								$g^\prime=c_j((2,0,0,0))\lambda_1+c_j((1,1,0,0))\lambda_2+c_j((1,0,1,0))\lambda_3+c_j((1,0,0,1))\lambda_4$ holds for $\xi=1$.
					\item If $k=3$ and $g_j=6\sum_{\alpha\in A_3} c_j(\alpha)\frac{\lambda^{\alpha}}{\alpha!}$, then 
							$g^\prime_j=2\sum_{\mu\le\ell}c_j(e_\xi+e_\mu+e_\ell)\frac{\lambda^{e_\mu+e_\ell}}{(e_\mu+e_\ell)!}$; e.g., 
						\begin{align*} 
							g^\prime_j=&2\lambda_1\big(c_j((2,1,0,0))\lambda_2+c_j((2,0,1,0))\lambda_3+c_j((2,0,0,1))\lambda_4\big)
									\\& +2\big( c_j(1,1,1,0)\lambda_2\lambda_3+c_j(1,1,0,1)\lambda_2\lambda_4+c_j(1,0,1,1)\lambda_3\lambda_4\big)
									\\&+c_j((3,0,0,0))\lambda_1^2+c_j((1,2,0,0))\lambda_2^2+c_j((1,0,2,0))\lambda_3^2+c_j((1,0,0,2))\lambda_4^2
						\end{align*}holds for $\xi=1$.
				\end{enumerate}
				\end{example}
				Recall that the two tetrahedra $T_1$ and $T_2$ in \cref{fig:Appendix_neighbouringTetrahedra} share  the face 
				$F=\partial T_1\cap\partial T_2=\textup{conv}\{P_1,P_2,P_3\}$ with fixed unit normal $\nu_F$. 
				Suppose $g|_{T_j}:=g_j$ in \eqref{eq:gonT1T2} for $j=1,2$ and 
				define  $g^\prime$ in \eqref{eq:gprimeon T1T2} with respect to a fixed index $\xi\in\{1,2,3\}$. 
				Let $\ell, m$ denote the two distinct indices with $\{\xi,\ell,m\}=\{1,2,3\}$ and $E:=\textup{conv}\{P_\ell,P_m\}$.
				\begin{lemma}\label{lem:gprime_C1_onT1T2}
					Under the present notation  $g\in C^1(T_1\cup T_2)$ is equivalent to the three conditions 
					$$g^\prime\in C^1(T_1\cup T_2)\text{ as well as } g\text{ and }\frac{\partial g}{\partial\nu_F}\text{ are continuous at }E.$$
				\end{lemma}
					\begin{proof}
						Without loss of generality fix $\xi=1$ and $E=\textup{conv}\{P_2,P_3\}$ throughout this proof.
						\cref{lem:g_C1_onT1T2}
						characterizes $g\in C^1(T_1\cup T_2)$ in \eqref{eq:con1_C1T1T2}--\eqref{eq:con2_C1T1T2}; but it also 
						applies to 
						$g^\prime\in P_{k-1}(\{T_1,T_2\})$ and shows that $g^\prime \in C^1(T_1\cup T_2)$ is equivalent to 
						\eqref{eq:con1_gprime}--\eqref{eq:con2_gprime},
						\begin{align}
							&c_1\big(e_1+(\beta;0)\big)=c_2\big(e_1+(\beta;0)\big)\quad \text{for all }\beta
							\in B_{k-1}:=\{\alpha\in\mathbb{N}_0^{3}|\, |\alpha|=k-1\},
							\label{eq:con1_gprime}\\
							&\sum_{\mu=1}^4 c_1(e_1+(\beta;0)+e_\mu)\nu_F\cdot \nabla \lambda_\mu|_{T_1}
								=\sum_{\mu=1}^4c_2(e_1+(\beta;0)+e_\mu)\nu_F\cdot \nabla \lambda_\mu|_{T_2}\quad\text{ for all }\beta\in B_{k-2}.
								\label{eq:con2_gprime}
						\end{align}
						Observe that the conditions \eqref{eq:con1_gprime}--\eqref{eq:con2_gprime} concern $e_1+(\beta;0)$ for $\beta\in B_{k-1}$ and  
						$\beta\in B_{k-2}$ and are included in 
						\eqref{eq:con1_C1T1T2}--\eqref{eq:con2_C1T1T2}. In other words \eqref{eq:con1_C1T1T2}--\eqref{eq:con2_C1T1T2} imply  
						\eqref{eq:con1_gprime}--\eqref{eq:con2_gprime}.
						\\[1em]\textit{\glqq $\Rightarrow$\grqq }
						The first implication assumes  $g\in C^1(T_1\cup T_2)$ and \cref{lem:g_C1_onT1T2} guarantees  
						\eqref{eq:con1_C1T1T2}--\eqref{eq:con2_C1T1T2} and so  
						 \eqref{eq:con1_gprime}--\eqref{eq:con2_gprime}. Recall from the very beginning of this proof that 
						 \eqref{eq:con1_gprime}--\eqref{eq:con2_gprime}   imply $g^\prime\in C^1(T_1\cup T_2)$.  
						Finally $g \in C^1(T_1\cup T_2)$ and $\nabla g\in C^0(T_1\cup T_2;\mathbb{R}^3)$ imply that  $g$ and $\frac{\partial g}{\partial \nu_F}$ 
						are continuous along $E=\textup{conv}\{P_2,P_3\}$.
						\\[1em]\textit{\glqq $\Leftarrow$\grqq } 
						The converse implication that assumes $g^\prime\in C^1(T_1\cup T_2)$, whence  
						\eqref{eq:con1_gprime}--\eqref{eq:con2_gprime}, and that $g$ and 
						$\frac{\partial g}{\partial \nu_F}$ are 
						continuous along $E\subset \{\lambda_1=0=\lambda_4\}$. 
						First \eqref{eq:con1_gprime} implies some conditions of \eqref{eq:con1_C1T1T2} and the remaining conditions read 
						$$c_1(0,\beta_2,\beta_3,0)=c_2(0,\beta_2,\beta_3,0)\text{ for all }\beta_2,\beta_3\in\mathbb{N}_0\text{ with }\beta_2+\beta_3=k.$$ 
						The remaining conditions concern $\beta\in B_{k}$ with $\beta_1=0$. 
						On the other hand, the continuity of $g$ along $E$ implies $g_1|_E=g_2|_E$. In terms of \eqref{eq:gonT1T2}, this reads 
						\begin{align}
							0=g_1|_E-g_2|_E
								=k!\sum_{\substack{\beta\in B_k\\\beta_1=0}}
									\Big(c_1\big((\beta;0)\big)-c_2\big((\beta;0)\big)\Big)
									\frac{\lambda^{\beta}}{\beta!}
									\quad \text{ at }x=\lambda_2P_2+\lambda_3P_3\in E. \label{eq:proof_cond1_gprimeC1}
						\end{align}
						Since the Bernstein polynomials 
						$\big({\lambda^\beta}/\beta!|\,\beta\in B_k, \beta_1=0)$ 
						form a basis of $P_k(E)$, \eqref{eq:proof_cond1_gprimeC1} implies the remaining condition in \eqref{eq:con1_C1T1T2}.   
						Consequently, \eqref{eq:con1_gprime} and $g|_E\in C^0(E)$ prove $g\in C^0(T_1\cup T_2)$. 
						Second, \eqref{eq:con2_gprime} implies the conditions of \eqref{eq:con2_C1T1T2} except those for $\beta\in B_{k-1}$ with 
						$\beta_1=0$. 
						Recall that \eqref{eq:gonT1T2} and \cref{rem:derivativeBernsteinBasis} show  						
						\begin{align*}
							\nu_F\cdot\nabla g_j|_E
								=k!\sum_{\substack{\beta\in B_{k-1}\\\beta_1=0}}\frac{\lambda^{\beta}}{\beta!} 
								\sum_{\mu=1}^4 c_j\big((\beta;0)+e_\mu\big) (\nu_F\cdot\nabla\lambda_\mu|_{T_j})
									\text{ at }x=\lambda_2P_2+\lambda_3P_3\in E 
						\end{align*}
						for $j=1,2$. Deduce that $\frac{\partial (g_1-g_2)}{\partial\nu_F}=\nu_F\cdot \nabla (g_1-g_2)=0$ in $E$ implies 
						$c_1\big((\beta;0)\big)=c_2\big((\beta;0)\big)$ for all $\beta\in B_{k-1}$ with 
						$\beta_1=0$. 
						 Consequently, 
						\eqref{eq:con1_gprime}--\eqref{eq:con2_gprime} and $g|_E$, $\frac{\partial g}{\partial \nu_F}\Big|_E\in C^0(E)$ imply 
						\eqref{eq:con1_C1T1T2}--\eqref{eq:con2_C1T1T2}. 
					\end{proof}
		
		\subsection{Piecewise quadratic polynomials in a $\textit{WF}$ partition of a tetrahedron}\label{sec:HCT_Triang} 
		Given a tetrahedron $K=\textup{conv}\{Q_1,Q_2,Q_3,Q_4\}\subset\mathbb{R}^3$ of a positive volume, select five center points 
		$c_K,\,c_1,\,c_2,\,c_3,\,c_4$  in $K$ as follows
			\begin{enumerate}[label=(\alph*)]
				\item\label{item:CTcon1} $c_K\in \textup{int}(K)$ (i.e. $\textup{dist}(c_K,\partial K)>0$)
				\item\label{item:CTcon2} $c_{m}\in\mathrm{relint}(F_{m})$
								 (i.e. $c_{m}=\sum_{\mu=1}^4\lambda_{\mu}Q_\mu$ for $\lambda_{1}+\dots+\lambda_{4}=1$ and $\lambda_{m}=0$ and all 
								 other $\lambda_{j}>0$) for the face $F_{m}$ of $K$ that is opposite to the vertex $Q_{m}$ for ${m}=1,2,3,4$.
			\end{enumerate}
			\begin{definition}[$\textit{WF}$ 
								partition in $3$D]\label{def:CTtriangulation}
				Given $c_K,\,c_1,\,c_2,\,c_3,\, c_4$ with \ref{item:CTcon1}--\ref{item:CTcon2} 
				in the tetrahedron $K=\allowbreak\textup{conv}\{Q_1{, } Q_2{, }Q_3{, }Q_4\}$,  let the regular triangulation 
				$\textit{WF3D}(K)=\widehat{\mathcal{T}}$ 
				be the set of $12$ 
				subtetrahedra obtained as the convex hull $\textup{conv}\{Q_j,Q_k,c_K,c_m\}$ 
				of two distinct vertices $Q_j$ and $Q_k$, $1\le j<k\le 4$, 
				the center point 
				$c_K$ of the tetrahedron, and  the center $c_m\in\mathrm{relint}(F_m)$ of the face $F_m$ of $K$ opposite to the vertex $Q_m$ for 
				$m\in \{1,2,3,4\}\setminus\{j,k\}$.
			\end{definition}
		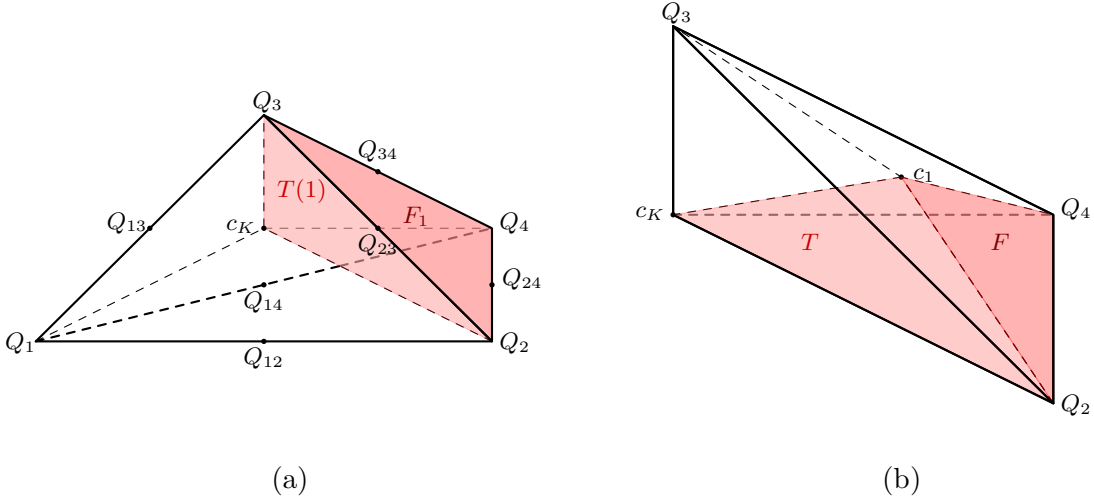
\begin{figure}
				\begin{minipage}{0.49\textwidth}\begin{center}
					\begin{tikzpicture}[x=30mm,y=30mm]
							\draw(0,1.33)node[shape=coordinate] (P) {P};
		  					\draw(0,0)node[shape=coordinate,label={[xshift=-2mm, yshift=-3mm] $Q_1$}] (Q1) {Q1};
		  					\draw(2,0)node[shape=coordinate,label={[xshift=3mm, yshift=-3mm] $Q_2$}] (Q2) {Q2};
		  					\draw(1,1)node[shape=coordinate,label={[xshift=0.5mm, yshift=-1mm] $Q_3$}] (Q3) {Q3};
		  					\draw(2,0.5)node[shape=coordinate,label={[xshift=3mm, yshift=-2mm] $Q_4$}] (Q4) {Q4};
		  					\draw(1.0,0.5)node[shape=coordinate,label={[xshift=-3mm, yshift=-2mm] $c_K$}] (cK) {cK};
		  					\draw(1.3333,0.3333)node[shape=coordinate] (Q1Q4) {Q1Q4};			  					
		  					\draw[dashed, thick,gray] (Q1Q4)--(Q4);
		  					\draw[dashed, thick] (Q1)--(Q1Q4);
		  					\draw[dashed] (Q2)--(cK);
		  					\draw[dashed, gray] (cK)--(Q4);
		  					\draw[dashed] (Q1)--(cK)--(Q3);
		  					\filldraw[red, opacity=0.3](Q2)--(Q4)--(Q3)--cycle;
		  					\filldraw[red, opacity=0.2](Q2)--(cK)--(Q3)--cycle;
		  					\draw [thick] (Q1)--(Q2)--(Q3)--cycle;
		  					\draw [thick] (Q2)--(Q4)--(Q3)--cycle;	
				
		  					\draw(1,0)node[shape=coordinate,label={[xshift=0mm, yshift=-5mm] $Q_{12}$}] (Q12) {Q12};	
		  					\draw(0.5,0.5)node[shape=coordinate,label={[xshift=-3mm, yshift=-2mm] $Q_{13}$}] (Q13) {Q13};
		  					\draw(1,0.25)node[shape=coordinate,label={[xshift=0mm, yshift=-5mm] $Q_{14}$}] (Q14) {Q14};
		  					\draw(1.5,0.5)node[shape=coordinate,label={[xshift=0mm, yshift=-5mm] $Q_{23}$}] (Q23) {Q23};
		  					\draw(2,0.25)node[shape=coordinate,label={[xshift=4mm, yshift=-2mm] $Q_{24}$}] (Q24) {Q24};
		  					\draw(1.5,0.75)node[shape=coordinate,label={[xshift=0mm, yshift=0mm] $Q_{34}$}] (Q34) {Q34};
		  					\fill (Q12) circle (1pt);
		  					\fill (Q13) circle (1pt);
		  					\fill (Q14) circle (1pt);
		  					\fill (Q23) circle (1pt);
		  					\fill (Q24) circle (1pt);
		  					\fill (Q34) circle (1pt);
		  						
		  					\draw(1.7,0.5)node[shape=coordinate,label={[xshift=-1mm, yshift=-1mm] \textcolor{black!50!red}{$F_1$}}] (c1) {c1};
		  					\draw(1.0,0.5)node[shape=coordinate,label={[xshift=5mm, yshift=2mm] \textcolor{black!20!red}{$T(1)$}}] (cK) {cK};
		  					
		  					\fill (cK) circle (1pt);
					\end{tikzpicture}\end{center}
				\end{minipage}
				\begin{minipage}{0.49\textwidth}\begin{center}
					\begin{tikzpicture}[x=50mm,y=50mm]
		  					\draw(2,0)node[shape=coordinate,label={[xshift=3mm, yshift=-3mm] $Q_2$}] (Q2) {Q2};
		  					\draw(1,1)node[shape=coordinate,label={[xshift=0.5mm, yshift=-1mm] $Q_3$}] (Q3) {Q3};
		  					\draw(2,0.5)node[shape=coordinate,label={[xshift=3mm, yshift=-2mm] $Q_4$}] (Q4) {Q4};
		  					
		  					\draw(1.0,0.5)node[shape=coordinate,label={[xshift=-3mm, yshift=-2mm] $c_K$}] (cK) {cK};
		  					\draw[thick] (Q2)--(Q3)--(Q4)--cycle;
		  					\draw[thick] (Q2)--(cK)--(Q3);
		  					\draw[thick,dashed,gray] (cK)--(Q4);
		  				
		  					\draw(1.6,0.6)node[shape=coordinate,label={[xshift=3mm, yshift=-2mm] $c_1$}] (c1) {c1};
		  					\fill (c1) circle (1pt);
		  					\draw[dashed] (Q2)--(c1)--(Q4);
		  					\draw[dashed] (Q2)--(c1)--(Q3);
		  					\draw[dashed] (c1)--(cK);  
		  					\fill (cK) circle (1pt);
		  					\filldraw[red, opacity=0.3](Q2)--(c1)--(Q4)--cycle;
		  					\filldraw[red, opacity=0.2](cK)--(c1)--(Q2)--cycle;
		  					\draw(1.0,0.5)node[shape=coordinate,label={[xshift=18mm, yshift=-6mm] \textcolor{black!20!red}{$T$}}] (cK) {cK};
		  					\draw(1.7,0.5)node[shape=coordinate,label={[xshift=8mm, yshift=-6mm] \textcolor{black!50!red}{$F$}}] (c1) {c1};
		  					\draw[thick] (Q2)--(Q3)--(Q4)--cycle;
		  					\draw[thick] (Q2)--(cK)--(Q3);
					\end{tikzpicture}\end{center}
				\end{minipage}	
				\begin{center}
				\begin{tabular}{p{0.49\textwidth}p{0.49\textwidth}}
				 \centering (a)&\centering (b)
				\end{tabular}
				\end{center}			
				\caption{Illustration of the division of $K=\textup{conv}\{Q_1,\dots,Q_4\}$ with center $c_K\in\textup{int}(K)$ into 
							$T(1),\dots, T(4)$ in (a) with $T(1)=\textup{conv}(c_K,F_1)$ marked in red. 
							Further subdivision of $T(1)$ with center $c_1\in\textup{int}(F_1)$ of face $F_1$ in the $\textit{WF}$ 
							partition in (b)  
							with one subsimplex $T=\textup{conv}(c_k,F)\in \textit{WF3D}(K)$ 
							with face $F=\textup{conv}\{c_1,Q_2,Q_4\}\subset F_1$ marked in red. 
							}\label{fig:Appendix_CT_Triang}
			\end{figure}
		Each face $F_{m}\in\mathcal{F}(K)$ is a triangle partitioned  in $\textit{WF3D}(K)$ 
		by connecting the center 
		$c_{m}\in\mathrm{relint}(F_{m})$ with its vertices $\mathcal{V}(F_m)$  as in the $2$D 
		$HCT$ partition of the triangle $F_m$. In  fact there are two steps in the partition illustrated in \cref{fig:Appendix_CT_Triang}. 
		First, the four tetrahedra $T({m}):=\textup{conv}(c_K,F_{m})$ for ${m}=1,2,3,4$ 
		 partition $K$ as displayed in \cref{fig:Appendix_CT_Triang}.a. 
		Second, $\widehat{\mathcal{T}}(c_{m}):=\{T\in\widehat{ \mathcal{T}}(K)\,|\,c_{m} \text{ vertex of }T\}$ partitions each of the tetrahedra 
		$T({1})$, $T(2)$, $T(3)$, $T(4)$  in three subtetrahedra; \cref{fig:Appendix_CT_Triang}.b illustrates the partition of $T(1)$. 
		The following \cref{lem:C1Tm_k=012} adopts the notation of a $\textit{WF}$ 
		triangulation $\widehat{\mathcal{T}}=\textit{WF3D}(K)$ 
		and considers each 
		$T({m})=\textup{conv}(c_K,F_{m})$ for ${m}\in\{1,\dots,4\}$ with its decomposition 
		$\widehat{\mathcal{T}}(T({m}))=\{T\in\widehat{\mathcal{T}}\,|\,T\subset T(m)\}=:\widehat{\mathcal{T}}(c_{m})$ into three subtetrahedra. 
			\begin{lemma}\label{lem:C1Tm_k=012}
				Any $f\in P_k(\widehat{\mathcal{T}}(c_{m}))$ for $k=0,1,2$ and $m=1,2,3,4$ satisfies $f\in C^1(T(m))$ if and only 
				if $f\in P_k(T(m))$ is a global polynomial in $T(m)$.
			\end{lemma}
				\begin{proof}
					Without loss of generality let $m=1$ and $T(1)=\textup{conv}\{c_K,Q_2,Q_3,Q_4\}$ as displayed in \cref{fig:Appendix_CT_Triang}.b.
					The assertion is trivial for $k=0$ and immediate for $k=1$ (because continuity means 
					$\nabla f(c_1)=\nabla f|_T(c_1)$ for all $T\in \widehat{\mathcal{T}}(c_1)$ and leads to a global constant vector $\nabla f$). 
					Suppose $k=2$ and $f\in C^1(T(1))\cap P_2(\widehat{\mathcal{T}}(c_1))$. 
					The function values $f(z)$ at the vertices $z\in \{c_K,Q_2,Q_3,Q_4\}$ of $T(1)$ 
					and at the edge midpoints $z\in \{Q_{jk}:=(Q_j+Q_k)/2:\,2\le j<k\le 4\}\cup\{ (c_K+Q_j)/2:\,2\le j\le 4\}$ of $T(1)$ 
					define the quadratic Lagrange	interpolation $I_2f\in P_2(T(1))$ in the tetrahedron $T(1)$. 
					Subtract this (global) quadratic Lagrange interpolation  $I_2f\in P_2(T(1))$ from $f$. The difference $(f-I_2f)(z)=0$ 
					vanishes for each vertex and each edge midpoint of $T(1)$.  
					Since the assertion $f\in P_2(T(1))$ is equivalent to $f-I_2f\equiv 0$,  
					we may and will assume without loss of generality that $f$ vanishes at these $10$ distinct points on the boundary 
					$\partial T(1)$ and conclude $f\equiv 0$ below.
					
					Identify any subtetrahedron  
					$T=\textup{conv}\{c_1,c_K, Q_j,Q_k\}\in\widehat{\mathcal{T}}(c_1)$ for distinct $2\le j<k\le 4$ 
					with $(c_1,c_K,Q_j,Q_k)$. The Bernstein basis of $f|_T\in P_2(T)$ leads to ordinates $(c_T(\alpha)|\alpha\in A_2)$ in  
					$$f|_T=2\sum_{\alpha\in A_2}c_T(\alpha)\lambda^\alpha/\alpha!\text{ at }
					x=\lambda_1c_1+\lambda_2c_K+\lambda_3Q_j+\lambda_4Q_k\in T\in\widehat{\mathcal{T}}(c_1).$$ 
					Since $f|_T$ is a quadratic function in each $T\in\widehat{\mathcal{T}}(c_1)$,
					it is in particular a quadratic function on each face 
					$F_{jk}:=\textup{conv}\{c_K,Q_j,Q_k\}$ for $2\le j<k\le 4$  in $T$ opposite to $c_1$. 
					The quadratic function $f|_{F_{jk}}\in P_2(F_{jk})$  
					vanishes at the vertices $c_K, Q_j, Q_k$  and at the midpoints $Q_{jk}, (c_K+Q_j)/2$, and $(c_K+Q_k)/2$ of the edges of the triangle $F_{jk}$. 
					Hence 
					 	$$(f|_T)\big|_{F_{jk}}
					 		=2\sum_{\substack{\alpha \in A_2\\\alpha_1=0}}c_T(\alpha)\lambda^\alpha/\alpha!=0 \quad
					 		\text{at }x=\lambda_2c_K+\lambda_3Q_j+\lambda_4Q_k\in F_{jk}$$
					vanishes 
					and  $c_T(0;\beta)=0$ for all $\beta\in\mathbb{N}_0^3$ with $|\beta|=2$ and $T\in\widehat{\mathcal{T}}(c_1)$.
					This leads to 
						 \begin{align}
							 f|_T=2\lambda_1\sum_{\beta\in A_1}\frac{c_T(e_1+\beta)}{(1+\beta_1)}\frac{\lambda^\beta}{\beta!}
							 \quad\text{at }x=\lambda_1c_1+\lambda_2c_K+\lambda_3Q_j+\lambda_4Q_k\in 
							 T\in\widehat{\mathcal{T}}(c_1).
							 \label{eq:proof_lem_C1Tm_k=012_fT}
						\end{align} 
					The continuity conditions in \cref{lem:gprime_C1_onT1T2} with $\xi=1$ in \eqref{eq:gprimeon T1T2} 
					motivate the piecewise affine function 
						\begin{align}
							f^\prime|_T:&=\sum_{\beta\in A_1}c_T(e_1+\beta)\frac{\lambda^\beta}{\beta!}=\sum_{\mu=1}^4 c_T(e_1+e_\mu)\lambda_\mu
							\quad\text{at }x\in T\in\widehat{\mathcal{T}}(c_1)
							 \label{eq:proof_lem_C1Tm_k=012_fprime}
						\end{align} 
						with ordinates $\big(c_T(e_1+e_\mu):\, \mu=1,\dots, 4\big)$ of $f|_T$ for each $T\in \widehat{\mathcal{T}}(c_1)$.
					Since $f\in C^1(T(1))$ and any two tetrahedra $T_+,\, T_-\in \widehat{\mathcal{T}}(c_1)$ share a face 
					$F=\partial T_+\cap\partial T_-=\textup{conv}\{ c_1,c_K,Q_j\}$  for some $2\le j\le 4$ with common vertex $c_1$, 
					\cref{lem:gprime_C1_onT1T2} shows $f^\prime\in C^1(T_+\cup T_-)$. 
					The repeated application of  \cref{lem:gprime_C1_onT1T2} leads to $f^\prime\in C^1(T(1))$. 
					Since $f^\prime \in P_1(\widehat{\mathcal{T}}(c_1))$ is piecewise affine, the already proven assertion of this lemma for $k=1$ implies 
					that $f^\prime\in C^1(T(1))\cap P_1(\widehat{\mathcal{T}}(c_1))=P_1(T(1))$  is one global affine function in $\mathbb{R}^3$. 
					Let $\varphi_1\in S^1(\widehat{\mathcal{T}}(c_1)):=P_1(\widehat{\mathcal{T}}(c_1))\cap C(T(1))$  denote 
					the first-order nodal basis function with $\varphi_1(c_1)=1$ and 
					$\varphi_1(c_K)=0=\varphi_1(Q_j)$ for any $2\le j\le 4$. In other words  $\varphi_1|_T=\lambda_1$ restricted to any 
					 $T\in \widehat{\mathcal{T}}(c_1)$ is the first barycentric coordinate $\lambda_1$ in $T$.
					A comparison of $f$ in \eqref{eq:proof_lem_C1Tm_k=012_fT} with $f^\prime$ in \eqref{eq:proof_lem_C1Tm_k=012_fprime} reveals 
					$$f=2f^\prime\varphi_1-\gamma{\varphi_1^2}\quad\text{in }C^1(T(1)),$$ 
					where $\gamma|_T:=c_T(2e_1)$ for any $T\in \widehat{\mathcal{T}}(c_1)$ 
					defines the piecewise constant function $\gamma\in P_0(\widehat{\mathcal{T}}(c_1))$. 
					The piecewise quadratic $f$ in $C^1$ and the piecewise affine 
					$\nabla f=2(f^\prime-\gamma\varphi_1)\nabla\varphi_1+2\varphi_1\nabla f^\prime$ are  
					continuous in $T(1)$. In particular $\nabla f$ is continuous 
					at the boundary $\partial T(1)\setminus\partial K:=\bigcup_{1\le j<k\le 4} F_{jk}\subset\{\varphi_1=0\}$ 
					for $F_{jk}:=\textup{conv}\{c_K,Q_j,Q_k\}$
					without the face $F_1\in\mathcal{F}(K)$, where $\varphi_1=0$ vanishes. 
					Hence $f^\prime\nabla \varphi_1$ is continuous there as well. 
					Since $\nabla\varphi_1$ jumps at $c_K$ and $Q_j\in\partial T(1)$, this implies 
					$f^\prime(Q_j)=0=f^\prime(c_K)$ for any $2\le j\le 4$. 
					Thus the affine function $f^\prime\in P_1(T(1))$ vanishes and $f=\gamma\varphi_1^2\in C^1(K)$ leads to $\gamma=0$ as well. 
					This shows that $f\in C^1(T(1))\cap P_2(\widehat{\mathcal{T}}(c_1))$ with $f|_{F_{jk}}=0$ for all $1\le j<k\le 4$ implies 
					$f\equiv 0$. 
				\end{proof}
			\begin{lemma}\label{lem:C1_k=012}
				Let $\widehat{\mathcal{T}}=\textit{WF3D}(K)$ 
				be a $\textit{WF}$ 
				partition of $K$ and $k=0,1,2$. Then $f\in P_k(\widehat{\mathcal{T}})$ satisfies 
				$f\in C^1(K)$
				 if and only if $f\in P_k(K)$ is a global polynomial in $K$.
			\end{lemma}
				\begin{proof}
					The proof is similar to that of \cref{lem:C1Tm_k=012}. 
					The assertion is trivial for $k=0$ and immediate for $k=1$. 
					Suppose $k=2$ and $f\in C^1(K)\cap P_2(\widehat{\mathcal{T}})$. 
					In particular  $f|_{T(m)}\in C^1(T(m))\cap P_2(\widehat{\mathcal{T}}(c_m))$  
					for any $m=1,2,3,4$ and \cref{lem:C1Tm_k=012} shows $f|_{T(m)}\in P_2(T(m))$. 
					Hence $f\in P_2(\widetilde{\mathcal{T}})$ 
					holds for $\widetilde{\mathcal{T}}:=\{T(1),\dots, T(4)\}$.
					The $5$ vertices in the triangulation $\widetilde{\mathcal{T}}$ are $c_K,Q_1,Q_2,Q_3,Q_4$ and 
					each subtetrahedra $T(m)=\textup{conv}\{c_K,Q_j,Q_k,Q_\ell\}\in\widetilde{\mathcal{T}}$ 
					for $1\le j<k<\ell\le 4$ and $\{j,k,\ell, m\}=\{1,2,3,4\}$ will be identified with $(c_K,Q_j,Q_k,Q_\ell)$.
					The Bernstein basis of $f|_{T(m)}\in P_2(T(m))$ leads to ordinates $(c_m(\alpha)|\alpha\in A_2)$ in  
					$$f|_{T(m)}=2\sum_{\alpha\in A_2}c_m(\alpha)\lambda^\alpha/\alpha!\text{ at }
					x=\lambda_1c_K+\lambda_2Q_j+\lambda_3Q_k+\lambda_4Q_\ell\in T(m)\in\widetilde{\mathcal{T}}$$ for $m=1,\dots,4$. 
					The subtraction of the global quadratic Lagrange interpolation on $K$ shows that, without loss of generality,  
					$f$ vanishes at the vertices $Q_1,\dots, Q_4$ and 
					the edge 
					midpoints  $Q_{jk}:=(Q_j+Q_k)/2$ for $1\le j<k\le 4$.	
					Since $f$ is a quadratic function on $\partial K=F_1\cup\dots\cup F_4=\{\lambda_1=0\}$,  
					 $$f|_{F_m}=2\sum_{\substack{\alpha\in A_2\\\alpha_1=0}}
					 			c_m(\alpha){\lambda^{\alpha}}/{\alpha!} =0 \quad\text{at }x=\lambda_2Q_j+\lambda_3Q_k+\lambda_4Q_\ell \in F_m$$
					 vanishes. Hence $c_m(0;\beta)=0$ for all $\beta\in \mathbb{N}_0^3$ with $|\beta|=2$  and $\ell=1,\dots,4$.
					 Consequently, 
						 \begin{align*}
						 	f|_{T(m)}=2\lambda_1\sum_{\beta\in A_1}\frac{c_m(e_1+\beta)}{(1+\beta_1)}\frac{\lambda^\beta}{\beta!}
						 	\quad\text{ at }x=\lambda_1c_K+\lambda_2Q_j+\lambda_3Q_k+\lambda_4Q_\ell\in T(m)\in\widetilde{\mathcal{T}}.
						\end{align*}
					Since all $T(m)\in\widetilde{\mathcal{T}}$ are face connected with common vertex $c_K$, 
					define $f^\prime$ for $\xi=1$ as in \eqref{eq:gprimeon T1T2}. Repeated applications of  \cref{lem:gprime_C1_onT1T2} 
					show   
					$f^\prime \in C^1(K)\cap P_1(\widehat{\mathcal{T}})$ and therefore
					$f^\prime \in P_1(K)$ by the proven assertion of this lemma. 
					The comparison of $f$ and  $f^\prime$ concludes the proof as that of \cref{lem:C1Tm_k=012}. Since the arguments apply verbatim, 
					further details are omitted.  
				\end{proof}
	
	\subsection{The $\textit{WF}$ finite element}\label{sec:HCT3D=FE}
		Given a tetrahedron $K$ and its $\textit{WF}$ 
		partition $\widehat{\mathcal{T}}=\textit{WF3D}(K)$ 
		from \cref{def:CTtriangulation}, this section 
		characterizes
		$$\textit{WF}(K)
		:=P_3(\widehat{\mathcal{T}})\cap C^1(K)$$
		as a $28$-dimensional vector space such that the $28$ degrees of freedom (dof) $L_1,\dots,L_{28}$ in \cref{tab:CT_3D_dof} are linear 
		independent and the triple $\big(K,\textit{WF}(K)
		,(L_1,\dots,L_{28})\big)$ forms a finite element in the sense of Ciarlet. 
		To describe the degrees of freedom from \cref{tab:CT_3D_dof} for the tetrahedron $K=\textup{conv}\{Q_1,\dots,Q_4\}$, determine,
		for each edge $E_{jk}=\textup{conv}\{Q_j,Q_k\}$  of $K$ with $1\le j<k\le 4$ and with midpoint $Q_{jk}=(Q_j+Q_k)/2$  a unit tangential vector 
		$\tau_{{jk}}\parallel Q_j-Q_k$ and two additional unit vectors $\nu_{jk,1},\nu_{jk,2}\in \mathbb{R}^3$ such that 
		$\textup{span}\{\tau_{{jk}},\nu_{jk,1},\nu_{jk,2}\}=\mathbb{R}^3$. Without loss of generality, assume that 
		$(\tau_{{jk}},\nu_{jk,1},\nu_{jk,2})$ is an orthonormal basis of $\mathbb{R}^3$. 
		The degrees of freedom with the enumeration in \cref{tab:CT_3D_dof} are the evaluation of a function $f\in C^1(K)$ and its gradient 
		$\nabla f$ at the vertices $Q_1,\dots, Q_4$ of $K$ and  the evaluation of the gradient 
		$\nabla f(Q_{jk})\cdot s$ in the direction $s=\nu_{jk,1}$ and $s=\nu_{jk,2}$ at each edge midpoint $Q_{jk}$.
			\begin{table}
			 	\begin{center}
			 		\setlength{\fboxsep}{0.5ex}
					\fbox{\parbox{0.94\linewidth}{
						\begin{align*}
						\begin{array}{ll}
							L_\mu(f)=f(Q_\mu) &\text{for } \mu=1,\dots, 4\\
					 		L_{4\ell+\mu}(f)=\nabla f(Q_\mu)\cdot e_\ell&\text{for }\mu=1,\dots, 4,\ \ell=1,2,3 \\
							L_{20-2^{4-j}+2k+m}(f)=\nabla f(Q_{jk})\cdot \nu_{jk,m}&\text{for } j=1,2,3, \  k=j+1,\dots,4,\  m=1,2
						\end{array}
						\end{align*}}
						}
				\caption{The $28$ degrees of freedom  $L_1,\dots, L_{28}$ of $\textit{WF}(K)$ 
						for any $f\in C^1(K)$.}	
				\label{tab:CT_3D_dof}
				 \end{center}
			\end{table}
		Throughout this section,  the four vertices $Q_1,\dots, Q_4\in\mathcal{V}$ of the tetrahedron $K$  are opposite to the faces 
		$F_1,\dots, F_4\in \mathcal{F}$.
		The  \cref{alg:CT3D} below with input 
		$x_1,\dots, x_{28}\in\mathbb{R}$ for the $28$ degrees of freedom computes $f\in P_3(\widehat{\mathcal{T}})$ in terms of ordinates 
		$(c_T(\alpha)|\alpha\in A_3,T\in\widehat{\mathcal{T}})$ in 
			\begin{align}
				f|_T=6\sum_{\alpha\in A_3}c_T(\alpha)\frac{\lambda^\alpha}{\alpha!}\quad
				\text{at }x=\sum_{\mu=1}^4\lambda_\mu P_{\sigma(T,\mu)}\in T\in\widehat{\mathcal{T}}	\label{eq:finAlg}
			\end{align}
		with the global enumeration $P_1,\dots,P_9$ of all the vertices in $\widehat{\mathcal{T}}$ with $P_1=c_K$. 
		For each $T\in\widehat{\mathcal{T}}$ let $P_{\sigma(T,1)}=P_1=c_K$, 
		$P_{\sigma(T,2)}=c_m$ for the unique $m\in\{1,2,3,4\}$ with $\partial T\cap\partial K\subset F_m$, and the remaining vertices 
		$P_{\sigma(T,3)},P_{\sigma(T,4)}\in \{Q_1,\dots,Q_4\}\setminus\{Q_m\}$. 
		The algorithm is contained in \cite{WF87}. 
			\begin{breakablealgorithm}
				\caption{Worsey-Farin-CT3D }\label{alg:CT3D}
				 \algorithmicrequire{ 	$x_1,\dots,x_{28}\in\mathbb{R}$, a tetrahedron
										 $K=\textup{conv}\{Q_1,Q_2,Q_3,Q_4\}$ 
										with $\textit{WF}$ 
										partition $\widehat{\mathcal{T}}=\textit{WF3D}(K)$ 
										and a local  enumeration so that 
										$T\equiv (P_{\sigma(T,1)},\dots,P_{\sigma(T,4)})$ with $P_{\sigma(T,1)}=c_K$, $P_{\sigma(T,2)}=c_m$, 
										and $P_{\sigma(T,3)},P_{\sigma(T,4)}\in \{Q_1,\dots,Q_4\}\setminus\{Q_m\}$ 
										for each 
										$T\in\widehat{\mathcal{T}}(c_m)$ with $m\in\{1,\dots,4\}$
										}
				\begin{enumerate}[label=(\alph*), ref=\alph*]
					\item\label{item:Alg_CT_preprocessing} 
						A preprocessing computes $f(Q_{jk})$ and $\nabla f(Q_{jk})$ with $Q_{jk}:=(Q_j+Q_k)/2$  for any
						$j,k=1,\dots,4$, note  $Q_{jj}\equiv Q_j$, from the data $x_1,\dots,x_{28}$ representing $L_1(f),\dots, L_{28}(f)$ 
						(cf. \cref{rem:preprocessing} below for further details). 
					\item\label{item:Alg_CT_gP2} 
						For all $\ell=1,\dots,4$ and $1\le j<k\le 4$, 
							\begin{align*}
								g(Q_\ell)&:=f(Q_\ell)-\frac{1}{3} (Q_\ell-c_K)\cdot\nabla f(Q_\ell),\\
								g(Q_{jk})&:=\frac{1}{3}(c_K-Q_{jk})\cdot \nabla f(Q_{jk})+f(Q_{jk}).
							\end{align*}
						Let $g\in P_2(K)$ denote  the quadratic Lagrange interpolation of 
						this data at the vertices $(Q_1,\dots,Q_4)$ and edge midpoints $(Q_{jk}:\, 1\le j<k\le 4)$ of $K$ and compute the ordinates
						$(c_T(\alpha)|\alpha\in A_3, \alpha_1\ge 1, T\in\widehat{\mathcal{T}})$ of 
							\begin{align*}
								g|_T=2\sum_{\beta\in A_2}c_T(e_1+\beta)\frac{\lambda^\beta}{\beta!}\quad
								\text{at }x=\sum_{\mu=1}^4\lambda_\mu P_{\sigma(T,\mu)}\in T\in\widehat{\mathcal{T}}.
							\end{align*}
					\item\label{item:Alg_CT_HCT}
						For all $m=1,\dots,4$ let $f_m\in HCT(F_m):=P_3(\widehat{\mathcal{F}}(F_m))\cap C^1(F_m)$ denote the $2$D $HCT$ finite 
						element interpolation of the input data (resp. data from (\ref{item:Alg_CT_preprocessing})) on the face $F_m$ of $K$ 
						with unit normal $\nu_{F_m}$.
						The degrees for freedom of $HCT(F_m)$ are the  
						 evaluation of $f$ and the tangential derivatives $\nabla f\times \nu_{F_m}$ 
						at the vertices  of the triangle $F_m$ as well as for each edge $E\in\mathcal{E}(F_m)$ with unit tangent vector $\tau_E$ of $F_m$ 
						the directional derivative in $\langle F_m\rangle$ perpendicular to the edge $E$ ($\nabla f\cdot(\tau_E\times\nu_{F_m})$)
						evaluated in the edge midpoint. 
						Compute coefficients $(c_T(\alpha)|\alpha\in A_3,\alpha_1=0, T\in\widehat{\mathcal{T}})$ such that 
							\begin{align*}
								f_m|_{F_m\cap \partial T}=6\sum_{\substack{\alpha\in A_3\\\alpha_1=0}}c_T(\alpha)\frac{\lambda^\alpha}{\alpha!}\quad
								\text{at }x=\sum_{\mu=2}^4 \lambda_\mu P_{\sigma(T,\mu)}\in F_m\cap\partial T\text{ for } T\in\widehat{\mathcal{T}}.
							\end{align*}
				\end{enumerate}
				 \algorithmicensure{  $(c_T(\alpha)|\,\alpha\in A_3,T\in\widehat{\mathcal{T}})$ and $f\in P_3(\widehat{\mathcal{T}})$ with
					$f|_T=6\sum_{\alpha\in A_3}c_T(\alpha){\lambda^\alpha}/{\alpha!}$ 
					 at $x=\sum_{\mu=1}^4\lambda_\mu P_{\sigma(T,\mu)}\in T\in\widehat{\mathcal{T}}$.
					}
			\end{breakablealgorithm}
			\FloatBarrier
			\begin{remark} [preprocessing in \cref{alg:CT3D}.\ref{item:Alg_CT_preprocessing}]\label{rem:preprocessing}
			Algorithm~\ref{alg:CT3D} aims at the computation of some global $C^1$ function $f$ with prescribed values $L_\mu(f)=x_\mu$ for $\mu=1,\dots,28$. The step 
			(\ref{item:Alg_CT_preprocessing}) provides $f(Q_{jk})$ and $\nabla f(Q_{jk})$ for any $Q_{jk}$ as follows. \\
			\emph{Step~a.1 for $j=k$.} 
			Then $Q_{jj}=Q_j$ is a vertex and $x_j$ directly represents $f(Q_j)=x_j$ and $x_{4\ell+j}$ represents 
			$\frac{\partial f}{\partial x_\ell}(Q_j)=x_{4\ell+j}$ for $\ell=1,2,3$. \\
			\emph{Step~a.2 for $j<k$.} Then $Q_{jk}$ is the midpoint of the edge $E=\textup{conv}\{Q_j,Q_k\}$ from $Q_j$ to $Q_k$ 
			with unit tangent vector $\tau_{jk}$ 
			in the tetrahedron $K$. 
			We consider $f|_E\in P_3(E)$ as a one-dimensional cubic polynomial that satisfies the following interpolation conditions. 
			The function values $f(Q_\ell)=x_\ell$ and 
			the one-dimensional derivative $f^\prime(Q_\ell)=\nabla f(Q_{\ell})\cdot\tau_{jk}$ at the vertices $Q_j$ and $Q_k$ (for $\ell=j,k$) 
			are prescribed in Step~a.1. 
			The four one-dimensional Hermite interpolation conditions determine a unique cubic polynomial $f|_E\in P_3(E)$ along the edge $E$.  
			This polynomial provides the values $f(Q_{jk})$ and 
			$f^\prime(Q_{jk})=\nabla f(Q_{jk})\cdot \tau_{jk}$. The final components of $\nabla f( Q_{jk})$ are prescribed by 
			$x_{20-2^{4-j}+2k+m}=\nabla f(Q_{jk})\cdot\nu_{jk,m}$ for $m=1,2$ and the fixed vectors $\nu_{jk,m}$   
			with $\textup{span}\{\tau_{{jk}},\nu_{jk,1},\nu_{jk,2}\}=\mathbb{R}^3$.    
			\end{remark}
			\begin{remark}[feasability of \cref{alg:CT3D}]\label{rem:feasbility_CT3D}
				The data in  \cref{alg:CT3D}.\ref{item:Alg_CT_gP2} define a unique $g\in P_2(K)$, which has a unique representation in the Bernstein 
				basis of each $T\in\widehat{\mathcal{T}}$. 
			    Since the $2$D $HCT$ finite element  is  a finite element in the sense of Ciarlet \cite[Thm.~6.1.2]{Ciarlet78}, 
			    $F_m\cap\partial T=:F\in\widehat{\mathcal{F}}(F_m)$ and 
				$f_m|_F\in P_3(F)$ allow for unique ordinates in \cref{alg:CT3D}.\ref{item:Alg_CT_HCT} 
				(cf. \cite{Me12} for details on the implementation of a reduced 
				HCT element in $2$D). 
			\end{remark}
			\begin{remark}[parameter dependence of \cref{alg:CT3D}]\label{rem:Alg_continuity}
				Suppose the input $x_1,\dots,$ $x_{28}\in\mathbb{R}$ of \cref{alg:CT3D} is fixed, while the center points $c_K\in\textup{int}(K)$ and 
				$c_m\in \textup{relint}(F_m)$ for $m=1,2,3,4$ may vary, i.e., the $\textit{WF}$ 
				partition $\widehat{\mathcal{T}}=\textit{WF3D}(K)$ 
				differs. 
				The point is that the output $f\in P_3(\widehat{\mathcal{T}})\cap C^1(\Omega)$ depends continuously on these five parameters 
				$c_K,c_1,\dots,c_4\in\mathbb{R}^{3}$.
				The preprocessing in \cref{alg:CT3D}.\ref{item:Alg_CT_preprocessing} depends only on the edges of $K$ (see \cref{rem:preprocessing}) 
				and is independent of the subtriangulation $\widehat{\mathcal{T}}$. 
				The values of $g\in P_2(K)$ in the vertices and edge midpoints in \cref{alg:CT3D}.\ref{item:Alg_CT_gP2} depend directly on the choice of $c_K$. 
				This dependence in $g$ is affine, whence continuous. 
				The $2$D $HCT$ function $f_m$ on the face $F_m$ of $K$ for $m=1,\dots,4$ in \cref{alg:CT3D}.\ref{item:Alg_CT_HCT} 
				depends on the choice of the face center $c_m$. The coefficients of the function $f_m\in P_3(\widehat{\mathcal{T}})$ for a given $c_m$ 
				are obtained through the solution of a linear system with an invertible matrix \cite[p.345, l.26ff]{Ciarlet78}; this is already exploited in 
				\cite[p.346, l.1ff]{Ciarlet78}. 	
				Hence the dependence of the ordinates $\big(c_T(\alpha):\,\alpha \in A_3, T\in\widehat{\mathcal{T}}\big)$ 
				 on the points $c_K,c_1,\dots,c_4$ that define the subtriangulation $\widehat{\mathcal{T}}$ is continuous.\hfill $\Box$ 
			\end{remark}
			The $\textit{WF}$ 
			finite element is a finite element in the sense of Ciarlet as announced in \cite{NeilanWu2019}, 
			Since \cite{WF87,NeilanWu2019} do not contain the proofs, this is carried out here for convenient reading.
			\begin{theorem}[$\textit{WF}$ 
							finite element]\label{thm:CT_FECiarlet}
				\begin{enumerate}[label=(\alph*), ref=\alph*]
					\item \label{item:CT_dof} For any input $x_1,\dots,x_{28}\in \mathbb{R}$, the output $f\in C^1(K)\cap P_3(\widehat{\mathcal{T}})$ 
							of \cref{alg:CT3D} is 
							continuously differentiable and 
							\begin{align}
								L_\mu f=x_\mu\quad\text{ for }\mu=1,\dots,28. \label{eq:falg_CT3_dof}
							\end{align}
					\item \label{item:CT_uniqueness} There exist at most one function $f\in C^1(K)\cap P_3(\widehat{\mathcal{T}})$ with 
							\eqref{eq:falg_CT3_dof}.
					\item \label{item:CT_Ciarlet}The triple $(K,C^1(K)\cap P_3(\widehat{\mathcal{T}}),(L_1,\dots,L_{28}))$ 
							is a finite element in the sense of Ciarlet. 
				\end{enumerate}
			\end{theorem}
				\begin{proofof}\textit{(\ref{item:CT_dof}).} This proof is divided in two steps. 
				\begin{enumerate}[label=\textit{Step \arabic*}, ref=\textit{Step \arabic*}, wide]
				\item \textit{($L_\mu f=x_\mu$ for $\mu=1,\dots,28$).}
					The output $f\in P_3(\widehat{\mathcal{T}})$ is a piecewise cubic polynomial. To prove  \eqref{eq:falg_CT3_dof} 
					fix one $T\in\mathcal{T}$, abbreviate $f_T:=f|_T$, 
					and check $L_\mu f_T=x_\mu$ 
					for any degree of freedom $L_\mu$, which can be evaluated for $f_T$.
					
					Suppose  $T:=\textup{conv}\{c_K,c_m,Q_j,Q_k\} \in\widehat{\mathcal{T}}$ for $1\le j <k\le 4$ and 
					$m\in\{1,\dots, 4\}\setminus\{j,k\}$ and identify $T$ with $(c_K,c_m,Q_j,Q_k)$.
					The output $f$ on $T\in\widehat{\mathcal{T}}$ reads 
						\begin{align*}
							f_T:=f|_T=6\sum_{\alpha\in A_3}c_T(\alpha){\lambda^\alpha}/{\alpha!}\quad\text{ at }
							x=\lambda_1c_K+\lambda_2c_m+\lambda_3 Q_j+\lambda_4 Q_k\in T.
						\end{align*}
					We can evaluate $10$ degrees of freedom from \cref{tab:CT_3D_dof} in $T:=\textup{conv}\{c_K,c_m,Q_j,Q_k\}$, 
					namely $L_\mu$, $L_{4+\mu}$, $L_{8+\mu}$, $L_{12+\mu}$ for $\mu=j$ and $\mu=k$,  
					as well as $L_{20-2^{4-j}+2k+1}$ and $L_{20-2^{4-j}+2k+2}$.
					\begin{proofof}\textit{\eqref{eq:falg_CT3_dof} for $\mu=j,k$. }
					Since  the $2$D HCT finite element in \cref{alg:CT3D}.\ref{item:Alg_CT_HCT} 
					interpolates the nodal values exactly, \cref{exp:k=3_m=4}.\ref{item:fPj} shows 
					\begin{align*}
						L_j(f_T)&=f_T(Q_j)=c_T(3e_3)=f_m|_T(Q_j)=x_j.
					\end{align*}
					The proof of $L_k(f_T)=x_k$ follows by symmetry of $j$ and $k$. 
					Since $L_j(f_T)=x_j$ for any $T\in\widehat{\mathcal{T}}$, the output function $f \in P_3(\widehat{\mathcal{T}})$ 
					is single valued and continuous at every vertex $Q_j\in \mathcal{V}$ of $K$ for $j=1,\dots, 4$   with $f(Q_j)=x_j$.
					\end{proofof}
					\begin{proofof}\textit{continuity at $Q_{jk}$. }					
					\cref{rem:preprocessing} explains that the data  $x_1,\dots,x_{28}$ 
					allow the computation of  a value $x_{jk}$ in the edge midpoint in \cref{alg:CT3D}.\ref{item:Alg_CT_preprocessing}. 
					\cref{alg:CT3D}.\ref{item:Alg_CT_HCT} interpolates $x_{jk}$ such that $f_m|_{F_m}(Q_{jk})=x_{jk}$. 
					For any $T\in\widehat{\mathcal{T}}(T(m))$ identified with $(c_K,c_m,Q_j,Q_k)$ as above, 
					\cref{exp:k=3_m=4}.\ref{item:fPkell} shows 
 					\begin{align}
						f_T(Q_{jk})&=\frac{1}{8}\big(f_T(Q_k)+3c_T(0,0,2,1)+3c_T(0,0,1,2)+f_T(Q_j)\big)\notag\\
										&=\frac{1}{8}\big(x_k+3c_T(0,0,2,1)+3c_T(0,0,1,2)+x_j\big)= f_m|_{F_m}(Q_{jk})=x_{jk}. \label{eq:proof_xjk} 
 					\end{align}
 					The result $f_T(Q_{jk})=x_{jk}$ is independent of $m\in\{1,2,3,4\}\setminus\{j,k\}$ and $T\in\widehat{\mathcal{T}}(T(m))$. 
 					Hence the output function $f \in P_3(\widehat{\mathcal{T}})$ 
					is single valued  and continuous at every edge midpoint $Q_{jk}=(Q_j+Q_k)/2$ of $K$ for $1\le j<k\le 4$   with $f(Q_{jk})=x_{jk}$.
					\end{proofof}
					\begin{proofof}\textit{\eqref{eq:falg_CT3_dof} for $\mu=4+j,8+j,12+j,4+k,8+k,12+k $. }	
					Abbreviate $X_j:=(x_{4+j},x_{8+j},x_{12+j})^\top$, recall $T\in\widehat{\mathcal{T}}(T(m))$ is identified with $(c_K,c_m,Q_j,Q_k)$.  
					Observe that the gradient  $\nabla f_T(Q_j)\in\mathbb{R}^3$ is uniquely determined by its scalar product with  
					$Q_j-Q_k$, $Q_j-c_m$, and $Q_j-c_K$. 
					The $2$D $HCT$ interpolation in \cref{alg:CT3D}.\ref{item:Alg_CT_HCT}  interpolates the tangential derivatives correctly.
					In other words $\nabla f_m|_{F_m}(Q_j)\times\nu_{F_m}=X_j\times\nu_{F_m}$ for the unit normal $\nu_{F_m}$ of $F_m$. 
					 The first two  directions $Q_j-Q_k$ and $Q_j-c_m$ are parallel to the plane  $\langle F_m\rangle$  and 
					 \cref{exp:k=3_m=4}.\ref{item:nablafPj_1} shows 
					\begin{align*}
						(Q_j-Q_k)\cdot\nabla f_T(Q_j)&
						=3\big(f_T(Q_j)-c_T(2e_3+e_4)\big)=3\big(x_j-c_T(2e_3+e_4)\big)\\&
						=(Q_j-Q_k)\cdot\nabla f_m|_{F_m}(Q_j)=(Q_j-Q_k)\cdot X_j,\\
						(Q_j-c_m)\cdot\nabla f_T(Q_j)
						&=3\big(f_T(Q_j)-c_T(2e_3+e_2)\big)=3\big(x_j-c_T(2e_3+e_2)\big)\\&=(Q_j-c_m)\cdot\nabla f_m|_{F_m}(Q_j)=(Q_j-c_m)\cdot X_j.
					\end{align*}
					The third direction $Q_j-c_K$ is not parallel to the plane $\langle F_m\rangle$.  
					\cref{exp:k=3_m=4}.\ref{item:nablafPj_1} and 
					\ref{exp:k=2_m=4}.\ref{item:gPj} imply
						$$
							(Q_j-c_K)\cdot\nabla f_T(Q_j)=3\big(f_T(Q_j)-c_T(2e_3+e_1)\big)=3\big(x_j-g(Q_j)\big).
						$$
					The choice of the nodal value 
					$g(Q_j)\equiv x_j-\frac{1}{3} (Q_j-c_K)\cdot X_j$ in \cref{alg:CT3D}.\ref{item:Alg_CT_gP2} 
					implies $$(Q_j-c_K)\cdot\nabla f_T(Q_j)=(Q_j-c_K)\cdot X_j.$$ This 
					concludes the proof of $\nabla f_T(Q_j)=X_j=(x_{4+j},x_{8+j},x_{12+j})^\top$. 
					The proof of $\nabla f_T(Q_k)=(x_{4+k},x_{8+k},x_{12+k})^\top$ follows by symmetry of $j$ and $k$. 					
					Hence the gradient of the output function $f \in P_3(\widehat{\mathcal{T}})$ 
					is single valued and continuous at every vertex $Q_j\in \mathcal{V}$ of $K$ for $j=1,\dots, 4$   with $\nabla f(Q_j)=X_j$.
					\end{proofof}
 					\begin{proofof}{\eqref{eq:falg_CT3_dof} for $\mu=20-2^{4-j}+2k+1,20-2^{4-j}+2k+2$. }
					The remaining degrees of freedom $L_{20-2^{4-j}+2k+1}$ and $L_{20-2^{4-j}+2k+2}$ 
					concern the gradient $\nabla f_T(Q_{jk})$ in the edge midpoint $Q_{jk}$. 
					 \cref{rem:preprocessing} explains how  $x_1,\dots,x_{28}$ determine a unique value $X_{jk}$ in 
					  \cref{alg:CT3D}.\ref{item:Alg_CT_preprocessing} that represents $\nabla f(Q_{jk})$ 
					 in \cref{alg:CT3D}.\ref{item:Alg_CT_gP2}--\ref{item:Alg_CT_HCT}. Recall $f(Q_{jk})=x_{jk}$ from \eqref{eq:proof_xjk} above. 
					 We want to verify that $\nabla f_T(Q_{jk})=X_{jk}$ for $T\in\widehat{\mathcal{T}}(T(m))$ identified with $(c_K,c_m,Q_j,Q_k)$. 
					The  gradient $\nabla f_T(Q_{jk})$ of the output of \cref{alg:CT3D} is uniquely determined by its scalar product with 
					$c_m-Q_{jk}$, $Q_j-Q_{jk}$, and $c_K-Q_{jk}$.	
					In the directions $c_m-Q_{jk}$ and $Q_j-Q_{jk}$,  parallel to the plane $\langle F_m\rangle$,  the gradient is determined by 
					 $f_m$ in \cref{alg:CT3D}.\ref{item:Alg_CT_HCT}.  
					\cref{exp:k=3_m=4}.\ref{item:nablafPkell}--\ref{item:nablafPkell_2} lead to 
					\begin{align*}
						(c_m-Q_{jk})\cdot\nabla f_T(Q_{jk})&=\frac{3}{4}\big(c_T(0,1,2,0)+2c_T(0,1,1,1)+c_T(0,1,0,2)\big)-3x_{jk}
							\\&=(c_m-Q_{jk})\cdot\nabla f_m|_{F_m}(Q_{jk})=(c_m-Q_{jk})\cdot X_{jk}. \\
						(Q_j-Q_{jk})\cdot\nabla f(Q_{jk})|_T&=\frac{3}{8}\big(x_j+c_T(2e_\gamma+e_\kappa)-c_T(2e_\kappa+e_\gamma)-x_k)
							\\&=(Q_j-Q_{jk})\cdot\nabla f_m|_{F_m}(Q_{jk})=(Q_j-Q_{jk})\cdot X_{jk}.
					\end{align*}
					For the remaining direction $c_K-Q_{jk}$, \cref{exp:k=3_m=4}.\ref{item:nablafPkell} and 
					\ref{exp:k=2_m=4}.\ref{item:gPkell} show
						\begin{align*}
							(c_K-Q_{jk})\cdot\nabla f_T(Q_{jk})&=\frac{3}{4}\big(c_T(e_1+2e_3)+2c_T(e_1+e_3+e_4)+c_T(e_1+2e_4)\big)-3f_T(Q_{jk})\\
								&=3 g(Q_{jk})-3x_{jk}. 
						\end{align*}
					Since the quadratic polynomial $g\in P_2(K)$ in \cref{alg:CT3D}.\ref{item:Alg_CT_gP2}  assumes the value 
					$g(Q_{jk})\equiv\frac{1}{3}(c_K-Q_{jk})\cdot X_{jk}+x_{jk}$ at $Q_{jk}$, it follows $(c_K-Q_{jk})\cdot\nabla f_T(Q_{jk})=(c_K-Q_{jk})\cdot X_{jk}$.
					The last three displayed identities  
					ensure $\nabla f_T(Q_{jk})=X_{jk}$ for any $T\in\widehat{\mathcal{T}}(Q_{jk})$.  
					Hence the gradient of the output function $f \in P_3(\widehat{\mathcal{T}})$ 
					is indeed single valued and continuous at every edge midpoint $Q_{jk}=(Q_j+Q_k)/2$ of $K$ for $1\le j<k\le 4$   with $\nabla f(Q_{jk})=X_{jk}$. 
					In particular, $ L_{20-2^{4-j}+2k+1} (f_T)=\nabla f_T(Q_{jk})\cdot\nu_{jk,1}=x_{20-2^{4-j}+2k+1}$ and 
					$ L_{20-2^{4-j}+2k+2} (f_T)=\nabla f_T(Q_{jk})\cdot\nu_{jk,2}=x_{20-2^{4-j}+2k+2}$. 
					This concludes the proof of \eqref{eq:falg_CT3_dof}.
	 				\end{proofof}
				\item \textit{($f\in C^1(K)$).}  It suffices to show $f|_{T_1\cup T_2}\in C^1(T_1\cup T_2)$ 
				for any two tetrahedra $T_1,T_2\in\widehat{\mathcal{T}}$ with common face 
				$F:=\partial T_1\cap \partial T_2$.
				There are two types of common faces in $\widehat{\mathcal{T}}$.
				\begin{enumerate}[label=\textit{Case \arabic*.}, wide]
						\item 
						Assume $T_1,T_2\in \widehat{\mathcal{T}}(c_m)$ for $m=1,\dots, 4$.
						In other words, let  $F=\textup{conv}\{c_K,c_m,Q_j\}$ and identify 
						$T_1$ with $(c_K,c_m,Q_j,Q_k)$ and $T_2$ with $(c_K,c_m,Q_j,Q_\ell)$ for $\{j,k,\ell,m\}=\{1,2,3,4\}$.
						The definition of $g\in P_2(K)\subset C^1(K)$ in \cref{alg:CT3D}.\ref{item:Alg_CT_gP2} 
						ensures $g|_{T_1\cup T_2}=f^\prime|_{T_1\cup T_2}$ for 
						$\xi=1$ in \eqref{eq:gprimeon T1T2}. 
						Hence by \cref{lem:gprime_C1_onT1T2}, $f|_{T_1\cup T_2}\in C^1(T_1\cup T_2)$ is equivalent to the 
						continuity of $f$ and $\nu_F\cdot\nabla f $ along $E:=\textup{conv}\{c_m,Q_j\}$. 
						It remains to verify  (i) continuity of $f$ along $E$ and (ii) continuity of  $\nu_F\cdot\nabla f $ along $E$. 
						\\\textit{Proof of (i). }
						Abbreviate  $c_\mu(\alpha):=c_{T_\mu}(\alpha)$ for $\mu=1,2$.  
						For $\alpha\in A_3$ with 
						$\alpha_1=0=\alpha_4$, the continuity of $f_m$ from \cref{alg:CT3D}.\ref{item:Alg_CT_HCT} shows $c_1(\alpha)=c_2(\alpha)$.
						For $\alpha_1\not=0=\alpha_4$, the continuity of $g$ from  
						\cref{alg:CT3D}.\ref{item:Alg_CT_gP2} shows $c_1(\alpha)=c_2(\alpha)$.  
						In particular, \eqref{eq:con1_C1T1T2} implies that $f|_E$ is continuous at $E$.  \hfill $\Box$
						 \\\textit{Proof of (ii). }
						 Note, since the ordinates $c_1(e_2+\alpha)=c_1(e_2+\alpha)$ are equal for all $\alpha\in A_2$ with $\alpha_4=0$, 
						\cref{lem:g_C1_onT1T2} proves the continuity of the piecewise quadratic function $g^{\prime\prime}\in P_2(\{T_1,T_2\})$
						defined by 
						\begin{align}
								g^{\prime\prime}|_T:=2\sum_{\beta\in A_2}c_T(e_2+\beta){\lambda^\beta}/{\beta!} \quad \text{ in }T\in\{T_1,T_2\};
								\label{eq:gprimeprime}
						\end{align}
						 $g^{\prime\prime}\in C(T_1\cup T_2)$.  
						In order to prove the continuity of $\nu_F\cdot\nabla f$ along $E$, recall that $f|_{F_m}=f_m\in C^1(F_m)$  on the face 			
						$F_m:=\textup{conv}\{Q_\ell, Q_j,Q_k\}\supset E$ of $K$ from \cref{alg:CT3D}.\ref{item:Alg_CT_HCT}. 
						Hence  the tangential derivative 
						$\nu_{F_m}\times\nabla f$ is continuous along 
						$F_m$ and in particular along $E$. 
						It remains to prove the continuity for one fixed direction  
						$\tau:=c_m-c_K\in \mathbb{R}^3\setminus\{0\}$ oblique to $\langle F_m\rangle$. 
						In $T\in\{T_1,T_2\}$,  \cref{rem:derivativeBernsteinBasis} shows  
						\begin{align*}
						\tau\cdot\nabla f|_T=6\sum_{\beta\in A_2}
							\sum_{\kappa=1}^4 c_T(\beta+e_\kappa)\frac{\lambda^\beta}{\beta!}(\tau\cdot\nabla\lambda_\kappa).
						\end{align*}
						Since $\tau\cdot\nabla\lambda_\kappa=\lambda_\kappa(c_m)-\lambda_\kappa(c_K)=\delta_{2\kappa}-\delta_{1\kappa}$, this is 
						\begin{align*}
						\tau\cdot\nabla f|_T=6\sum_{\beta\in A_2}\big(c_T(\beta+e_2)-c_T(\beta+e_1)\big)\frac{\lambda^\beta}{\beta!}
						=3g^{\prime\prime}|_T-3g|_T
						\end{align*}
						with $g^{\prime\prime}$ from \eqref{eq:gprimeprime} and $g$ from \cref{alg:CT3D}.\ref{item:Alg_CT_gP2} . 
						Since $g^{\prime\prime}$ and $g$ are continuous in $T_1\cup T_2$,  $\nabla f|_{F_m}$ is continuous in particular
						along $E$.	\hfill$\Box$ 		
						\item Assume $T_1,T_2\in \widehat{\mathcal{T}}$ share a face $F=\textup{conv}\{c_K,Q_j,Q_k\}$ and 
						identify $T_1$ with $(c_K,$ $Q_j,$ $Q_k,$ $c_m)$ and $T_2$ with $(c_K,Q_j,Q_k,c_\ell)$ for $\{j,k,\ell,m\}=\{1,2,3,4\}$. 
						The definition of $g\in P_2(K)\subset C^1(K)$ in \cref{alg:CT3D}.\ref{item:Alg_CT_gP2} ensures 
						$g|_{T_1\cup T_2}=f^\prime|_{T_1\cup T_2}$ for 
						$\xi=1$ in \eqref{eq:gprimeon T1T2}. 
						Hence by \cref{lem:gprime_C1_onT1T2}, $f|_{T_1\cup T_2}\in C^1(T_1\cup T_2)$ is equivalent to the 
						continuity of $f$ and $\nu_F\cdot\nabla f $ along $E:=\textup{conv}\{Q_j,Q_k\}$. 
						It remains to verify  (i) continuity of $f$ along $E$ and (ii) continuity of  $\nu_F\cdot\nabla f $ along $E$. 
						\\\textit{Proof of (i). }
						Since $f|_E$ is uniquely determined in 
						\cref{alg:CT3D}.\ref{item:Alg_CT_preprocessing}, 			 
						the equality of $f_m(x)=f_\ell(x)$ at $x=\lambda_2Q_j+\lambda_3Q_k\in E$
						 proves $c_1(\alpha)=c_2(\alpha)$  for $\alpha\in A_3$ with $\alpha_1=0=\alpha_4$. The equality of the ordinates shows the 
						 continuity of $f$ along $E$. \hfill $\Box$
						 \\\textit{Proof of (ii). }
						Since (each component of) $\nabla f|_{T_1}$ and $\nabla f|_{T_2}$ 
						is a quadratic polynomial along $E$, which 
					 	coincides at the three points $Q_j$, $Q_{jk}$, $Q_k$, it follows $\nabla f|_{T_1}=\nabla f|_{T_2}$ on $E$.  
					 	In particular $\nu_F\cdot\nabla f $ is continuous along along $E$.
						\end{enumerate} 
				\end{enumerate}
				\end{proofof}
				
				\begin{proofof}\textit{(\ref{item:CT_uniqueness}).}
					Suppose two functions $f_1,f_2\in P_3(\widehat{\mathcal{T}})\cap C^1(K)$ satisfy $L_\mu(f_j)=x_\mu$ for $\mu=1,\dots,28$ and $j=1,2$. 
					Their difference $f:=f_2-f_1$ satisfies $L_\mu(f)=0$ for all $\mu=1,\dots,28$.  
					Identify $T:=\textup{conv}\{c_K,c_m,Q_j,Q_k\}\in \widehat{\mathcal{T}}$ 
					with $(c_K,c_m,Q_j,Q_k)$ for $1\le j<k\le 4$ and $m\in\{1,2,3,4\}\setminus\{j,k\}$. 
					Let $(c_T(\alpha)|\,\alpha\in A_3, T\in\widehat{\mathcal{T}})$ denote the coefficients of $f\in C^1(K)$ in 
					$$f|_T=6\sum_{\alpha\in A_3}c_T(\alpha){\lambda^\alpha}/{\alpha!} \quad
					\text{ at }x=\lambda_1c_K+\lambda_2c_m+\lambda_3Q_j+\lambda_4Q_k\in T.$$
					Since $L_1(f)=\dots=L_{28}(f)=0$, \cref{rem:feasbility_CT3D} implies $f(Q_\ell)=f(Q_{jk})=0$ and 
					$\nabla f(Q_\ell)=\nabla f(Q_{jk})=0$ for all $\ell=1,\dots,4$ and $1\le j<k\le 4$. 
					The restriction $f|_F$ of $f$ to a face $F$ of $K$ is a $2$D HCT finite element function with vanishing associated degrees of freedom. 
					Hence $f|_F\equiv 0$. 
					Consequently, 
					$$0=f_T|_{\partial K}=\sum_{\substack{\alpha\in A_3\\ \alpha_1=0}} c_T(\alpha)\frac{\lambda^\alpha}{\alpha!} 
					\quad\text{at }x=\lambda_2 c_m+\lambda_3Q_j+\lambda_4Q_k\in T\cap \partial K$$ 
					vanishes. This implies $c_T(\alpha)=0$ for all $\alpha\in A_3 $ with $\alpha_1=0$. 
					For $\alpha_1\ge 0$, $c_T(\alpha)$ is related to $f^\prime \in P_2(\widehat{\mathcal{T}})$ with 
					\begin{align}
					 	f^\prime|_T:=2\sum_{\beta\in A_2}c_T(e_1+\beta){\lambda^\beta}/{\beta!}\quad\text{ at }
					 	x=\lambda_1c_K+\lambda_2c_m+\lambda_3Q_j+\lambda_4Q_k\in T.
					 	\label{eq:fprime_uniqueness}
					\end{align}
					\cref{lem:gprime_C1_onT1T2} applies to any two $T_+,T_-\in\widehat{\mathcal{T}}$ that share a face 
					$T_+\cap T_-=F\in\widehat{\mathcal{F}}$. Hence \eqref{eq:fprime_uniqueness} defines  $f^\prime \in C^1(K)\cap P_2(\widehat{\mathcal{T}})$ and 
					\cref{lem:C1_k=012} shows $f^\prime \in P_2(K)$. 
					For any $T\in\mathcal{T}$ with vertex $Q_\ell\in \mathcal{V}(T)$ and
					$(\lambda_1,\dots,\lambda_4)=e_\kappa$ for $\kappa\in\{3,4\}$,  \eqref{eq:nablaf(Pj)} in \cref{exp:k=3_m=4} reads
					\begin{align*}
						0=\nabla f(Q_\ell)
							=3\sum_{\mu=1}^4c_T(2e_\kappa+e_\mu)\nabla\lambda_\mu=3 c_T(2e_\kappa+e_1)\nabla \lambda_1
					\end{align*}
					with $c_T(\alpha)=0$ for $\alpha_1=0$ in the last equality.
					It follows with \cref{exp:k=2_m=4}.\ref{item:gPj} that $0=c_T(e_1+2e_\kappa)=f^\prime (Q_\ell)$ for any $1\le \ell\le 4$.
					Since for any $T\in\widehat{\mathcal{T}}$ with vertices $Q_j$ and $Q_k$ the edge midpoint $Q_{jk}$ has barycentric coordinates 
					$(\lambda_1,\dots,\lambda_4)=e_3/2+e_4/2$, \eqref{eq:nablaf(Pjk)} in  \cref{exp:k=3_m=4} implies 
					\begin{align*}
						0&=\nabla f(Q_{jk})=\frac{3}{4}	\sum_{\mu=1}^4\big(c_T(2e_3+e_\mu)+2c_T(e_3+e_4+e_\mu)+c_T(2e_4+e_\mu)\big)\nabla\lambda_\mu
						\\&=\frac{3}{4}\big(c_T(2e_3+e_1)+2c_T(e_3+e_4+e_1)+c_T(2e_4+e_1)\big)\nabla\lambda_1
					\end{align*}
					with $c_T(\alpha)=0$ for $\alpha_1=0$ in the last equality. 
					From \cref{exp:k=2_m=4}.\ref{item:gPkell} follows 
					\begin{align*}
						0=c_T(2e_3+e_1)+2c_T(e_3+e_4+e_1)+c_T(2e_4+e_1)
								=4f^\prime(Q_{jk})|_T
					\end{align*}
					for any $T\in\widehat{\mathcal{T}}$ with vertices $Q_j$ and $Q_k$ for $1\le j<k\le 4$. 
					This shows that $f^\prime \in P_2(K)$ vanishes at all Lagrange points for a quadratic polynomial in $K$, whence 
					$f^\prime\equiv 0$.
					It follows $c_T(\alpha)=0$ for all $T\in\widehat{\mathcal{T}}$ and $\alpha\in A_3$; in other words $f\equiv 0$. 
				\end{proofof}
				\begin{proofof}\textit{(\ref{item:CT_Ciarlet}).}
					This follows directly from (\ref{item:CT_dof}), which proves that \cref{alg:CT3D}  allows the construction of 
					$f\in P_3(\widehat{\mathcal{T}})\cap C^1(K)$ for all input data $x_1,\dots,x_{28}$, 
					and (\ref{item:CT_uniqueness}), which proves the uniqueness of the 
					function with the values  $x_1,\dots,x_{28}$ in the  degrees of freedom $L_1,\dots, L_{28}$ in \cref{tab:CT_3D_dof}.
				\end{proofof}
	\subsection{The $\textit{WF}$ on adjacent tetrahedra} 
		This section is devoted to  the compatibility of $\textit{WF}$ 
		for two tetrahedra with a 
		common face. The point of departure is a first Lemma on one tetrahedron. 
		Let  $\widehat{\mathcal{T}}$  denote 
		the $\textit{WF}$ 
		partition of a simplex $K$ with vertices $Q_1,\dots,Q_4$ and opposite faces $F_1,\dots,F_4$. 
		Recall, $c_m\in\mathrm{relint}(F_m)$ and $c_K\in\mathrm{relint}(K)$ from the $\textit{WF}$ 
		partition in \cref{def:CTtriangulation}, $T(m)=\textup{conv}(c_K, F_m)$, and 
		abbreviate $\tau:=c_m-c_K\in\mathbb{R}^3\setminus\{0\}$ as well as $\widehat{\mathcal{T}}(c_m):=\widehat{\mathcal{T}}(T(m))$ for one fixed 
		$m \in \{1,\dots,4\}$.
			\begin{lemma}\label{lem:fP3_Tprime} 
				Let $f\in P_3(\widehat{\mathcal{T}}(c_m))$  and $g:={\partial f}/{\partial \tau}\in P_2(\widehat{\mathcal{T}}(c_m))$.
				Then $f\in C^1(T(m))$ if and only if $f|_{F_m}\in C^1(F_m)$ and $g\in P_2(T(m))$.
			\end{lemma}
				\begin{proof}
					Without loss of generality let $m=4$, $\tau:=c_4-c_K$, and $T(4)=\textup{conv}(c_K,F_4)$. 
					\\[1em]\textit{\glqq $\Rightarrow$\grqq } 
					The continuity of $f\in C^1(T(4))$ and its derivatives imply $f|_{F_4}\in C^1(F_4)$. 
					\cref{lem:C1Tm_k=012} leads to $P_2(\widehat{\mathcal{T}}(c_4))\cap C^1(T(4))=P_2(T(4))$.  
					Since $g\in P_2(\widehat{\mathcal{T}}(c_4))\cap C^0(T(4))$, it suffices to show $g\in C^1(T(4))$. 
					Fix the common side $F=\textup{conv}\{c_K,c_4,Q_j\}=\partial T_1\cap\partial T_2$ of two neighbouring tetrahedra 
					$T_1,T_2\in \widehat{\mathcal{T}}(c_4)$ identified with $T_1=(c_K,c_4,Q_j,Q_k)$ 
					and $T_2=(c_K,c_4,Q_j,Q_\ell)$ for $\{j,k,\ell\}=\{1,2,3\}$. 	
					Since $f|_{T_j}=6\sum_{\alpha\in A_{3}}c_j(\alpha)\lambda^\alpha/\alpha!$ for $j=1,2$  defines $f\in C^1(T(m))$ in $T_j$, 
					\cref{lem:gprime_C1_onT1T2} shows that  $f^\prime_\xi|_{T_j}:=2\sum_{\beta\in A_{k-1}}c_j(e_\xi+\beta)\lambda^\beta/\beta!$ 
					for $j=1,2$ 
					in \eqref{eq:gprimeon T1T2} forms $f^\prime_\xi\in C^1(T_1\cup T_2)$ for any $\xi\in \{1,2,3\}$. 
					Notice that \eqref{eq:derivativeBernstein} implies 
					$\tau\cdot\nabla f|_{T_j}=6\sum_{\beta\in A_2}\sum_{\mu=1}^4c_j(\beta+e_\mu)(\tau\cdot\nabla\lambda_\mu)\lambda^\beta/\beta!$ 
					with 
					$\tau\cdot\nabla\lambda_\mu=\lambda_\mu(c_4)-\lambda_\mu(c_K)=\delta_{2\mu}-\delta_{1\mu}$.  Hence a.e. in $T_j$ it holds  
					\begin{align*}
					 	g|_{T_j}=\tau\cdot\nabla f|_{T_j}=6\sum_{\beta\in A_2}\big(c_j(\beta+e_2)-c_j(\beta+e_1)\big)\lambda^\beta/\beta!
					 	=3f^{\prime}_2|_{T_j}-3f^\prime_1|_{T_j}\text{ for }j=1,2.
					\end{align*}
					This and $f_1^\prime, f_2^\prime \in C^1(T_1\cup T_2)$ show  $g\in C^1(T_1\cup T_2)$. 
					Since $\widehat{\mathcal{T}}(c_4)$ is side-connected, it follows successively 
					$g\in C^1(T(4))\cap P_2(\widehat{\mathcal{T}}(c_4))$. 
					\cref{lem:C1Tm_k=012} implies $g\in P_2(T(4))$. \hfill $\Box$
					\\[1em]\textit{\glqq $\Leftarrow$\grqq }
					Suppose $g\in P_2(T(4))$ and $f|_{F_4}\in C^1(F_4)$. 
					A one-dimensional integration shows for all $x\in F_m$ and $t\in\mathbb{R}$ with $x+t\tau\in T(4)$ that 
					\begin{align}
						f(x+t\tau)=f(x)+\int_0^t g(x+s\tau)\,\textup{d}s. \label{eq:proof_fP3_Tprime}
					\end{align}
					(Notice that $0<s<t$ and $x\in F_4$ with $x+t\tau\in T$ for some tetrahedron $T\in\widehat{\mathcal{T}}$ 
					imply hat $x+s\tau$ belongs to the same tetrahedron by definition of $\tau$.)
					The assumptions $g\in P_2(T(4))$ and $f|_{F_4}\in C^1(F_4)$ show in combination with \eqref{eq:proof_fP3_Tprime} that $f\in C(T(4))$ and that  
					$\sigma \cdot\nabla f$ is continuous for all 
					$\sigma\in \textup{span}\{\tau, Q_2-Q_1,Q_3-Q_1\}\in\mathbb{R}^3$ (in any direction parallel to the 
					plane $\langle F_4\rangle$  and in the oblique direction 
					$\tau \not \in \textup{span}\{ Q_2-Q_1,Q_3-Q_1\}\parallel\langle F_4\rangle$). 
					This implies $\nabla f$ in $C^0(T(4);\mathbb{R}^3)$ and concludes the proof. 
				\end{proof}
				The main result of this section requires some notation on two tetrahedra with a common side. 
				Suppose that $K_+,K_-\in\mathcal{T}$ are two neighbouring tetrahedra in a regular triangulation $\mathcal{T}$ 
				that share a common face $F=\partial K_+\cap\partial K_-$. Suppose that the center points $c_F$, $c_{K_+}$, and $c_{K_-}$ 
				in the $\textit{WF}$ 
				partition $\widehat{\mathcal{T}}(K_\pm)$ of \cref{def:CTtriangulation} belong to one straight line. 
				In other words $c_{K_+}-c_F\parallel c_{K_-}-c_F$. Suppose the degrees of freedom for some piecewise $\textit{WF}$ 
				(finite element) function 
				$f\in \textit{WF} (\{K_+,K_-\})$ 
				with $f_{\pm}:=f|_{K_{\pm}}\in \textit{WF}(K_{\pm}):=C^1(K_{\pm})\cap P_3(\widehat{\mathcal{T}}(K_{\pm}))$ 
				coincide at 
				$F=\textup{conv}\{Q_1,Q_2,Q_3\}$. 
				The latter means that 
				\begin{align}
					f_+(Q_\ell)=f_-(Q_\ell), \  \nabla f_+(Q_\ell)=\nabla f_-(Q_\ell),
					\text{ and }
					\nabla f_+(Q_{jk})\times \tau_{jk}=\nabla f_-(Q_{jk})\times\tau_{jk}
				\label{eq:commondof_K+K-} 
				\end{align}{ for all }$1\le \ell\le 3, \ 1\le j<k\le 3$.
				(Recall $Q_{jk}:=(Q_j+Q_k)/2$ is the midpoint of the edge $E_{jk}:=\textup{conv}\{Q_j,Q_k\}$ with 
				unit tangent vector $\tau_{jk}$ of fixed orientation, e.g., $\tau_{jk}=(Q_k-Q_j)/|E_{jk}|$, in 
				the triangle $F$.)   
			\begin{theorem}\label{thm:CT_K1K2} 
				 Under the present notation of $18$ coinciding degrees of freedom \eqref{eq:commondof_K+K-}, $f_{\pm}$ form a $C^1$ function 
				 $f\in C^1(K_+\cup K_-)$ on the patch $\omega(F):=\textup{int}(K_+\cup K_-)$ of $F$.  
			\end{theorem}
			\begin{proof}
				Recall that the $2D$ $HCT$ function $f_m|_{F_m}\in HCT(F_m)$ in \cref{alg:CT3D}.\ref{item:Alg_CT_HCT} 
				is computed from twelve degrees of 
				freedom at $F_m$. 
				For $f_+$ and $f_-$ these twelve degrees of freedom coincide along $F_m\equiv F$ 
				owing to the assumption \eqref{eq:commondof_K+K-}.
				Hence the interpolating $2$D $HCT$ function $f_{m\pm}$  in \cref{alg:CT3D}.\ref{item:Alg_CT_HCT}  
				is the same for each of the two tetrahedra $K_{\pm}$, i.e., $f_{m+}|_{K_+}=f_{m-}|_{K_-}=:f_m|_F\in HCT(F)$.  
				\cref{alg:CT3D}.\ref{item:Alg_CT_HCT} guarantees 
				that $f_m|_F$ defines the $\textit{WF}$ 
				function $f_{\pm}|_F=f|_{F_m}$ at $F$. 
				It follows $f_+|_F=f_-|_F=f_m|_F$ and so $f\in C^0(\omega(F))$.
				This and $f_{\pm}\in C^1(K_{\pm})$ also imply that the tangential derivatives 
				$\nu_F\times \nabla f_+(x)=\nu_F\times \nabla f_-(x)$ coincide  at any $x\in F$ 
				(owing to the Hadamard jump condition with unit normal vector $\nu_F$ of $F$).
				To conclude the proof of $f\in C^1(\omega (F))$, it remains to verify $\tau \cdot \nabla f_+(x)=\tau \cdot \nabla f_-(x)$ 
				for the direction					
				\begin{align}
					\tau:=\frac{c_{K_+}-c_F}{|c_{K_+}-c_F|}=-\frac{c_{K_-}-c_F}{|c_{K_-}-c_F|}\in\mathbb{R}^3\setminus\{0\} \label{eq:def_tau_CTK1K2}
				\end{align}
				oblique to $\langle F\rangle$. The equality in \eqref{eq:def_tau_CTK1K2} results  
				from the assumption $c_{K_+}-c_F\parallel c_{K_-}-c_F$.
				Since $f_{\pm}\in C^1(K_{\pm})$, \cref{lem:fP3_Tprime} implies that 
				$g_+:=\tau\cdot \nabla f_+\in P_2(K_+(m))$ and $g_-:= \tau\cdot \nabla f_-\in P_2(K_-(m))$ are quadratic polynomials on 
				$K_{\pm}(m):=\textup{conv}(F,c_{\pm})$, in particular $g_{\pm}|_F=\tau\cdot\nabla f_{\pm}|_F\in P_2(F)$. 
				The continuity of $f_{\pm}$ at $F$ and the equality \eqref{eq:commondof_K+K-} of the $18$ degrees of freedom at $F$ 
				imply that the quadratic 
				polynomials $g_+\in P_2(K_+(m))$ and $g_-\in P_2(K_-(m))$ coincide at $Q_1$, $Q_2$, $Q_3$, $Q_{12}$, 
				$Q_{13}$, and $Q_{23}\in F$ (see \cref{rem:preprocessing} for the edge midpoints). 
				The latter are the six Lagrange points of $P_2(F)$ in the triangle $F$ and uniquely determine $g_{\pm}|_F\in P_2(F)$.
				Consequently, $g_+|_F=g_-|_F\in P_2(F)$. This shows that $\tau\cdot\nabla f$ is continuous at $F$ and concludes the proof. 
			\end{proof}
			\FloatBarrier
	\subsection{The shape-regularity of the $\textit{WF}$ partition} 
			Let $\mathbb{T}$ denote the set of uniformly shape-regular triangulation $\mathcal{T}$ of a bounded polyhedral 
			Lipschitz domain $\Omega\subset\mathbb{R}^3$ into tetrahedra: There exists a global constant 
			$C_{\mathrm{sr}}>0$ with  
			$|T|^{1/3}\le h_T\le  C_{\mathrm{sr}}|T|^{1/3}$
			for any tetrahedra $T\in\mathcal{T}\in\mathbb{T}$ with diameter $h_T$ and volume $|T|$. 
			Throughout this section, let $\mathcal{V}$ denote set of vertices, $\mathcal{F}$ the set of faces, 
			and $\mathcal{E}$ the set of edges in the triangulation $\mathcal{T}\in\mathbb{T}$. 
			For any $E\in\mathcal{E}$ fix one unit tangent vector $\tau_E\parallel E$ and $\nu_{E,1},\,\nu_{E,2}\in\mathbb{R}^3$ such that
							 $(\tau_E,\,\nu_{E,1},\,\nu_{E,2})$ form an orthonormal basis of $\mathbb{R}^3$.
		
			\begin{definition}[global $\textit{WF}$ 
								partition]\label{def:CT_triang_points}
					For any $T\in\mathcal{T}$, let $c_T$ denote the midpoint of the incircle; 
					for any interior face $F:=\partial T_+\cap\partial T_-\in\mathcal{F}(\Omega)$, 
					let $c_F$ be the intersection of $F$ with the straight line through 
					$c_{T_+}$ and $c_{T_-}$ of the neighbouring tetrahedra 
					$T_{\pm}\in\mathcal{T}$;
					for any boundary face $F\in\mathcal{F}(\partial\Omega)$, let $c_F:=\textup{mid}(F)$ denote the center of gravity.
					Let $\widehat{\mathcal{T}}$ denote the subtriangulation of ${\mathcal{T}}$, where any tetrahedron $T\in\mathcal{T}$ is partitioned  
					in $\widehat{\mathcal{T}}(T)=\textit{WF3D}(T)$ 
					according to \cref{def:CTtriangulation}.
			\end{definition}
			The center points $c_T$ and $c_F$ from \cref{def:CT_triang_points} satisfy 
		 	the conditions \ref{item:CTcon1}--\ref{item:CTcon2} in Section~\ref{sec:HCT_Triang} according to the following \cref{thm:shapeRegularity3D}. 
		 	In particular, the global $\textit{WF}$ 
		 	partition in \cref{def:CT_triang_points} is well-defined.
		 	
			\begin{theorem}[shape regularity of $\widehat{\mathcal{T}}$]\label{thm:shapeRegularity3D}
				There exists a constant $\varepsilon>0$ (that exclusively depends on $\mathbb{T}$), 
				such that the points selected in \cref{def:CT_triang_points} satisfy 
				$\textup{dist}(c_T,\partial T)\ge \varepsilon h_T$ for all $T\in\mathcal{T}\in\mathbb{T}$ 
				and  $\textup{dist}(c_F,\partial F)\ge \varepsilon h_F$ for all $F\in\mathcal{F}$.
				In particular, the global $\textit{WF}$ 
				partition  $\widehat{\mathcal{T}}$ in \cref{def:CT_triang_points} is uniformly shape-regular 
				with a shape-regularity constant that is bounded in terms of $C_{\mathrm{sr}}$ in the definition of $\mathbb{T}$. 		 	
		 	\end{theorem}	
		 	\begin{figure}
				 	\begin{center}
				 	\begin{tabular}{ccc}
				  	\scalebox{0.9}{\begin{tikzpicture}[x=30mm,y=30mm]
				 							\draw(0,0)node[shape=coordinate,label={[xshift=2mm, yshift=-2mm] $m$}] (m) {m};
				 							\draw(0,0)node[shape=coordinate,label={[xshift=15mm, yshift=-17mm] $B(m,\varrho)$}](a) {a};
				 							\draw (0,0) circle (1);
				 							\draw(-2,0)node[shape=coordinate,label={[xshift=-4mm, yshift=-2mm] $\tau_E$}] (tau) {tau};
				 							\draw(-2,0)  circle (3pt);
				 							\draw[dashed] (m)--(tau);	
				 							\draw[dashed] (0,0)--++([shift=(120:1)]0,0);
				 							\draw[dashed] (0,0)--++([shift=(240:1)]0,0);			
				 							\draw(-0.5,0.866)node[shape=coordinate,label={[xshift=3mm, yshift=-3mm] $m_1$}] (m1) {m1};	
				 							\draw(-0.5,-0.866)node[shape=coordinate,label={[xshift=3mm, yshift=-2mm] $m_2$}] (m2) {m2};	
				 							\draw[thick] (tau)--(m1);
				 							\draw[thick] (tau)--(m2);
				 							\draw[->, rotate around={270:(-0.5,0.866)}] (m1)--(-0.875,0.6495);
				 							\draw[->, rotate around={90:(-0.5,-0.866)}] (m2)--(-0.875,-0.6495);
				 							\draw[thick, rotate around={180:(-0.5,0.866)}] (m1)--(-1.25,0.433);
				 							\draw[thick, rotate around={180:(-0.5,-0.866)}] (m2)--(-1.25,-0.433);
				 							\filldraw[white, opacity=0.5] (-2,0)--([shift=(-30:0.3)]-2,0) arc (-30:30:0.3) -- cycle;
				 							\draw[thick] (tau)--(m1);
				 							\draw[thick] (tau)--(m2);
				 							\draw ([shift=(-30:0.3)]-2,0) arc (-30:30:0.3);
				 							\draw(-2,0)node[shape=coordinate,label={[xshift=5mm, yshift=-3mm] $\alpha(E)$}] (alpha) {alpha};
				 							\draw ([shift=(-30:0.6)]-2,0) arc (-30:0:0.6);
				 							\draw(-2,0)node[shape=coordinate,label={[xshift=13mm, yshift=-7mm] $\nicefrac{\alpha(E)}{2}$}] (alpha2) 
				 							{alpha2};
				 							\draw(-0.6,1.0392) node[shape=coordinate,label={[xshift=2mm, yshift=-1mm] $\nu_{F_1}$}] (n1) {n1};	
				 							\draw(-0.2,1.0392) node[shape=coordinate,label={[xshift=3mm, yshift=2mm] $\langle F_1\rangle $}] (n12) {n12};	
				 							\draw(-0.6,-1.0392) node[shape=coordinate,label={[xshift=3mm, yshift=-3mm] $\nu_{F_2}$}] (n2) {n2};	
				 							\draw(-0.2,-1.0392) node[shape=coordinate,label={[xshift=3mm, yshift=-8mm] $\langle F_2\rangle$}]  (n22) 
				 							{n22};	
				 							  \pic[draw, fill=red,fill opacity=0.2, angle radius=4mm,"$\cdot$" opacity=1] {angle=tau--m1--m};
				 							  \pic[draw, fill=red,fill opacity=0.2, angle radius=4mm,"$\cdot$" opacity=1] {angle=n12--m1--n1};
				 							  \pic[draw, fill=red,fill opacity=0.2,angle radius=4mm,"$\cdot$" opacity=1] {angle=m--m2--tau};
				 							    \pic[draw, fill=red,fill opacity=0.2, angle radius=4mm,"$\cdot$" opacity=1] {angle=n2--m2--n22}; 
				 							 
				 							 \draw [decorate,decoration={brace,amplitude=10pt,raise=4pt},yshift=0pt]
													(tau) -- (m1) node [black,midway,yshift=7mm, xshift=-1mm] {$d(E)$};
											\draw [decorate,decoration={brace,amplitude=10pt,mirror, raise=4pt},yshift=0pt]
													(tau) -- (m2) node [black,midway,yshift=-7mm, xshift=-1mm] {$d(E)$};
											\draw [decorate,decoration={brace,amplitude=10pt,raise=9pt},yshift=0pt]
													(m) -- (m2) node [black,midway,yshift=-4mm, xshift=7mm] {$\varrho$};
						\end{tikzpicture}}
					 	&\phantom{xx}&
					 	{\begin{tikzpicture}[x=49mm,y=49mm]
											\draw(1,1)node[shape=coordinate, label={[xshift=0.5mm, yshift=-1mm] $\phantom{P}$}] (Q3) {Q3};
						  					\draw (0,0)--(1,0)--(0,1)--cycle;
						  					\draw(0,0)node[shape=coordinate,label={[xshift=2mm, yshift=-0.6mm] $F$}] (P1) {P1};
						  					 \filldraw[white!80!red] (0.11,0.11)--(0.75,0.11)--(0.11,0.75)--cycle; 
						  					 \draw[black!30!red](0.11,0.11)--(0.75,0.11)--(0.11,0.75)--cycle;
						  					\draw(0.11,0.11)node[shape=coordinate,label={[xshift=2mm, yshift=-0.6mm] 
						  					\textcolor{black!30!red}{$F^\prime$}}] 
						  					(P11) {P11};
						  					\draw (0,-0.2) node [shape=coordinate] (P0) {P0};
						 \end{tikzpicture}}
						 \\
						 (a)&&(b) 	
				 	\end{tabular} 
					 \end{center}
					 \caption{ Notation in the profile of the intersection $\langle F_1\rangle\cap \langle F_2\rangle$ 
					 			of two planes through the faces $F_1,F_2\in \mathcal{F}(T)$ 
					 			with common edge $E$ indicated by $\tau_E$ as $\odot$ perpendicular to the plane depicted in (a) and  of  $F^\prime$ 
							   from \eqref{eq:Formregularity_defFprime} (in red) in (b).}\label{fig:Formregularity_Profile}
				 \end{figure}
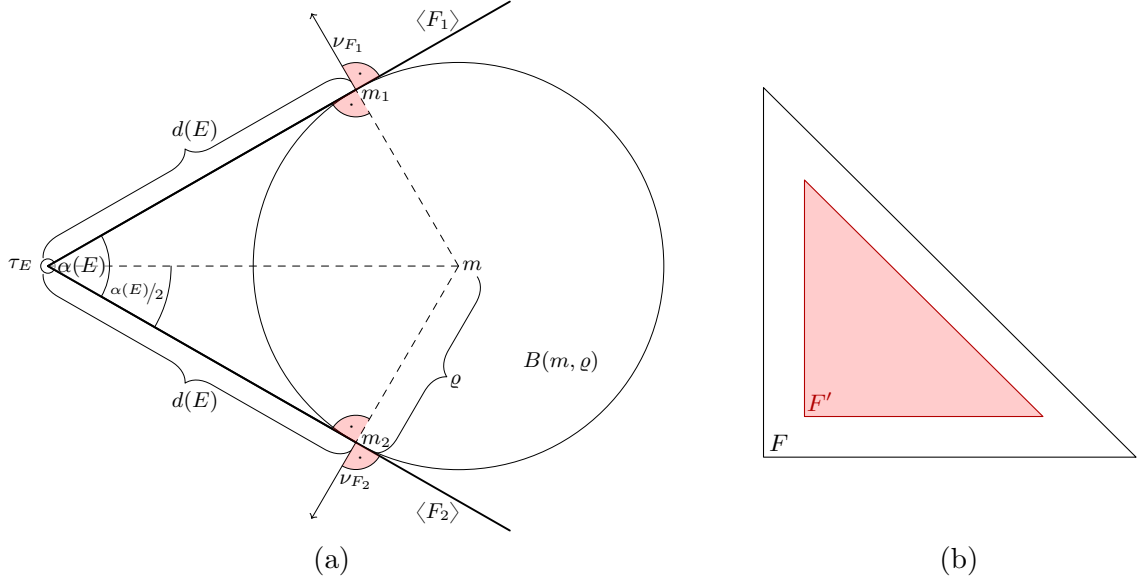
			\begin{proof}
			\begin{enumerate}[label=\textit{Step \arabic* --},wide]
			\item \textit{Preparation. }
			First consider one tetrahedron $T\in\mathcal{T}$ with largest ball $B(m,\varrho)\subset T$ contained in $T$  with midpoint $m$ and radius $\varrho$ 
			and a common edge $E=\partial F_1\cap \partial F_2\in\mathcal{E}(T)$ of  two of its faces 
			$F_1,F_2\in\mathcal{F}(T)$, $F_1\not = F_2$.
			
			 \cref{fig:Formregularity_Profile}.a illustrates the following notation. 
			 Let $\alpha(E)$ denote the dihedral angle of the distinct planes $\langle F_1\rangle$ and $\langle F_2\rangle$ through the faces 
			 $F_1,F_2\in\mathcal{F}(T)$ with common edge $E$.
			 The two intersecting planes touch the ball $\overline{B}(m,\varrho)$ at the two points 
			 $m_j=\langle F_j\rangle\cap \partial B(m,\varrho)\in\textup{int}(F_j)$ for $j=1,2$.
			 The plane $\langle m,m_1,m_2\rangle$ through the points $m,\,m_1,$ and $m_2$ is perpendicular to the unit tangent vector $\tau_E$ of the edge $E$.
			 The  distance $d(E)$ of $m_j$ to the straight line $\langle E\rangle$ through $E$ coincides  for $j=1$ and $j=2$, 
			 $d(E):=\textup{dist}(m_1,\langle E\rangle)=\textup{dist}(m_2,\langle E\rangle)$,  
			 according to elementary geometry in a deltoid. 
			 The shape regularity of a tetrahedron allows many equivalent characterizations e.g. in terms of angles as discussed in \cite{BKK08}. 
			 The lower and upper bounds of $\alpha(E)$ from \cite[Eq.(4)]{BKK08} lead in 
			 $\varrho=d(E)\tan\big({\alpha(E)}/{2}\big)$ to the equivalence
				\begin{align}
					d(E)\approx\varrho\approx h_T\approx h_F\quad\text{for all }E\in\mathcal{E}(T). \label{eq:Formregularity_varrho}
				\end{align}
			  (Recall that the equivalence constants in the notation $\approx$ depend exclusively on $C_{\mathrm{sr}}$.) 
			  This means that the point $m_F\in \partial B(m,\varrho)\cap F$  for the interior side $F$ belongs to the subtriangle
			  	\begin{align}
			  		F^\prime:=\{x\in F|\,\textup{dist}(E,x)\ge d(E)\text{ for all }E\in \mathcal{E}(F)\}\ni m_F \label{eq:Formregularity_defFprime}
				\end{align}	  		
			  depicted in \cref{fig:Formregularity_Profile}.b. This and \eqref{eq:Formregularity_varrho} will be employed in the sequel. 	
			\item\textit{Distance conditions. } 
			Suppose that $F=\partial T_1\cap\partial T_2$ is a common side of two distinct neighbouring tetrahedra $T_1,\,T_2\in\mathcal{T}$.
			Then each $T_j$ has a maximal incircle $B(m_j,\varrho_j)\subset T_j$ with midpoint $m_j$ and radius $\varrho_j$ and the touching point 
			$m_{F,j}\in F_j^\prime$ from \eqref{eq:Formregularity_defFprime}. 
			The shape regularity \cite{BKK08} means $h_{T_1}\approx\varrho_1\approx h_F\approx \varrho_2\approx h_{T_2}$.
			Hence the choice $m_j=c_{T_j}$ in \cref{def:CT_triang_points} 
			guarantees $\varepsilon h_{T_j}\le \textup{dist}(c_{T_j},\partial T)=\varrho_j\approx h_{T_j}$ for $j=1,2$ with a constant $\varepsilon\approx 1$. 
			
			Let $d_j(E)$ replaces the distance $d(E)$ for $T$ replaced by $T_j$ for $j=1,2$ in \eqref{eq:Formregularity_defFprime} to define $F^\prime_j$. 
			Then  $F^\prime_1\subseteq F^\prime_2$ or  $F^\prime_2\subseteq F^\prime_1$ holds. 
			Let 
			$$F^\prime:=F^\prime_1\cup F^\prime_2= \big\{x\in F|\,\textup{dist}(E,x)\ge \min\{d_1(E),d_2(E)\}\text{ for all }E\in\mathcal{E}(F)\big\}$$ 
			denote the bigger set. 
			The face center $c_F$ in \cref{def:CT_triang_points} lies in the convex hull of the two touching points 
			$m_{F,1}=\big(\partial B(m_1,\varrho_1)\cap F\big)\in F^\prime_1\subseteq F^{\prime}$ 
			and $m_{F,2}=\big(\partial B(m_2,\varrho_2)\cap F\big)\in F^\prime_2\subseteq F^{\prime}$. 
			Since the triangle $F^{\prime}$ depicted in \cref{fig:Formregularity_Profile}.b is convex, $c_F\in F^{\prime}$. 
			This and \eqref{eq:Formregularity_varrho} (with $h_F\approx |F|^{1/2}\approx\varrho_1\approx\varrho_2$ from shape regularity \cite{BKK08}) prove that 
			$\textup{dist}(c_F,\partial F)\ge \min_{E\in\mathcal{E}(F)}\min_{j=1,2}d_j(E)\approx h_F$.
			This leads to $\varepsilon\approx 1$ and $\varepsilon h_F\le \textup{dist}(c_F,\partial F)$.     
			\item\textit{Shape regularity of $\widehat{\mathcal{T}}$. }
			Any tetrahedron $\widehat{T}=\textup{conv}\{c_T,c_F,Q_j,Q_k\}\in\widehat{\mathcal{T}}$ 
			is contained in one $ T\in\mathcal{T}$ with center point $c_T$ of the maximal incircle $B(c_T,\varrho)\subset T$ in $T$  
			and face $F\in\mathcal{F}(T)$ such that 
			$c_F\in \mathrm{relint}(F)$ is the point $c_F\in F^\prime$ from \cref{def:CT_triang_points}  
			and  $E=\textup{conv}\{Q_j,Q_k\}\in\mathcal{E}(F)$ is an edge of $F$.
			Then 
			$$h_{\widehat{T}}:=\textup{diam}(\widehat{T})=\max\{|c_T-c_F|, |c_T-Q_j|, |c_T-Q_k|, |c_F-Q_j|,|c_F -Q_k|,|Q_j-Q_k|\}\le h_T.$$ 
			By construction and shape regularity of $\mathcal{T}$ with \eqref{eq:Formregularity_varrho}--\eqref{eq:Formregularity_defFprime} 
			it follows $h_T\approx \varrho\le \min\{|c_T-c_F|, |c_T-Q_j|, |c_T-Q_k|\}$, $h_T\approx d(E)\le \min\{|c_F-Q_j|,|c_F -Q_k|\}$, 
			and $h_T\approx h_E=|Q_j-Q_k|$. In other words $h_{\widehat{T}}\approx h_T$. 
			The equivalence $h_{\widehat{T}}\approx h_T\approx \varrho\approx |F|^{1/2}$ and the observation $\widehat{T}=\textup{conv}(F,c_T)$ 
			imply  
			$$|\widehat{T}|=\frac{|F|\varrho}{3}\approx h_T^3\approx h_{\widehat{T}}^3.$$ 
			This guarantees the shape regularity of $\widehat{\mathcal{T}}$ 
			\cite[Eq.(1)]{BKK08} and the shape regularity constant solely depends on $\mathbb{T}$.  
			\end{enumerate}
		\end{proof}
		\subsection{The scaling of the $\textit{WF}$ basis functions} 
		Throughout this section, $\widehat{\mathcal{T}}$ is the $\textit{WF}$ 
		partition  of a triangulation $\mathcal{T}\in\mathbb{T}$ as in 
		\cref{def:CT_triang_points} with set of vertices $\mathcal{V}$, set of faces $\mathcal{F}$, and set of edges 
		$\mathcal{E}$.  For any edge $E\in\mathcal{E}$
		fix the orientation of a unit tangent vector $\tau_E$ along $E$ as well as $\nu_{E,1},\,\nu_{E,2}\in\mathbb{R}^3$ such that
		$(\tau_E,\,\nu_{E,1},\,\nu_{E,2})$ form an orthonormal basis of $\mathbb{R}^3$. 	Set
		\begin{align}
				\textit{WF}(\mathcal{T})				
					:=P_3(\widehat{\mathcal{T}})\cap H^2_0({\Omega}).\label{eq:def_HCT}
			\end{align}
			\begin{definition}[global degrees of freedom]\label{def:CT_dof_global}
				The global $M=4|\mathcal{V}(\Omega)|+2|\mathcal{E}(\Omega)|$ 
				degrees of freedom are $\big(L_{z,\mu}|\,z\in\mathcal{V}(\Omega),\mu=1,\dots,4\big)$ 
				and $\big(L_{E,\mu}|\,E\in\mathcal{E}(\Omega),\mu=1,2\big)$ defined,  for all 
				$f\in C^1(\overline{\Omega})$, by 
					\begin{align}
					\begin{split}
						&L_{z,1}f:=f(z), \, L_{z,\kappa+1}:=\frac{\partial f}{\partial x_\kappa}(z)\quad \text{for }\kappa=1,2,3;\\
						&L_{E,\mu}f:=\frac{\partial f}{\partial \nu_{E,\mu}}(\text{mid}(E))\qquad \qquad \text{for }\mu=1,2.
					\end{split}\label{eq:global_dof}
					\end{align}
					These $M$ degrees of freedom are also enumerated as $L_1,\dots, L_M$. 
					Note, that each degree of freedom $L_\ell $ for $\ell=1,\dots, M$ is a directional derivative of order 
					$\mu_\ell=0$ or $\mu_\ell=1$ 
					at an interior node 
					$z_\ell \in \mathcal{N}(\Omega)
						:=\mathcal{V}(\Omega)\cup \{\textup{mid}(E):\, E\in\mathcal{E}(\Omega)\}$.
			\end{definition}	
				For each node $z_\ell\in\mathcal{N}(\Omega)$ let 
				$\mathcal{T}(z_\ell):=\{T\in\mathcal{T}|\, z_\ell\in\partial T\}$ denote the set of adjacent simplices and let  
				$\omega_\ell:=\textup{int}(\bigcup\mathcal{T}(z_\ell))$ denote the nodal patch with  
				volume $h_\ell^{3}:=|\omega_\ell|$. 
			 
		 \cref{thm:Scaling_CT_fine} below extends the analysis for a quasi-affine family of finite elements in \cite[\S 6.1]{Ciarlet78} to $3$D 
		 and implies 
		the correct scaling of the dual basis functions   for $\textit{WF}({\mathcal{T}})$. 
		The point is that  \cref{thm:shapeRegularity3D} lead to uniform bounds in (\ref{item:CT_basis_scaling}). 
		Recall that $a \lesssim b$ abbreviates $a\le Cb$ with a generic constant $C$ that only 
		depends the shape-regularity constant $C_{\mathrm{sr}}$ of $\mathcal{T}\in\mathbb{T}$. 
		(The constants below depend on $\varepsilon >0$ from \cref{thm:shapeRegularity3D} which depends solely on $\mathbb{T}$.) 
			\begin{theorem}[scaling of $\textit{WF}$ 
							basis functions]\label{thm:Scaling_CT_fine}
			\begin{enumerate}[label=(\alph*), ref=\alph*, wide]
			\item\label{item:CT_basis_exists}
					There exists a unique nodal basis $(\varphi_1,\dots,\varphi_M)$ of $\textit{WF}(\mathcal{T})$ 
					which satisfies 
					$L_k(\varphi_\ell)=\delta_{k\ell}\text{ for } k,\ell=1,\dots, M$.  
					\begin{notation}	
						If $L_\ell=L_{z,\mu}$ for $z\in\mathcal{V},\,\mu=1,\dots,4$, relabel 
						$\varphi_{z,\mu}:=\varphi_\ell$. 
						If $L_\ell=L_{E,\mu}$ for $E\in\mathcal{E},\,\mu=1,2$, relabel 
						$\varphi_{E,\mu}:=\varphi_\ell$. 
					\end{notation}
			\item\label{item:CT_basis_scaling}
				Suppose the degree of freedom $L_\ell$ for $\ell=1,\dots,M$ is a directional derivative of order 
				$\mu_\ell\in\{0,1\}$ at the node $z_\ell \in \mathcal{N}(\Omega)$ with nodal patch 
				$\omega_\ell$  of volume $h_\ell^{3}:=|\omega_\ell|$. Then  
				 $$
				 	\Vert \varphi_{\ell}\Vert_{H^s(\Omega)}	\lesssim h_\ell^{3/2+\mu_\ell-s} \qquad\text{ for all }s=0,1,2.
				 $$
			\end{enumerate}		
			\end{theorem}
			\begin{proofof}\textit{(\ref{item:CT_basis_exists}). }
			Consider one global degree of freedom $L_\ell$ with $\ell=1,\dots,M$ from \cref{def:CT_dof_global} associated with the node 
			$z_\ell\in\mathcal{N}(\Omega)$. 
		 	For any $K\in\mathcal{T}(z_\ell)$ the duality translates in input data $(x_1,\dots,x_{28})\in\mathbb{R}^{28}$ for \cref{alg:CT3D} 
		 	as follows: 
			Suppose $L_\ell=L_{E,\mu}$ for some edge 
			$E=\textup{conv}\{Q_j,Q_k\}\in\mathcal{E}(K)$ in $K$ and $\mu\in\{1,2\}$, i.e., $z_\ell=\textup{mid}(E)$. 
			Recall $(\tau_E,\nu_{E,1},\nu_{E,2})$ is an orthonormal basis of $\mathbb{R}^3$ and 
			$(\tau_{jk},\nu_{jk,1},\nu_{jk,2})$ for $E=E_{jk}$ denotes another basis  of $\mathbb{R}^3$ in \cref{sec:HCT3D=FE}.   
			Since $\tau_{jk}\parallel\tau_E$, it holds 
			$\nu_{E,\mu}=\sum_{m=1}^2\alpha_m \nu_{jk,m}$ for some $\alpha_1,\alpha_2\in \mathbb{R}$. 
			Let $x_{{20-2^{4-j}+2k+m}}=\alpha_m$ for $m=1,2$ (for the enumeration of \cref{tab:CT_3D_dof}) . 
			\cref{alg:CT3D} with input data
			 $(0,  \dots, 0, x_{{20-2^{4-j}+2k+1}},x_{{20-2^{4-j}+2k+2}},0,\dots,0)\in\mathbb{R}^{28}$  (and $\widehat{\mathcal{T}}(K)$)  
			 computes a function $\varphi_{\ell,K}\in \textit{WF}(K):=P_3(\widehat{\mathcal{T}}(K))\cap C^1(K)$ 
			 with 
			\begin{align*}
				L_\ell(\varphi_{\ell,K})&=\nu_{E,\mu} \cdot \nabla \varphi_{\ell,K}(\textup{mid}(E))
				=\alpha_1\nu_{jk,1}\cdot  \nabla \varphi_{\ell,K}(Q_{jk})+
				\alpha_2\nu_{jk,2}\cdot  \nabla \varphi_{\ell,K}(Q_{jk})
				\\&=
				\alpha_1 x_{{20-2^{4-j}+2k+1}}  +\alpha_2  x_{{20-2^{4-j}+2k+2}} 
				=\alpha_1^2+\alpha_2^2=|\nu_{E,\mu}|^2=1
			\end{align*} 
			which   vanishes $L_k(\varphi_{\ell,k})=0$ in all other degrees of freedom $k=1,\dots, M$ with $k\not = \ell$. 
			(For the second normal vector $\nu_{E,\kappa}=\sum_{m=1,2}\beta_m \nu_{jk,m}$ 
			for some $\beta_1,\beta_2\in \mathbb{R}$, $\kappa\in\{1,2\}\setminus\{\mu\}$ holds 
			$0=\nu_{E,\kappa}\cdot\nu_{E,\mu}=\alpha_1\beta_1+\alpha_2\beta_2$ and whence 
			$L_{E,\kappa}(\varphi_{\ell,K})=\sum_{m=1}^2 \beta_1 x_{{20-2^{4-j}+2k+m}}=0$.) 
			 Similarly suppose $L_\ell=L_{z,\mu}$ for some vertex $z=Q_j\in\mathcal{V}(K)$ in $K$ and $\mu\in\{1,\dots, 4\}$. 
			 Set  $x_{4(\mu-1)+j}=1$. 
			 \cref{alg:CT3D} with input data $(0,  \dots, 0, x_{4(\mu-1)+j},0,\dots,0)\in\mathbb{R}^{28}$   
			 (and $\widehat{\mathcal{T}}(K)$)  
			  computes a function 
			 $\varphi_{\ell,K}\in \textit{WF}(K):=P_3(\widehat{\mathcal{T}}(K))\cap C^1(K)$  
			 with 
			\begin{align*}
				L_\ell(\varphi_{\ell,K})
				&=\varphi_{\ell,K}(Q_j)=x_{j}=1\qquad\text{ for }\mu=1,
				\\L_\ell(\varphi_{\ell,K})
				&=e_{\mu-1}\cdot\varphi_{\ell,K}(Q_j)=x_{4(\mu-1)+j}=1\qquad\text{ for }\mu=2,3,4,
			\end{align*}
			which   vanishes $L_k(\varphi_{\ell,k})=0$ in all other degrees of freedom $k=1,\dots, M$ with $k\not =\ell$. 
			Let $\varphi_\ell\in P_3(\widehat{\mathcal{T}})$ denote the global function locally defined as 
			$\varphi_\ell|_K:=\varphi_{\ell,K}$	for any $K\in\mathcal{T}(z_\ell)$ (and $\varphi_\ell|_K=0$ otherwise). 	
			In summary the functions $\varphi_1,\dots, \varphi_M \in P_3(\widehat{\mathcal{T}})$ satisfy the duality property in 
			(\ref{item:CT_basis_exists}).
			For any two tetrahedra $K_\pm$ with common side $F=	\partial K_+\cap \partial K_-\in\mathcal{F}(\Omega)$ 
			the respective center points $c_{K_+}\in\textup{int}(K_+)$, $c_{K_-}\in\textup{int}(K_-)$, and $c_{F}\in\textup{relint}(F)$ in 
			\cref{def:CT_triang_points} lie on a straight line and the values in the degrees of freedom in  \eqref{eq:commondof_K+K-} 
			along $F$ coincide. 
			Hence  \cref{thm:CT_K1K2} guarantees $C^1$ conformity of $\varphi_\ell|_{K_+\cup K_-}$ 
			on the side patch $\overline{\omega}(F)=K_+\cup K_-$. 
			Successive application of this argument leads to $\varphi_\ell\in C^1(\Omega)$ with support 
			$\textup{supp}(\varphi_\ell)=\overline{\omega_\ell}$ in the nodal patch of $z_\ell\in\mathcal{N}(\Omega)$. 
			In particular this implies vanishing boundary conditions. 
			Hence $\varphi_\ell\in P_3(\widehat{\mathcal{T}})\cap H^2_0(\Omega) =\textit{WF}(\mathcal{T})$ 
			is a $\textit{WF}$ 
			function. 
			For any $f\in \textit{WF}(\mathcal{T})$ 
			holds $f|_K\in \textit{WF}(K)$ 
			and  along  any interior face 
			the values in the degrees of freedom  (cf. \eqref{eq:commondof_K+K-}) coincide. 
			Hence the definition of  $\varphi_\ell|_K=\varphi_{\ell,K}$ for any $K\in\mathcal{T}$  and \cref{thm:CT_FECiarlet}
			imply $f\in \textup{span}\{\varphi_\ell:\, \ell=1,\dots, M\}$.
			This concludes the proof of (\ref{item:CT_basis_exists}). 
			\end{proofof}
				The proof of \cref{thm:Scaling_CT_fine}.\ref{item:CT_basis_scaling} follows the methodology of  
				\cite[Thm.~6.1.3]{Ciarlet78} with an affine equivalent alternative finite element.
				Consider a tetrahedron $K=\textup{conv}\{Q_1,Q_2,Q_3, Q_4\}$ 
				with vertices $Q_1,\dots, Q_4\in\mathcal{V}(K)$ 
				and edges $E_{jk}:=\textup{conv}\{Q_j,Q_k\}\in\mathcal{E}(K)$ 
				with midpoints $Q_{jk}:=(Q_j+Q_k)/2$ for $1\le j<k\le 4$. 
				\begin{definition}[alternative set of linear functionals]\label{def:alternativeFE}
					Define $28=16+12$ linear functionals 
					$ \mathcal{L}(K):=\big(L_{jk}|\,j,k=1,\dots, 4\big)\cup\big(L_{abc} |\,\{a,b,c,d\}=\{1,2,3,4\},a<b\big)$ on 
					$K=\textup{conv}\{Q_1,Q_2,$ $Q_3, Q_4\}$, 
					for $j,k=1,\dots, 4$ and all pairwise distinct $a,b,c\in\{1,2,3,4\}$ with $a<b$, by  
					\begin{align}
					\begin{split}
						L_{jk}:=\begin{cases}
										\delta_{Q_j}\quad&\text{if } j=k\\
										(Q_k-Q_j)\cdot\delta_{Q_j}\nabla &\text{if } j\not =k
									\end{cases}
						\qquad\text{and}\qquad
						L_{abc}:=(Q_c-Q_{ab})\cdot \delta_{Q_{ab}}\nabla.
					\end{split}\label{eq:dof_competing_FE}
					\end{align}
					(Recall  the point evaluation $\delta_P(f)=f(P)$ and $\delta_P\nabla(f)=\nabla f(P)$ 
					for any $P\in\overline{\Omega}$ and 
					$f\in C^1(\overline{\Omega})$.) 
				\end{definition} 
				\begin{lemma}[local basis]\label{lem:properties_alternativeFE}
					There exits a constant $C_a>0$ (that exclusively depends on $\mathbb{T}$), such that 
					for any tetrahedron $K\in\mathcal{T}\in\mathbb{T}$ with $\textit{WF}$ 
					partition $\widehat{\mathcal{T}}(K)$, there exists a basis 
					$\big(\varphi_{jk}:\,j,k=1,\dots, 4\big)\cup\big(\varphi_{abc}:\,\{a,b,c,d\}=\{1,2,3,4\},a<b\big)$ of 
					$\textit{WF}(K)$ 
					that satisfies the following. 
					The $28$ functions are dual to the linear functionals in 
					$\mathcal{L}(K)$ from \eqref{eq:dof_competing_FE}, in that 
					\begin{align}
						L_{jk}(\varphi_{\ell m})=\delta_{j\ell}\delta_{\ell m},\quad L_{abc}(\varphi_{def})=\delta_{ad}\delta_{be}\delta_{cf}, 
						\quad \text{and}\quad L_{jk}(\varphi_{abc})=0=L_{abc}(\varphi_{jk})	 \label{eq:duality_alternative}	
					\end{align}					  
					 for all $j,k,\ell,m=1,\dots, 4$,  
					all pairwise distinct $a,b,c\in\{1,2,3,4\}$ with $a<b$, and all pairwise distinct $d,e,f\in\{1,2,3,4\}$ with $d<e$,  
					and satisfy, for any $s=0,1,2$, that 
					\begin{align}
									&\sum_{j,k=1}^4\Vert {\varphi}_{jk}\Vert_{H^s(K)}
								+\sum_{\substack{\{a,b,c,d\}=\{1,2,3,4\}\\ a<b}}\Vert {\varphi}_{abc}\Vert_{H^s(K)}\le C_a h_K^{3/2-s}. 
								\label{eq:scaling_alternativ}
					\end{align}
				\end{lemma}
				\begin{proof}
				\begin{enumerate}[label=\textit{Step \arabic* --}, ref=Step \arabic*,wide] 
				\item\textit{linear independence. }
					For each $K:=\textup{conv}\{Q_1,Q_2,Q_3,Q_4\}\in\mathcal{T}\in\mathbb{T}$ 
					the linear functionals in $\mathcal{L}(K)$ from \cref{def:alternativeFE} are linearly independent. 
					The proof concerns $f\in \textit{WF}(K) 
					:=P_3(\widehat{\mathcal{T}}(K))\cap C^1(K)$ with $L_{jk}f=0=L_{abc}f$ 
					for all indices $j,k$ and $a,b,c$.  
					Hence  \eqref{eq:dof_competing_FE} implies that $f$ vanishes at the vertices $Q_1,\dots,Q_4\in\mathcal{V}(K)$ of $K$.
					Since for fixed $j\in\{1,2,3,4\}$ the three vectors $Q_j-Q_k\in\mathbb{R}^3$ 
					for $k\in\{1,2,3,4\}\setminus\{j\}$ are linearly independent, 
					 \eqref{eq:dof_competing_FE} implies that $\nabla f$ vanishes at the vertices $Q_1,\dots,Q_4\in\mathcal{V}(K)$ of $K$ as well. 
					 Moreover the tangential derivative $\tau_{ab}\cdot\nabla f$ vanishes along each edge 
					 $E_{ab}:=\textup{conv}\{Q_a,Q_b\}\in\mathcal{E}(K)$ of $K$ 
					 with unit tangent vector $\tau_{ab}$ for $1\le a<b\le 4$. 
					 Since $\textup{span}\{\tau_{ab}, Q_c-Q_{ab},Q_d-Q_{ab}\}=\mathbb{R}^3$ 
					 for all  $\{a,b,c,d\}=\{1,2,3,4\}$, \eqref{eq:dof_competing_FE} implies that 
					 $\nabla f(Q_{ab})$ vanishes in all midpoints $Q_{ab}:=(Q_a+Q_b)/2$ of edges 
					 $E_{ab}\in\mathcal{E}(K)$ of $K$ as well.
					We conclude that $L_1f=\dots=L_{28}f=0$ for the $\textit{WF}$ 
					degrees of freedom from \cref{tab:CT_3D_dof} 
					and so $f\equiv 0$ by \cref{thm:CT_FECiarlet}. 
					Hence $(K,\textit{WF}(K) 
					,\mathcal{L}(K))$ is a finite element in the sense of Ciarlet. 
					Consequently, there exist nodal basis functions 
					${\varphi}_{jk}\in \textit{WF}(K)$ 
					and  ${\varphi}_{abc}\in \textit{WF}(K)$ 
					with the required duality properties 
					\eqref{eq:scaling_alternativ}.
				\item\label{step:scaling_Tref}\textit{nodal basis in the reference tetrahedron $\widehat{T}$. }
					Let $\Phi: \widehat{T}\to K$ denote an affine diffeomorphism from  the reference tetrahedron  
					$\widehat{T}:=\textup{conv}\{0,e_1,e_2, e_3\}$ onto $K$. 
					The center points in the $\textit{WF}$ 
					partition $\widehat{\mathcal{T}}$ in  \cref{def:CT_triang_points} 
					are input parameter in \cref{alg:CT3D}. 
					 For each $K\in\mathcal{T}$ the center points $c_K\in \textup{int}(K)$ and $c_m\in \textup{relint}(F_m)$ for $m=1,\dots,4$ 
					of the faces  $F_1,\dots, F_4\in\mathcal{F}(K)$ of $K$  are mapped to five points in $\widehat{T}$. 
					Namely the point 
					$\widehat{c}_0:=\Phi^{-1}(c_K)$ and $\widehat{c}_m:=\Phi^{-1}(c_m)\in\widehat{F}_m=\Phi^{-1}(F_m)$ for $m=1,\dots,4$. 
					For each set of points $S:=(\widehat{c}_0,\widehat{c}_1,\dots, \widehat{c}_4)\in\mathbb{R}^{15}$ 
					suppressed in the notation, there 
					exists a set  of nodal basis functions 
					$\widehat{\varphi}_{jk}$ and $\widehat{\varphi}_{abc}$ for $j,k=1,\dots, 4$ and for all $\{a,b,c,d\}=\{1,2,3,4\}$ 
					with $a<b$ associated with $\widehat{\mathcal{L}}:=\mathcal{L}(\widehat{T})$ in \eqref{eq:dof_competing_FE}. 
					\cref{thm:shapeRegularity3D} shows that the parameters 
					$S:=(\widehat{c}_0,\widehat{c}_1,\dots, \widehat{c}_4)\in\mathbb{R}^{15}$ belong to a 
					compact set ${\mathcal{C}}$ of points in $\mathbb{R}^{15}$ such that  the distance to the respective relative boundaries is positive.
					Note that for each fixed parameter $S\in{\mathcal{C}}$ 
					the algorithm \cref{alg:CT3D} (with adapted preprocessing in step (\ref{item:Alg_CT_preprocessing})) maps
					 input data $x_{jk}$ and $x_{abc}$ onto a function $f\in \textit{WF}(\mathcal{T})$ 
					with $L_{jk}(f)=x_{jk}$ and $L_{abc}(f)=x_{abc}$ from 
					$\widehat{\mathcal{L}}:=\big(L_{jk}|\,j,k=1,\dots, 4\big)\cup\big(L_{abc} |\,\{a,b,c,d\}=\{1,2,3,4\},a<b\big)$.
					\cref{rem:Alg_continuity} points out that the evaluation of this interpolation operator depends continuously on the parameters in the compact 
					set ${\mathcal{C}}$. Hence the norms of the nodal basis functions 
					$\widehat{\varphi}_{jk}$ and $\widehat{\varphi}_{abc}$ for any $j,k=1,\dots, 4$ and for all $\{a,b,c,d\}=\{1,2,3,4\}$ with $a<b$ 
					for the resulting $\textit{WF}$ 
					partition of  the reference tetrahedron $\widehat{T}$  
					remain 
					bounded by a universal 
					constant $C_1>0$ as in \cite[Eq. (6.1.30)]{Ciarlet78} for $2$D,
					\begin{align}
								\sum_{j,k=1}^4\Vert \widehat{\varphi}_{jk}\Vert_{H^s(\widehat{T})}
								+\sum_{\substack{\{a,b,c,d\}=\{1,2,3,4\}\\ a<b}}\Vert \widehat{\varphi}_{abc}\Vert_{H^s(\widehat{T})}\le C_1. 
								\label{eq:proof_scalingprep_1}
					\end{align}  
					Note that the constant $C_1$  depends exclusively on ${\mathcal{C}}$ and so on $\varepsilon>0$ from \cref{thm:shapeRegularity3D}. 
				\item \textit{affine equivalence. }
					Recall the affine diffeomorphism 
					$\Phi: \widehat{T}\to K$  from  the reference tetrahedron  $\widehat{T}:=\textup{conv}\{0,e_1,e_2, e_3\}$ onto $K$ from \ref{step:scaling_Tref}. 
					The point is that the alternative finite element  $(K,\textit{WF}(K) 
					,\mathcal{L}(K))$ is affine equivalent to 
					$(\widehat{T}, \textit{WF}(\mathcal{T}) 
					,\widehat{\mathcal{L}})$ in that  
					$\widehat{\mathcal{L}}=\mathcal{L}(K)\circ \Phi$; cf. \cite{Ciarlet78,BS08} for the concept of affine equivalence.
					In particular the local nodal basis $\varphi_{jk}\in \textit{WF}(K)$ 
					and $\varphi_{abc}\in \textit{WF}(K)$ 
					associated with 
					\eqref{eq:dof_competing_FE} and the local nodal basis functions 
					$\widehat{\varphi}_{jk}\in \textit{WF}(\mathcal{T})$ 
					and  $\widehat{\varphi}_{abc}\in \textit{WF}(\mathcal{T})$ 
					on the reference tetrahedron 
					$\widehat{T}$ satisfy 
					$${\varphi_{jk}}=\widehat{\varphi}_{jk}\circ\Phi^{-1}\quad\text{and}\quad {\varphi_{abc}}=\widehat{\varphi}_{abc}\circ\Phi^{-1}$$ 
					for any 
					$j,k=1,\dots, 4$ and all pairwise distinct $a,b,c\in\{1,2,3,4\}$ with $a<b$. 
					Standard scaling properties in \cite[Thm.~3.1.2--3.1.3]{Ciarlet78} lead to a 
					universal constant $C_2$ (that depends only on the reference tetrahedron $\widehat{T}$) such that for any $s=0,1,2$ 
					\begin{align*}
						\Vert\varphi_{abc}\Vert_{H^s(K)}
									\le C_2|K|^{1/2}\varrho_K^{-s}\Vert \widehat{\varphi}_{abc}\Vert_{H^s(\widehat{T})}
					\quad\text{and}\quad
						\Vert\varphi_{jk}\Vert_{H^s(K)}
									\le C_2|K|^{1/2}\varrho_K^{-s}\Vert \widehat{\varphi}_{jk}\Vert_{H^s(\widehat{T})}
					\end{align*}
					with the radius $\varrho_K$ of the biggest incircle in $K$. The  shape regularity of $K\in\mathcal{T}\in\mathbb{T}$ 
					guarantees $\varrho_K\approx h_K\approx |K|^{1/3}$ and so 
					\begin{align*}
						|K|^{1/2}\varrho_K^{-s}\le C_3 h_K^{3/2-s}.
					\end{align*}
					The combination of the last two displayed estimates with \eqref{eq:proof_scalingprep_1}
					concludes the proof of \eqref{eq:scaling_alternativ} with $C_a=C_1C_2C_3$.
					\end{enumerate}
				\end{proof}
					
					The point in the proof of  \cref{thm:Scaling_CT_fine}.\ref{item:CT_basis_scaling} below 
					is that each (global) nodal basis function $\varphi_\ell$ for $\ell=1,\dots, M$ 
					restricted to one tetrahedron 
					$K\subset\textup{supp}(\varphi_\ell)$ is identical with its interpolation in the local basis from \cref{lem:properties_alternativeFE}, i.e., 
					\begin{align}
						\varphi_\ell|_K
						=\sum_{j,k=1}^4L_{jk}(\varphi_\ell)\varphi_{jk}+\sum_{\substack{\{a,b,c,d\}=\{1,2,3,4\}\\ a<b}}L_{abc}(\varphi_\ell)\varphi_{abc}.
						\label{eq:scaling_interpolation}
					\end{align}
					The interpolation \eqref{eq:scaling_interpolation} and a triangle inequality show 
					\begin{align}
						\Vert \varphi_\ell\Vert_{L^2(K)}
						\le\sum_{j,k=1}^4\vert L_{jk}(\varphi_\ell)\vert \Vert\varphi_{jk}\Vert_{L^2(K)}
							+\sum_{\substack{\{a,b,c,d\}=\{1,2,3,4\}\\ a<b}}\vert L_{abc}(\varphi_\ell)\vert \Vert\varphi_{abc}\Vert_{L^2(K)}.
							\label{eq:scalingsum2beanalyzed}
					\end{align} 
					A careful analysis of the contributions $|L_{jk}(\varphi_\ell)|$ and $|L_{abc}(\varphi_\ell)|$ on the right-hand side 
					of \eqref{eq:scalingsum2beanalyzed} for each of the three types of global nodal basis functions from \cref{def:CT_dof_global} 
					and the application of the scaling in \eqref{eq:scaling_alternativ} conclude the proof.
 				
 				\begin{proofof}\textit{\cref{thm:Scaling_CT_fine}.\ref{item:CT_basis_scaling} for $\varphi_\ell=\varphi_{E,\mu}$ for $E\in\mathcal{E}(\Omega)$.}
							We discuss the scaling first for the basis function related to the global degree of freedom 
							\begin{align*}
								L_{E,\mu}:=\nu_{E,\mu}\cdot \delta_{\textup{mid}(E)}\nabla \quad \text{ for }\mu=1,2 
											\text{ and }E\in\mathcal{E}.
							\end{align*}
							Consider one tetrahedron $K:=\textup{conv}\{Q_1,\dots,Q_4\}\in\mathcal{T}(E)$ with $\textit{WF}$ 
							partition $\widehat{\mathcal{T}}(K)$.
							Suppose, without loss of generality, $E=\textup{conv}\{Q_1,Q_2\}=E_{12}$. 
							Note $\varphi_\ell|_K\in \textit{WF}(K)$ 
							satisfies $L_{E,\mu} \varphi_\ell=1$ and all the other degrees of freedom of $\varphi_\ell$ 
							in \cref{def:CT_dof_global} vanish.
							In particular, $\varphi_\ell$ and $\nabla \varphi_\ell$ vanish at $Q_1,\dots, Q_4$. This implies in 	\eqref{eq:scalingsum2beanalyzed}, that  						
							$$ L_{jk}(\varphi_\ell)=0 \quad \text{for all }j,k=1,\dots,4.$$
							Along any edge $E_{jk}\in \mathcal{E}(K)$ with unit tangent vector $\tau_{jk}$ for $1\le j<k\le 4$ 
							we consider $\varphi_\ell|_{E_{jk}}\in P_3(E_{jk})$ as a one-dimensional cubic polynomial with four 
							vanishing degrees of freedom. 
							Hence the one-dimensional derivative $\varphi_\ell^\prime=\tau_{jk}\cdot \nabla \varphi_\ell=0$ vanishes. 
							For any edge $E_{jk}\in\mathcal{E}(K)\setminus\{E_{12}\}$ in $K$ besides $E_{12}$ the gradient 
							$\nabla \varphi_\ell$ in the midpoint $Q_{jk}$
							vanishes in the direction $\nu_{jk,1}$ and $\nu_{jk,2}$. In other words $\nabla \varphi_\ell(Q_{jk})=0$ for all 
							$(j,k)\not =(1,2)$ with $1\le j<k\le4$. 									This implies  
							$$ L_{abc}(\varphi_\ell)=0 \quad\text{ for all }(a,b,c)\not=(1,2,3)\text{ or }  (a,b,c)\not=(1,2,4).$$
							Consequently, \eqref{eq:scaling_interpolation} reads
							\begin{align*}
								\varphi_\ell
									= L_{123}(\varphi_\ell)\varphi_{123}+L_{124}(\varphi_\ell)\varphi_{124}
								=(Q_3-Q_{12})\cdot \nabla \varphi_\ell(Q_{12})\varphi_{123}+(Q_4-Q_{12})\cdot \nabla \varphi_\ell(Q_{12})\varphi_{124}.
							\end{align*}
							Recall $\tau_{12} \cdot \nabla \varphi_\ell(Q_{12})\equiv 0$ for the unit tangent vector $\tau_{12}\equiv \tau_E$ 
							of the edge $E_{12}\equiv E $ and 	$L_{E,\kappa}(\varphi_\ell)=\nu_{E,\kappa}\cdot\nabla \varphi_\ell(Q_{12})=0$ 
							for the normal $\nu_{E,\kappa}$ with $\kappa\in \{1,2\}\setminus\{\mu\}$. 
							Since $L_\ell(\varphi_\ell)=\nu_{E,\mu}\cdot \nabla \varphi_\ell(Q_{12})=1$, this implies for $k=3,4$ 
							\begin{align*}
								|L_{12k}(\varphi_\ell)|
									&= \big|(Q_k-Q_{12})\cdot \big((\nu_{E,\mu}\cdot \nabla \varphi_\ell(Q_{12}))\nu_{E,\mu}
									+(\nu_{E,\kappa}\cdot \nabla \varphi_\ell(Q_{12}))\nu_{E,\kappa}
																+(\tau_E \cdot \nabla \varphi_\ell(Q_{12}))\tau_{E}\big)\big|
									\\& = \big|(Q_k-Q_{12})\cdot \nu_{E,\mu}|\le h_K.
							\end{align*}
							For  \eqref{eq:scalingsum2beanalyzed} this implies that 
							\begin{align*}
								\Vert \varphi_\ell\Vert_{H^s(K)}=
								\Vert \varphi_{E,\kappa}\Vert_{H^s(K)}
									\le h_K\sum_{k=3,4} \Vert\varphi_{12\ell}\Vert_{H^s(K)}
									\le C_a h_K^{5/2-s}
							\end{align*}
							with the upper bound from \cref{lem:properties_alternativeFE} in the last step. 
							The shape regularity of $\mathcal{T}\in\mathbb{T}$ implies $|K|\approx h_K^3$ and the boundedness of the number of tetrahedra 
							$|\mathcal{T}(E)|$ in the support of $\varphi_{E,\kappa}$. This concludes the proof.  
					\end{proofof}
					\begin{proofof}\textit{\cref{thm:Scaling_CT_fine}.\ref{item:CT_basis_scaling} for $\varphi_\ell=\varphi_{z,1}$ for $z\in\mathcal{V}(\Omega)$. }   
							Consider the basis function associated with  
							\begin{align*}
								L_{z,1}:=\delta_{z} \quad \text{ for }z\in\mathcal{V}
							\end{align*}
							and  $K:=\textup{conv}\{Q_1,\dots,Q_4\}\in\mathcal{T}(z)$ with $\textit{WF}$  
							partition $\widehat{\mathcal{T}}(K)$. 
							Suppose, without of loss of generality, $z= Q_1$. 
							Note $\varphi_\ell|_K\in \textit{WF}(K)$ 
							satisfies $L_{z,1}(\varphi_\ell)=\varphi_\ell(Q_1)=1=L_{11}(\varphi_\ell)$ 
							and all the other degrees of freedom in \cref{def:CT_dof_global} 
							of $\varphi_\ell$ vanish. 
							Since $\varphi_\ell$  vanishes at $Q_2,Q_3,Q_4$ and $\nabla \varphi_\ell$ vanish at $Q_1,Q_2,Q_3,Q_4$, \eqref{eq:dof_competing_FE} implies 
								$$ L_{kk}(\varphi_\ell)=0\text\quad\text{and}\quad L_{jk}(\varphi_\ell)=0\quad \text{ for } 1\le j<k \le 4.$$
							For any edge $E\in\mathcal{E}$ with midpoint $\textup{mid}(E)$ holds 
							$L_{E,1}(\varphi_\ell)=\nu_{E,1}\cdot\nabla f(\textup{mid}(E)) = 0=\nu_{E,2}\cdot\nabla \varphi_\ell(\textup{mid}(E))
							=L_{E,2}(\varphi_\ell)$.
							This implies for all pairwise distinct $a,b,c\in\{1,2,3,4\}$ with $a<b$ 
							and edges $E_{ab}\in\mathcal{E}(K)$ with unit tangent vector $\tau_{ab}$ and midpoint $Q_{ab}$, that 
								$$ L_{abc}(\varphi_\ell)= (Q_c-Q_{ab})\cdot \nabla \varphi_\ell(Q_{ab})
											= (Q_c-Q_{ab})\cdot (\tau_{ab}\cdot \nabla \varphi_\ell(Q_{ab})) \tau_{ab}.$$
							Along any $E_{jk}\in \mathcal{E}(K)$ for $2\le j<k\le 4$ in $K$ without the vertex $Q_1$ with unit tangent vector $\tau_{jk}$, 
							consider $\varphi_\ell|_{E_{jk}}\in P_3(E_{jk})$ as a one-dimensional cubic polynomial with vanishing degrees of freedom. 
							Hence the one-dimensional derivative $\varphi_\ell^\prime=\tau_{jk}\cdot \nabla \varphi_\ell=0$ vanishes. 
							In other words, 
								$$L_{abc}(\varphi_\ell)= (Q_c-Q_{ab})\cdot (\tau_{ab}\cdot \nabla \varphi_\ell(Q_{ab})) \tau_{ab} =0 \quad\text{ for } 2\le a<b\le 4$$
							and \eqref{eq:scaling_interpolation} reads
							\begin{align*}
								\varphi_\ell	
									=L_{11}(f)\varphi_{11}+ \sum_{\substack{2\le b, c\le 4\\ b\not=c}}L_{1bc}\varphi_{1bc}
								=\varphi_{11}+ \sum_{\substack{2\le b, c\le 4\\ b\not=c}}
									(Q_c-Q_{1b})\cdot \big((\tau_{{1b}}\cdot \nabla \varphi_\ell(Q_{1b}))\tau_{{1b}}\big)\varphi_{1bc}.	
							\end{align*}
							Along any edge $E=E_{1b}\in\mathcal{E}(K)$  for $b=2,3,4$ with unit tangent vector $\tau_E\equiv\tau_{1b}$ 
							consider $\varphi_\ell|_E\in P_3(E)$ as one-dimensional function uniquely determined by the given values of $\varphi_\ell$ 
							and the tangential derivative $\tau_E\cdot\nabla \varphi_\ell$ in the endpoints $Q_1$ and $Q_b$ of $E$. 
							The parametrisation $f_\tau(s)$ for $s\in (0,h)$ of $\varphi_\ell|_E$ with $f_\tau(0)=\varphi_\ell(Q_1)=1$, 
							$f_\tau (h)=\varphi_\ell(Q_b)=0$, and 
							$f^\prime_\tau=\tau_E\cdot\nabla \varphi_\ell$ vanishing in both endpoints, 
							reads $$f_\tau(s)=h^{-3}(s-h)^2(2s+h).$$ 
							Hence,  
							$|\tau_E\cdot\nabla \varphi_\ell(Q_{1b})|=|f_\tau^\prime(h/2)|=1.5 h^{-1}\le C_4 h_K^{-1}$.					
							In other words, for all $b=2,3,4$,  
							$$|L_{1bc}(\varphi_\ell)|= |(Q_c-Q_{1b})\cdot (\tau_{{1b}}\cdot\nabla \varphi_\ell(Q_{1b}))\tau_{{1b}}|
								\le |(Q_c-Q_{1b})\cdot \tau_{{1b}}|C_4 h_K^{-1}\le C_4.
								$$	
							For  \eqref{eq:scalingsum2beanalyzed} this implies that 
							\begin{align*}
								\Vert \varphi_\ell \Vert_{H^s(K)}
									=\Vert \varphi_{z,1}\Vert_{H^s(K)}
									= \Vert \varphi_{11}\Vert_{H^s(K)}+ C_4 \sum_{\substack{2\le b, c\le 4\\ b\not=c}}\Vert\varphi_{1bc}\Vert_{H^s(K)}
									\le \max\{1,C_4\} C_a h_K^{3/2-s}
							\end{align*}
							with the upper bound from \cref{lem:properties_alternativeFE} in the last step. 	
							The shape regularity of $\mathcal{T}\in\mathbb{T}$ implies $|K|\approx h_K^3$ and the boundedness of the number of tetrahedra 
							$|\mathcal{T}(z)|$ in the support of $\varphi_{z,1}$. This concludes concludes the proof.
						\end{proofof}
					\begin{proofof}\textit{\cref{thm:Scaling_CT_fine}.\ref{item:CT_basis_scaling} for $\varphi_\ell=\varphi_{z,\kappa+1}$ for $z\in\mathcal{V}(\Omega)$. } 
							Finally, regard the basis function associated with 
							\begin{align*}
								L_{z,\kappa+1}:=e_\kappa\cdot \delta_{z}\nabla \quad \text{ for }\kappa=1,2,3, \text{ and }z\in\mathcal{V}
							\end{align*}
							and consider  $K:=\textup{conv}\{Q_1,\dots,Q_4\}\in\mathcal{T}(z)$ with $\textit{WF}$ 
							partition $\widehat{\mathcal{T}}(K)$. 
							Suppose, without of loss of generality, $z=Q_1$. 
							Note $\varphi_\ell|_K\in \textit{WF}(K)$ 
							satisfies $L_{z,\kappa+1}\varphi_\ell=\partial \varphi_\ell/\partial x_\kappa=1$ 
							and all the other degrees of freedom of $\varphi_\ell$ vanish.
							Since $\varphi_\ell$  vanishes at $Q_1,Q_2,Q_3,Q_4$ and $\nabla \varphi_\ell$ vanish at $Q_2,Q_3,Q_4$, 
							the definition in \eqref{eq:dof_competing_FE} shows 
								$$ L_{jj}(\varphi_\ell)=0\quad\text{ for } j=1,\dots,4\quad
								\text{and}\quad L_{jk}(\varphi_\ell)=0\quad \text{ for } 2\le j<k \le 4.$$
							For any edge $E\in\mathcal{E}$ with midpoint $\textup{mid}(E)$ holds 
							$L_{E,1}(\varphi_\ell)=\nu_{E,1}\cdot\nabla \varphi_\ell(\textup{mid}(E)) = 0
							=\nu_{E,2}\cdot\nabla \varphi_\ell(\textup{mid}(E))=L_{E,2}(\varphi_\ell)$.
							This implies all pairwise distinct $a,b,c\in\{1,2,3,4\}$ with $a<b$ 
							and edges $E_{ab}\in\mathcal{E}(K)$ with unit tangent vector $\tau_{ab}$ and midpoint $Q_{ab}$, that 
								$$ L_{abc}(f)= (Q_c-Q_{ab})\cdot \nabla \varphi_\ell(Q_{ab})
											= (Q_c-Q_{ab})\cdot (\tau_{ab}\cdot \nabla \varphi_\ell(Q_{ab})) \tau_{ab}.$$
							Along any $E_{jk}\in \mathcal{E}(K)$ for $2\le j<k\le 4$ in $K$ without the vertex $Q_1$ with unit tangent vector $\tau_{jk}$, 
							consider $\varphi_\ell|_{E_{jk}}\in P_3(E_{jk})$ as a one-dimensional cubic polynomial with vanishing degrees of freedom. 
							Hence the one-dimensional derivative $\varphi_\ell^\prime=\tau_{jk}\cdot \nabla \varphi_\ell=0$ vanishes. 
							In other words, 
								$$L_{abc}(\varphi_\ell)= (Q_c-Q_{ab})\cdot (\tau_{ab}\cdot \nabla \varphi_\ell(Q_{ab})) \tau_{ab} =0 \quad\text{ for } 2\le a<b\le 4$$
							and \eqref{eq:scaling_interpolation} reads
							\begin{align*}
								\varphi_\ell
								&=\sum_{k=2}^4 L_{1k}(\varphi_\ell)\varphi_{1k}+ \sum_{\substack{2\le b, c\le 4\\ b\not=c}}L_{1bc}(\varphi_\ell)\varphi_{1bc}
								\\&=\sum_{k=2}^4 (Q_k-Q_1)\cdot\nabla\varphi_\ell(Q_1) \varphi_{1k}+
								\sum_{\substack{2\le b, c\le 4\\ b\not=c}} 
									(Q_c-Q_{1b})\cdot \big((\tau_{{1b}}\cdot \nabla \varphi_\ell(Q_{1b}))\tau_{{1b}}\big)\varphi_{1bc}.		
							\end{align*}
							Since $\partial \varphi_\ell/\partial x_{\ell}(Q_1)=0$ for $\ell\in \{1,2,3\}\setminus\{\kappa\}$ 
							and $\partial \varphi_\ell/\partial x_{\kappa}(Q_1)=L_{z,\kappa+1}(\varphi_\ell)=1$, it holds for $k=2,3,4$
							$$ |L_{1k}(\varphi_\ell)|=|(Q_k-Q_1)\cdot\nabla \varphi_\ell(Q_1)|
											=|(Q_k-Q_1)\cdot e_\kappa \partial \varphi_\ell/\partial x_{\kappa}(Q_1)|
											=|(Q_k-Q_1)\cdot e_\kappa|\le h_K. $$
							Along any edge $E=E_{1b}\in\mathcal{E}(K)$  for $b=2,3,4$ with unit tangent vector $\tau_E=\tau_{1b}$ 
							consider $\varphi_\ell|_E\in P_3(E)$ as one-dimensional function uniquely determined by the given values of $\varphi_\ell$ 
							and the tangential derivative $\tau_E\cdot\nabla \varphi_\ell$ in the endpoints $Q_1$ and $Q_b$ of $E$. 
							Let $c:=e_\kappa\cdot\tau_E\le 1$. 
							The parametrisation $f_\tau(s)$ for $s\in (0,h)$ of $f|_E$ with $f_\tau(0)=\varphi_\ell(Q_1)=0$, $f_\tau(h)=\varphi_\ell(Q_b)=0$, and 
							$f^\prime_\tau=\tau\cdot\nabla \varphi_\ell$ with $f^\prime(0)=c$ and $f^\prime(h)=0$, 
							reads $$f_\tau(s)=c(s-h)^2s/h^2.$$ Hence 
							$|\tau_{1b}\cdot\nabla \varphi_\ell(Q_{1b})|=|f^\prime_\tau (h/2)|= c/4\le 1$ and 
							$$|L_{1bc}(f)|= |(Q_c-Q_{1b})\cdot (\tau_{{1b}}\cdot\nabla \varphi_\ell(Q_{1b}))\tau_{{1b}}|
								\le |(Q_c-Q_{1b})\cdot \tau_{{1b}}|\le h_K.
								$$	
							For  \eqref{eq:scalingsum2beanalyzed} this implies that 
							\begin{align*}
								\Vert \varphi_\ell\Vert_{H^s(K)}
									=\Vert \varphi_{z,\kappa+1}\Vert_{H^s(K)}
									\le h_K \sum_{k=2}^4 \Vert {\varphi}_{1k}\Vert_{H^s(K)} 
										+ h_K \sum_{\substack{2\le b, c\le 4\\ b\not=c}}\Vert\varphi_{1bc}\Vert_{H^s(K)}
									\le C_a h_K^{5/2-s}
							\end{align*}
							with the upper bound from \cref{lem:properties_alternativeFE} in the last step. 
							The shape regularity of $\mathcal{T}\in\mathbb{T}$ implies $|K|\approx h_K^3$ and the boundedness of the number of tetrahedra 
							$|\mathcal{T}(z)|$ in the support of $\varphi_{z,\kappa+1}$. This  concludes the proof.
					\end{proofof}		 
\FloatBarrier
\newpage
 \section{Conforming companion} \label{sec:Appendix_Companion}
 \subsection{Overview}
 		Given a regular triangulation $\mathcal{T}\in\mathbb{T}$ and the Morley finite element space $M(\mathcal{T})$, 
 		this sections contains the details of the definition of the conforming companion operator $J_M:M(\mathcal{T})\to V:=H^2_0(\Omega)$ 
 		in $3$D and the proof of Theorem~3.1. 
 		\begin{extratheorem}[properties of $J_M$]
			 	There exists a constant $M_2\approx 1$ (that exclusively depends on $\mathbb{T}$) 
			 	and a conforming companion $J_Mv_M\in V:=H^2_0(\Omega)$ for any $v_M\in M(\mathcal{T})$ with 
				 \begin{enumerate}[label=(\alph*), ref=\alph* ,wide]
			 		\item\label{item:App_J_rightInverse} 
			 				$J_M$ is a right inverse to the interpolation $I_M$ in that
			 					$
			 						I_M\circ J_M =\textup{id} \text{ in  } {M(\mathcal{T})}, 
			 					$
					\item\label{item:App_J_approx} 
						$
							\Vert h_{\mathcal{T}}^{-{2}} (1-J_M)v_{M}\Vert_{L^2(\Omega)}
							+\vvvert (1-J_M)v_{M}\vvvert_{\mathrm{pw}}
							\le {M_2} \min_{v\in V}\vvvert v_{M}-v\vvvert_{\mathrm{pw}}
						$,			
			 		\item\label{item:App_J_orthogonality}
						the orthogonality	$(1-J_M)(M(\mathcal{T}))\perp P_{2}(\mathcal{T})$  
						holds in $L^2(\Omega)$.		
			 	\end{enumerate}
			\end{extratheorem}
 		For the triangulation $\mathcal{T}\in\mathbb{T}$ with set edges $\mathcal{E}$ 
		fix one unit tangent vector $\tau_E\parallel E$ and choose two unit vectors $\nu_{E,1},\,\nu_{E,2}\in\mathbb{R}^3$ for any edge 
		$E\in\mathcal{E}$ 
		such that $(\tau_E,\,\nu_{E,1},\,\nu_{E,2})$ form an orthonormal basis of $\mathbb{R}^3$, and define the $\textit{WF}$ 
		partition $\widehat{\mathcal{T}}$ as in \cref{def:CT_triang_points}.
		Consider the conforming subspace with homogeneous boundary conditions as in \eqref{eq:def_HCT} 
		\begin{align*}
				\textit{WF} (\mathcal{T})
					:=P_3(\widehat{\mathcal{T}})\cap H^2_0(\Omega). 
		\end{align*}
		The design in Sections~\ref{sec:Morley_J1}--\ref{sec:Morley_J4} below contains four steps. Given any Morley function 
		$v_M\in M(\mathcal{T})$, first the enrichment $J_1v_M\in \textit{WF}(\mathcal{T})$ 
		is computed via averaging of nodes. 
		The global degrees of freedom for the Morley finite element read  
		\begin{align}
			L_E(f):= \intmean_{E} f\,\textup{d}s 
			\text{ and }  L_F(f):=\intmean_{F} \nabla f\cdot\nu_F\,\textup{d}\sigma\text{ for  any } E\in \mathcal{E},\, F\in \mathcal{F},
			\text{ and }	f\in H^2(\Omega).\label{eq:Morley_dof}
		\end{align}	
		Note, for a Morley function $v_M\in M(\mathcal{T})$ the integral means in \eqref{eq:Morley_dof} are single valued, 
		in that for any edge $E\in\mathcal{E}$ and any adjacent tetrahedron 
		$T\in\mathcal{T}(E)$ holds 
		$\intmean_E v_M|_T\,\textup{d}s=L_E(v_M)$ and 
		for any face $F\in\mathcal{F}$ and any adjacent tetrahedron 
		$T\in\mathcal{T}(F)$ holds 
		$\intmean_{F} \nabla v_M|_T\cdot\nu_F\,\textup{d}\sigma=L_F(v_M)$. 
		Two separate corrections are necessary to guarantee 
			$L_E(v_M)=\intmean_E v_M\,\textup{d}s=L_E(J_Mv_M)$ for any $E\in\mathcal{E}$ and 
			$L_F(v_M)=\intmean_F\nabla v_M\cdot \nu_F\,\textup{d}\sigma=L_F(J_Mv_M)$ for any $F\in\mathcal{F}$ in $3$D. 
		The integral means over the edges are corrected in $J_2$. The integral means of the 
		normal derivatives along the faces are corrected in $J_3$. This guarantees that $J_3$ is a right inverse of $I_M$ 
		in that $I_MJ_3v_M=v_M$. Finally the correction $J_4v_M$ 
		is designed such that its $L^2$ projection  in $P_2(\mathcal{T})$ equals  
		$v_M\in M(\mathcal{T})$, $\Pi_2J_Mv_M=v_M$. 
		Note that in an abstract point of view \cref{sec:Morley_J1} concerns nodal values for vertices $z\in\mathcal{V}$ and midpoints 
		$\textup{mid}(E)$ of edges $E\in\mathcal{E}$ ($0$-simplices), 
		\cref{sec:Morley_J2} concerns the integral mean for each edge $E\in\mathcal{E}$ ($1$-simplex), 
		\cref{sec:Morley_J3}  concerns one integral mean for each face $F\in\mathcal{F}$ ($2$-simplex), and finally 
		\cref{sec:Morley_J4} concerns the volume contribution for each tetrahedron $T\in\mathcal{T}$ ($3$-simplex).
	    The corrections in the last three steps are independent of the choice of the conforming space $\textit{WF}(\mathcal{T})$ 
		and possible in any space dimension; cf. \cite{Gal15, VZ19} for $n=2$ and 
		Theorem~3.1.\ref{item:App_J_rightInverse}--\ref{item:App_J_approx}.
		
	\subsection{Design of $J_1$}\label{sec:Morley_J1}
		Determine the values for the degrees of freedom of $\textit{WF}(\mathcal{T})$ 
		by averaging the respective values of the Morley function on all 
			adjacent simplices $T\in\mathcal{T}$ and assure homogeneous boundary conditions to define 
			$J_1:M(\mathcal{T})\to \textit{WF}(\mathcal{T})\subset V$ 
			as follows.
		\begin{definition}[enrichment $J_1$]\label{def:J1}		
			For any vertex $z\in\mathcal{V}$ with nodal patch $\omega(z)=\textup{int}(\bigcup \mathcal{T}(z))$,
			$\mathcal{T}(z):=\{T\in\mathcal{T}|\, z\in\mathcal{V}(T)\}$, of $|\mathcal{T}(z)|$ many tetrahedra with the vertex $z$ and 
			for any multi-index $\alpha\in\mathbb{N}_0^3$ with 
			$0\le|\alpha|\le 1$, set 
				\begin{align}
					(D^\alpha J_1v_M)(z)&:=
										\frac{1}{|\mathcal{T}(z)|}\sum_{T\in\mathcal{T}(z)}(D^\alpha v_M|_T)(z)\quad 
										\text{ for }z\in \mathcal{V}(\Omega)
										 \label{eq:def_J1_z}
				\end{align}
				$\big(\text{and } =0 \text{ for } z\in \mathcal{V}(\partial\Omega)$ owing to the homogeneous boundary conditions in $V\big )$.
			For any edge $E\in\mathcal{E}$ with midpoint $\textup{mid}(E)$, edge patch $\omega(E)=\textup{int}(\bigcup \mathcal{T}(E))$,  
			$\mathcal{T}(E):=\{T\in\mathcal{T}|\, E\in\mathcal{E}(T)\}$, of $|\mathcal{T}(E)|$ many tetrahedra with common edge $E$, 
			and unit tangential vector $\tau_E$ with a fixed orientation, set
				\begin{align}
					\tau_E\times\nabla J_1v_M(\textup{mid}(E))&:=
										\frac{1}{|\mathcal{T}(E)|}\sum_{T\in\mathcal{T}(E)}\tau_E\times
										(\nabla v_M|_T)(\textup{mid}(E))\text{ for }
										E\in \mathcal{E}(\Omega)	\label{eq:def_J1_E}		
				\end{align}
			$\big(\text{and } =0 \text{ for }E\in \mathcal{E}(\partial \Omega)$ owing to the homogeneous boundary conditions in $V\big )$. 
			That means equality in the 		
			directions  $\{\nu_{E,1},\nu_{E,2}\}$ perpendicular  to the tangential unit vector $\tau_E$.
		\end{definition}
			The enrichment operator $J_1$ satisfies the approximation condition as shown in \cref{lem:prop_J1} below. 
			We aim to provide explicit constants whenever possible. To this end we need in the proof of \cref{lem:prop_J1} 
			the following \cref{lem:fundamentalInequality} as addition to the 
			related observations in \cite[Lem.~B--C]{CP18}, which may be of independent interest. Similar results are contained in 
			\cite[Lem~4.2]{CH17}. 
	\begin{lemma}[inequality]\label{lem:fundamentalInequality}
		Any $x=(x_1,\dots,x_m)\in\mathbb{R}^m$, $m\in\mathbb{N}$ satisfies 
		\begin{align}
			2\Big(1-\cos\Big(\frac{\pi}{2m+1}\Big)\Big)|x|^2
			\le \sum_{j=1}^{m-1}(x_{j+1}-x_j)^2+\min\{x_1^2,\dots, x_m^2\}. \label{eq:inequality}
		\end{align}
		The factor $2\Big(1-\cos\Big(\frac{\pi}{2m+1}\Big)\Big)$ in \eqref{eq:inequality} 
		is sharp in that there exists some $x\in \mathbb{R}^m\setminus\{0\}$ with equality in \eqref{eq:inequality}.
	\end{lemma}
		\begin{proof}
		   	Abbreviate 
		   	$\varrho(x)^2:=\sum_{j=1}^{m-1}(x_{j+1}-x_j)^2
		   	+\min\{x_1^2,\dots, x_m^2\}$ 
		   	for any $x=(x_1,\dots,x_m)\in\mathbb{R}^m\setminus\{0\}$, let $S:=\{x\in\mathbb{R}^m|\,|x|=1\}$ denote the unit sphere in $\mathbb{R}^m$,
		   	and $S_+:=\{x\in S|\, 0\le x_j\text{ for all }j=1,\dots, m\}$ its positive half.
		   	Then $\varrho((|x_1|,\dots, |x_m|))\le \varrho(x)$ follows from $(|x_{j+1}|-|x_j|)^2\le |x_{j+1}-x_j|^2$ for $j=1,\dots, m-1$.
		   	This implies the equality after the definition of 
		   	\begin{align*}
				\Lambda:=\min_{x\in S}{\varrho(x)^2}=\min_{x\in S_+}\varrho(x)^2.
			\end{align*}
			Lemma B in the appendix of \cite{CP18} shows that a minimum of the sum $\sum_{j=1}^{m-1}(x_{j+1}-x_j)^2$ 
			is attained for ordered coefficients in the vector $x\in S_+$, hence
			\begin{align*}
				\Lambda=\min_{\substack{x\in S_+\\0\le x_1\le x_2\le\dots\le x_m}}{\varrho(x)^2}=x\cdot B x
			\end{align*} 	
			with the tridiagonal matrix
			$$
				B:=
					\begin{pmatrix}
						2 	&-1 		& 0 		 &\cdots		 &0\\
						-1	&2		&-1		 &\ddots 	 &\vdots\\
						0 	&\ddots&\ddots &\ddots 	 &0\\
						\vdots&   \ddots   &  -1 		&2& -1\\
						0  &	\cdots		&			 0&-1	 		 &1
					\end{pmatrix}.
			$$ 
			The positive eigenvalues $0<\lambda_k=2\big(1-\cos({(2k-1)\pi}/({2m+1}))\big)$ for $k=1,\dots,m$ of the $B$ are computed 
			in \cite[Thm.3.2.viii]{YC08},
			 hence from the Rayleigh quotient follows $\lambda_1|x|^2\le x\cdot B x$ for any $x\in\mathbb{R}^m$ and 
			 $\Lambda=\lambda_1=2\big(1-\cos\big(\pi/(2m+1)\big)\big)$.
			 The eigenspace $E(\lambda_1)$ is spanned by 
			 $x=\big(\sin(\pi/(2m+1)),\sin(2\pi/(2m+1)),\dots, \sin(m\pi/(2m+1))\big)\in\mathbb{R}^m\setminus\{0\}$.
			 Since $0\le \sin(k\pi/(2m+1))\le\sin((k+1)\pi/(2m+1))$ for any $k=1,\dots,m-1$,  
			 the normed eigenvector $v:=x/|x|\in E(\lambda_1)\cap S_+$ 
			 has ordered coefficients and leads to the equality $\Lambda=\varrho(v)^2$; in other words
			 $\Lambda=\lambda_1$ is the optimal constant.  
		\end{proof}
			\begin{lemma}[properties of $J_1$]\label{lem:prop_J1}
				There exists a constant $C_1\approx 1$ (that exclusively depends on $\mathbb{T}$) 
			 	such that  the enrichment operator $J_1v_M\in \textit{WF}(\mathcal{T})$ 
			 	for any $v_M\in M(\mathcal{T})$, $\mathcal{T}\in\mathbb{T}$, satisfies
				$$			\Vert h_{\mathcal{T}}^{-{2}} (1-J_1)v_{M}\Vert_{L^2(\Omega)}
							+\vvvert (1-J_1)v_{M}\vvvert_{\mathrm{pw}}
							\le {C_1} \min_{v\in V}\vvvert v_{M}-v\vvvert_{\mathrm{pw}}.$$
						
			\end{lemma}
			\begin{proof}
				This proof consist of three steps. For \ref{item:proof_J1jump}--\ref{item:proof_jumpfine} fix one tetrahedron $K\in\mathcal{T}$. 
				\begin{enumerate}[label=\textit{Step \arabic*.},ref={Step \arabic*},wide]
				\item \label{item:proof_J1jump}\textit{Proof of $ h_K^{-4}\Vert v_M-J_1v_M\Vert_{L^2(K)}^2	
				 				\lesssim \sum_{z\in\mathcal{V}(K)}\sum_{F\in\mathcal{F}(z)}h_F\Vert [D^2 v_M]_F\times 
				 				\nu_F\Vert_{L^2(F)}^2$.}
				 The $\textit{WF}$ 
				interpolation on $K$  preserves the enrichment $J_1v_M\in \textit{WF}(K)$ 
				and quadratic polynomials like the Morley function 
				$v_M|_K\in P_2(K)\subset \textit{WF}(K)$. 
				Recall the nodal  basis functions $\varphi_{z,1}|_K$, $\varphi_{z,\kappa+1}|_K$, and $\varphi_{E,\mu}|_K$  for any 
				$z\in\mathcal{V}(K)$, $E\in\mathcal{E}(K)$, $\kappa=1,2,3$, 
				$\mu=1,2$, and the degrees of freedom of $\textit{WF}(K)$ 
				from \cref{def:CT_dof_global}. Then the $\textit{WF}$ 
				interpolation reads 
				\begin{align*}
				 (v_M-J_1v_M)|_K
				 		=&\sum_{z\in\mathcal{V}(K)}\Big((v_M-J_1v_M)(z)\varphi_{z,1}+\sum_{\kappa=1}^3 e_{\kappa}\cdot
				 		\nabla (v_M-Jv_M)(z)\varphi_{z,\kappa+1}\Big)\\
				 			&+\sum_{E\in\mathcal{E}(K)}\sum_{\mu=1}^2\nu_{E,\mu}\cdot\nabla(v_M-J_1v_M)
				 			(\textup{mid}(E))\varphi_{E,\mu}.			
				\end{align*}
				The triangle inequality and 
				the scaling of the basis functions 
				$\Vert \varphi_{z,1}\Vert_{L^2(K)}\lesssim h_K^{3/2}$, $\Vert\varphi_{z,\kappa+1}\Vert_{L^2(K)}\lesssim h_K^{5/2}$, 
				and $\Vert\varphi_{E,\mu}\Vert_{L^2(K)}\lesssim h_K^{5/2}$ in \cref{thm:Scaling_CT_fine} reveal for a constant $c_1\lesssim 1$ 
				\begin{align}
					c_1 ^{-1}\Vert v_M-J_1v_M\Vert_{L^2(K)}^2
						\le &\,h_K^{3}\sum_{z\in\mathcal{V}(K)}\Big(\vert(v_M|_K-J_1v_M)(z)\vert^2
									+h_K^2\sum_{\kappa=1}^3 \vert e_{\kappa}\cdot\nabla (v_M|_K-J_1v_M)(z)\vert^2\Big)\notag\\
				 			&+h_K^{5}\sum_{E\in\mathcal{E}(K)}\sum_{\mu=1}^{2}\vert\nu_{E,\mu}\cdot
				 			\nabla(v_M|_K-J_1v_M)(\textup{mid}(E))\vert^2.\label{eq:proof_Approx_vM-J1vM}
				\end{align}
				In \cref{def:J1}  the given local values $L_j(v_M|_T)$ of $v_M$ for every $\textit{WF}$ 
				degree of freedom
				$L_j$ for $j=1,\dots, 28$ and every tetrahedron $T\in\mathcal{T}$ are averaged (and homogeneous boundary conditions ensured).
				This in mind, suppose $z\in\mathcal{V}(K)$ and analyse the term $|v_M(z)-J_1v_M(z)|$ first. 
				There are three cases of interest that lead to the constant $C_z\lesssim 1$ in \eqref{eq:vM-J1vMz} below.  
				The vertex $z\in\mathcal{V}(K)$ could be 
				\begin{enumerate}[label=\textit{Case \arabic*.}, ref=case~\arabic*, wide]
					\item\label{case1:zinpartialOmega_noF} 
				a boundary vertex $z\in\mathcal{V}(\partial\Omega)$, but $K$ has no boundary side that contains $z$ in that
				$\mathcal{F}(K)\cap \mathcal{F}(z)\cap\mathcal{F}(\partial\Omega)=\emptyset$,
					\item\label{case2:zinpartialOmega_capF} 
				 a boundary vertex 	$z\in\mathcal{V}(\partial \Omega)\cap F$ 
				that belongs to a boundary side $F\in\mathcal{F}(\partial\Omega)\cap\mathcal{F}(K)$ of $K$, or 
					\item\label{case3:zinintOmega_capF}  an interior vertex  $z\in\mathcal{V}(\Omega)$.
				 \end{enumerate}			
				In \ref{case1:zinpartialOmega_noF} holds $J_1v_M(z)=0$ due to the homogeneous boundary conditions in $V$.  
				Then \cref{lem:fundamentalInequality}  applies and contains the optimal constant.
				Choose a face-connected subset $\{T_1,\dots,T_k\}\subset \mathcal{T}(z)$ with $2\le k\le |\mathcal{T}(z)|=:J_z\lesssim 1$, $T_k=K$, 
		 		$z\in F_j:=T_j\cap T_{j+1}\in\mathcal{F}$ for any $j=1,\dots, k-1$, and 
		 		boundary simplex $T_1$ such that $z\in F_0\in\mathcal{F}(T_1)\cap\mathcal{F}(\partial\Omega)$.
		 		Set $x_j:=v_M|_{T_j}(z)-J_1v_M(z)=v_M|_{T_j}(z)\in\mathbb{R}$ for $j=1,\dots, k$ and $x:=(x_1,\dots,x_k)\in\mathbb{R}^k$. 
		 		If $x\in\mathbb{R}^k\setminus\{0\}$,  \cref{lem:fundamentalInequality} shows 
		 				\begin{align*}
		 					2(1-\cos(\pi/(2k+1)))|x|^2\le \sum_{j=1}^k(x_{j+1}-x_j)^2+\min\{x_1^2,\dots,x_k^2\}
		 																	\le \sum_{j=1}^k(x_{j+1}-x_j)^2+x_1^2. 
		 				\end{align*}
		 		The jump $[v]_F\in L^1(F)$ of a piecewise Lipschitz continuous function $v$ across 
				$F\in\mathcal{F}$ reads $[v]_F:=v|_{T_+}-v|_{T_-}$ for an interior side $F=\partial T_+\cap \partial T_-\in\mathcal{F}(\Omega)$ 
				that belongs to the neighbouring simplicies $T_+,T_-\in \mathcal{T}$, which are labelled such that 
				the unit outward normal along the boundary of $T_+$ (resp. $T_-$) satisfies $\nu_{T_+}|_F=\nu_F$ (resp. $\nu_{T_-}|_F=-\nu_F$),
				while $[v]_F:=v$ along $F\in \mathcal{F}(\partial \Omega)$. 
				Hence the side-connectivity of $\{T_1,\dots,T_k\}\subset \mathcal{T}(z)$ implies, $\text{for any }j=1,\dots,k-1$, that 
				\begin{align*}
						|x_{j+1}-x_j|=\big|v_M|_{T_j}(z)-v_M|_{T_{j+1}}(z)\big|=\big|[v_M]_{F_j}(z)]\big| \text{ and }
						|x_1|=\big|v_M|_{T_1}(z)\big|=\big|[v_M]_{F_0}(z)\big|.
				\end{align*}						 		
		 		We conclude (with the Euclidean norm $|\bullet|$ of $x\in\mathbb{R}^k$) in \ref{case1:zinpartialOmega_noF}, 
		 				\begin{align}
		 					\big|v_M|_{K}(z)-J_1v_M(z)\big|^2	&\le |x|^2\le \frac{1}{2(1-\cos(\pi/(2k+1))}\sum_{j=0}^k|[v_M]_{F_j}(z)|^2.
		 						\label{eq:vM-J1vM_case1}  
		 				\end{align}
		 		In \ref{case2:zinpartialOmega_capF}, for the boundary node $z\in\mathcal{V}(\partial \Omega)\cap F$ 
				that belongs to a boundary side $F\in\mathcal{F}(\partial\Omega)\cap\mathcal{F}(K)$ of $K$, 
				holds $J_1v_M(z)=0$ and the jump definition with $[v]_F:=v$ 
		 		along any boundary side, in particular along $F\in \mathcal{F}(\partial \Omega)\cap\mathcal{F}(K)$, shows 
		 			\begin{align}
		 					\big|v_M|_{K}(z)-J_1v_M(z)\big|=\big|v_M|_{K}(z)\big|=\big|[v_M]_{F}(z)\big|. \label{eq:vM-J1vM_case2}
		 			\end{align}
		 		In \ref{case3:zinintOmega_capF}, for the interior vertex $z\in\mathcal{V}(\Omega)$, holds
				$J_1v_M(z)=J_z^{-1}\sum_{T\in\mathcal{T}(z)}v_M|_T(z)$ with the abbreviation $J_z:=|\mathcal{T}(z)|\lesssim 1$. 
				This is the standard averaging and an analogue argumentation with explicit constant is provided Step 3 in the proof 
				of Theorem 6.1 in \cite {CP18}. For convenience we repeat the main arguments in the present notation. 
				Choose an enumeration $\{T_1,\dots, T_{J_z}\}$ of $\mathcal{T}(z)$ such that the values 
				$x_j:=v_M|_{T_j}(z)-J_1v_M(z)\in\mathbb{R}$ 
				for $j=1,\dots, J_z$ are ordered in the sense that $x_1\le x_2\le \dots\le x_{J_z}$. 
				The sum $\sum_{j=1} ^{J_z} x_j=0$ vanishes by definition of $J_1v_M(z)$.  
				Since 	$\big|v_M|_{K}(z)-J_1v_M(z)\big|\le\max_{T\in\mathcal{T}(z)}\big|v_M|_{T}-J_1v_M(z)\big|=\max_{j=1,\dots,J_z}|x_j|$, 
				 Lemma~C in \cite{CP18} implies 
				(in the same abstract notation as \cref{lem:fundamentalInequality} above) that 
				\begin{align*}
					\big|v_M|_{K}(z)-J_1v_M(z)\big|^2\le\frac{(J_z-1)(2J_z-1)}{6J_z} {\sum_{j=1}^{J_z-1}|x_{j+1}-x_j|^2}. 
				\end{align*}
				It remains a reordering to arrange for jump contributions on the right-hand side. 
				To this end let 
				$\mathcal{J}:=\big\{\{\alpha,\beta\}:\,T_\alpha,T_\beta\in\mathcal{T}(z)\text{ and } \partial T_\alpha\cap \partial T_\beta\in\mathcal{F}(z)\}$ 
				denote unordered index pairs  of all tetrahedra in $\mathcal{T}(z)$ 
				 which share a side in $\mathcal{F}(z):=\{F\in\mathcal{F}|\, z\in\mathcal{V}(z)\}$. 
				 Since  $\mathcal{T}(z)$ is side-connected,  $\mathcal{J}$ is connected, in the sense that 
				for all $\alpha,\beta\in\{1,\dots,J_z\}$ and $\alpha\not = \beta$ there are $k\in\mathbb{N}$ pairs 
				$\{\alpha_1,\alpha_2\},\,\{\alpha_2,\alpha_3\},\,\dots,\,\{\alpha_{k},\alpha_{k+1}\}\in\mathcal{J}$ with 
				$\alpha_1=\alpha$ and $\alpha_{k+1}=\beta$. 
				Therefore some arguments in the language of graph theory in Lemma~B in \cite{CP18} imply that 
				\begin{align*}
					 			\sum_{j=1}^{J_z-1}|x_{j+1}-x_j|^2
								\le\sum_{\{\alpha,\beta\}\in\mathcal{J}} |x_{\alpha}-x_\beta|^2
								=\sum_{\{\alpha,\beta\}\in\mathcal{J}}\big|(v_M|_{T_{\alpha}})(z)-(v_M|_{T_\beta})(z)\big|^2. 			
				\end{align*}
				The combination of the last two displayed estimates with the jump definition shows in \ref{case3:zinintOmega_capF} that 
		 		\begin{align}
		 			\big|v_M|_K(z)-J_1v_M(z)\big|\le \frac{(J_z-1)(2J_z-1)}{6J_z}
		 									\sum_{F\in\mathcal{F}(z)}|[v_M]_F(z)|^2. \label{eq:vM-J1vM_case3}
		 		\end{align}
				
		 		Abbreviate 
				$$
					C_z:=\max\bigg\{1,\frac{(J_z-1)(2J_z-1)}{6J_z},\frac{1}{2(1-\cos(\frac{\pi}{2J_z+1}))}\bigg\}\lesssim 1.
				$$
				The combination of \eqref{eq:vM-J1vM_case1}--\eqref{eq:vM-J1vM_case3} implies for any vertex $z\in\mathcal{V}(K)$ of 
				$K\in\mathcal{T}$ that 
		 		\begin{align}
		 			\big|v_M(z)|_K-J_1v_M(z)\big|\le 
		 									C_z\sum_{F\in\mathcal{F}(z)}|[v_M]_F(z)|^2. \label{eq:vM-J1vMz}
		 		\end{align}
		 		
		 		Since $J_1$ is defined by nodal averaging for all degrees of freedom $L_1,\dots, L_{28}$ of $\textit{WF}$ 
		 		(with homogeneous boundary conditions), 
				the argumentation for $|v_M|_K(z)-J_1v_M(z)|$ with \ref{case1:zinpartialOmega_noF}--\ref{case3:zinintOmega_capF} 
				applies verbatim for  $\vert e_{\kappa}\cdot\nabla (v_M|_K-J_1v_M)(z)\vert$ and 
				$\vert\nu_{E,\mu}\cdot\nabla(v_M|_K-J_1v_M)(\textup{mid}(E))\vert$. 
				Let $\mathcal{N}:=\mathcal{V}\cup\{\textup{mid}(E)|\,E\in\mathcal{E}\}$ denote the set of all vertices and edge midpoints.
				Let $\mathcal{T}(z):=\{T\in\mathcal{T}|\, z\in\partial T\}$ denote  the 
				set of all tetrahedra with node $z\in\mathcal{N}$ and 
				$\mathcal{F}(z):=\{F\in\mathcal{F}|\, z\in\partial F\}$ set of all faces with node $z\in\mathcal{N}$. 
				In abstract notation 
		 		let $J_A$ denote a general averaging operator; 
		 		for any $k\in\mathbb{N}_0$, a piecewise polynomial $v\in P_k(\mathcal{T})$ of degree at most $k$, 
		 		the tetrahedron $K\in\mathcal{T}$, and 
		 		the node $z\in\mathcal{N}(K)$  set 
		 		$J_A v(z):= \vert\mathcal{T}(z)\vert^{-1}\sum_{T\in\mathcal{T}(z)}v|_T(z)$ for each interior node 
		 		$z\in \mathcal{N}(\Omega)$ and 
		 		$J_A v(z):= 0$ else.
		 		Given this notation, 
		 		$(D^\alpha J_1v_M)(z)=J_A(D^\alpha v_M)(z)$ holds for $z\in\mathcal{V}$, $\alpha\in\mathbb{N}_0^3$ with 
				$0\le|\alpha|\le 1$, as well as, $\tau_E\times\nabla J_1v_M(\textup{mid}(E))=J_A(\tau_E\times\nabla v_M)(\textup{mid}(E))$ for 
				$E\in\mathcal{E}$. In the argumentation above we could replace $J_1$ by $J_A$ and $z\in\mathcal{V}$ by 
				$z\in\mathcal{N}$ to obtain \eqref{eq:vM-J1vMz}  
				in the general setting. Hence, 
				\begin{align*}
						\big\vert e_{\kappa}\cdot\nabla (v_M|_K-J_1v_M)(z)\big\vert&\le C_z\sum_{F\in\mathcal{F}(z)}|[e_\kappa\cdot\nabla v_M]_F(z)|^2,
						\\
							\big\vert\nu_{E,\mu}\cdot\nabla(v_M|_K-J_1v_M)(\textup{mid}(E))\big\vert
							&\le C_E\sum_{F\in\mathcal{F}(E)}\vert [\nu_{E,\mu}\cdot\nabla v_M]_F(\textup{mid}(E))\vert^2
				\end{align*}								
					with 
				$$C_E:=\max\bigg\{1,\frac{(J_E-1)(2J_E-1)}{6J_E},\frac{1}{2(1-\cos(\frac{\pi}{2J_E+1}))}\bigg\}\lesssim 1$$ 
				for any edge $E\in\mathcal{E}$ with set of adjacent faces 
				$\mathcal{F}(E):=\{F\in\mathcal{F}:E\in\mathcal{E}(F)\}=\mathcal{F}(\textup{mid}(E))$ 
				of cardinality $J_E=|\mathcal{F}(E)|\lesssim 1$. 
				
				Hence \eqref{eq:proof_Approx_vM-J1vM} simplifies to 
								 \begin{align}
										c_1^{-1} \Vert v_M-J_1v_M\Vert_{L^2(K)}^2
												\le \ &h_K^{3}\sum_{z\in\mathcal{V}(K)}C_z\sum_{F\in\mathcal{F}(z)}
												\Big(|[v_M]_F(z)|^2
									+ h_K^2\vert [\nabla v_M]_F(z)\vert^2\Big)\notag\\
										&+h_K^{5}\sum_{E\in\mathcal{E}(K)}C_E\sum_{F\in\mathcal{F}(E)}
				 					\vert [\nabla v_M]_F(\textup{mid}(E))\vert^2\label{eq:proof_Approx_||vm-J1vM||}
								 \end{align}	
				(The shape regularity of $\mathcal{T}\in\mathbb{T}$ guarantees that the number of tetrahedra $J_z$ and $J_E$ in the patches, i.e., 
				$C_z$ and $C_E$, 
				are uniformly bounded for any 
				$z\in\mathcal{V}$ and $E\in\mathcal{E}$.)
				Let $S$ denote either the set of vertices $S=\mathcal{V}$ or the set of edge midpoints 
				$S=\{\textup{mid}(E)|\,E\in\mathcal{E}\}=:\mathcal{M}$.
				For a fixed polynomial degree $k\in\mathbb{N}_0$ and any $z\in\mathcal{N}=\mathcal{V}\cup\mathcal{M}$,
				there exists a constant $C(k,S)>0$, such that the following inverse estimate holds 
			 		\begin{align}
			 			|F|\,\big|[v]_F(z)\big|^2\le C(k,S) \Vert [v]_F\Vert_{L^2(F)}^2\quad 
			 			\text{ for any }v\in P_k(\mathcal{T})\text{ and }F\in\mathcal{F}(z). 	
			 			\label{eq:proof_Approx_jump}
			 		\end{align}
			 	In the cases of interest the optimal constant $C(k,S)>0$ can be computed as in \cite[Lem.~D]{CP18}. 
			 	The inverse mass matrix for the basis functions dual to the Lagrange interpolation 
			 	points and a Cauchy-Schwarz inequality prove $C(1,\mathcal{V})=n^2=9$ and $C(2,\mathcal{V})=36$ for any vertex 
			 	$z\in\mathcal{V}$ as well as 
			 	$C(1,\mathcal{M})=n(n-1)/2=3$ and $C(2,\mathcal{M})=39/4$  for any edge midpoint $z=\textup{mid}(E)\in\mathcal{M}$.
				Moreover, an inverse estimate for  $[v_M]_F\in P_2(F)$  ensures 
				$c_{\mathrm{inv}}^{-1}\Vert [v_M]_F\Vert_{L^2(F)}\le h_F \Vert [\nabla v_M]_F\Vert_{L^2(F)}\le h_K \Vert [\nabla v_M]_F\Vert_{L^2(F)} $. 
				In combination with \eqref{eq:proof_Approx_jump} this recasts \eqref{eq:proof_Approx_||vm-J1vM||} as 
			\begin{align*}
				c_1^{-1} \Vert v_M-J_1v_M\Vert_{L^2(K)}^2
												\le \  &h_K^{5}\sum_{z\in\mathcal{V}(K)}C_z\sum_{F\in\mathcal{F}(z)}
												\Big( c_{\mathrm{inv}}^2 36 \Vert [\nabla v_M]_F\Vert_{L^2(F)}^2/|F|
									+ 9 \Vert [\nabla v_M]_F\Vert_{L^2(F)}^2/|F|\Big)\notag\\
										&+h_K^{5}\sum_{E\in\mathcal{E}(K)}C_E\sum_{F\in\mathcal{F}(E)}3 \Vert [\nabla v_M]_F\Vert^2_{L^2(F)}/|F|
			\end{align*}
				Since for any edge $E=\textup{conv}\{P_j,P_k\} \in\mathcal{E}$ with vertices $P_j,P_k\in \mathcal{V}(E)$, it  
				holds  $\mathcal{F}(E)\subset \mathcal{F}(P_j)\cup\mathcal{F}(P_k)$, this simplifies to 			
							 \begin{align*}	
				 				h_K^{-5}\Vert v_M-J_1v_M\Vert_{L^2(K)}^2	
				 				\lesssim \sum_{z\in\mathcal{V}(K)}\sum_{F\in\mathcal{F}(z)}\Vert [\nabla v_M]_F\Vert_{L^2(F)}^2/|F|.
							\end{align*}					
				Recall that the integral mean of the gradient $\intmean_F[\nabla v_M]_F\,\textup{d}\sigma$ 
				vanishes for any side $F\in\mathcal{F}$ from \cite[Lem.~4]{MX06} as utilized in \cite[Thm.~2.1.b]{CP_Part1}.  
				Hence the Poincar\'e inequality with Payne-Weinberger constant \cite{PW60,Bebendorf2003} shows
				\begin{align*}
					\Vert [\nabla v_M]_F\Vert_{L^2(F)}\le h_F\pi^{-1}\Vert [D^2 v_M]_F\times \nu_F\Vert_{L^2(F)}
				\end{align*}
				with the tangential components $[D^2v_M]_F\times \nu_F$ of the jump $[D^2v_M]_F$ along any edge $F\in\mathcal{F}$ and  
					the (piecewise) Hessian $D^2$. This concludes the first step, since 
					$h_F\approx h_K\approx |F|^{1/2}$ holds for any $K\in\mathcal{T}$ and $F\in\mathcal{F}(K)$.
				
				\item \label{item:proof_jumpfine}\textit{Proof of $h_F \Vert [D^2 v_M]_F\times \nu_F\Vert_{L^2(F)}^2
							\lesssim\min_{v\in H^2(\omega(F))}\Vert  D^2_{\mathrm{pw}}(v_M-v)\Vert_{L^2(\omega(F))}^2.$}
				This step is comparable to the $2$D  analysis in \cite[Prop.~2.3]{Gal15} but demands the observations on the $\textup{Curl}$ operator 
				in $3$D from \cite[Thm~3.2]{CBJ01}. 
				For $\boldsymbol{\Psi}=(\Psi_1,\Psi_2,\Psi_2)^\top\in H^1(\Omega)^{3\times 3}$ the $\textup{Curl}$ is applied row-wise, 
				\begin{align*}
					\textup{Curl}(\boldsymbol{\Psi}):=(\nabla \times\Psi_1,\nabla\times \Psi_2,\nabla \times \Psi_3)^\top. 
				\end{align*}
				First, suppose $F\in\mathcal{F}(\Omega)=\partial T_+\cap \partial T_-$ is an interior side 
				and for any vertex $z\in\mathcal{V}(T)$ let $\phi_z\in S^1(\mathcal{F}):=P_1(\mathcal{T})\cap C(\Omega)$ denote the hat-function with 
				$\phi_z(P)=\delta_{zP}$ for any $P\in\mathcal{V}$. 
				Define the face bubble $b_F:=60\prod_{z\in \mathcal{V}(F)} \phi_z \in H^1_0(\omega(F))$ that vanishes on the boundary of 
				the face patch $\omega(F):=\textup{int}(T_+\cup T_-)$. Then $\int _F b_F\,\textup{d}\sigma=|F|$ and 
				$\Vert b_F\Vert_{L^\infty(\omega(E))}=20/3$. Given 
				$[D^2 v_M]_F\in \mathbb{R}^{3\times 3}$, 
				set 
				$\Psi_F:= b_F([D^2 v_M]_F\times \nu_F)\in H^1_0(\omega(F);\mathbb{R}^{3\times 3})$.  
				Since $[D^2 v_M]_F\in \mathbb{R}^{3\times 3}$ is constant,  
				\begin{align*}
					 \Vert [D^2 v_M]_F\times \nu_F\Vert_{L^2(F)}^2=\Vert b_F^{1/2} [D^2 v_M]_F\times \nu_F\Vert_{L^2(F)}^2. 
				\end{align*}
				For any $v\in H^2(\omega(E))$ the tangential jump $[D^2 v]_F\times \nu_F$ vanishes. Hence an integration by parts 
				shows
				\begin{align*}
					\Vert b_F^{1/2} [D^2 v_M]_F\times \nu_F\Vert_{L^2(F)}^2
						&=\int_F  \big([D^2 v_M]_F\times \nu_F\big) \Psi\,\textup{d}\sigma
						 =\sum_{T\in \omega(F)}\int_{\partial T} \big((D^2 v_M - v)\times \nu_F\big):\Psi\,\textup{d}\sigma	
						\\& =\big(\textup{Curl}\Psi_F,D^2_{\mathrm{pw}}(v_M-v)\big)_{L^2(\omega(F))}.	
				\end{align*}
				The Cauchy-Schwarz inequality, an inverse estimate, and the boundedness of $\Vert b_F\Vert_{L^\infty(\omega(E))}$ 
				prove that the last term is bounded by
				\begin{align*}
					\Vert \textup{Curl}\Psi_F\Vert_{L^2(\omega(F))} \Vert D^2_{\mathrm{pw}}(v_M-v)\Vert_{L^2(\omega(F))}
						&\lesssim h_F^{-1/2} \Vert [D^2 v_M]_F\times \nu_F\Vert_{L^2(F)}\Vert D^2_{\mathrm{pw}}
						(v_M-v)\Vert_{L^2(\omega(F))}.
				\end{align*}
				This implies for any interior side $F\in\mathcal{F}(\Omega)$ that 
					\begin{align}
						h_F \Vert [D^2 v_M]_F\times \nu_F\Vert_{L^2(F)}^2\lesssim\min_{v\in H^2(\omega(F))}\Vert  D^2_{\mathrm{pw}}
						(v_M-v)\Vert_{L^2(\omega(F))}^2.  \label{eq:jumpJ1_interioredge}
					\end{align}
				Second, suppose $F\in\mathcal{F}(\partial \Omega)$ is a boundary face with face patch $\omega(F)=T_+$. 
				The face bubble  $b_F:=60\prod_{z\in \mathcal{V}(F)} \phi_z \in H^1_0(\omega(F))$ vanishes 
				on all faces $G\in\mathcal{F}(T_+)\setminus\{F\}$ of $T_+$ besides $F$. Hence the argumentation above holds verbatim and shows 
				\eqref{eq:jumpJ1_interioredge} for $F\in\mathcal{F}(\partial\Omega)$. 
				
				\item \textit{Conclusion of the proof.} 
					The combination of \ref{item:proof_J1jump}--\ref{item:proof_jumpfine} implies for any $K\in\mathcal{T}$ 
					\begin{align*}
						 h_K^{-4}\Vert v_M-J_1v_M\Vert_{L^2(K)}^2
						 	\lesssim \sum_{F\in\mathcal{F}(z)}\min_{v\in H^2(\omega(F))}\Vert  D^2_{\mathrm{pw}}(v_M-v)\Vert_{L^2(\omega(F))}^2.
						  \label{eq:local_bestapprox}
					\end{align*}
					The sum over all tetrahedra $K\in\mathcal{T}$ with  the finite overlap of the face patches $(\omega(F): \, F\in\mathcal{F})$ 
					concludes the proof of the approximation property 
					\begin{align*}
					 	\Vert h_{\mathcal{T}}^{-2}(1-J_1)v_M\Vert_{L^2(\Omega)}\lesssim \min_{v\in H^2_0(\Omega)}
					 	\vvvert v-v_M\vvvert_{\mathrm{pw}}.
					\end{align*}
			  		The combination with the inverse estimate
			  		 $\vvvert (1-J_1)v_M\vvvert_{\mathrm{pw}}\lesssim \Vert h_{\mathcal{T}}^{-2}(1-J_1)v_M\Vert_{L^2(\Omega)}$ 
			  		for the piecewise polynomials in $\textit{WF}(\mathcal{T})\subset P_3(\widehat{\mathcal{T}})$ 
			  		concludes the proof of \cref{lem:prop_J1}. 
			  		\end{enumerate}
			\end{proof}
	\subsection{Design of $J_2$}	\label{sec:Morley_J2} 		
	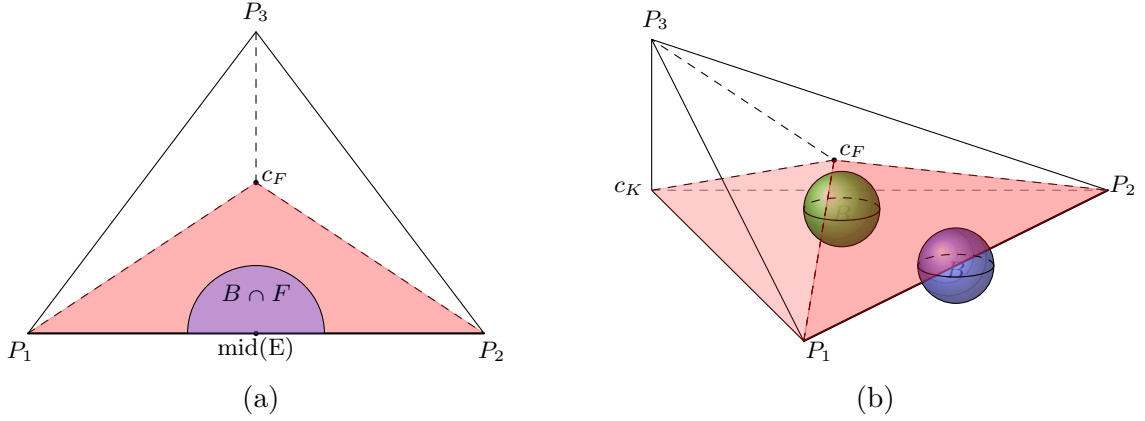
\begin{figure}
		\begin{center}
		\begin{tabular}{ccc}
			{\begin{tikzpicture}[x=30mm,y=20mm]
					  					\draw(0,0)node[shape=coordinate,label={[xshift=-1mm, yshift=-5mm] $P_1$}] (P1) {P1};
					  					\draw(1,1)node[shape=coordinate,label={[xshift=2.5mm, yshift=-1.5mm] $c_F$}] (P3) {P3};
					  					\draw(2,0)node[shape=coordinate,label={[xshift=1mm, yshift=-5mm] $P_2$ }] (P2) {P2};
					  					\draw(1,2)node[shape=coordinate,label={[xshift=0mm, yshift=0mm] $P_3$}] (P4) {P4};
					  					\draw(1,0) node[shape=coordinate,label={[xshift=0mm, yshift=-5mm] $\textup{mid(E)}$ }] (mE) {mE};
					  					\fill (P3) circle (1pt);
					  					\fill (mE) circle (1pt);
					  					\draw[dashed] (P1)--(P3)--(P2);
					  					\draw[dashed] (P4)--(P3);
					  					\filldraw[red, opacity=0.3](P1)--(P2)--(P3)--cycle;
					  					\draw (P1)--(P4)--(P2);
					  					\draw[thick] (P1)--(P2);
					  	\pic[fill=blue!60,fill opacity=0.4,angle radius=9mm,"\textcolor{blue!10!black}
					  	{\footnotesize $B\cap F$}" opacity=1] 			
					  					{angle=P2--mE--P1};		
					  					\draw (0.7,0) arc(180:0:0.9cm);
					\end{tikzpicture} } &\phantom{x}&
					\begin{tikzpicture}[x=1cm,y=1cm]
		  					\draw(6,0)node[shape=coordinate,label={[xshift=2mm, yshift=-4mm] $P_1$}] (Q2) {Q2};
		  					\draw(4,4)node[shape=coordinate,label={[xshift=0.5mm, yshift=0mm] $P_3$}] (Q3) {Q3};
		  					\draw(10,2)node[shape=coordinate,label={[xshift=2mm, yshift=-2mm] $P_2$}] (Q4) {Q4};
		  					\draw(8,1)node[shape=coordinate, label={[xshift=0mm, yshift=-3mm]\textcolor{black!90!blue}{${B}$}}] (ME) {ME};
		  					\draw(4,2)node[shape=coordinate,label={[xshift=-3mm, yshift=-2mm] $c_K$}] (cK) {cK};
		  					\draw (Q2)--(Q3)--(Q4)--cycle;
		  					\draw (Q2)--(cK)--(Q3);
		  					\draw[dashed,gray] (cK)--(Q4);
		  					\draw[thick](Q2)--(Q4);
		  					\shade[ball color = blue!60, opacity = 0.4] (ME) circle (0.5cm);
  							\draw (ME) circle (0.5cm);
  							\draw(8.5,1)node[shape=coordinate] (E) {E};
  							\draw[dashed] (E) arc(0:180:0.5cm and 0.15cm);
  							\draw (E) arc(360:180:0.5cm and 0.15cm);
 							\draw(6.5,1.75)node[shape=coordinate, label={[xshift=0mm, yshift=-3mm]\textcolor{black!90!green}{$\widehat{B}$}}] (cE) {cE};
 							\shade[ball color = green!60, opacity = 0.6] (cE) circle (0.5cm);
  							\draw (cE) circle (0.5cm);	 
  							\draw(7,1.75)node[shape=coordinate] (cE2) {cE2};
  							\draw[dashed] (cE2) arc(0:180:0.5cm and 0.15cm);
  							\draw (cE2) arc(360:180:0.5cm and 0.15cm);
		  					\draw(6.4,2.4)node[shape=coordinate,label={[xshift=2.5mm, yshift=-1mm] $c_F$}] (c1) {c1};
		  					\fill (c1) circle (1pt);
		  					\draw[dashed] (Q2)--(c1)--(Q4);
		  					\draw[dashed] (Q2)--(c1)--(Q3);
		  					\draw[dashed] (c1)--(cK);  
		  					
		  					\filldraw[red, opacity=0.3](Q2)--(c1)--(Q4)--cycle;
		  					\filldraw[red, opacity=0.2](cK)--(c1)--(Q2)--cycle;		
					\end{tikzpicture}\\ (a)& & (b)
					\end{tabular}
					\caption{Subtetrahedron $T(4)=\textup{conv}(c_K,F_4)$ of $K\in\mathcal{T}$ with face 
								$F=F_4=\textup{conv}\{P_1,P_2,P_3\}\in\mathcal{F}(K)$ and 
								edge $E=\textup{conv}\{P_1,P_2\}\in\mathcal{E}(K)$ and dashed $\textit{WF}$ 
								 triangulation $\widehat{\mathcal{T}}(K)$ with $T=\textup{conv}(c_K,c_F,E)\in\widehat{\mathcal{T}}$ 
								marked in red, the support $\textup{supp}(\xi_E)=B$ of $\xi_E$ from \eqref{eq:def_xiE} in blue and 
								the disjoint ball $\widehat{B}\subset T$ in green. }				
								\label{fig:Appendix_inverseestimate}	
		\end{center}
		\end{figure}
					This step corrects the integral mean over each internal edge $E\in\mathcal{E}$ such that 
					$\intmean_E J_2 v_M\,\textup{d}s=\intmean_E v_M\,\textup{d}s.$
				 	For any edge $E\in\mathcal{E}(\Omega)$ define the ball $B:=B(\textup{mid}(E),R_E)\subset\widehat{\omega}(E)$ of radius 
				 	$R_E>0$ with midpoint $\textup{mid}(E)$. The radius $R_E$ is chosen maximal such that 
					\begin{enumerate}[label=(\roman*)]
					\item 
						$B$ belongs to the edge patch
				 		$\widehat{\omega}(E):=\textup{int}\big(\bigcup_{T\in \widehat{\mathcal{T}}(E)} T\big)$ 
				 		of $E$ in the $\textit{WF}$ 
				 		triangulation $\widehat{\mathcal{T}}$ and 
				 	\item 
				 		in any adjacent $\textit{WF}$-subtetrahedra 
				 		${T}\in\widehat{\mathcal{T}}(E)$ exists a point $\widehat{c}_{{T}}\in {T} $ such that 
				 		$\widehat{B}:=B(\widehat{c}_{{T}},R_E)\subset \mathrm{int}({T})$ is a 
				 		disjunct ball $B\cap\widehat{B}=\emptyset$ of the same size as $B$ and lies in the interior of ${T}$. 
					\end{enumerate}	
					\cref{fig:Appendix_inverseestimate}.b illustrates the disjoint balls $B$ and $\widehat{B}$. 								 	
				 	The shape regularity of $\mathcal{T}\in \mathbb{T}$ and $\widehat{\mathcal{T}}$ in \cref{thm:shapeRegularity3D} 
				 	guarantees  $h_T\approx R_E\approx h_E$. 
				 	Set for any $E\in\mathcal{E}(\Omega)$
				 	\begin{align}
				 		\xi_E(x):=\begin{cases}
				 				\frac{|E|}{R_E}\Big(1-3\frac{|y|^2}{R_E^2}+2\frac{|y|^3}{R_E^3}\Big) \quad &\text{ for }  x\in \overline{B} 
				 				\text{ and }y:=x-\textup{mid}(E),\\ 
				 				0\quad &\text{ for }x\not \in \overline{B}.
				 		\end{cases}\label{eq:def_xiE}
				 	\end{align} 
				 	By construction holds $\xi_E\in C^1(\mathbb{R}^3)\cap  H^2_0(B)$ 
				 	with $\intmean_E \xi_D\,\textup{d}s=\delta_{ED}$ for any $E,D\in\mathcal{E}$ 
				   and the  support $\textup{supp}(\xi_E)=B\subset \widehat{\omega}(E)$. 
				   \begin{definition}[$J_2$]
					For any $v_M\in M(\mathcal{T})$ set 
								\begin{align}
									J_2(v_M):= J_1 v_M+
										\sum_{E\in\mathcal{E}(\Omega)}\Big(\intmean_E (v_M-J_1v_M)\,\textup{d}s\Big)\xi_E\in V.
										\label{eq:def_J2_M_1}
								\end{align}
					\end{definition}
				
			\begin{lemma}[properties of $J_2$]\label{lem:prop_J2}
				There exists a constant $C_2\approx 1$ (that exclusively depends on $\mathbb{T}$) 
			 	such that  the companion $J_2v_M\in  V$  
			 	for any $v_M\in M(\mathcal{T})$ satisfies
			 	 \begin{enumerate}[label=(\alph*), ref=\alph* ,wide]
			 		\item\label{item:App_J2_rightInverse} 
			 				 $L_E (J_2v_M)=L_E(v_M)$ for any $E\in\mathcal{E}(\Omega),$
					\item\label{item:App_J2_approx} 
						$			\Vert h_{\mathcal{T}}^{-{2}} (1-J_2)v_{M}\Vert_{L^2(\Omega)}
							+\vvvert (1-J_2)v_{M}\vvvert_{\mathrm{pw}}
							\le {C_2} \min_{v\in V}\vvvert v_{M}-v\vvvert_{\mathrm{pw}}.$
			 	\end{enumerate}
			\end{lemma}
				\begin{proofof}\textit{(\ref{item:App_J2_rightInverse}). }			
					This holds by construction, since for any $v_M\in M(\mathcal{T})$ and $E\in\mathcal{E}(\Omega)$ 
					\begin{align*}
						L_E(J_2v_M)&=\intmean_E J_2 v_M\,\textup{d}s
									=\intmean_EJ_1 v_M\,\textup{d}s+\sum_{D\in\mathcal{E}(\Omega)}\Bigg(\intmean_D(v_M-J_1v_M)\,\textup{d}s\Bigg)
									\intmean_D\xi_E\textup{d}s
									\\&=	\intmean_EJ_1 v_M\,\textup{d}s+\intmean_E(v_M-J_1v_M)\,\textup{d}s =\intmean_E v_M\,\textup{d}s=L_E(v_M).
					\end{align*}					 
				\end{proofof}
				\begin{proofof}\textit{(\ref{item:App_J2_approx}). }			
					Recall the discrete trace inequality 
					$\Vert v\Vert_{L^2(F)}^2\le {|F|}/{|T|}\Vert v\Vert_{L^2(T)}(\Vert v\Vert_{L^2(T)} 
								+{2h_T}/{m}\Vert \nabla v\Vert_{L^2(T)})$ 
					 for any $v\in H^1(T)$ and any $m$-simplex $T$ with $(m-1)$-subsimplex $F$, $m\ge 2$. 
					(This is a direct consequence from the discrete trace identity \cite[Eq. (2.5)]{CP_Part1} (see \cite{CGR12,CH17}).)
					 For any $g\in \textit{WF}(\mathcal{T})$ 
					 and any edge $E\in\mathcal{E}$ of a side $F\in\mathcal{F}(E)$,
					  the  trace inequality and a Cauchy-Schwarz inequality show 
					 \begin{align*}
					 	\intmean_E g\,\textup{d}s
					 		&\le |E|^{-1/2}\Vert g\Vert_{L^2(E)}
					 		\le |E|^{-1/2} \sqrt{\frac{|E|}{|F|}\Vert g\Vert_{L^2(F)}(\Vert g\Vert_{L^2(F)}+{h_F}\Vert \nabla g\Vert_{L^2(F)})}\\
					 		&\le |F|^{-1/2} \sqrt{1+c_{\mathrm{inv}}}\Vert g\Vert_{L^2(F)}
					 \end{align*} with an inverse estimate $\Vert \nabla g\Vert_{L^2(F)}\le c_{\mathrm{inv}} h_F^{-1}\Vert g\Vert_{L^2(F)}$ in the last 
					 step. 
					The trace inequality for any side $F\in\mathcal{F}$ of a tetrahedron $K\in\mathcal{T}(F)$ reveals, for any 
					$g\in\textit{WF}(\mathcal{T})$, 
					that
					\begin{align}
						 \frac{|K|}{|F|}\Vert g\Vert_{L^2(F)}^2\le\Vert g\Vert_{L^2(K)}(\Vert v\Vert_{L^2(K)} 
								+{2h_K}/{3}\Vert \nabla v\Vert_{L^2(K)}) \le (1+3C_{\textup{inv}}/2)\Vert v\Vert_{L^2(K)} \label{eq:F_traceIneq}			
					\end{align}	with an inverse estimate $\Vert \nabla g\Vert_{L^2(K)}\le C_{\mathrm{inv}} h_K^{-1}\Vert g\Vert_{L^2(K)}$ in the last 
					step. The combination of the two displayed inequalities shows, for $g=v_M-J_1v_M\in  \textit{WF}(\mathcal{T})$ 
					and any edge $E\in\mathcal{E}$ of a tetrahedron $K\in\mathcal{T}(E)$, that 
					\begin{align*}
						\intmean_E (v_M-J_1v_M)\,\textup{d}s\le \sqrt{1+c_{\mathrm{inv}}}\sqrt{1+3C_{\mathrm{inv}}/2}|K|^{-1/2}\Vert v_M-J_1v_M\Vert_{L^2(K)}.
					\end{align*}
					Since straightforward computations reveal  
					$\Vert\xi_E\Vert_{L^2(K)}=|E|\sqrt{38\pi R_E/315}\approx h_K^{3/2}$ for any $E\in\mathcal{E}(K)$ and $\xi_E$ in \eqref{eq:def_xiE}, 
					the last displayed estimate and a triangle inequality prove 
						\begin{align*}
							\Vert v_M&-J_2v_M\Vert_{L^2(K)}-\Vert v_M-J_1 v_M\Vert_{L^2(K)}
									\le \sum_{E\in\mathcal{E}(K)\cap\mathcal{E}(\Omega)}
										\Big|\intmean_E (v_M-J_1v_M)\,\textup{d}s\Big|\Vert\xi_E\Vert_{L^2(K)}
									\\&\le 6\frac{|E|\sqrt{38\pi R_E}}{\sqrt{315}}\frac{\sqrt{1+c_{\mathrm{inv}}}\sqrt{1+3C_{\mathrm{inv}}/2}}{|K|^{1/2}} 
											\Vert v_M-J_1v_M\Vert_{L^2(K)}
									\lesssim \Vert v_M-J_1v_M\Vert_{L^2(K)}
						\end{align*} for any $K\in\mathcal{T}$ with the shape-regularity $h_K\approx |E|\approx R_E\approx|K|^{1/3}$  in the last step. 
					The  inverse estimate 
					$\vert v_M-J_2v_M\vert_{H^2(K)}\le\Vert h_K^{-2} (v_M-J_2v_M)\Vert_{L^2(K)}$, 
					a summation over all $K\in\mathcal{T}$, and \cref{lem:prop_J1} conclude the proof. 	
				\end{proofof}
			The following \cref{lem:inverseEstimate} provides more insight into the inverse estimates for functions of the type $\xi_E$ in  
			\eqref{eq:def_xiE}. 
		In particular standard inverse estimates hold for the operator $J_M\equiv J_4$ in \cref{def:JM} below. 
	\begin{lemma}[inverse estimates]\label{lem:inverseEstimate}
	  There exists a constant $C\approx1$ (that exclusively depends on $\mathbb{T}$), 
	  such that for any tetrahedron $G\in\mathcal{T}\in\mathbb{T}$ 
	  or triangle $G\in\mathcal{F}$,
	  \begin{align*}
	  	\vert v_M-J_Mv_M|_{H^1(G)}\le C\textup{diam}(G)^{-1}\Vert v_M-J_Mv_M\Vert_{L^2(G)}.
	  \end{align*}
	  Notice that the $H^1$-seminorm concerns the tangential gradient $\nu_F\times \nabla v_M$ for a face $G=F\in\mathcal{F}$ 
	  with unit normal $\nu_F$ of fixed orientation and the full gradient $\nabla v_M$  for a tetrahedra $G=T\in\mathcal{T}$.
	\end{lemma}
	\begin{proof}
			The difference $v_M-J_Mv_M$ can be rewritten locally for any $G\in \mathcal{T}$ or $G\in\mathcal{F}$ as
			\begin{align*}
				(v_M-J_Mv_M)|_G=p_2+g_{\textit{WF}}+p_7+p_{10}+\sum_{E\in\mathcal{E}(G)}\alpha_E\xi_E
			\end{align*}
			with polynomials $p_2\in P_2(G)$, $p_7\in P_7(G)$, $p_{10}\in P_{10}(G)$, a $\textit{WF}$ 
			function $g_{\textit{WF}}\in P_3(\widehat{\mathcal{T}})$, 
			$\alpha_E\in\mathbb{R}$, and  $\xi_E$ in \eqref{eq:def_xiE} with 
			$\textup{supp}(\xi_E)=B(\textup{mid}(E),R_E)=:B$ for any edge $E\in\mathcal{E}(G)$. 
			\begin{enumerate}[label=\textit{Case \arabic*.}, wide ]
			\item 
					Suppose $G=F\in\mathcal{F}(T)$ is a face and without loss of generality regard $F\subset\mathbb{R}^2$.  
					The key argument is that the support $\textup{supp}(\xi_E)\subset \widehat{\omega}(E)$ of $\xi_E$ lies in the edge patch of 
					$E\in\mathcal{E}$ in $\widehat{\mathcal{T}}$. 
					In other words the $\textit{WF}$ 
					function $g_{\textit{WF}}|_{\overline{\omega}}\in P_3(\overline{\omega})$ is a cubic polynomial in 
					$\overline{\omega}:=\textup{supp}(\xi_E)\cap G$ 
					and the intersection $\omega=B\cap F\subset\mathbb{R}^2$ is the half ball in 
					$\mathbb{R}^2$ illustrated in \cref{fig:Appendix_inverseestimate}.a.
					Hence, it remains to show, for all $p_{10}\in P_{10}(\omega)$  and $\alpha_E\in \mathbb{R}$, that 
					\begin{align*}
						|p_{10}+\alpha_E\xi_E|_{H^1(\omega)}\lesssim h_E^{-1}\Vert p_{10}+\alpha_E\xi_E\Vert_{L^2(\omega)}.
					\end{align*}
					Let $H_+:=\mathbb{R}\times (0,\infty)$ and consider the re-scaling $\varphi: B(0,R_E)\to B(0,1)$, 
					$x\mapsto \varphi(x)=x/R_E:=\vartheta$ with $ |\textup{det}D\varphi|=R_E^{-2}$ and abbreviate $\hat{f}(x):=f(xR_E)$ for $|x|<1$, 
					such 
					that  $f(x)=\hat{f}(x/R_E)=\hat {f}(\vartheta)$ and 
					$\nabla_{\vartheta}\hat{f}(\vartheta)=R_E\nabla_x f(x)$.
					A substitution shows
					\begin{align*}
						\Vert \hat{f}\Vert_{L^2(B(0,1)\cap H_+)}^2
						&=R_E^{-2} \,\Vert {f}\Vert_{L^2(B(0,R_E)\cap H_+)}^2, \\
						\vert\hat{f} \vert_{H^1(B(0,1)\cap H_+)}^2
							&=\int _{B(0,1)\cap H_+}|\nabla_y \hat{f}(y)|^2\,\textup{d}y
							=\int _{B(0,R_E)\cap H_+}|\nabla_\vartheta\hat{f}(\vartheta)|^2|\textup{det}D\varphi|\,\textup{d}x
						\\&	=\int _{B(0,R_E)\cap H_+}R_E^2|\nabla_x f(x)|^2|\textup{det}D\varphi|\,\textup{d}x
							=\vert {f} \vert_{H^1(B(0,R_E))\cap H_+)}^2.
					\end{align*}
					For radius one there exists a constant $C_1>0$, such that for all $p_{10}\in P_{10}$ and $\alpha_E\in \mathbb{R}$,
					\begin{align}
						|p_{10}+\alpha_E\xi_E|_{H^1(B(0,1)\cap H_+)}\le C_1\Vert p_{10}+\alpha_E\xi_E\Vert_{L^2(B(0,1)\cap H_+)},\label{eq:proof_inverseestimate_face}
					\end{align}
					since $\Vert p_{10}+\alpha_E\xi_E\Vert_{L^2(B(0,1)\cap H_+)}=0$ on the right-hand side of \eqref{eq:proof_inverseestimate_face} 
					implies on the left-hand side $|p_{10}+\alpha_E\xi_E|_{H^1(B(0,1)\cap H_+)}=0$. 
					The substitution shows for $\omega= B(0,R_E)\cap H^+$, 
					\begin{align*}
						|p_{10}+\alpha_E\xi_E|_{H^1(\omega)}\le C_1 R_E^{-1}\Vert p_{10}+\alpha_E\xi_E\Vert_{L^2(\omega)}
					\end{align*}
					and $R_E\approx h_E\approx h_F\approx h_T$ concludes the proof since $C_1$ is independent of any mesh-size factor. 
			\item 
					For a tetrahedron, the shape regularity of the  $\textit{WF}$ 
					triangulation $\widehat{\mathcal{T}}$ as well as $\mathcal{T}\in\mathbb{T}$ 
					allow the restriction to 
					$G=T\in \widehat{\mathcal{T}}$. 
					Hence, it remains to show for all $p_{10}\in P_{10}(T)$ and $\alpha_E\in \mathbb{R}$,
					\begin{align*}
						|p_{10}+\alpha_E\xi_E|_{H^1(T)}\lesssim h_T^{-1}\Vert p_{10}+\alpha_E\xi_E\Vert_{L^2(T)}.
					\end{align*}
					The key argument in this case is the existence of $\widehat{B}=B(\widehat{c}_T, R_E)$ with $\widehat{B}\cap B=\emptyset$. 
					A triangle inequality, the inverse estimate for $p_{10}\in P_{10}(T)$ with constant $c_{10}>0$ and for $\xi_E$ 
					in \eqref{eq:def_xiE} with directly calculated constant show 
					$$\vert p_{10}+\alpha_E\xi_E\vert_{H^1(T)}
								\le \frac{6\sqrt{6}}{\sqrt{19}R_E}\,\Vert \alpha_E\xi_E\Vert_{L^2(T)}+ h_T^{-1}c_k\Vert p_{10}\Vert_{L^2(T)}.$$
					Since the shape regularity ensures $R_E\approx h_T$, a triangle inequality leads to 
					\begin{align*}
						\vert p_{10}+\alpha_E\xi_E\vert_{H^1(T)}&
								\lesssim h_T^{-1} \big(\Vert \alpha_E \xi_E\Vert_{L^2(T)}+ \Vert p_{10}\Vert_{L^2(T)}\big)
								\le h_T^{-1} \big(\Vert \alpha_E \xi_E+p_{10}\Vert_{L^2(T)}+ 2\Vert p_{10}\Vert_{L^2(T)}\big).
					\end{align*}	
					The shape regularity guarantees further that there  exists a ball $\widetilde{B}:=B(\widehat{c}_T,\widetilde{R})$ with radius 
					$\widetilde{R}\approx h_T\approx R_E$ that contains $T$, 
					$\widehat{B}\subset T\subset\widetilde{B}$. 
					Standard scaling arguments for 
					polynomials show  for the polynomial extension $p_{10}\in P_{10}(\widetilde{B})$ (not-relabelled), 
					that 
					\begin{align*}
						\Vert p_{10}\Vert_{L^2(T)}\le \Vert p_{10}\Vert_{L^2(\widetilde{B})}\lesssim\Vert p_{10}\Vert_{L^2(\widehat{B})}.
					\end{align*}				 
					The constant in the last step is bounded because the radius quotient $\widetilde{R}/R_E\lesssim 1$ is bounded. 
					Since the support of $\xi_E$ satisfies $\textup{supp}(\xi_E)\cap \widehat{B}=\emptyset$, the estimate
					$\Vert p_{10}\Vert_{L^2(\widehat{B})}= \Vert p_{10}+\alpha_E \xi_E\Vert_{L^2(\widehat{B})}\le  \Vert p_{10}+\alpha_E \xi_E\Vert_{L^2(T)}$
					concludes the proof. 
			\end{enumerate}
	\end{proof}			
		\subsection{Design of $J_3$}	\label{sec:Morley_J3} 
					This step corrects the integral mean of the normal derivative over each internal face $F\in\mathcal{F}(\Omega)$ such that 
					$\intmean_F \nu_F\cdot\nabla J_3 v_M\,\textup{d}s=\intmean_F\nu_F\cdot \nabla v_M\,\textup{d}s.$ 
					Suppose $T=\textup{conv}\{z_1,\dots,z_{4}\}\in\mathcal{T}$ denotes a tetrahedron with face  
					${F}=\textup{conv}\{z_1,z_2, z_3\}\in{\mathcal{F}}(T)$ opposite to the 
					vertex $z_4\in\mathcal{V}(T)$. Recall the barycentric coordinate ${\lambda}_k\in P_1({T})$ in $T$ associated with 
					$z_k\in\mathcal{V}(T)$ for $k=1,\dots, 4$ 
					and set 
					\begin{align*}
						\zeta_{F}|_T
							:=\frac{7!}{2} (\nu_T\cdot\nu_{{F}})\textup{dist}(z_{4},{F}) ({\lambda}_1\lambda_2{\lambda}_3)^2\, {\lambda}_{4} 
							\in P_{7}(T) 
					\end{align*}
					for $T\in\mathcal{T}(F)$. This defines $\zeta_{{F}}\in P_{7}(\mathcal{T}(F))\cap W^{2,\infty}_0(\omega(F))$ 
					with $\intmean_{{G}} \nabla \zeta_{{F}}\cdot\nu_{G}\,\textup{d}s =\delta_{{G}{F}}$ for all sides 
				 	$F,{G}\in{\mathcal{F}}$ in the face 
				 	patch $ \omega_{{F}}:=\textup{int}\big(\bigcup_{T\in \mathcal{T}({{F}})} T\big)$; cf. \cite[Prop.2.6]{Gal15} for $2$D. 
				 	\begin{definition}[$J_3$] 
					For any $v_M\in M(\mathcal{T})$ set 
						\begin{align}
								 J_3(v_M):= J_2 		
								 				v_M+\sum_{F\in\mathcal{F}(\Omega)}\Big(\intmean_{{F}}\nabla(v_M-J_2v_M)
								 					\cdot\nu_{{F}}\,\textup{d}s\Big)\zeta_{{F}}
								 	\in V.\label{eq:def_J2_M}
						\end{align}
					\end{definition}
					\begin{lemma}[properties of $J_3$]\label{lem:prop_J3}
						There exists a constant $C_3\approx 1$ (that exclusively depends on $\mathbb{T}$) 
					 	such that  the companion $J_3v_M\in V$  for any $v_M\in M(\mathcal{T})$ satisfies
					 	 \begin{enumerate}[label=(\alph*), ref=\alph* ,wide]
					 		\item\label{item:App_J3_rightInverse} 
					 				$J_3$ is a right inverse to the interpolation $I_M$ in that
					 					$
					 						I_M\circ J_3 =\textup{id} \text{ in  } {M(\mathcal{T})}, 
					 					$
							\item\label{item:App_J3_approx} 
								$			\Vert h_{\mathcal{T}}^{-{2}} (1-J_3)v_{M}\Vert_{L^2(\Omega)}
									+\vvvert (1-J_3)v_{M}\vvvert_{\mathrm{pw}}
									\le {C_3} \min_{v\in V}\vvvert v_{M}-v\vvvert_{\mathrm{pw}}.$
					 	\end{enumerate}
					\end{lemma}
						\begin{proofof}\textit{(\ref{item:App_J3_rightInverse}). }			
							Since $\zeta_F|_{\partial T}\equiv 0$, it holds $L_E(J_3v_M)=L_E(J_2v_M)$ for any $E\in\mathcal{E}$. 
							Hence \cref{lem:prop_J2}.\ref{item:App_J2_rightInverse} implies $L_E(J_3v_M)=L_E(v_M)$ and 
							by construction holds $L_F(J_3v_M)=L_F(v_M)$ for any $F\in\mathcal{F}(\Omega)$. 
							This proves that $J_3$ satisfies the right-inverse property (\ref{item:App_J3_rightInverse}), since 
							$$
								I_M(v):=\sum_{F\in\mathcal{F}(\Omega)}\intmean_F\nabla v\cdot\nu_F\,\textup{d}\sigma\, \phi_F
								+\sum_{E\in\mathcal{E}(\Omega)}\intmean_E v\,\textup{d}s\,\phi_E\quad \text{for any } v\in V.
							$$
						\end{proofof}
						\begin{proofof}\textit{(\ref{item:App_J3_approx}). }			
							The inverse estimate in \cref{lem:inverseEstimate} guarantees in particular for $g:=v_M-J_2v_M$,  
							$\Vert\nabla g\Vert_{L^2(F)}\le h_F^{-1} \tilde{c}_{\mathrm{inv}}\Vert g\Vert_{L^2(F)}$. 
							Moreover the discrete trace inequality \eqref{eq:F_traceIneq} holds 
							with a constant $\tilde{C}_{\mathrm{inv}}$ from \cref{lem:inverseEstimate}. 
							For any  $F\in\mathcal{F}$, $K\in\mathcal{T}(F)$, and  $v_M\in M(\mathcal{T})$, this and a Cauchy-Schwarz inequality imply						
							\begin{align*}
								\intmean_{{F}}\nabla(v_M-J_2v_M)\cdot\nu_{{F}}\,\textup{d}s	
									&\le |F|^{-1/2} \Vert \nabla(v_M-J_2v_M)\Vert_{L^2(F)}	
									\le |F|^{-1/2}h_F^{-1}	\tilde{c}_{\mathrm{inv}}\Vert v_M-J_2v_M\Vert_{L^2(F)}		
								\\&	\le |K|^{-1/2}h_F^{-1}	\tilde{c}_{\mathrm{inv}} \sqrt{1+2	\tilde{C}_{\mathrm{inv}}/3}\Vert v_M-J_2v_M\Vert_{L^2(K)}.	
							\end{align*}
							Straightforward computations reveal  that 
							$\Vert \zeta_F\Vert_{L^2(K)}=6 h_K|K|^{1/2}/\sqrt{12155}\approx h_K^{5/2}$ holds for any side $F\in\mathcal{F}$ of $K\in\mathcal{T}(F)$.
							In combination with the last displayed estimate and a triangle inequality, this proves
							\begin{align*}
								\Vert v_M&-J_3v_M\Vert_{L^2(K)}-\Vert v_M-J_2 v_M\Vert_{L^2(K)}
									\le \sum_{F\in\mathcal{F}(K)\cap\mathcal{F}(\Omega)}
										\Big|\intmean_F (v_M-J_2v_M)\,\textup{d}s\Big|\Vert\zeta_F\Vert_{L^2(K)}
									\\&\le 4\frac{ 6 h_K}{\sqrt{12155}}\frac{\tilde{c}_{\mathrm{inv}} \sqrt{1+2	\tilde{C}_{\mathrm{inv}}/3}}{h_F}	 
											\Vert v_M-J_2v_M\Vert_{L^2(K)}
									\lesssim \Vert v_M-J_2v_M\Vert_{L^2(K)}
							\end{align*} for any $K\in\mathcal{T}$ with the shape-regularity $h_K\approx h_F$ in the last step. 
							The  inverse estimate 
							$\vert v_M-J_3v_M\vert_{H^2(K)}\lesssim\Vert h_K^2 (v_M-J_3v_M)\Vert_{L^2(K)}$ from \cref{lem:inverseEstimate}, 
							a summation over all $K\in\mathcal{T}$, and \cref{lem:prop_J2}.\ref{item:App_J2_approx} conclude the proof. 	
						\end{proofof}	
		\subsection{Design of $J_4=:J_M$} \label{sec:Morley_J4}
			In order to assure the $L^2$ orthogonality of $v_M-J_Mv_M$ onto $P_2(\mathcal{T})$, 
			we propose a local correction for each tetrahedron   $T\in\mathcal{T}$ 
			with the product of a squared volume bubble-function and a 
			Riesz representation in $P_2(T)$. 
			For any 
			$T:=\textup{conv}\{z_1,z_2,z_3,z_4\}\in\mathcal{T}$,  let $\lambda_k\in P_1(T)$ denote 
			the barycentric coordinate  in $T$  for $k=1,\dots,4$ 
			associated with the 
			vertex $z_k\in \mathcal{V}(T)$.
			Let 
			\begin{align}
				b_T:=4^8\prod_{k=1}^4\lambda_k^2\in P_{8}(T)\cap H^2_0(T)\subset V\label{eq:def_bT}
			\end{align}			 
			denote the squared volume bubble that satisfies $\Vert b_T\Vert_{L^\infty(T)}=1$.
			Let $v_T\in P_2(T)$ denote the Riesz representation of the linear functional 
	  		\begin{align*}
	  			\ell_T:P_2(T)\to \mathbb{R}, \qquad w_T\mapsto \int_T(v_M-J_3v_M)w_T\,\textup{d}x
	  		\end{align*}
	  		in the Hilbert space $P_2({T})$ endowed with the weighted $L^2$ scalar product $(b_T\bullet,\bullet)_{L^2(T)}$.
			The Riesz representation $v_T\in P_2(T)$ satisfies for any $v_M\in M(\mathcal{T})$ 
			\begin{align}
				(v_M-J_3v_M,w_T)_{L^2(T)}=(b_Tv_T,w_T)_{L^2(T)}\quad\text{for all }w_T\in P_2(T). \label{eq:J4_RieszRep}
			\end{align}
			It can be computed  as 
			 $v_T=\sum_{\ell=1}^{10}\Big(\int_T (v_{M}-J_3v_{M})\Phi_{T,\ell}\,\textup{d}x\Big)\Phi_{T,\ell}b_T$ 		
			 with a basis   $\Phi_{T,1}, \dots, $ $\Phi_{T,10}$  of $P_2(T)$ that satisfies 
			$\int_T \Phi_{T,\ell}\Phi_{T,k}b_T\,\textup{d}x=\delta_{\ell,k}$ for any $\ell,\, k=1,\dots, 10$.
			\begin{definition}[$J_4\equiv J_M$]\label{def:JM}				
				Define $J_4 v_M\in V$ for any $v_M\in\mathcal{M}(T)$ by
			\begin{align}
				J_Mv_M := J_4 v_M:=J_3v_{M}+\sum_{T\in\mathcal{T}}v_Tb_T.\label{eq:def_J3_M}
			\end{align}
			\end{definition} 
		\begin{lemma}[properties of $J_4\equiv J_M$]
			For any $v_M\in\mathcal{M}(T)$ the conforming companion $J_4v_M=J_M v_M$ satisfies Theorem~3.1.
		\end{lemma}
		\begin{proofof}\textit{Theorem~3.1.\ref{item:App_J_rightInverse}.}
					Since  $(v_Tb_T)|_{\partial T}\equiv 0 \equiv (D(v_Tb_T))|_{\partial T}$  vanishes along the boundary of any $T\in\mathcal{T}$, 
					$L_E(J_Mv_M)=L_E(J_3v_M)$ and $L_F(J_Mv_M)=L_F(J_3v_M)$ holds for any $E\in\mathcal{E}$ and $F\in\mathcal{F}$. 
					Hence $J_M\equiv J_4$ satisfies the right inverse property according to \cref{lem:prop_J3}.\ref{item:App_J3_rightInverse}.
		\end{proofof}
		\begin{proofof}\textit{Theorem~3.1.\ref{item:App_J_approx}.}
			For $w_T=v_T\in P_2(T)$, \eqref{eq:J4_RieszRep} and  a Cauchy-Schwarz inequality imply
			\begin{align}
				(b_Tv_T,v_T)_{L^2(T)}=(v_M-J_3v_M,v_T)_{L^2(T)}\le \Vert v_M-J_3v_M\Vert_{L^2(T)}\Vert v_T\Vert_{L^2(T)}.\label{eq:J4_approx_1}
			\end{align}
			The volume bubble $b_T$ in \eqref{eq:def_bT} generates an inverse estimate (cf. \cite[\S 3.6]{Verf2013})
			\begin{align}
				C_{\mathrm{inv}}\Vert w_T\Vert_{L^2(T)}\le \Vert b_T^{1/2}w_T\Vert_{L^2(T)}\le\Vert w_T\Vert_{L^2(T)}\quad\text{for all }w_T\in P_2(T)
				\label{eq:J4_approx_2}
			\end{align}
			with a universal constant $C_{\mathrm{inv}}\approx 1$ and the bound $\Vert b_T^2\Vert_{L^{\infty}(T)}=1$ in the last step. 
			The combination of \eqref{eq:J4_approx_1}--\eqref{eq:J4_approx_2} proves
			\begin{align*}
					C_{\mathrm{inv}}^2\Vert v_T\Vert_{L^2(T)}^2\le \Vert b_T^{1/2}v_T\Vert_{L^2(T)}\le \Vert v_M-J_3v_M\Vert_{L^2(T)}\Vert v_T\Vert_{L^2(T)}.
			\end{align*}
			Hence $C_{\mathrm{inv}}^2\Vert v_T\Vert_{L^2(T)}\le \Vert v_M-J_3v_M\Vert_{L^2(T)}$. 
			The definition \eqref{eq:def_J3_M} and $\Vert b_T^2\Vert_{L^{\infty}(T)}=1$ imply 
			\begin{align*}
				\Vert J_Mv_M-J_3v_M\Vert_{L^2(T)}=\Vert v_T b_T\Vert_{L^2(T)}\le \Vert v_T\Vert_{L^2(T)}\le C_{\mathrm{inv}}^{-2} \Vert v_M-J_3v_M\Vert_{L^2(T)}.
			\end{align*}	
			This holds for all $T\in\mathcal{T}$ and so a triangle inequality implies 
			\begin{align*}
				\Vert h_{\mathcal{T}}^{-2}(v_M-J_Mv_M)\Vert_{L^2(\Omega)}
					&=\Vert h_{\mathcal{T}}^{-2}(v_M-J_3v_M)\Vert_{L^2(\Omega)}+\Vert h_{\mathcal{T}}^{-2}(J_Mv_M-J_3v_M)\Vert_{L^2(\Omega)}
					\\&\le (1+C_{\mathrm{inv}}^{-2} )\Vert h_{\mathcal{T}}^{-2}(v_M-J_3v_M)\Vert_{L^2(\Omega)}.
			\end{align*}					 
			The combination with \cref{lem:prop_J3}.\ref{item:App_J3_approx} concludes the proof of 
			\begin{align*}
				\Vert h_{\mathcal{T}}^{-2}(v_M-J_Mv_M)\Vert_{L^2(\Omega)}\le (1+C_{\mathrm{inv}}^{-2} )C_3 \min_{v\in V}\vvvert v_{M}-v\vvvert_{\mathrm{pw}}.
			\end{align*}
			This and the inverse estimate $\vvvert v_M-J_Mv_M\vvvert_{\mathrm{pw}}\le \Vert h_{\mathcal{T}}^{-2}(v_M-J_Mv_M)\Vert_{L^2(\Omega)}$ from 
			\cref{lem:inverseEstimate} conclude the proof of (\ref{item:App_J_approx}). 
		\end{proofof}
		\begin{proofof}\textit{Theorem~3.1.\ref{item:App_J_orthogonality}.}
			The definition \eqref{eq:def_J3_M} and the identity \eqref{eq:J4_RieszRep} imply 
			\begin{align*}
				(v_M-J_Mv_M,w_T)_{L^2(T)}=(v_M-J_3v_M,w_T)_{L^2(T)}-(b_Tv_T,w_T)_{L^2(T)}=0\quad \text{for all }w_T\in P_2(T).
			\end{align*}
			This proves the claim. 
		\end{proofof}
		\begin{remark}[general boundary conditions] 
			The construction in Subsections~\ref{sec:Morley_J1}--\ref{sec:Morley_J4} works for more general boundary conditions of 
			$M(\mathcal{T})$. 
			In \eqref{eq:def_J1_z}--\eqref{eq:def_J1_E} the vanishing boundary values can be replaced by the same averaging  formula as 
			in the interior degrees of freedom. 
			This simplifies the proof of \cref{lem:prop_J1} in that \eqref{eq:vM-J1vM_case3} holds directly for any $z\in\mathcal{V}$ etc. 
		 	The corrections in \eqref{eq:def_J2_M_1}--\eqref{eq:def_J2_M} are possible for boundary edges and sides as well and guarantee 
		 	the right inverse property. \cref{sec:Morley_J4} is independent of boundary conditions. \hfill $\Box$
		\end{remark}
	
		\begin{footnotesize}
		\bibliographystyle{alpha}
\bibliography{BibGLB}
		\end{footnotesize}